\documentclass[a4paper,11pt]{article}
\usepackage{epsfig}
\usepackage{xcolor}
\usepackage{pgfplots}
\usetikzlibrary{patterns}
\usepackage{tikz-cd}
\usepackage{graphics}
\usepackage{graphicx}
\usepackage{amsmath}
\usepackage{amssymb}
\usepackage{amsthm}
\usepackage{authblk}
\usepackage{multirow}
\usepackage{longtable}
\usepackage{lscape}
\usepackage{booktabs}
\usepackage{lineno}
\usepackage{caption}
\usepackage{subcaption}
\newtheorem{thm}{Theorem}
\newtheorem{rem}{Remark}

\newtheorem{prop}[thm]{Proposition}
\newtheorem{lem}[thm]{Lemma}
\newtheorem{cor}[thm]{Corollary}

\usepackage{epstopdf}
\usepackage[margin=1in]{geometry}
\allowdisplaybreaks

\usepackage{natbib}
\setcitestyle{authoryear,open={(},close={)}}
\usepackage{datetime}

\makeatletter
\newcommand{\pushright}[1]{\ifmeasuring@#1\else\omit\hfill$\displaystyle#1$\fi\ignorespaces}
\newcommand{\pushleft}[1]{\ifmeasuring@#1\else\omit$\displaystyle#1$\hfill\fi\ignorespaces}
\makeatother

\begin{document}
\title{Two-Stage Robust Optimization with Decision Dependent Uncertainty}
\author {Bo Zeng and Wei Wang}
\date{}
\maketitle

\vspace{-10pt}

\renewcommand{\baselinestretch}{1.4}

\abstract{The type of decision dependent uncertainties (DDUs) imposes a great challenge in decision making, while existing methodologies are not sufficient to  support many real practices.  In this paper, we present a systematic study to handle this challenge in two-stage robust optimization~(RO). Our main contributions include three sophisticated  variants of column-and-constraint generation method to exactly compute DDU-based two-stage RO. By a novel application of core  concepts of linear programming, we provide rigorous analyses on their computational behaviors. Interestingly, in terms of the iteration complexity of those algorithms, DDU-based two-stage RO is not more demanding than its decision independent uncertainty (DIU) based counterpart. It is worth highlighting a counterintuitive  discovery that converting a DIU set into a DDU set by making use of ``deep knowledge'' and then computing the resulting DDU-based formulation may lead to a significant improvement. Indeed, as shown in this paper, in addition to capturing the actual dependence existing in the real world, DDU is a powerful and flexible tool to represent and leverage analytical properties or simply domain expertise to achieve a strong solution capacity. So, we believe it will open a new direction to solve large-scale DIU- or DDU-based RO. Other important results include basic  structural properties for two-stage RO, an approximation scheme to deal with mixed integer recourse, and a couple of enhancement techniques for the developed algorithms, as well as an organized numerical study to help us appreciate all algorithms and enhancement techniques' computational performances.}

\section{Introduction}
\vspace{-2pt}
Addressing critical uncertainties by analytical approaches is always a central theme in all decision making related disciplines. It has been well recognized that, if uncertainties are handled inappropriately or insufficiently, the derived decision could be either infeasible or with a poor performance once implemented in practice. Hence, many research efforts have been devoted to characterizing uncertainties, incorporating them by extending deterministic decision making models, and developing  efficient solution methodologies, as well as investigating important real issues affected by serious uncertainties.

In the literature, there are three popular schemes to capture uncertainties within an optimization problem. One is stochastic programming (SP), where the uncertainty can be described by a distribution or can be closely approximated by it. One is robust optimization (RO), which uses a set to include all concerned possibilities (hereafter referred to as the \textit{uncertainty set}).  The third one is distributionally robust optimization (DRO), which can be treated as an integration of SP and RO. It assumes that the uncertainty follows a distribution, while that distribution could be anyone in a set of distributions. Currently, all those schemes, along with their algorithms, have been widely adopted in different applications according to  the decision making context and  the nature of underlying uncertainty.

Traditionally, the description of uncertainty is static, i.e., it is not affected by the choice of decision. For the corresponding decision making problem, it means that the distribution, the uncertainty set, or the set of distributions, is fixed for SP, RO and DRO respectively. Such descriptive method is appropriate for many random factors (and their direct consequences) occurring in the natural world, such as wind speed and natural disasters, whose uncertain behaviors are independent of decision making.
Also, if the decision maker has little control on but just passively responds to a random factor existing in a societal system, that factor can also be considered as static in her decision making model. A typical example is the classical newsvendor model, where the demand uncertainty of a perishable product is often described by a fixed distribution regardless of the order quantity of the newsvendor. In the literature, static or fixed uncertainty is  commonly known as \textit{decision-independent uncertainty} (DIU).

Nevertheless, it is often the case in real life that some random factor is substantially affected by the choice of decision, which is therefore referred to as \textit{decision-dependent uncertainty} (DDU). For example,  the cost of one product is  uncertain until it is actually produced~\citep{jonsbraaten1998class}. Hence, a production decision serves not only an instruction to produce but also an investment to refine the information on the production cost. Another obvious situation happens in system maintenance, where components' reliabilities or failure rates change with respect to maintenance decisions \citep{kobbacy2008complex,zhu2021multicomponent}. Similarly, when protecting critical infrastructure systems, hardening some component can convert it from vulnerable to  very robust to future destruction \citep{brown2006defending}.  Indeed, one interesting phenomenon, i.e., \textit{induced demand}, has been observed, which  captures the impact of capacity expansion decision on traffic demand. As more future traffic demand, which can only be estimated, could be generated when highway's capacity is increased, the system planner should take such decision induced demand into account when making any expansion improvement to transportation infrastructure. Other examples could be found in network design and facility location problems \citep{ahmed2000strategic}, offshore oil and gas planning \citep{goel2004stochastic}, the R\&D project portfolio management
\citep{solak2010optimization}, and those listed in \citep{apap2017models}. We also note in the current literature that decision-dependent uncertainty is often called endogenous uncertainty.

\vspace{-5pt}
\subsection{Optimization with DDUs}
To make sound decisions under the aforementioned circumstances, the impact of decisions on the concerned uncertainties should be taken into account. Nevertheless, compared to DIUs that are well-investigated with numerous publications, DDUs are essentially more complex and challenging, and the related research is much less.

 The DDU issue is initially investigated within an SP model \citep{jonsbraaten1998class}, where the underlying scenario tree of the production cost changes with respect to decisions on what items and when they will be produced. Since then, a non-trivial amount of studies adopt SP to model and deal with DDUs in different applications, e.g.,  infrastructure network design \citep{bhuiyan2020stochastic}, power system capacity expansion \citep{zhan2016generation} and operating room scheduling \citep{hooshmand2018adapting}. Nevertheless, the methodological study on this topic remains limited, especially for the case where parameters of the distribution are affected by decisions.  Up to now, there are just a handful of related papers available in the literature, e.g., \citet{apap2017models}, \citet{hellemo2018decision},  \citet{motamed2021multistage}, \citet{pantuso2021node} and related references therein.

As for RO and DRO, which are rather new optimization schemes, we note that only a few studies consider DDUs. For the software partitioning problem that minimizes the total run time,  \citet{spacey2012robust} formulate a single-stage DDU-based RO model where the DDU set captures uncertain decision-dependent execution orders and calling frequencies.
Robust combinatorial optimization problems with a structured DDU polytope  imposed on costs are studied in  \citet{poss2014robust}. A general type of DDU-based robust linear optimization model is considered in \citet{nohadani2018optimization}, where the model is proven to be NP-complete. General DDU set construction and modeling methods are investigated in  \citet{lappas2018robust} and \citet{feng2021multistage}. It is noted that single stage DDU-based RO can often be reformulated as monolithic models  computable by professional mixed integer program (MIP) solvers. To solve DDU-based two-stage and multistage formulations, decision rule \citep{zhang2020unified} and parametric programming based computational methods \citep{avraamidou2020adjustable} have been developed. It is interesting to note a couple of studies presented in \citet{vayanos2020robust,vayanos2020active} that adopt DDU sets to capture uncertain information whose availability  is decision-dependent, a clear demonstration of active learning. Similarly, DRO with DDUs has also been studied  recently \citep{noyan2018distributionally,luo2020distributionally,ryu2019nurse,yu2020multistage,doan2021distributionally,basciftci2021distributionally,feng2021multistage}. One assumption often made is that a DDU set has a finite decision-independent support. It is generally valid and beneficial for DRO, as assigning a zero probability to a particular element equivalently excludes it from the support.
Hence, decision-dependence is reflected by varying probability mass functions only.

Specific to research on DDU-based RO, we note that more
systematic studies are definitely needed.  The most desired research is the development of exact and efficient algorithms to handle general problems, as they are virtually missing when practical two-stage and multistage ones need to be solved. On the  one hand, those algorithms are certainly fundamental in theoretical and computational perspectives. On the other hand, they are of an essential modeling value when constructing a DDU set in decision making, especially by a data-driven approach on raw data. Note that two-way interactions between the decision maker and the random factor are introduced once a DDU set is considered. If they cannot be accurately observed and evaluated, we might result in a wrong or inappropriate DDU set. Hence, those exact algorithms, instead of approximation ones, are critical for validating or refining DDU sets, as well as the decision making problem. Undoubtedly,  rigorous theoretical studies, e.g., convergence and computational complexity, on the developed  algorithms should be carried out to ensure their performances and to support their applications.

We mention that structures of DDU sets and connections between DIU and DDU sets should also be closely examined.  A deeper understanding on them could help us improve  modeling methodologies and resulting formulations, develop computational methods, and strengthen our solution capacity. Note that different from DIU set that mainly contains information of randomness,  DDU set presents significant knowledge that links decision and randomness. Taking one step further, we can treat it as a tool to actively convey our knowledge into RO decision making, if hidden connection between the ``concerned'' uncertainty and our decision can be revealed and analytically represented. This strategy or direction, we believe, should be of great theoretical and practical values in future research. A preliminary study is presented in this paper and results  are deemed very promising.  


%
%
%
%
%
%
%
%
%
%
%
%
%
%
%
%
%
%
%
%
%
%
\vspace{-5pt}
\subsection{Main Results and Contributions of the Paper}
To address the insufficiency existing in the current literature,  our investigation and analyses have made three original and fundamental contributions, along with several new or important developments. We list them in the following with brief descriptions.

\noindent $(i)$ \textit{Strong exact algorithms: Fast computation with a flexible framework.} \ We develop three sophisticated variants of the column-and-constraint generation (C\&CG) method to solve DDU-based  two-stage RO exactly. The second variant (and possibly third one) demonstrates a superior solution capacity to handle practical instances with ordinary DDU sets that are subject to left- and right-hand-side (LHS and RHS) decision dependence, respectively. They also provide a general framework into which structure or application oriented customizations can be readily incorporated. We also note that they naturally extend and complement previous research work in a unified manner, and remarkably broaden our capability to handle difficult uncertainty issues in decision making.

\noindent $(ii)$ \textit{Rigorous complexity analyses: A novel showcase of core linear programming concepts.} We present mathematical proofs and analyses on  convergence and computational complexity for all developed algorithms.  Note that their  iteration complexities generalize and are comparable to those developed for DIU set, indicating two-stage RO with DDU is not more demanding than that with DIU theoretically. It is worth noting that our derivations, together with algorithm development, showcase several core linear programming (LP) concepts, including optimality conditions, basis, reduced cost, and projection, which help us resolve fundamental challenges in a new and more sophisticated optimization paradigm. Indeed,  those new analyses and results ensure the applicability and performance of developed algorithms to handle more complex DDU sets or recourse problems.\\
\noindent $(iii)$ \textit{Converting DIU to DDU: A counterintuitive discovery with a great potential. 
} When studying decision making subject to DIU, we generally do not think of any DDU related methodologies. Nevertheless,  as a matter of fact, we are often able to obtain structural properties or strong insights of a DIU set that characterize or connect worst case scenarios with respect to the first stage decision. Mathematically, such ``deep knowledge'', if proven theoretically, helps us convert the original DIU set into a DDU set, and thus build an equivalent DDU-based reformulation. Surprisingly, compared to directly solving the original DIU-based formulation, computing the new reformulation could  converge in notably fewer iterations, generate solutions with better quality, and reduce our solution time significantly. Even if that ``deep knowledge'' is rather heuristic,  computing the DDU-based formulation leads to quantifiable approximation solutions.  We believe that this discovery is encouraging and of a great value for future research. It unveils the value of structural properties or domain expertise (or \textit{knowledge} in a more general sense) in solving practical-scale DIU-based models, and points out a new direction to develop powerful RO methodologies. Note that even for an original DDU set, we can revise it to convey deep knowledge to reduce its complexity so that a better computational performance can be achieved.




\noindent $(iv)$ \textit{Other critical or noteworthy results:} \\
\indent $(a)$ A set of basic structural properties for two-stage RO with DDU are derived and analyzed. Also, a few reformulation strategies that convert structured DDU-based two-stage RO into equivalent DIU-based or single-stage counterparts are presented;

\indent $(b)$ An approximation scheme to handle RO with mixed integer recourse problem is presented. Since it generates both lower and upper bounds to the RO model, the quality of approximation can be easily accessed, which is often very high;

 \indent$(c)$ Two strategies for algorithm enhancements are designed. One is to take advantage of the model structure to minimize the impact of multiple optimal solutions on cutting sets, and the other one is to generate  cutting sets with Pareto optimality.

\vspace{-5pt}
\subsection{Organization of the Paper and Basic Notations}
 The remainder of this paper is organized as follows. Section \ref{sect_general_formulation} introduces the general mathematical formulation of two-stage RO with DDU, presents a set of basic properties, and provides a few reformulation strategies that convert structured DDU-based two-stage RO into equivalent DIU-based or single-stage counterparts. Section  \ref{sect_algorithms} presents three sophisticated variants of C\&CG method  to exactly solve two-stage RO with DDU. For each of them, it provides the mathematical foundation, detailed operations, and analyses of convergence and iteration complexity. With those algorithms, Section \ref{sect_more} studies enhancements and extensions to further improve our solution capacity on RO problems. Section  \ref{sect_numerical} demonstrates the performances of those algorithms and techniques on instances of robust facility location problems considering different types of uncertainty sets.

\textit{Notations}: Throughout this paper, we use bold lower and uppercase letters to denote vectors and matrices, respectively.  
A vector $\mathbf x$'s $i$th component is denoted by $x_i$.  Most sets are denoted by special calligraphic uppercase letters, while blackboard bold  uppercase letters are reserved for several well-recognized sets.  Scalars are marked in regular font. Unless otherwise noted, all vectors are column vectors, and will be labeled with superscript $\intercal$ if transposed.  Parameter $M$ is introduced to specifically denote a sufficiently large number, and $\mathbf 0$ and $\mathbf 1$ represent a vector with all its entries being 0s and 1s, respectively.  
 Also, following the convention, the optimal value of an infeasible minimization (maximization, respectively) problem is set to $+\infty$ (-$\infty$, respectively).

\vspace{-5pt}
\section{DDU-Based Two-Stage Robust Optimization}
\label{sect_general_formulation}
\vspace{-3pt}
In this section, we first introduce the general mathematical formulation of two-stage RO with DDU formally. Then, we derive a set of basic properties of this RO model. Finally, we present several reformulation strategies for structured DDU sets so that the original DDU-based formulation can be converted into one computable by existing methods.

\vspace{-5pt}\subsection{The general formulations}
In a two-stage decision making procedure, the decision maker initially determines the value of the first stage decision variable $\mathbf x$ before the materialization of random factor $\mathbf u$. Then, after the uncertainty is cleared, she has an opportunity to make a recourse decision for mitigation, which, however, is restricted by her choice of $\mathbf x$  and $\mathbf{u}$'s realization. Specifically, let $\mathbf x=(\mathbf x_c,\mathbf x_d)$ denote the first stage decision variable vector with $\mathbf x_c$ and $\mathbf x_d$ representing collectively its continuous   and discrete components, respectively. Similarly, $\mathbf y=(\mathbf y_c,\mathbf y_d)$ and $\mathbf{u}=(\mathbf{u}_c, \mathbf{u}_d)$ denote the recourse decision variable vector and the uncertainty variable vector, respectively, both of which could contain continuous and discrete variables. The general mathematical formulation of  two-stage RO with DDU is
\begin{eqnarray}
\label{eq_2RO}
\mathbf{2-Stg \ RO}: \ \ \ \mathit{w}^*& = & \min_{\mathbf x\in \mathcal{X}}\mathbf c_1\mathbf x+\max_{\mathbf u\in\mathcal U(\mathbf x)} \ \min_{\mathbf y\in\mathcal Y(\mathbf x,\mathbf u)}\mathbf c_2\mathbf y
\end{eqnarray}
where
\begin{align}
\mathcal{X} =\{\mathbf{x}\in \mathbb{Z}^{m_x}_+\times \mathbb{R}^{n_x}_+: \mathbf{Ax}\geq \mathbf{b}\}
\end{align}
is the feasible set of variables $\mathbf x$,
$\mathcal U(\mathbf x)$ is the uncertainty set as in the following,
\begin{align}
\label{eq_uncer_set}
\mathcal U(\mathbf x)=\left\{\mathbf u\in\mathbb Z^{m_u}_+\times \mathbb{R}^{n_u}_+:\mathbf F(\mathbf x)\mathbf u \leq\mathbf h+\mathbf{Gx}\right\},
\end{align}
and $\mathcal Y(\mathbf x,\mathbf u)$ is the feasible set of the recourse problem as in the next.
\begin{align}
\label{eq_recourse_set}
\mathcal Y(\mathbf x,\mathbf u)=\left\{\mathbf y\in\mathbb Z^{m_y}_+\times \mathbb{R}^{n_y}_+:\mathbf B_2\mathbf y\geq\mathbf d-\mathbf B_1\mathbf x-\mathbf E\mathbf u\right\}
\end{align}
Coefficient vectors $\mathbf{c}_1$, $\mathbf{c}_2$ (both are row vectors), $\mathbf{b}$, $\mathbf h$, $\mathbf d$, and constraint matrices $\mathbf{A}$, $\mathbf F(\mathbf x)$ (its column-wise representation is $[F_1(\mathbf  x),\dots,F_{m_u+n_u}(\mathbf{x})]$), $\mathbf{G}$, $\mathbf B_1$, $\mathbf B_2$, $\mathbf E$ are all with appropriate dimensions. Regarding $\mathcal U(\mathbf x)$, the DDU set defined by a point-to-set map in \eqref{eq_uncer_set}, we differentiate two types of decision-dependence: RHS dependence and LHS dependence. The former one has $\mathbf x$ appeared in RHS of \eqref{eq_uncer_set} only, while the latter one has $\mathbf x$ in \eqref{eq_uncer_set}'s LHS  only. Actually, if $\mathcal{U}(\mathbf x)$  have both RHS and LHS dependence, the dependence can be converted into LHS one by appending $-\mathbf {Gx}$ as a column to $\mathbf{F(x)}$ and extending $\mathbf u$ with one more dimension that takes value $1$.

The fully concatenated  RO problem in \eqref{eq_2RO} subsumes the classical DIU-based two-stage RO model as a special case  \citep{ben2004adjustable,bertsimas2010optimality,Bertsimas.2010}. By denoting the DIU set by $\mathcal{U}^0$, $\mathbf{2-Stg \ RO}$ reduces to the DIU-based RO if we set $\mathcal U(\mathbf x) = \mathcal{U}^0$ for all $\mathbf x\in \mathcal{X}$.  As we may vary sets $\mathcal{X}$, $\mathcal U(\mathbf x)$, or $\mathcal Y(\mathbf x,\mathbf u)$  in \eqref{eq_2RO} in the remainder of this paper, let $w(\mathcal{X}, \mathcal{U}(\mathbf x), \mathcal{Y}(\mathbf{x},\mathbf{u}))$ denote the optimal value of \eqref{eq_2RO} subject to those sets.

Using reformulation techniques to derive a new formulation is a key and frequently used strategy in the study of RO. We say that, in the context of two-stage RO, two formulations are equivalent if they share the same optimal value, and one's optimal first stage solution is also optimal to the other one, and vice versa. In the following, we present the epigraph formulation
of DDU-based two-stage RO. Clearly, it is equivalent to  \eqref{eq_2RO}.
\begin{subequations}
	\begin{align}
		\mathit{w}^* =  \ \min&\quad\mathbf c_1\mathbf x+\eta \label{eq_2RO_epigraph-1} \\
		\mbox{s.t.} & \quad \mathbf x\in \mathcal{X} \label{eq_2RO_epigraph-2} \\
		& \ \ \ \eta\geq \{\mathbf c_2\mathbf y: \mathbf y\in\mathcal Y(\mathbf x,\mathbf u)\} \
		\forall \mathbf u\in\mathcal U(\mathbf x)
		\label{eq_2RO_epigraph-3}
	\end{align}
\end{subequations}

 Next, we present three very mild assumptions that practically do not impose restriction. \\
\noindent $(\textit {A1})$ For any $\mathbf{x}\in \mathcal{X}$, $\mathcal{U}(\mathbf x)\neq \emptyset$; \\
\noindent  $(\textit{A2})$ $\mathcal{U}(\mathbf x)$ is a bounded set, i.e., for any given $\mathbf{x}\in \mathcal{X}$, $u(\mathbf x)_j<\infty$ \ $\forall j$;\\
\noindent $(\textit{A3})$ The next MIP  has a finite optimal value.
  \begin{eqnarray}
  \label{eq_MIP_RO}
  \min\{\mathbf{c_1x+c_2y}: \mathbf{x}\in \mathcal{X}, \mathbf{u}\in \mathcal{U}(\mathbf x), \mathbf{y}\in \mathcal{Y}(\mathbf x, \mathbf u)\}
  \end{eqnarray}

Assumption (\textit{A1}) is introduced to substantiate the two-stage decision making framework. Note that if  $\mathcal{U}(\mathbf x^0) =\emptyset$ for some $\mathbf{x}^0$, the recourse problem is not defined and \eqref{eq_2RO} is trivially unbounded.   The next result eliminates such $\mathbf x^0$ from consideration.


\begin{lem}
\label{lem_replication}
Let $\mathbf u'$ be an independent replicate of $\mathbf u$. If the two-stage RO formulation in \eqref{eq_2RO} is non-trivial, it is equivalent to 
\begin{eqnarray}
\label{eq_eliminate_emptyU}
 \min_{\mathbf x\in\mathcal X, \mathbf{u}'\in \mathcal{U}(\mathbf x)}\mathbf c_1\mathbf x+\max_{\mathbf u\in\mathcal U(\mathbf x)} \ \min_{\mathbf y\in\mathcal Y(\mathbf x,\mathbf u)}\mathbf c_2\mathbf y,
\end{eqnarray}
which satisfies assumption \textit{(A1)}. \hfill $\square$
\end{lem}
Note that \eqref{eq_eliminate_emptyU} actually is in the form of \eqref{eq_2RO}.  
Regarding assumption (\textit{A2}), it basically holds as the practical random factor under consideration is bounded in general. Nevertheless, it may not be straightforward that $\mathcal{U}(\mathbf x)$ satisfies that assumption for any $\mathbf x\in \mathcal{X}$. If it concerns us, we can simply augment $\mathcal{U}(\mathbf x)$ by including $\mathbf u\leq \mathbf{\overline u}$, where $\mathbf{\overline u}$ is a reasonable upper bound for that random factor. 

We say $\mathbf x^0\in \mathcal{X}$ is infeasible if the recourse problem is infeasible for some $\mathbf u\in \mathcal{U}(\mathbf x^0)$, and otherwise it is feasible. Note that the whole $\mathbf{2-Stg \ RO}$ is infeasible if no first stage decision in $\mathcal{X}$ is feasible. As in the following, the optimization problem defined in assumption (\textit{A3}) helps us detect its infeasibility.
\begin{lem}
\label{lem_feasibility}
Formulation \eqref{eq_MIP_RO} is a relaxation to \eqref{eq_2RO}. And if \eqref{eq_MIP_RO} is infeasible, so is \eqref{eq_2RO}. \hfill $\square$
\end{lem}
Note that once the infeasibility is observed, correction or revision is probably  necessary for \eqref{eq_2RO}, especially for the dependence reflected in $\mathcal U{\mathbf (x)}$, if it is built for a real problem.   
Hence, in the remainder of this paper, unless otherwise noted, those three assumptions hold. In the next subsections, we present a set of structural properties of $\mathbf{2-Stg \ RO}$, and discuss some special cases that can be solved by existing methods after reformulation.


\vspace{-2pt}
\subsection{Basic Structural Properties}
The next result directly follows from the facts that $(i)$ formulation \eqref{eq_2RO} reduces to DIU-based two-stage RO model if $\mathcal U(\mathbf x) = \mathcal{U}^0$ for all $\mathbf x$, and $(ii)$ two-stage RO with DIU, even in its simplest linear form with pure continuous decisions, is NP-hard \citep{ben2004adjustable}.

\begin{thm}
Formulation $\mathbf{2-Stg \ RO}$ defined in \eqref{eq_2RO} is NP-hard, even if $m_x=m_y=m_u=0$, i.e., only continuous variables are involved. \hfill $\square$
\end{thm}

Next, we discuss several results regarding $\mathbf{2-Stg \ RO}$ in \eqref{eq_2RO} by varying the underlying sets, which help us derive different relaxations.  
\begin{prop}
\label{prop_relax_1}
\begin{description}
  \item[$(i)$] Consider two sets $\mathcal{X}^1$ and $\mathcal{X}^2$ such that  $\mathcal{X}^1\subseteq \mathcal{X}^2$. We have
      $$\mathit{w}(\mathcal X^1, \mathcal{U}(\mathbf x), \mathcal Y(\mathbf x,\mathbf u))
\geq \mathit{w}(\mathcal X^2, \mathcal{U}(\mathbf x), \mathcal Y(\mathbf x,\mathbf u)).$$
  \item[$(ii)$] Consider two sets $\mathcal{U}^1(\mathbf x)\subseteq \mathcal{U}^2(\mathbf x)$ $\forall \mathbf{x}\in \mathcal{X}$. We have
      $$\mathit{w}(\mathcal X, \mathcal{U}^2(\mathbf x), \mathcal Y(\mathbf x,\mathbf u))\geq \mathit{w}(\mathcal X, \mathcal{U}^1(\mathbf x), \mathcal Y(\mathbf x,\mathbf u)).$$
  \item[$(iii)$]  Consider two sets $\mathcal{Y}^1(\mathbf x, \mathbf u)\subseteq \mathcal{Y}^2(\mathbf x, \mathbf u)$ $\forall \mathbf{x}\in \mathcal{X}, \mathbf{u}\in \mathcal{U}(\mathbf x)$. We have
  \[
  	\mathit{w}(\mathcal X, \mathcal{U}(\mathbf x), \mathcal Y^1(\mathbf x,\mathbf u))
  	\geq \mathit{w}(\mathcal X, \mathcal{U}(\mathbf x), \mathcal Y^2(\mathbf x,\mathbf u)). \pushQED{\qed}\qedhere
  \] 
\end{description} 
\end{prop}
The aforementioned results can be proven easily, given that the RO problems in RHS are relaxations to ones in LHS. Next result simply follows from Proposition \ref{prop_relax_1}.  Note that by definition $\mathit{w}^*=\mathit{w}(\mathcal{X}, \mathcal{U}(\mathbf x), \mathcal{Y}(\mathbf x,\mathbf u))$.

\begin{cor}
\label{cor_3_level_relaxation}
  $(i)$ Let $\mathcal{X}_r$, $\mathcal{U}_r(\mathbf x)$ and $\mathcal{Y}_r(\mathbf x,\mathbf u)$ denote the continuous relaxations of $\mathcal{X}$, $\mathcal{U}(\mathbf x)$ and $\mathcal{Y}(\mathbf x,\mathbf u)$, respectively. We have
  \begin{eqnarray*}
  \mathit{w}(\mathcal{X}_r,\mathcal{U}(\mathbf x),\mathcal{Y}_r(\mathbf x,\mathbf u))\leq \mathit{w}(\mathcal{X},\mathcal{U}(\mathbf x),\mathcal{Y}_r(\mathbf x,\mathbf u))\leq \mathit{w}^*\leq
  \mathit{w}(\mathcal{X},\mathcal{U}_r(\mathbf x),\mathcal{Y}(\mathbf x,\mathbf u)).
  \end{eqnarray*}
  $(ii)$ Consider sets $\mathcal{\hat X}\subseteq \mathcal{X}$, $\mathcal{\hat U}(\mathbf x)\subseteq
  \mathcal{U}(\mathbf x)$ for all $\mathbf x\in \mathcal{X}$, and $\mathcal{\hat Y}(\mathbf x,\mathbf u)\subseteq \mathcal{Y}(\mathbf x,\mathbf u)$
  for all $\mathbf x\in \mathcal{X}$ and $\mathbf u\in \mathcal U(\mathbf x)$. We have
  \begin{eqnarray*}
  \mathit{w}(\mathcal{X},\mathcal{\hat U}(\mathbf x),\mathcal{Y}_r(\mathbf x,\mathbf u))\leq \mathit{w}^* \leq \mathit{w}(\mathcal{X},\mathcal{U}(\mathbf x),\mathcal{\hat Y}(\mathbf x,\mathbf u))\leq
  \mathit{w}(\mathcal{\hat X},\mathcal{U}(\mathbf x),\mathcal{\hat Y}(\mathbf x,\mathbf u)).
  \end{eqnarray*}
  $(iii)$ Consider two sets $\mathcal{U}^1(\mathbf x)$ and $\mathcal{U}^2(\mathbf x)$ such that
  $\mathcal{U}^1(\mathbf x)\cup \mathcal{U}^2(\mathbf x)\subseteq \mathcal{U}(\mathbf x)$ for all
  $\mathbf{x}\in \mathcal{X}$, and set $\mathcal{U}^3(\mathbf x)$
  such that  $\mathcal{U}(\mathbf x)\subseteq \mathcal{U}^3(\mathbf x)$
  for all $\mathbf x\in \mathcal{X}$. We have
    \begin{align*}
  	&\max\{\mathit{w}(\mathcal{X},\mathcal{U}^1(\mathbf x),\mathcal{Y}(\mathbf x,\mathbf u)), \mathit{w}(\mathcal{X},\mathcal{U}^2(\mathbf x),\mathcal{Y}(\mathbf x,\mathbf u))\}  \leq \mathit{w}(\mathcal{X},\mathcal{U}^1(\mathbf x)\cup \mathcal{U}^2(\mathbf x),\mathcal{Y}(\mathbf x,\mathbf u)) \leq \\ & \  \ \mathit{w}^*
  	\leq \mathit{w}(\mathcal{X},
  	\mathcal{U}^3(\mathbf x),\mathcal{Y}(\mathbf x,\mathbf u)). \pushQED{\qed}\qedhere
  \end{align*}
\end{cor}

%
%
%
%
%
%
%


\begin{rem}
$(i)$ It is worth highlighting that the results in Proposition \ref{prop_relax_1} and Corollary \ref{cor_3_level_relaxation} do not depend on linear structure of $\mathbf{2-Stg \ RO}$ in \eqref{eq_2RO}. They are actually  valid for general two-stage RO that involve  mixed integer nonlinear program or mixed integer nonlinear set. \\
$(ii)$ Although the results of Corollary \ref{cor_3_level_relaxation} directly follow from Proposition \ref{prop_relax_1}, they provide a mathematical framework for us to construct strong approximations to $\mathbf{2-Stg \ RO}$ as shown in the next sections. For example, by considering the continuous relaxation or some restriction of the recourse problem, we can derive lower or upper bounds of the optimal value of \eqref{eq_2RO}. Similar results can be obtained by modifying the underlying uncertainty set. \\
$(iii)$ We  mention that $\mathcal{U}^1(\mathbf x)$ and $\mathcal{U}^2(\mathbf x)$ do not need to be disjoint, which offers us a great flexibility in constructing and analyzing those decision-dependent sets. Certainly, compared to static uncertainty set that can be easily expanded or split, a deep understanding on the interaction between $\mathbf x$ and $\mathcal{U}(\mathbf x)$ is necessary to perform those operations.
\end{rem}

%

Next, we consider the case where the recourse problem is an LP. For a given $\mathbf x\in\mathcal{X}$, note that by assumption (\textit{A3})  the dual problem of the recourse problem is always feasible. So, we have the following equivalence.
\begin{eqnarray}
	\label{eq_max_min}
	&&  \max_{\mathbf{u}\in \mathcal{U}(\mathbf x)}\min\{\mathbf{c_2y}: \mathbf y\in \mathbb{R}^{n_y}_+: \mathbf B_2\mathbf y\geq \mathbf d-\mathbf B_1\mathbf x-\mathbf E\mathbf u\}\\
	&& = \max\{(\mathbf d-\mathbf B_1\mathbf x-\mathbf E\mathbf u)^\intercal\boldsymbol\pi: \mathbf{u}\in
	\mathcal{U}(\mathbf x), \mathbf B_2^\intercal\boldsymbol\pi \leq \mathbf c_2^\intercal, \boldsymbol\pi\geq \mathbf{0}\}.  \label{eq_max_max}
\end{eqnarray}
Note that \eqref{eq_max_max} is a disjoint bilinear program over two independent sets. For a fixed $\boldsymbol\pi$, it reduces to the following MIP problem
$$\max\{(\mathbf d-\mathbf B_1\mathbf x-\mathbf E\mathbf u)^\intercal\boldsymbol\pi: \mathbf{u}\in \mathcal{U}(\mathbf x)\},$$  or its equivalent linear program  \begin{eqnarray}
\label{eq_LP_convexhull}
\max\{(\mathbf d-\mathbf B_1\mathbf x-\mathbf E\mathbf u)^\intercal\boldsymbol\pi: \mathbf{u}\in co(\mathcal{U}(\mathbf x))\},
\end{eqnarray}
where $co(\mathcal{U}(\mathbf x))$ denotes the convex hull of $\mathcal{U}(\mathbf x)$. According to assumption (\textit{A2}), $co(\mathcal{U}(\mathbf x))$ is a polytope. Hence, the next result simply follows, which 
characterizes optimal, i.e., the worst case, scenarios and generalizes a similar result in~\citet{zeng2013solving}. 

\begin{prop}
	\label{prop_extremepoint}
	Suppose that $m_y=0$. For a given $\mathbf x\in \mathcal{X}$, we have \\
	$(i)$ there exists an optimal solution $\mathbf u^*$ to $\max_{\mathbf{u}\in \mathcal{U}(\mathbf x)} \min_{\mathbf y\in \mathcal{Y}(\mathbf x, \mathbf u)} \mathbf c_2\mathbf y$ that is
	an extreme point of $co(\mathcal{U}(\mathbf x))$; \\
	$(ii)$ $\max_{\mathbf{u}\in \mathcal{U}(\mathbf x)} \min_{\mathbf y\in \mathcal{Y}(\mathbf x, \mathbf u)} \mathbf c_2\mathbf y = \max_{\mathbf{u}\in co(\mathcal{U}(\mathbf x))} \min_{\mathbf y\in \mathcal{Y}(\mathbf x, \mathbf u)} \mathbf c_2\mathbf y.$ \hfill$\square$
\end{prop}

Proposition \ref{prop_extremepoint} is particularly useful when the uncertainty set reduces to a DIU one. Note that the convex hull of a DIU set has a fixed and finite set of extreme points. By enumerating those extreme points (i.e., scenarios) and including a replicate of the recourse problem for each of them, the next result simply follows.

\begin{cor}
	\label{cor_DIU_Omega}
  When $m_y=0$ and $\mathcal{U}(\mathbf x)=\mathcal{U}^0$ for $\mathbf x\in \mathcal{X}$,
  let $\mathcal{P}_{co(\mathcal{U}^0)}=\{\mathbf u^1,\dots,\mathbf u^{|\mathcal{P}_{co(\mathcal{U}^0)}|}\}$ be the set of extreme points of $co(\mathcal{U}^0)$.
  Then, $\mathbf{2-Stg \ RO}$ in \eqref{eq_2RO} is equivalent to the following single-level optimization problem.
	\begin{align}
	\ \ \ \mathit{w}^* = \min_{\mathbf x\in \mathcal{X}} \quad&\mathbf c_1\mathbf x+\eta\nonumber\\
	\begin{split}\mathrm{s.t.} \quad & \mathbf x\in \mathcal{X} \\
		& \eta\geq \mathbf c_2\mathbf y^{k}, \ k=1, \dots, |\mathcal{P}_{co(\mathcal{U}^0)}|
	\end{split}\label{eq_DIU_enumeration}\\	
	&\mathbf y^{k}\in\mathcal Y(\mathbf x,\mathbf u^{k}), \ k=1, \dots, |\mathcal{P}_{co(\mathcal{U}^0)}|\pushQED{\qed}\qedhere\nonumber
\end{align}
\end{cor}
Since the set of extreme points is generally huge, computing this extreme point based reformulation, which is a regular single-level MIP, might not be practically feasible. Nevertheless, a similar formulation built upon a small subset of $\mathcal{P}_{co(\mathcal{U}^0)}$ is a relaxation to \eqref{eq_DIU_enumeration}, and therefore provides a lower bound to its optimal value. Together with some upper bound technologies, this strategy should be of a key value to derive  exact solutions, which indeed is the foundation for the development of basic C\&CG method \citep{zeng2013solving}.

For a general DDU set, the aforementioned reformulation is definitely not applicable, noting that $co(\mathcal{U}(\mathbf x))$ and
its extreme point set change with respect to $\mathbf x$. Yet, it is interesting to note in Section \ref{sect_algorithms} that the philosophy still can be extended to support algorithm development.  Indeed, one nice property of the disjoint bilinear program in \eqref{eq_max_max}  is that once we fix $\mathbf u$ or $\boldsymbol\pi$, the whole program reduces to an MIP or LP that is easier to analyze. The next result takes advantage of this property to analyze \eqref{eq_max_max}.


\begin{cor}
\label{cor_maxmax_dual}
When $m_y=0$, let $\Pi=\{\mathbf B_2^\intercal\boldsymbol\pi \leq \mathbf c_2^\intercal, \  \boldsymbol\pi\geq \mathbf{0}\}$.
For the bilinear program in \eqref{eq_max_max}, when it has a finite optimal value, there exists an optimal solution $(\mathbf u^*,  \boldsymbol\pi^*)$ with $\boldsymbol \pi^*$ being an extreme point of $\Pi$.  When it is unbounded, there exists  $\mathbf u^*\in \mathcal{U}(\mathbf x)$ and $\boldsymbol\gamma^*$, which is an extreme ray of $\Pi$, such that $(\mathbf d-\mathbf B_1\mathbf x-\mathbf{Eu}^*)^\intercal\boldsymbol\gamma^* > 0$.\qed
\end{cor}
\begin{rem}
Unlike $\mathcal{U}(\mathbf x)$,  $\Pi$ is a fixed polyhedron independent of $\mathbf x$ and $\mathbf u$, which indicates the sets of extreme points and extreme rays are finite and fixed. By assumption (A3), $\Pi$ is clearly non-empty. Also, when \eqref{eq_max_max} is unbounded, $\mathbf u^*$ mentioned in Corollary \ref{cor_maxmax_dual} causes the recourse problem to be infeasible. According to definition, the given first stage decision $\mathbf x$ is infeasible to $\mathbf{2-Stg \ RO}$ in \eqref{eq_2RO}.
\end{rem}

\subsection{Reformulations of Special Cases of Two-Stage RO with DDU}
\label{subsect_specialDDU}
In this subsection, we consider two-stage RO with DDU in special structures. By using a few reformulation techniques,  either \eqref{eq_2RO} with those DDU  sets could be converted into DIU-based two-stage RO, or its two-stage structure would reduce to a single-stage one. Then, the resulting reformulations can be directly solved by using algorithms existing in the literature. Actually, those reformulation techniques complement the algorithm development presented in the next section, as most of these structured DDU-based ROs cannot be handled by our new algorithms. We note that the first two reformulation techniques have appeared in simple forms in the context of particular applications, and our work is to abstract them with more generality and rigorous analyses.

\subsubsection{Neutralization Reformulation}
Consider a two-stage RO with a DDU set that is in the form of
\begin{eqnarray}
\label{eq_US_binary}
\mathcal{U}(\mathbf x)=\{\mathbf {u}\in \{0,1\}^{m_u}: \mathbf{Fu}\leq \mathbf{h}, \mathbf u\leq \mathbf {g}(\mathbf x)\},
\end{eqnarray}
 and has the downward closedness property, i.e., if $\mathbf u^1\in \mathcal{U}(\mathbf x)$ and
 $\mathbf u^2\leq \mathbf u^1$,  it follows that $\mathbf u^2\in \mathcal{U}(\mathbf x)$. Actually, because $\mathbf u$ is binary, constraint $\mathbf u\leq \mathbf {g}(\mathbf x)$ naturally has this property. Hence, it is equivalent to say that constraint $\mathbf{Fu}\leq \mathbf{h}$ has this property.

This type of DDU set has been often used in reliability and security applications of critical infrastructures, e.g., \citet{brown2006defending}, \citet{cappanera2011optimal} and \citet{yuan2019cost}, where the uncertainty set represents the components' failures or unavailabilities that could be controlled by the first stage hardening or protection decisions.

Since $g_i(\mathbf x), i=1,\dots, m_u$ are completely determined in the first stage, without loss of generality, we assume that $g_i(\mathbf x)\in\{0,1\} \ \forall i$, and the first stage set $\mathcal{X}$ is augmented by a set of additional binary variables $\mathbf x'$ and constraints $\mathbf x'=\mathbf {g}(\mathbf x)$. Accordingly, we can replace the corresponding constraint in \eqref{eq_US_binary} by $\mathbf u\leq \mathbf x'$. Next, we present a technique that neutralizes  scenarios not in $\mathcal{U}(\mathbf x)$ to build a DIU-based RO reformulation.  
To simplify the exposition in this paper, we denote the
Hadamard product of vectors $\mathbf a$ and $\mathbf b$, which returns the element-wise multiplication of $\mathbf a$ and $\mathbf b$,  by  $\mathbf a\circ \mathbf b$.

\begin{prop}
\label{prop_reform_DIU}
Formulation $\mathbf{2-Stg \ RO}$ with  the DDU set defined in \eqref{eq_US_binary} is equivalent to the following DIU-based 2-stage RO
$$\min_{\mathbf x\in \mathcal{X}}\mathbf c_1\mathbf x+\max_{\mathbf u\in \mathcal{U}^0}\ \min_{\mathbf y\in\bar{\mathcal Y}(\mathbf x,\mathbf u)}\mathbf c_2\mathbf y$$
where $\mathcal{U}^0=\{\mathbf {u}\in \{0,1\}^{m_u}: \mathbf{Fu}\leq \mathbf{h}\}$ and
\begin{eqnarray*}
	\bar{\mathcal Y}(\mathbf x,\mathbf u)=\{\mathbf y\in\mathbb Z^{m_y}_+\times \mathbb{R}^{n_y}_+:\mathbf B_2\mathbf y\geq\mathbf d-\mathbf B_1\mathbf x-\mathbf E(\mathbf u\circ \mathbf x')\},
\end{eqnarray*}
or equivalently
\begin{eqnarray*}
(\mathbf y, \mathbf v)\in \bar{\mathcal Y}(\mathbf x,\mathbf u)&=&\{\mathbf y\in\mathbb Z^{m_y}_+\times \mathbb{R}^{n_y}_+, \mathbf v\in \mathbb{R}^{m_u}_+:\mathbf B_2\mathbf y\geq\mathbf d-\mathbf B_1\mathbf x-\mathbf E\mathbf v\\
&& \ \mathbf v\leq \mathbf x', \mathbf v\leq \mathbf u, \mathbf {v}\geq
\mathbf {x}'+\mathbf u-\mathbf{1}\}.
\end{eqnarray*}
\end{prop}
\begin{proof}
	See its proof in Appendix \ref{apd:proofS2}.
\end{proof}

\begin{rem}
$(i)$ Note that variables $\mathbf v$ introduced in linearizing the recourse problem are continuous. Hence, if the original recourse problem is an LP, introducing $\mathbf v$ does not change its nature. Actually, from the computational point of view, such a linearization is not necessary since both $\mathbf x'$ and $\mathbf u$ are given as parameters for the recourse problem. We can simply use the first set $\bar{\mathcal Y}(\mathbf x,\mathbf u)$ to derive optimal solutions.  \\
$(ii)$ It is interesting to observe the neutralization effect on scenarios that do not belong to $\mathcal U(\mathbf x)$. Through reformulation, those scenarios do not impose any new challenge to the recourse problem as they are projected to some scenarios within $\mathcal {U}(\mathbf x)$. We refer to this  technique as \textit{neutralization reformulation}, and believe that it could be extended to handle other DDU sets. Indeed, as in the following, the result in Proposition \ref{prop_reform_DIU} can be generalized to handle one class of mixed integer DDU sets that has not been investigated in the literature.
\end{rem}

Consider  a DDU set in the following form
\begin{equation}
\label{eq_US_MIP}
\begin{split}
\mathcal{U}(\mathbf x) = \{\mathbf {u}\in \{0,1\}^{m_u}\times \mathbb{R}^{n_u}_+: & \mathbf{Fu}\leq \mathbf{h}, u_i\leq x'_i, \ i=1,\dots, m_u, \\
 & u_i\leq u^0_i x'_i, \ i=m_u+1,\dots, m_u+n_u\},
\end{split}
\end{equation}
where $u^0_i$ is the upper bound parameter for continuous $u_i, \ i=m_u+1,\dots, m_u+n_u$, and $\mathbf x'=\{x'_1,\dots, x'_{m_u+n_u}\}$ are binary variables determined in the first stage. Again,  $\mathcal{U}(\mathbf x)$ has the downward closedness property if \eqref{eq_US_binary} has. Using a neutralization
reformulation similar to that in Proposition~\ref{prop_reform_DIU},  two-stage RO with \eqref{eq_US_MIP} can be converted into a DIU-based one.
\begin{cor}
Formulation $\mathbf{2-Stg \ RO}$ with the DDU set defined in \eqref{eq_US_MIP} is equivalent to the following DIU-based 2-stage RO
$$\min_{\mathbf x\in \mathcal{X}}\mathbf c_1\mathbf x+\max_{\mathbf u\in \mathcal{U}^0}\ \min_{\mathbf y\in \mathcal Y'(\mathbf x,\mathbf u)}\mathbf c_2\mathbf y$$
where $\mathcal{U}^0=\{\mathbf {u}\in \{0,1\}^{m_u}\times \mathbb{R}^{n_u}_+:
\mathbf{Fu}\leq \mathbf{h}, \ u_i\leq  u^0_i, i=m_u+1,\dots, m_u+n_u\}$ and
\begin{equation*}
\mathcal Y'(\mathbf x,\mathbf u) = \{\mathbf y\in\mathbb Z^{m_y}_+\times \mathbb{R}^{n_y}_+:\mathbf B_2\mathbf y\geq\mathbf d-\mathbf B_1\mathbf x-\mathbf E(\mathbf u\circ \mathbf x')\}. \pushQED{\qed}\qedhere
\end{equation*}
\end{cor}

\subsubsection{Normalization Reformulation}
Another special DDU set is in the form of
\begin{eqnarray}
\label{eq_US_hypercube}
\mathcal U(\mathbf x)=\{\mathbf u\in \mathbb R^{n_u}_+:   x_{li}\leq  u_i \leq  x_{hi} \ \forall i\},
\end{eqnarray}
i.e., a hypercube set with $ x_{li}$ and $ x_{hi}$ being the first stage decision variables. This type of uncertainty sets arises from an interesting application in power grids where the system operator determines the do-not-exceed (DNE) limits for the random renewable generation. As long as the  renewable generation is within the range defined by those limits, the grid can, with recourse operations, accommodate without sacrificing its reliability \citep{zhao2014variable,wang2016robust}. As noted in \citet{zhao2014variable}, although set \eqref{eq_US_hypercube} is subject to DNE limits determined in the first stage, it can be normalized into a \textit{unit hypercube}, a DIU set. Next, we generalize this result to regular two-stage RO with this type of DDU sets.

\begin{prop}
\label{prop_reform_DIU2}
Let $\mathbf x_l$ and $\mathbf x_h$ denote the subvectors consisting of the first stage variables involved in \eqref{eq_US_hypercube}, i.e., $ x_{li}$ and $ x_{hi}$, $i=1,\dots, n_u$, respectively.
Then two-stage RO with the DDU defined \eqref{eq_US_hypercube} is equivalent to
$$\min_{\mathbf x\in \mathcal{X}}\mathbf c_1\mathbf x+\max_{\mathbf u\in \mathcal{U}^0}\ \min_{\mathbf y\in\mathcal Y'(\mathbf x,\mathbf u)}\mathbf c_2\mathbf y$$
where $\mathcal{U}^0=[0,1]^{n_u}$ and
\[
\mathcal Y'(\mathbf x,\mathbf u)=\{\mathbf y\in\mathbb Z^{m_y}_+\times \mathbb{R}^{n_y}_+:\mathbf B_2\mathbf y\geq\mathbf d-\mathbf B_1\mathbf x-\mathbf E(\mathbf x_l+\mathbf u\circ(\mathbf{x}_ h-\mathbf{x}_l))\}. \pushQED{\qed}\qedhere
\]
\end{prop}
\begin{rem}
$(i)$ Note again that $\mathbf x_l$, $\mathbf x_h$ and $\mathbf u$ are given as parameters for the recourse problem. Hence, it is not necessary to perform linearization from the computational point of view. Moreover, if the recourse problem is an LP, $\mathcal{U}^0$ can be set to $\mathcal{U}^0=\{0,1\}^{n_u}$. \\
$(ii)$ Since the whole reformulation is based on the normalization of the hypercube set, we refer to this technique as \textit{normalization reformulation}. It actually is applicable to other simple sets, e.g., ball or ellipsoidal sets, where the first stage decisions only affect their centers and sizes. Also, it can be employed together with the neutralization one to handle DDU sets that have structures of both \eqref{eq_US_MIP} and \eqref{eq_US_hypercube}.
\end{rem}

\subsubsection{Order Switching Reformulation}
Sometimes the impact of the random factor is to change coefficients of recourse variables in the objective function. Assume that both variables in the uncertainty set and in the recourse problem are continuous, and the recourse problem in \eqref{eq_2RO} is in the form of
\begin{eqnarray}
\label{eq_switching_recourse}
\begin{split}
\min & \quad \mathbf{\hat c}_2(\mathbf u)\mathbf y \\
   \mbox{s.t.} &\quad \mathbf y\in \mathcal Y(\mathbf x) = \{\mathbf{B}_2\mathbf y\geq \mathbf d-\mathbf B_1\mathbf x, \mathbf y\in \mathbb{R}^{n_y}_+\}
\end{split}
\end{eqnarray}
where $\mathbf{\hat c}_2(\mathbf u)= (\hat {\mathbf E}\mathbf u)^\intercal+ \mathbf c_2$ returns a row vector with a dimension matching with that of $\mathbf y$. Clearly, feasible set $\mathcal Y(\mathbf x)$ is independent of the random factor. Two-stage RO with such a structure, which is not in the standard form of \eqref{eq_2RO}, can be significantly simplified by the minimax theorem \citep{sion1958general}. We first present a result regarding switching the order between the uncertainty set and the recourse problem. Note that assumptions (\textit{A1-A3}) still hold, and (\textit{A3}) indicates that
$\exists \mathbf{x}\in \mathcal{X}$ such that $\mathcal Y(\mathbf x)\neq \emptyset .$

%

\begin{lem}
\label{lem_switching}
For $\mathbf x\in \mathcal{X}$ such that $\mathcal Y(\mathbf x)$ is non-empty, we have
$$\max_{\mathbf u \in \mathcal U(\mathbf x)} \min_{\mathbf y\in \mathcal Y(\mathbf x)} \mathbf{ \hat c}_2(\mathbf u) \mathbf y = \min_{\mathbf y \in \mathcal Y(\mathbf x)} \max_{\mathbf u\in \mathcal U(\mathbf x)} \mathbf {\hat c}_2(\mathbf u) \mathbf y. $$
\end{lem}
\begin{proof}
See its proof in Appendix \ref{apd:proofS2}.
\end{proof}

Consequently, two-stage RO with this type of recourse problem can be greatly simplified.
\begin{prop}
\label{prop_RO_switching}
  Two-stage RO with the DDU set defined in \eqref{eq_uncer_set} and the recourse problem defined in \eqref{eq_switching_recourse} is equivalent to
\begin{eqnarray}
\label{eq_equiv_switching}
\min_{\mathbf x\in \mathcal{X}, \mathbf y\in \mathcal Y(\mathbf x)} \mathbf c_1\mathbf x + \mathbf c_2\mathbf y+\max_{\mathbf u\in \mathcal{U}(\mathbf x)} (\hat{\mathbf{E}}\mathbf u)^\intercal\mathbf y,
\end{eqnarray}
which is further equivalent to the following single-level optimization problem.
\begin{eqnarray*}
\min && \mathbf c_1\mathbf x+\mathbf c_2\mathbf y+ \mathbf{(h+Gx)}^\intercal\boldsymbol\lambda\\
\mathrm{s.t.} && \mathbf x\in \mathcal{X}, \mathbf y\in \mathcal Y(\mathbf x)\\
&& \mathbf F(\mathbf x)^\intercal\boldsymbol\lambda  - \mathbf{\hat E}^\intercal\mathbf y \geq \mathbf 0\\
&& \boldsymbol\lambda \geq \mathbf{0}
\end{eqnarray*}
\end{prop}
\begin{proof}
	See its proof in Appendix \ref{apd:proofS2}.
\end{proof}
\begin{rem}
$(i)$ Formulation \eqref{eq_equiv_switching} is actually a single-stage RO with DDU set $\mathcal{U}(\mathbf x)$ that appears in the objective function only.  The single-level equivalence is a bilinear  mixed integer program. When the size is manageable, it now can be readily solved by some professional solver. When $\mathbf x$ variables are binary or can be represented by their binary expansions, that bilinear program can be linearized into a mixed integer linear program that could be handled by professional solvers with a greater scalability.\\
$(ii)$ Since both the bi- and single-level equivalences are based on switching the order between the uncertainty set and the recourse problem, we refer to this technique as order-switching reformulation. We mention that when either $\mathcal U(\mathbf x)$ or $\mathcal Y(\mathbf x)$ is a mixed integer set, order switching does not ensure the equality in Lemma \ref{lem_switching}. Indeed, the min-max formulation in \eqref{eq_equiv_switching} only yields an upper bound to the optimal value of the original two-stage RO, according to the max-min inequality \citep{calafiore2014optimization}.
\end{rem}

As mentioned, the presented reformulation techniques (and the associated solution techniques) are restricted to  those specially structured DDU-based two-stage RO. In the next section, we focus on algorithm development to directly compute more general DDU-based RO problems, without depending on any  ``DDU to DIU'' reformulations.

\section{Algorithms for Regular Two-stage RO with DDU}
\label{sect_algorithms}
In this section, we study regular $\mathbf{2-Stg \ RO}$ with continuous DDU set and recourse problem, i.e., with $m_u=m_y=0$. For this type of problems, there is no efficient and general algorithm appeared in the literature yet. In light of the philosophy behind  basic C\&CG and existing variants \citep{zeng2013solving,zhao2011exact,zeng2014solving}, we develop three new variants that are able to derive exact solutions for two-stage RO with DDU sets.
It is worth mentioning that many results presented in this section not only generalize previous ones for two-stage RO with DIU, but also shed critical lights on understanding the fundamental structures of DDU-based two-stage RO.

\subsection{The First Variant: Benders $+$ C\&CG}
We first present a  single-level equivalent reformulation for $\mathbf{2-Stg \ RO}$ in \eqref{eq_2RO} that supports our algorithm development. Then, we describe the complete operations of this variant, along with definitions of sub- and master problems. Finally, we analyze the algorithm's convergence issue and complexity regarding the number of iterations upon termination.

\subsubsection{A Natural Single-Level Reformulation by Enumeration}
As shown in the development of Benders decomposition and C\&CG methods, it is necessary to have an equivalent reformulation, which often is of a simple but large-scale structure, to establish the mathematical foundation for decomposition algorithms. To achieve this, we first consider the following bilevel reformulation of the original tri-level model in \eqref{eq_2RO}, which can be simply obtained by dualizing the recourse problem as shown before Proposition \ref{prop_extremepoint}.
\begin{eqnarray}
\label{eq_2RO_NC_bilevel}
\min_{\mathbf x\in \mathcal{X}} \mathbf{c}_1\mathbf x+ \max\{(\mathbf d-\mathbf B_1\mathbf x-\mathbf E\mathbf u)^\intercal\boldsymbol\pi: \mathbf{u}\in
    \mathcal{U}(\mathbf x), \mathbf B_2^\intercal\boldsymbol\pi \leq \mathbf c_2^\intercal, \boldsymbol\pi\geq \mathbf{0}\}.
\end{eqnarray}
Note that this bilevel model has a lower-level problem that is a complex disjoint bilinear program. By Proposition \ref{prop_extremepoint}, Corollary \ref{cor_maxmax_dual} and their remarks, we can derive, through enumeration, a much simpler but large-scale linear bilevel reformulation. Recall that $\Pi$ represents the polyhedron of the dual problem for the recourse problem.

\begin{thm}
\label{thm_PI_enumeration}
Let $\Sigma_{\Pi}=(\mathcal{P}_{\Pi}, \mathcal{R}_{\Pi})$ with $\mathcal{P}_{\Pi}$ and $\mathcal{R}_{\Pi}$ being the sets of extreme points and extreme rays of $\Pi$, respectively. Then, $\mathbf{2-Stg \ RO}$ in \eqref{eq_2RO} (and  the equivalence in \eqref{eq_2RO_NC_bilevel}) is equivalent to a bilevel linear optimization program as in the following.
\begin{subequations}
\label{eq_2stgRO_PI}
\begin{align}
\mathbf{2-Stg \ RO(\Sigma_{\Pi})}: \   w^*=\min \ & \ \mathbf{c}_1\mathbf x+ \eta\\
\mathrm{s.t.} \ & \ \mathbf x\in \mathcal{X}\\
            &  \bigg\{\eta\geq  (\mathbf d-\mathbf B_1\mathbf x)^\intercal\boldsymbol\pi + \max_{\mathbf u\in \mathcal{U}(\mathbf x)} \{(-\mathbf E\mathbf u)^\intercal\boldsymbol\pi\}\bigg\} \forall \boldsymbol\pi\in \mathcal{P}_{\Pi} \label{eq_enu_optimality}\\
            &  \bigg\{(\mathbf d-\mathbf B_1\mathbf x)^\intercal\boldsymbol\gamma 
            +\max_{\mathbf v\in \mathcal{U}(\mathbf x)}(-\mathbf E\mathbf v)^\intercal\boldsymbol\gamma \leq 0\bigg\} \forall \boldsymbol\gamma\in \mathcal{R}_{\Pi} \label{eq_enu_feasibility}
\end{align}
\end{subequations}
\end{thm}
 \begin{proof}
	See its proof in Appendix \ref{apd:proofS3}.
\end{proof}

Note in \eqref{eq_enu_feasibility} that variable $\mathbf v$ is an alias of $\mathbf u$. 
Obviously, this bilevel optimization model has many lower-level LPs in \eqref{eq_enu_optimality} and \eqref{eq_enu_feasibility}. By directly making use of LP's KKT or primal-dual optimality conditions, $\mathbf{2-Stg \ RO(\Sigma_{\Pi})}$ can be further converted into a single-level optimization problem. In the following, we present such one based on  KKT conditions. 

Let $\mathcal{P}_{\Pi}$ be $\{\boldsymbol\pi^1,\dots, \boldsymbol\pi^{K_p}\}$ with $K_p=|\mathcal{P}_{\Pi}|$, and  $\mathcal{R}_{\Pi}$ be $\{\boldsymbol\gamma^1,\dots, \boldsymbol\gamma^{K_r}\}$ with $K_r=|\mathcal{R}_{\Pi}|$, respectively.  The LPs in \eqref{eq_enu_optimality} and \eqref{eq_enu_feasibility} can be generalized into the following one.
\begin{eqnarray}
	\label{eq_LP_parametric}
	\mathbf{LP}(\mathbf x, \boldsymbol\beta): \ \ \max\{
	(-\mathbf{Eu})^\intercal\boldsymbol\beta: \ \mathbf u\in \mathcal{U}(\mathbf x)\}.
\end{eqnarray}
When $\boldsymbol\beta=\boldsymbol\pi^k$, we use $\mathcal{OU}(\mathbf x,\boldsymbol\pi^k)$ to denote the optimal solution set of $\mathbf{LP}(\mathbf x, \boldsymbol\pi^k)$, which can be defined by its KKT conditions, i.e.,
\begin{eqnarray}
\label{eq_KKT_point}
\mathcal{OU}(\mathbf x, \boldsymbol\pi^k) =  \left\{\begin{array}{l}
			 \mathbf F(\mathbf x)\mathbf u^k\leq \mathbf h+\mathbf{Gx}\\
             \mathbf F(\mathbf x)^\intercal\boldsymbol\lambda^k\geq - \mathbf E^\intercal\boldsymbol\pi^k \\
             \boldsymbol\lambda^k \circ (\mathbf h+\mathbf{Gx}-\mathbf F(\mathbf x)\mathbf u^k)=\mathbf{0}\\
             \mathbf u^k \circ (\mathbf F(\mathbf x)^\intercal\boldsymbol\lambda^k+ \mathbf E^\intercal\boldsymbol\pi^k) = \mathbf 0\\
             \mathbf u^k\geq \mathbf{0}, \boldsymbol\lambda^k\geq \mathbf{0},
            \end{array}\right\}
\end{eqnarray}
where $\boldsymbol\lambda^k$ denotes dual variables of constraints in $\mathcal U(\mathbf x)$. The first and the second constraints are constraints of primal and dual problems, respectively, and the third and the forth ones are complementarity constraints. Similarly, when $\boldsymbol\beta=\boldsymbol\gamma^l$, let $\mathcal{OV}(\mathbf x,\boldsymbol\gamma^l)$ denote the optimal solution set of $\mathbf{LP}(\mathbf x, \boldsymbol\gamma^l)$. It can be simply obtained through replacing $\mathbf u^k$ by $\mathbf v^l$ and $\boldsymbol\lambda^k$ by $\boldsymbol\zeta^l$ in \eqref{eq_KKT_point}, with $\boldsymbol\zeta^l$ denoting the corresponding dual variables. We note that another primal-dual optimality condition based representation for $\mathcal{OU}(\mathbf x,\boldsymbol\pi^k)$ and $\mathcal{OV}(\mathbf x,\boldsymbol\gamma^l)$ can be seen in \eqref{eq_dual_point} in  Appendix \ref{apd:proofS3}, which maybe computationally more friendly if $\mathbf x$ is binary.

By enumerating those KKT conditions based sets, the next result simply follows.

\begin{cor}
\label{cor_KKT_BendersRe}
Formulation $\mathbf{2-Stg \ RO(\Sigma_{\Pi})}$ is equivalent to the following single-level optimization problem.
\begin{align}
	\mathbf{2-Stg \ RO(\Sigma^K_{\Pi})}: \ \min \quad & \mathbf{c}_1\mathbf x+ \eta\nonumber\\		
	\begin{split}
		\mathrm{s.t.} \quad & \mathbf x\in \mathcal{X}\\
		& \eta\geq (\boldsymbol\pi^k)^\intercal\mathbf d-(\boldsymbol\pi^k)^\intercal\mathbf B_1\mathbf x - (\boldsymbol\pi^k)^\intercal\mathbf E\mathbf u^k, \ k=1,\dots,K_p\\
		& (\mathbf u^k, \boldsymbol\lambda^k) \in \mathcal{OU}(\mathbf x, \boldsymbol\pi^k), \ k=1,\dots, K_p\\
		& (\boldsymbol\gamma^l)^\intercal\mathbf d-(\boldsymbol\gamma^l)^\intercal\mathbf B_1\mathbf x - (\boldsymbol\gamma^l)^\intercal\mathbf E\mathbf v^l \leq 0, \  l=1,\dots, K_r
	\end{split} \label{eqn_2stgRO_Pi_KKT}\\
	& (\mathbf v^l, \boldsymbol\zeta^l)\in \mathcal{OV}(\mathbf x, \boldsymbol\gamma^l), \ l=1,\dots, K_r \pushQED{\qed}\qedhere \nonumber
\end{align} 

\end{cor}

Certainly, it would be desirable if we can directly compute the single-level formulation
$\mathbf{2-Stg \ RO(\Sigma^K_{\Pi})}$. Nevertheless, similar to any enumeration based reformulations, the numbers of variables and constraints are enormous. A subset of them should provide a relaxation and therefore yield a lower bound to $\mathbf{2-Stg \ RO(\Sigma^K_{\Pi})}$.

\begin{cor}
	\label{cor_Pi_partial_KKT}
	Let $\Sigma_{\hat {\Pi}}=(\mathcal{\hat P}_{\Pi}, \mathcal{\hat R}_{\Pi})$ with $\mathcal{\hat P}_{\Pi}\subseteq \mathcal{P}_{\Pi}$ and $\mathcal{\hat R}_{\Pi}\subseteq \mathcal{R}_{\Pi}$, and denote the formulation~\eqref{eqn_2stgRO_Pi_KKT} defined with respect to $\Sigma_{\hat\Pi}$ by $\mathbf{2-Stg \ RO(\Sigma^K_{\hat\Pi})}$. Then,  $\mathbf{2-Stg \ RO(\Sigma^K_{\hat\Pi})}$ is a relaxation to $\mathbf{2-Stg \ RO(\Sigma^K_{\Pi})}$ (and $\mathbf{2-Stg \ RO}$  in \eqref{eq_2RO}). And the optimal value of $\mathbf{2-Stg \ RO(\Sigma^K_{\hat\Pi})}$ is smaller than or equal to $w^*$.    \hfill $\square$
\end{cor}
It is clear that the strength of the aforementioned lower bound depends on the components in the subset. To achieve a better strength, one strategy could be  augmenting subset $\Sigma_{\hat {\Pi}}$ iteratively as needed, which is actually the basic idea of the algorithm development in this subsection. As $\max-\min$ (or $\min-\max$) bilevel optimization and complementarity constraints  (e.g., those in $\mathcal{OU}$) are heavily involved in  algorithm development in this paper, we discuss briefly an MIP-based solution method before presenting algorithm operations.

\subsubsection{Computing $\max-\min$ Bilevel Linear Optimization Problem}
 In this paper, all  $\max-\min$ bilevel linear optimization problems, unless otherwise stated, are in the form of  
 \begin{align}
 \label{eq_bilevel_general}
 	v^*=\max_{\mathbf x\in \mathcal X}  \min \{\mathbf c_{\mathbf y}\mathbf y:  \mathbf B_{\mathbf y}\mathbf y\geq \mathbf d-\mathbf{B}_{\mathbf x}\mathbf x, \mathbf{y}\geq \mathbf{0}\},	
 \end{align}
where $\mathcal{X}$ is a non-empty bounded mixed integer set.

 
For the special case where $v^*=-\infty$, i.e., the lower-level problem is unbounded for all $\mathbf x\in \mathcal X$, it can be easily detected by checking whether  $\{\mathbf{B}^\intercal_{\mathbf y}\boldsymbol{\lambda}\leq \mathbf c^\intercal_{\mathbf y}, \boldsymbol{\lambda}\geq \mathbf 0\}$ is an empty set. Note that according to assumption (\textit{A3}), it does not occur to any $\max-\min$ problem in this paper. For another special case where $v^*=\infty$, there exists $\mathbf x\in \mathcal{X}$ to which the lower-level problem is infeasible. To check if this is the case, we construct an extended formulation of~\eqref{eq_bilevel_general} similar to that in \citet{lee2014robust} as in the following.
  \begin{align*}
  	v^*_f=\max_{\mathbf x\in \mathcal X} \min \{\mathbf{1}^\intercal\tilde{\mathbf y}: \mathbf{B}_{\mathbf y}+\tilde{\mathbf y}\geq \mathbf d-\mathbf B_{\mathbf x}\mathbf x, \mathbf y\geq \mathbf 0, \tilde{\mathbf y}\geq \mathbf 0 \}
  \end{align*}
  Clearly, with artificial variable $\tilde{\mathbf y}$, the extended formulation is feasible and bounded for all $\mathbf x\in \mathcal{X}$. The next result, which can be proven easily by contradiction, shows that its optimal value helps us perform that check.  \begin{prop}
  	\label{prop_generalBilevel_feasibility}
  	For the $\max-\min$ problem in \eqref{eq_bilevel_general}, the lower-level problem is feasible for all $\mathbf x\in \mathcal{X}$ if and only if the optimal value of its extended formulation is zero, i.e., $v^*_f=0$. Otherwise, i.e., $v^*_f>0$,  an optimal solution to its extended formulation renders the lower-level problem of \eqref{eq_bilevel_general} infeasible. \hfill$\square$
  \end{prop}

   In the case where \eqref{eq_bilevel_general} has a finite optimal value,  we can replace the lower-level problem by its optimality conditions. For example, if KKT conditions are employed, \eqref{eq_bilevel_general} can be converted into the following single-level formulation
   \begin{align}
   	\label{eq_bilevel_general_KKT}
   	\begin{split}
   	v^*=  \max\quad&\mathbf c_{\mathbf y}\mathbf y\\
   	 \mathrm{s.t.}\quad&\mathbf x\in \mathcal X, \ \ \ \mathbf B_{\mathbf y}\mathbf y\geq \mathbf d-\mathbf{B}_{\mathbf x}\mathbf x, \ \ \ \mathbf B_{\mathbf y}^\intercal\boldsymbol\lambda\leq \mathbf	c^\intercal_{\mathbf y}\\
   	 & \boldsymbol\lambda\circ(\mathbf B_{\mathbf y}\mathbf y-\mathbf d+\mathbf{B}_{\mathbf x}\mathbf x) = \mathbf 0, \ \ \ \mathbf y\circ (\mathbf	c^\intercal_{\mathbf y}- \mathbf B_{\mathbf y}^\intercal\boldsymbol\lambda)=\mathbf 0\\
   	 & \mathbf{y}\geq \mathbf{0}, \ \boldsymbol\lambda\geq \mathbf 0
   \end{split}
   \end{align}
   where $\boldsymbol\lambda$ denotes the dual variables of the lower-level problem.  Constraints in the second row are complementarity constraints that have appeared in  $\mathcal{OU}$ and $\mathcal{OV}$. Those constraints actually can be linearized, converting \eqref{eq_bilevel_general_KKT} into an MIP. Specifically, consider the $j$-th constraint in the  first set of complementarity constraints in \eqref{eq_bilevel_general_KKT}, i.e.,
   $$\lambda_j (\mathbf d-\mathbf{B}_{\mathbf x}\mathbf x-\mathbf{B}_{\mathbf y}\mathbf y)_j=0.$$
   Let $\delta_j$ be a binary variable and recall that $M$ denotes a sufficiently large number. Then, this complementarity constraint can be replaced by the next two linear constraints.
   \begin{equation*}
   	\lambda_j\leq M\delta_j,  \  (\mathbf{B}_{\mathbf y}\mathbf y - \mathbf d+\mathbf{B}_{\mathbf x}\mathbf x)_j\leq M(1-\delta_j)
   \end{equation*}
   Note that when $\delta_j=0$, together with the nonnegativity of $\lambda_j$, we have
   $\lambda_j=0$. When $\delta_j=1$, together with the primal constraint $(\mathbf B_{\mathbf y}\mathbf y)_j\geq (\mathbf d-\mathbf{B}_{\mathbf x}\mathbf x)_j$, we have $(\mathbf B_{\mathbf y}\mathbf y - \mathbf d-\mathbf{B}_{\mathbf x}\mathbf x)_j=0$. Hence, such binary variables and linear constraints help us achieve the same effect as those complementarity constraints.
   Applying this technique to every complementarity constraint, e.g., those in \eqref{eq_bilevel_general_KKT}, or in $\mathcal{OU}$ and  $\mathcal{OV}$ in \eqref{eqn_2stgRO_Pi_KKT}, a linear MIP formulation is obtained.

  \begin{rem}
  $(i)$
  Certainly we can apply the strong duality to convert the $\max-\min$ formulation in \eqref{eq_bilevel_general} into a single-level bilinear one, like the one in \eqref{eq_max_max}. As noted earlier, it is a model computable (with a rather restricted capacity) by some professional solvers. Indeed, this strategy works very well when every bilinear term is a product between a binary variable and a continuous variable, since it simply can be linearized. For other cases, the KKT based reformulation (in its  linearized form) often has a better computational performance. \\
  $(ii)$ As the aforementioned solution strategies are generally applicable, we do not describe detailed implementation to any specific $\max-\min$ or $\min-\max$ formulation presented in the remainder of this paper, and simply assume that it can be solved.
  \end{rem}

  Next, we describe a decomposition algorithm based on Corollaries \ref{cor_KKT_BendersRe} and \ref{cor_Pi_partial_KKT} to solve $\mathbf{2-Stg \ RO}$, which dynamically introduces new variables and constraints to represent optimality conditions and valid inequalities using the C\&CG framework.

\subsubsection{Dynamically Generating Optimality Conditions by C\&CG}
\label{subsubsec_Variant1}
The basic principle behind C\&CG, a master-subproblem algorithm framework, is to iteratively generate both new variables and constraints to strengthen a relaxation. The generated variables and constraints are mainly in the form of a replicate for the recourse problem \citep{zeng2013solving}, or a replicate for the dual problem and/or  KKT conditions for some substructure \citep{zhao2011exact,zeng2014solving}. According to the structure presented in $\mathbf{2-Stg \ RO(\Sigma^K_{\Pi})}$, it is rather clear that the majority of new variables and constraints appear in the form of KKT conditions presented in $\mathcal{OU}$ and $\mathcal{OV}$.

Before presenting the detailed algorithm, we first define three subproblems for a given $\mathbf x^*$. Our first subproblem is designed to detect the feasibility of $\mathbf x^*$.
\begin{eqnarray}
\label{eq_SP1}
\mathbf{SP1}: \ \eta_f (\mathbf x^*) = \max_{\mathbf u\in \mathcal U(\mathbf x^*)}  & \min
\{\mathbf{1}^\intercal\tilde{\mathbf y}: \mathbf{B}_2\mathbf y+\tilde{\mathbf y}\geq \mathbf d-\mathbf B_1\mathbf x^*-\mathbf E\mathbf u, \mathbf y\geq \mathbf 0, \tilde{\mathbf y}\geq
\mathbf 0 \}
\end{eqnarray}
 By definition, $\mathbf x^*$ is feasible if  the recourse problem is feasible for all scenarios in $\mathcal{U}(\mathbf x^*)$.  According to Proposition \ref{prop_generalBilevel_feasibility}, solving $\mathbf{SP1}$ provides a certification regarding $\mathbf{x}^*$'s feasibility. 
\begin{cor}
A given $\mathbf x^*$ is feasible to $\mathbf{2-Stg \ RO}$ in \eqref{eq_2RO} and its equivalences if and only if $\eta_f(\mathbf x^*)=0$. \hfill$\square$
\end{cor}

On the one hand, when $\eta_f(\mathbf x^*)=0$, i.e., the recourse problem is feasible for all $\mathbf u\in \mathcal{U}(\mathbf x^*)$, we compute the second subproblem to evaluate the worst case performance of $\mathbf x^*$ by deriving the worst case scenario $\mathbf u^*_s$ and its recourse cost.
\begin{eqnarray}
\label{eq_SP2}
	\mathbf{SP2}: \ \eta_s(\mathbf x^*)=\max_{\mathbf u\in \mathcal{U}(\mathbf x^*)} \min\{
\mathbf c_2\mathbf y: \mathbf y\in \mathcal{Y}(\mathbf x^*, \mathbf u)\}
\end{eqnarray}
As mentioned, $\mathbf{SP2}$ can be solved through replacing the minimization problem by its KKT conditions or its dual problem. For either computational method the optimal values for dual variables, denoted by $\boldsymbol\pi^*$, define an extreme point of $\Pi$. According to \eqref{eq_LP_parametric}, we have
\begin{eqnarray}
	\label{eq_maxmi-LP}
	\eta_s(\mathbf x^*) =(\mathbf d-\mathbf B_1\mathbf x^*)^\intercal\boldsymbol\pi^*+ \mathbf{LP}(\mathbf x^*, \boldsymbol\pi^*).
\end{eqnarray}

On the other hand,  when $\eta_f(\mathbf x^*)>0$, its optimal solution to \eqref{eq_SP1}, denoted by $\mathbf u^*_f$, clearly causes the recourse problem to be infeasible. If this is the case, we compute the third subproblem, which is the dual of the recourse problem for $\mathbf u^*_f$.
\begin{eqnarray}
\label{eq_SP3}
\mathbf{SP3}: \ \max\{(\mathbf d- \mathbf B_1 \mathbf x^*- \mathbf {E}\mathbf{u}^*_f)^\intercal\boldsymbol\pi: \boldsymbol\pi\in \Pi\}
\end{eqnarray}
Since $\mathbf{SP3}$ is actually unbounded with respect to $(\mathbf x^*, \mathbf u^*_f)$, solving it will identify an extreme ray of $\Pi$, denoted by $\boldsymbol\gamma^*$, along which the optimal value is infinite. We then set the corresponding worst case recourse cost  $\eta_s(\mathbf x^*)$ to $+\infty$ by convention.

With those subproblems defined, we next present the detailed operations of the first variant of C\&CG method for $\mathbf{2-Stg \ RO}$, including the master problem. Note that the master problem, according to Corollary \ref{cor_Pi_partial_KKT}, is a relaxation of \eqref{eq_2RO} and hence provides a lower bound (denoted by $\textit{LB}$). Also, a fixed $\mathbf x$ and $\eta_s(\mathbf x)$ together provide an upper bound (denoted by $\textit{UB}$). Once the difference between those bounds is not more than $T
\!O\!L$, which sets the optimality tolerance, the algorithm terminates.   \\

\noindent\textbf{Variant 1: Benders C\&CG}
\begin{description}
	\item[Step 1] Set $LB = -\infty$, $UB=+\infty$, $t=1$,  and $\mathcal{\hat P}_{\Pi}=\mathcal{\hat R}_{\Pi}=\emptyset$.
	\item[Step 2] Solve the following master problem.
	\begin{eqnarray*}
		\label{eq_master_1}
		\mathbf{MP1}: \ \underline w=\min_{\mathbf{x}, \eta} && \mathbf{c}_1\mathbf{x}+ \eta \notag \\
		\mbox{s.t.} && \mathbf{x} \in \mathcal{X} \notag\\
		&& \eta\geq \boldsymbol\pi^\intercal\mathbf d- \boldsymbol\pi^\intercal\mathbf B_1\mathbf x-\boldsymbol\pi^\intercal\mathbf E\mathbf u^{\boldsymbol\pi} \ \ \forall \boldsymbol\pi\in \mathcal{\hat P}_{\Pi}\\
		&& (\mathbf u^{\boldsymbol\pi},\boldsymbol\lambda^{\boldsymbol\pi})\in \mathcal{OU}(\mathbf x, \boldsymbol\pi) \ \ \forall \boldsymbol\pi\in \mathcal{\hat P}_{\Pi}\\
		&& \boldsymbol\gamma^\intercal\mathbf d-\boldsymbol\gamma^\intercal\mathbf B_1\mathbf x - \boldsymbol\gamma^\intercal\mathbf E\mathbf v^{\boldsymbol\gamma} \leq 0 \  \ \forall \boldsymbol\gamma\in \mathcal{\hat R}_{\Pi}\\
		&& (\mathbf v^{\boldsymbol\gamma}, \boldsymbol\zeta^{\boldsymbol\gamma})\in \mathcal{OV}(\mathbf x, \boldsymbol\gamma) \ \forall \boldsymbol\gamma\in \mathcal{\hat R}_{\Pi}
	\end{eqnarray*}
	If it is infeasible, report infeasibility of $\mathbf{2-Stg \ RO}$ in \eqref{eq_2RO} and terminate. Otherwise, derive its optimal solution $\mathbf x^*, \eta, (\mathbf u^{\boldsymbol\pi}, \boldsymbol\lambda^{\boldsymbol\pi}) \ \forall \boldsymbol\pi \in \mathcal{\hat P}_{\Pi}, (\mathbf v^{\boldsymbol\gamma}, \boldsymbol\zeta^{\boldsymbol\gamma}) \ \forall \boldsymbol\gamma\in \mathcal{\hat R}_{\Pi}$, and its optimal value $\underline w$. Update
         $LB= \underline w$.
	\item[Step 3] Solve subproblem  $\mathbf{SP1}$ in \eqref{eq_SP1} and derive optimal 	$\mathbf u^*_f$ and $\eta_f(\mathbf x^*)$.
	\item[Step 4] Cases based on $\eta_f(\mathbf x^*)$
    \begin{description}
		 \item[(Case A): $\eta_f(\mathbf x^*)=0$] \textrm{}\\ $(i)$ compute $\mathbf{SP2}$ in \eqref{eq_SP2} to derive $\eta_s(\mathbf x^*)$, $\mathbf u^*_s$ and corresponding extreme point $\boldsymbol\pi^*$ of $\Pi$; $(ii)$ update $\mathcal{\hat P}_{\Pi}= \mathcal{\hat P}_{\Pi}\cup\{\boldsymbol\pi^*\}$ and augment master problem $\mathbf{MP1}$ accordingly, i.e., create variables $\mathbf u^{\boldsymbol\pi^*}$ and $\boldsymbol\lambda^{\boldsymbol\pi^*}$, and add the following constraints to $\mathbf{MP1}$.
		 \begin{subequations}
		  \label{eq_BD_optimality}
		 \begin{align}
          & \eta \geq (\boldsymbol\pi^*)^\intercal\mathbf d-(\boldsymbol\pi^*)^\intercal\mathbf B_1\mathbf x - (\boldsymbol\pi^*)^\intercal\mathbf E\mathbf u^{\boldsymbol\pi^*}  \label{eq_optimality_BD1}\\
          & (\mathbf u^{\boldsymbol\pi^*}, \boldsymbol\lambda^{\boldsymbol\pi^*}) \in \mathcal{OU}(\mathbf x, \boldsymbol\pi^*) \label{eq_optimality_BD2}
          \end{align}
		 \end{subequations}
		\item[(Case B): $\eta_f(\mathbf x^*) > $ 0] \textrm{}\\
		$(i)$ compute  $\mathbf{SP3}$ in \eqref{eq_SP3} to derive an extreme ray $\boldsymbol\gamma^*$ of $\Pi$, and set $\eta_s(\mathbf x^*)=+\infty$; $(ii)$ update $\mathcal{\hat R}_{\Pi}= \mathcal{\hat R}_{\Pi}\cup\{\boldsymbol\gamma^*\}$ and augment master problem $\mathbf{MP1}$ accordingly, i.e., create variables $\mathbf v^{\boldsymbol\gamma^*}$ and $\boldsymbol\zeta^{\boldsymbol\gamma^*}$, and add the following constraints  to $\mathbf{MP1}$.
		 \begin{subequations}
		 \label{eq_BD_feasibility}
	     \begin{align}
          &   (\boldsymbol\gamma^*)^\intercal\mathbf d-(\boldsymbol\gamma^*)^\intercal\mathbf B_1\mathbf x - (\boldsymbol\gamma^*)^\intercal\mathbf E\mathbf v^{\boldsymbol\gamma^*} \leq 0 \label{eq_feasibility_BD1}\\
      &  (\mathbf v^{\boldsymbol\gamma^*}, \boldsymbol\zeta^{\boldsymbol\gamma^*})\in \mathcal{OV}(\mathbf x, \boldsymbol\gamma^*) \label{eq_feasibility_BD2}
      	\end{align}
	    \end{subequations}
        \end{description}
       \item[Step 5] Update $UB=\min\{UB, \mathbf{c}_1\mathbf{x}^* +\eta_s(\mathbf x^*)\}.$
	\item[Step 6] If $UB-LB\leq T\!O\!L$, return $\mathbf{x}^*$ and
	terminate. Otherwise, set $t=t+1$ and go to \textbf{Step~2}.~\hfill $\square$
\end{description}

\begin{rem}
$(i)$ Note that for every iteration, new variables and constraints are generated either in the form of \eqref{eq_BD_optimality} or \eqref{eq_BD_feasibility}. On the one hand, constraints in \eqref{eq_optimality_BD1} and \eqref{eq_feasibility_BD1}  are introduced that are closely analogous to the well-known cutting planes in Benders decomposition. On the other hand, similar to \cite{zeng2014solving}, new KKT conditions, including both variables and constraints, are dynamically generated to define sets $\mathcal{OU}$ and $\mathcal{OV}$ in \eqref{eq_optimality_BD2} and \eqref{eq_feasibility_BD2}, which are then used to populate \eqref{eq_optimality_BD1} and \eqref{eq_feasibility_BD1}. Since \eqref{eq_BD_optimality} or \eqref{eq_BD_feasibility} is introduced when $\mathbf x^*$ has a finite worst case performance or is just infeasible, respectively, \eqref{eq_BD_optimality} is referred to as the \textit{``optimality cutting set''}, and \eqref{eq_BD_feasibility} the \textit{``feasibility cutting set''}, respectively. As this variant demonstrates an integration of Benders and C\&CG, it is also referred to as ``Benders C\&CG'' method. \\
$(ii)$ Regarding KKT conditions based sets $\mathcal{OU}$ and $\mathcal{OV}$,  only $\mathbf u^{\pi^*}$ or $\mathbf v^{\gamma^*}$ are adopted in defining constraint \eqref{eq_optimality_BD1} or \eqref{eq_feasibility_BD1}.
Let $\mathcal{OU}_{\mathbf u}(\mathbf x,\boldsymbol\pi^*)=\{\mathbf u: (\mathbf u,\boldsymbol\lambda)\in \mathcal{OU}(\mathbf x,\boldsymbol\pi^*) \ \textrm{for some} \ \boldsymbol\lambda\}$ and $\mathcal{OV}_{\mathbf u}(\mathbf x,\boldsymbol\gamma^*)=\{\mathbf v: (\mathbf v, \boldsymbol\zeta)\in \mathcal{OV}(\mathbf x,\boldsymbol\gamma^*) \ \textrm{for some} \ \boldsymbol\zeta\}$, i.e.,  the projections of $\mathcal{OU}(\mathbf x,\boldsymbol\pi^*)$ and $\mathcal{OV}(\mathbf x,\boldsymbol\gamma^*)$ onto the subspace hosting $\mathbf u$, respectively. Equivalently, \eqref{eq_optimality_BD2} and \eqref{eq_feasibility_BD2} can be re-stated as $$\mathbf u^{\boldsymbol\pi^*}\in \mathcal{OU}_{\mathbf u}(\mathbf x,\boldsymbol\pi^*) \ \ \textrm{and} \ \  \mathbf v^{\boldsymbol\gamma^*}\in \mathcal{OV}_{\mathbf u}(\mathbf x,\boldsymbol \gamma^*).$$ We highlight that this projection based interpretation is important and it indicates the reduction  of Variant 1 to the classical Benders-dual cutting plane method. \\
$(iii)$ Note that  all operations and analyses regarding $\mathcal{OU}$ and $\mathcal{OV}$ presented in this paper still hold if we replace KKT conditions based representations by their primal-dual optimality condition based equivalences described in \eqref{eq_dual_point} in  Appendix~\ref{apd:proofS3}.    
\end{rem}

The classical Benders-dual cutting plane method is developed for two-stage RO with DIU set $\mathcal U^0$ \citep{terry2009robust,thiele2009robust,bertsimas2012adaptive,jiang2014two,zeng2013solving}. When $\mathcal U(\mathbf x)=\mathcal U^0$ for $\mathbf x\in \mathcal X$ and consider a given extreme point $\boldsymbol\pi$ or extreme ray $\boldsymbol\gamma$ of $\Pi$, $\mathcal{OU}_{\mathbf u}(\cdot,\boldsymbol\pi)$ or $\mathcal{OV}_{\mathbf u}(\cdot, \boldsymbol\gamma)$ reduces to a fixed set of points of $\mathcal {U}^0$ that are optimal with respect to $\boldsymbol\pi$ or $\boldsymbol\gamma$. Hence, we can simply fix variable $\mathbf u^{\boldsymbol\pi}$ or $\mathbf v^{\boldsymbol\gamma}$ to any point of $\mathcal{OU}_{\mathbf u}(\cdot,\boldsymbol\pi)$ or $\mathcal{OV}_{\mathbf u}(\cdot, \boldsymbol\gamma)$, and completely eliminate \eqref{eq_optimality_BD2} and \eqref{eq_feasibility_BD2} from the algorithm. Consequently,
an optimality (or feasibility, respectively) cutting set  reduces to an optimality (or feasibility, respectively) cut developed in the traditional Benders-dual cutting plane method.
\begin{prop}
	\label{prop_BD_reduction}
	When uncertainty set $\mathcal U(\mathbf x)=\mathcal U^0$ for all $\mathbf x\in \mathcal X$, Variant 1 reduces to Benders-dual cutting plane method. \hfill $\square$
\end{prop}

\subsubsection{Analysis of Convergence and Complexity}
In order to theoretically understand the behavior of Variant 1, we investigate  its convergence and the iteration complexity issues. 

\begin{thm}
\label{thm_V1_convergence}
Suppose $T\!O\!L=0$. When terminates, Variant 1 either reports that $\mathbf{2-Stg \ RO}$ in \eqref{eq_2RO} is infeasible,  or converges to its optimal value and exact solution. 
\end{thm}
\begin{proof}
		See its proof in Section \ref{apd:proofS3}.
\end{proof}

According to the proof of Theorem \ref{thm_V1_convergence}, we can easily bound the number of iterations of the algorithm.
\begin{cor}
	\label{cor_BDCCG_complexity}
	 The number of iterations of Variant 1 before termination is bounded by  $(|\mathcal{P}_{\Pi}|+|\mathcal{R}_{\Pi}|)$, i.e., $K_{p}+K_{r}$. Hence, the algorithm is of $O(\binom{n_y+\mu_y}{\mu_y})$ iteration complexity where  $\mu_y$ denotes the number of rows of matrix $\mathbf B_2$ in $\mathcal Y(\mathbf x, \mathbf u)$. \hfill$\square$
\end{cor}
\begin{rem}
 We mention that this computational complexity result only depends on $\Pi$, completely independent of $\mathcal{U}(\mathbf x)$. For the classical Benders-dual cutting plane method developed for two-stage RO with DIU set $\mathcal U^0$, it is shown in \citet{zeng2013solving} that its iteration complexity is $O((|\mathcal{P}_{\Pi}|+|\mathcal{R}_{\Pi}|)*|\mathcal{P}_{\mathcal U^0}|)$, where $\mathcal{P}_{\mathcal U^0}$ denotes the set of extreme points of $\mathcal{U}^0$. As Variant 1 reduces to  Benders-dual cutting plane method when the DDU set reduces to $\mathcal{U}^0$, Corollary \ref{cor_BDCCG_complexity} clearly  strengthens our previous understanding on the iteration complexity.
\end{rem}

When $\mathcal{X}$ is a finite discrete set, its carnality simply bounds the iteration complexity.
\begin{prop}
\label{prop_cardinality_repeat_BD}
If $\mathbf x^*$ is an optimal solution to $\mathbf{MP1}$ in both iterations $t_1$ and $t_2$ with $t_1<t_2$, we have $LB=UB$ in iteration $t_2$ and $\mathbf x^*$ is optimal to $\mathbf{2-Stg \ RO}$.  Also, if $\mathcal{X}$ is a finite discrete set, the number of iterations is bounded by the cardinality of $\mathcal{X}$, i.e., $|\mathcal{X}|$.  \hfill$\square$ 	
\end{prop}
\begin{rem}
	Note from proofs that those convergence and complexity results  hold regardless of RHS and LHS dependence. 
	Also, they are valid if the DDU set is a more complex  second order conic (SOC)  set. Hence, theoretically speaking, the whole algorithm can be implemented to compute related two-stage RO instances, as long as the master problem and subproblems can be exactly solved by the available oracle. 
\end{rem}


%

\subsection{The Second Variant: Parameterized C\&CG}
\label{subsect_Variant2}
Although the first variant is theoretically sound, we observe that  its performance is not satisfactory. Even for small-scale instances, it often fails to converge after many iterations.  As noted earlier, it is an extension of the less-effective Benders-dual cutting plane method. Given the great performance demonstrated by  basic C\&CG method in computing two-stage RO with DIU, we are motivated to develop more sophisticated variants to handle DDU-based models. We expect those variants should be strong if they reduce to basic C\&CG  when the involved uncertainty set is actually a DIU set. In particular, it would be very significant if new variants have similar iteration complexities when dealing with complex DDU sets as the basic one when dealing with  DIU sets. In this subsection and the following one, we develop and analyze two variants of C\&CG that have those desirable properties.

\subsubsection{A Less Straightforward Single-Level Reformulation}
Similar to the development of Variant 1, we first present a simple but large-scale reformulation to lay down the foundation of this new C\&CG variant. In the following, we derive a result regarding the $\max-\min$ substructure of $\mathbf{2-Stg \ RO}$ for a given $\mathbf x$, through making use of projections of sets $\mathcal{OU}$ and $\mathcal{OV}$. Recall that $\boldsymbol\pi$ and $\boldsymbol\gamma$ are extreme point and extreme ray of $\Pi$, respectively.
\begin{lem}
\label{lem:parametric:maxmin}
For a given $\mathbf x$, we have
\begin{align*}
\max_{\mathbf u\in \mathcal{U}(\mathbf x)}\min\{\mathbf c_2\mathbf y: \mathbf y\in \mathcal{Y}(\mathbf x, \mathbf u)\}=\max_{\mathbf u\in \ \mathcal{U^*}(\mathbf x)\bigcup \mathcal{V^*}(\mathbf x)}  \min\{\mathbf c_2\mathbf y: \mathbf y\in \mathcal{Y}(\mathbf x, \mathbf u)\} 	
\end{align*}
where $\mathcal{U^*}(\mathbf x)=\bigcup \limits_{k=1}^{K_p} \mathcal{OU}_{\mathbf u}(\mathbf x, \boldsymbol\pi^k)$ and $\mathcal{V^*}(\mathbf x)=\bigcup\limits_{l=1}^{K_r} \mathcal{OV}_{\mathbf u}(\mathbf x, \boldsymbol\gamma^l)$.
\end{lem}
\begin{proof}
	See its proof in Section \ref{apd:proofS3}.
\end{proof}

\begin{thm}
	Formulation $\mathbf{2-Stg \ RO}$ in \eqref{eq_2RO} (and its other  equivalences) is equivalent to a bilevel linear optimization program as in the following.
\begin{subequations}
	\label{eq_2stgRO_PI_OM}
	\begin{align}
		\mathbf{2\!-\!Stg \ RO(\Omega_{\Pi})}: \   w^*=\min \ & \ \mathbf{c}_1\mathbf x+ \eta   \\
		\mathrm{s.t.} \ & \ \mathbf x\in \mathcal{X}\\
		&  \bigg\{\eta\geq  \mathbf{c_2}\mathbf y, \ \mathbf y\in \mathcal Y(\mathbf x, \mathbf u)\notag \\
		&  \  \ \mathbf{u}\in \arg\max\{(-\mathbf E\mathbf u)^\intercal\boldsymbol\pi:\mathbf u\in \mathcal{U}(\mathbf x)\}\bigg\}  \ \forall \boldsymbol\pi\in \mathcal{P}_{\Pi} \label{eq_enu2_optimality}\\
		& \bigg\{ \mathbf y\in \mathcal Y(\mathbf x, \mathbf v)\notag \\
		&  \ \mathbf{v}\in \arg\max\{(-\mathbf {Ev})^\intercal\boldsymbol\gamma: \mathbf v\in \mathcal{U}(\mathbf x)\}\bigg\}  \ \forall \boldsymbol\gamma\in \mathcal{R}_{\Pi},\label{eq_enu2_feasibility}
 	\end{align}
\end{subequations}
which can be further converted equivalently into the next single-level optimization problem.
\begin{align}
	\mathbf{2\!-\!Stg \ RO(\Omega^K_{\Pi})}:    w^*=\min \quad  &  \mathbf{c}_1\mathbf x+ \eta \nonumber\\
	\begin{split}
		\mathrm{s.t.} \quad  & \mathbf x\in \mathcal{X}\\
		&  \eta\geq  \mathbf{c_2}\mathbf y^k, \mathbf y^k\geq \mathbf 0, k=1,\dots,K_p \\
		&\mathbf{B}_2\mathbf y^k\geq \mathbf{d-\mathbf B}_1\mathbf x-\mathbf{Eu}^k, k=1,\dots,K_p \\
		& (\mathbf{u}^k,\boldsymbol\lambda^k)\in
		\mathcal{OU}(\mathbf x,\boldsymbol\pi^k),  k=1,\dots, K_p\\
		&\mathbf{B}_2\mathbf y^l\geq \mathbf{d-\mathbf B}_1\mathbf x-\mathbf{Ev}^l, \ \mathbf y^l\geq \mathbf 0,  l=1,\dots,K_r
	\end{split}\label{eq_2stgRO_PIOmega} \\
	&  (\mathbf{v}^l,\boldsymbol\zeta^l)\in \mathcal{OV}(\mathbf x, \boldsymbol\gamma^l), l=1,\dots, K_r \pushQED{\qed}\qedhere \nonumber
\end{align}
\end{thm}
Again, as a result of enumeration, $\mathbf{2\!-\!Stg \ RO(\Omega_{\Pi})}$ and $\mathbf{2\!-\!Stg \ RO(\Omega^K_{\Pi})}$ are large-scale bilevel and nonlinear optimization problems, respectively.

\begin{rem}
	$(i)$ It is worth noting a critical observation regarding the structure of \eqref{eq_2stgRO_PIOmega}. Although $\mathcal{U}(\mathbf x)$ is not fixed,  by their definitions, $\mathcal{OU}$ and $\mathcal{OV}$ actually yield a parametric approach to characterize non-trivial scenarios
	 in  $\mathcal{U}(\mathbf x)$. Those scenarios, which are optimal to some $\boldsymbol\pi$ or $\boldsymbol\gamma$ of $\Pi$, change with respect to $\mathbf x$. Then, by enumerating $\boldsymbol\pi$ and $\boldsymbol\gamma$ of $\Pi$, which  is fixed and independent of $\mathbf x$, \eqref{eq_2stgRO_PIOmega} can fully capture the impact of the DDU set by considering its parametrically represented scenarios. Nevertheless, given that sets $\mathcal P_{\Pi}$ and $\mathcal R_{\Pi}$ are large-scale and $\mathcal{OU}$ and $\mathcal{OV}$ are complex,  directly computing \eqref{eq_2stgRO_PIOmega} is unrealistic. \\
	 $(ii)$ Since a replicate of recourse problem is introduced and associated to every $\mathcal{OU}$ or $\mathcal{OV}$,
	 $\mathbf{2\!-\!Stg \ RO(\Omega_{\Pi})}$ is in a form similar to that of \eqref{eq_DIU_enumeration} of Corollary \ref{cor_DIU_Omega}.  One might believe that $\mathbf{2\!-\!Stg \ RO(\Omega_{\Pi})}$ reduces to \eqref{eq_DIU_enumeration} when $\mathcal U(\mathbf x)=\mathcal U^0$ for all $\mathbf x$. Actually, it is not the typical case even if $\mathcal{OU}$ and $\mathcal{OV}$ are singletons for every $\boldsymbol\pi$ and $\boldsymbol\gamma$. Clearly, if they are singletons, they are some extreme points of $\mathcal{U}^0$. Nevertheless, the map from $\boldsymbol\pi$ and $\boldsymbol\gamma$ to extreme points of  $\mathcal{U}^0$ is not surjective or injective. Only a subset of 	extreme points of $\mathcal{U}^0$ will be identified by enumerating $\boldsymbol\pi$ and $\boldsymbol\gamma$.
\end{rem}
Next, we present a result analogous to  Corollary \ref{cor_Pi_partial_KKT}.

\begin{cor}
	\label{cor_Pi_partial_KKTOmega}
	 Consider set $\Sigma_{\hat {\Pi}}$ defined in Corollary \ref{cor_Pi_partial_KKT} and denote formulation~\eqref{eq_2stgRO_PIOmega} defined with respect to $\Sigma_{\hat {\Pi}}$ by $\mathbf{2-Stg \ RO(\Omega^K_{\hat\Pi})}$. Then,  $\mathbf{2-Stg \ RO(\Omega^K_{\hat\Pi})}$ is a relaxation to $\mathbf{2-Stg \ RO(\Omega^K_{\Pi})}$ (and
	 $\mathbf{2-Stg \ RO}$ in \eqref{eq_2RO}). And the optimal value of $\mathbf{2-Stg \ RO(\Omega^K_{\hat\Pi})}$ is smaller than or equal to $w^*$.    \hfill $\square$
\end{cor}

\subsubsection{Variant 2: Generate Parametric Recourse Problems by C\&CG}
As mentioned earlier, when $\mathbf{2-Stg \ RO}$ is with a DIU set,   basic C\&CG generates and incorporates a replicate of recourse problem in an iterative fashion for every identified critical scenario \citep{zeng2013solving}. With $\mathbf{2\!-\!Stg \ RO(\Omega_{\Pi})}$ (and its single-level reformulation $\mathbf{2\!-\!Stg \ RO(\Omega^K_{\Pi})}$), such strategy now is feasible to be implemented in a parametric fashion for DDU-based $\mathbf{2-Stg \ RO}$. We refer to this new procedure as Variant 2. As the fundamental differences between Variants 1 and 2 are in the definitions of dynamically generated cutting sets and the resulting master problems,  we, to minimize repeated descriptions, just specify the  modifications on  top of Variant 1 in the following.  Unless noted explicitly, all concepts and notations have already been defined in the development of Variant 1.

\begin{description}
	\item[$\bullet$ Master Problem] Throughout the algorithm, $\mathbf{MP1}$ is replaced by the next one.
\begin{eqnarray*}
	\label{eq_master_1}
	\mathbf{MP2}: \ \underline w=\min && \mathbf{c}_1\mathbf{x}+ \eta \notag \\
	\mbox{s.t.} && \mathbf{x} \in \mathcal{X} \notag\\
	&& \eta\geq \mathbf{c_2}\mathbf y^{\boldsymbol\pi} \ \forall \boldsymbol\pi\in \mathcal{\hat P}_{\Pi}\\
	&& \mathbf{B}_2\mathbf y^{\boldsymbol\pi}\geq \mathbf{d-\mathbf B}_1\mathbf x-\mathbf{Eu}^{\boldsymbol\pi} \ \forall \boldsymbol\pi\in \mathcal{\hat P}_{\Pi}\\
	&& (\mathbf{u}^{\boldsymbol\pi},\boldsymbol\lambda^{\boldsymbol\pi})\in
	\mathcal{OU}(\mathbf x,\boldsymbol\pi),  \ \  \mathbf y^{\boldsymbol\pi} \geq \mathbf 0 \ \ \forall \boldsymbol\pi\in \mathcal{\hat P}_{\Pi}\\
	&&\mathbf{B}_2\mathbf y^{\boldsymbol\gamma}\geq \mathbf{d-\mathbf B}_1\mathbf x-\mathbf{Ev}^{\boldsymbol\gamma} \ \ \forall \boldsymbol\gamma\in \mathcal{\hat R}_{\Pi}\\
	&&  (\mathbf{v}^{\boldsymbol\gamma},\boldsymbol\zeta^{\boldsymbol\gamma})\in \mathcal{OV}(\mathbf x, \boldsymbol\gamma), \ \ \mathbf y^{\boldsymbol\gamma}\geq \mathbf 0 \ \ \forall \boldsymbol\gamma\in \mathcal{\hat R}_{\Pi}
\end{eqnarray*}
\item[$\bullet$ Cutting Sets]: For every $\mathbf u^{\boldsymbol\pi}$ or $\mathbf v^{\boldsymbol\gamma}$ defined by sets $\mathcal{OU}$ and $\mathcal{OV}$, we generate and incorporate a replicate of recourse problem (including recourse variables and constraints) into $\mathbf{MP2}$. Specifically, in \textbf{Step 4 [Case A]} where $\eta_f(\mathbf x^*)=0$,  the optimality cutting set in \eqref{eq_BD_optimality} is replaced by the following one.
\begin{equation}
	\label{eq_CCG_optimality}
	\begin{split}
		& \eta\geq \mathbf{c_2}\mathbf y^{\boldsymbol\pi^*}, \ \mathbf{B}_2\mathbf y^{\boldsymbol\pi^*}\geq \mathbf{d-\mathbf B}_1\mathbf x-\mathbf{Eu}^{\boldsymbol\pi^*} \\
		& (\mathbf{u}^{\boldsymbol\pi^*},\boldsymbol\lambda^{\boldsymbol\pi^*})\in
		\mathcal{OU}(\mathbf x,\boldsymbol\pi^*), \ \mathbf y^{\boldsymbol\pi^*} \geq \mathbf 0
	\end{split}
\end{equation}
\indent In \textbf{Step 4 [Case B]} where $\eta_f(\mathbf x^*)>0$,  the feasibility cutting set in \eqref{eq_BD_feasibility} is replaced by the following one.
\begin{equation}
	\label{eq_CCG_feasibility}
	\begin{split}
		& \mathbf{B}_2\mathbf y^{\boldsymbol\gamma^*}\geq \mathbf{d-\mathbf B}_1\mathbf x-\mathbf{Ev}^{\boldsymbol\gamma^*} \\
		& (\mathbf{v}^{\boldsymbol\gamma^*},\boldsymbol\lambda^{\boldsymbol\gamma^*})\in
		\mathcal{OV}(\mathbf x, \boldsymbol\gamma^*), \ \mathbf y^{\boldsymbol\gamma^*} \geq \mathbf 0
	\end{split}
\end{equation}
\end{description}

Different from basic C\&CG and Variant 1, Variant 2 dynamically generates both optimality conditions characterizing non-trivial scenarios in a parametric way and their associated recourse problems. Hence, we also refer to this variant as \textit{parametric C\&CG} method.

Similar to Proposition \ref{prop_BD_reduction}, the next result follows directly from the fact that $\mathcal{OU}$ and $\mathcal{OV}$ reduce to fixed sets when $\mathcal{U}(\mathbf x)$ is DIU. 
By using any points in those sets to populate cutting sets (hence eliminating $\mathcal{OU}$ and $\mathcal{OV}$),  \eqref{eq_CCG_optimality} and \eqref{eq_CCG_feasibility}  reduce to cutting sets of  basic C\&CG exactly. Indeed, we can treat Variant 2 as an algebraic generalization of basic C\&CG. 

\begin{prop}
	\label{prop_PCCG_reduction}
	When uncertainty set $\mathcal U(\mathbf x)=\mathcal U^0$ for all $\mathbf x\in \mathcal X$, Variant 2 reduces to basic C\&CG method presented in \citet{zeng2013solving}. \hfill $\square$
\end{prop}

Analogous to the comparison between basic C\&CG and Benders-dual methods for~$\mathbf{2-Stg}$  $\mathbf{RO}$ with DIU, Variant 2 generally produces stronger lower bounds in its executions.

\begin{prop}
	\label{prop_V2betterLB}
	Assume that the same sets $\mathcal{\hat P}_{\Pi}\subseteq \mathcal{P}_{\Pi}$ and $\mathcal{\hat R}_{\Pi}\subseteq \mathcal{R}_{\Pi}$ have been included for defining sets $\mathcal{OU}$ and $\mathcal{OV}$  in both $\mathbf{MP1}$ and $\mathbf{MP2}$. Then, the optimal value of $\mathbf{MP1}$ is an underestimation of that of $\mathbf{MP2}$.
\end{prop}
\begin{proof} See its proof in Appendix \ref{apd:proofS3}.
\end{proof}	

Next, we would like to highlight a few features of Variant 2 in the following.
\begin{rem}
\label{rem_advantageCCG}
	This new variant presents a few advantages that have practical significance.  \\
($i$) \textbf{Unified cutting sets:} It has been observed in \citet{zeng2013solving} that basic C\&CG can unify both optimality and  feasibility cutting sets into the same form. Variant 2 also has this property. Given that $\mathcal{OV}_{\mathbf u}(\mathbf x, \boldsymbol\gamma^*)\in \mathcal{U}(\mathbf x)$ for an extreme ray $\boldsymbol\gamma^*$ of $\Pi$, it is feasible to augment~\eqref{eq_CCG_feasibility} with  $\eta\geq \mathbf{c}_2\mathbf{y}^{\boldsymbol\gamma^*}$, which is then of the same form as \eqref{eq_CCG_optimality}.  Hence, optimality and feasibility cutting sets are unified into the same structure, and we can simply generate variables and constraints as those in  \eqref{eq_CCG_optimality} for every identified extreme point or ray of $\Pi$. Because it might lead to stronger lower bounds without extra computational overhead,  unless noted explicitly, Variant 2 adopts the unified cutting set as its default. \\
($ii$) \textbf{Fast computational performance:} Compared to Variant 1, Variant 2 demonstrates a superior capacity in computing complex instances. As shown in our numerical studies presented in Section \ref{sect_numerical},  Variant 2 is generally a few orders of magnitude faster than Variant~1. We believe that such a huge difference is related to and generalizes what we have observed between Benders-dual and   basic C\&CG methods, given that Variants 1 and 2 reduce to them respectively when the involved uncertainty set is DIU. \\
($iii$) \textbf{Lower bound insensitive to big-M:} As the algorithm proceeds,  $\mathbf{MP2}$ could have many complementarity constraints (in the case of KKT conditions based) or bilinear constraints (in the case of primal-dual based) in $\mathcal{OU}$ and $\mathcal{OV}$ representations. As mentioned earlier, those nonlinear constraints can often be linearized using binary variables and big-M technique. We usually have some concern on the choice of big-M coefficient, since if it is not large sufficiently we will have a local optimal solution that disqualifies the corresponding value of $\mathbf{MP2}$ as a lower bound. Nevertheless, this is not the case. Regardless of the value of  big-M parameter, we have $\mathcal{OU}_{\mathbf u}(\mathbf x, \boldsymbol\pi)\in \mathcal{U}(\mathbf x)$ and $\mathcal{OV}_{\mathbf u}(\mathbf x, \boldsymbol\gamma)\in \mathcal{U}(\mathbf x)$ for all $\boldsymbol\pi$ and $\boldsymbol\gamma$, which renders the optimal value of $\mathbf{MP2}$ a valid lower bound to $\mathbf{2-Stg \ RO}$ all the time. Hence, it can be seen that $\mathbf{MP2}$ is less sensitive to the choice of big-M. 
\\
($iv$) \textbf{A flexible platform to include non-trivial scenarios parametrically:} Note that $\boldsymbol\pi$ and $\boldsymbol\gamma$ just serve as  criteria to introduce some critical scenarios parametrically, through $\mathcal{OU}$ and $\mathcal{OV}$. They do not directly affect $\eta$ or the feasibility of $\mathbf{MP2}$. So, in addition to exactly being extreme points or extreme rays of $\Pi$, they can be derived by domain expertise or heuristics, or even take values out of $\Pi$, which offers a great flexibility for us to incorporate any useful knowledge.  As demonstrated in Section \ref{subsect_Uniqueness_Pareto}, this flexibility can be utilized to achieve a strong result on the iteration complexity or a fast computation.
\end{rem}

\subsubsection{Analysis of Convergence and Complexity}
\label{subsect_pCCG_complexity}

According to  roles of sets $\mathcal{P}_{\Pi}$ and $\mathcal{R}_{\Pi}$ and by adopting the same idea presented in the proof of Theorem \ref{thm_V1_convergence}, the convergence result of Variant 2 can be obtained easily.
\begin{thm}
\label{thm_V2_converge}
	Suppose $T\!O\!L=0$.  When terminates, Variant 2 either reports that $\mathbf{2-Stg \ RO}$ in \eqref{eq_2RO} is infeasible, or converges to its optimal value and exact solution. ~\hfill$\square$
\end{thm}

Similar to Corollary \ref{cor_BDCCG_complexity} and Proposition \ref{prop_cardinality_repeat_BD}, next results simply follow.
\begin{cor}
	\label{cor_V2_complexity}
  The number of iterations of Variant 2 before termination is bounded by  $(|\mathcal{P}_{\Pi}|+|\mathcal{R}_{\Pi}|)$, i.e., $K_{p}+K_{r}$. Hence, the algorithm is of $O(\binom{n_y+\mu_y}{\mu_y})$ iteration complexity.~\hfill$\square$
\end{cor}

\begin{prop}
\label{prop_cardinality_repeat_CCG}
	If $\mathbf x^*$ is an optimal solution to $\mathbf{MP2}$ in both iterations $t_1$ and $t_2$ with $t_1<t_2$, we have $LB= UB$ in iteration $t_2$ and $\mathbf x^*$ is optimal to $\mathbf{2-Stg \ RO}$. Hence, if $\mathcal{X}$ is a finite discrete set, the number of iterations of Variant 2 is bounded by $|\mathcal{X}|$. ~\hfill$\square$   	
\end{prop}

\begin{rem}
$(i)$ As those convergence and complexity results
are primarily derived by making use of $\Pi$'s fixed polyhedral structure, they generally hold for $\mathcal{U}(\mathbf x)$ with both RHS and LHS dependence, as well as with more complex convex structure. \\ 
$(ii)$ Given the strong theoretical and computational performance of   basic C\&CG in computing DIU-based RO, it is desirable to theoretically prove that Variant 2 generalizes it in the context of DDU.  Nevertheless, regardless of the fact that Variant 2 reduces to  basic C\&CG if the underlying uncertainty set becomes DIU, such reduction connection is not reflected in the aforementioned complexity result. One essential reason is that, unlike a fixed DIU set, $\mathcal{U}(\mathbf x)$ changes with respect to $\mathbf x$ and its structure becomes indeterminate.  Hence, the classical extreme point/ray based analyses are not applicable to~$\mathcal{U}(\mathbf x)$.
\end{rem}

Nevertheless, we would like to mention that a core LP concept, i.e., ``basis'', can  be utilized as a powerful tool with regard to an indeterminate structure. Specifically when the DDU set is with RHS dependence, it helps us obtain more general and stronger iteration complexity for Variant 2. Indeed, it ensures the convergence and bounds the complexity even if the recourse problem is beyond the linear program exhibited in \eqref{eq_2RO}.

Let $\mathfrak{B}$ denote a basis of $\mathbf{LP}(\mathbf x, \boldsymbol\beta)$ defined in \eqref{eq_LP_parametric}, and $\mathfrak N$ contains all variables not in $\mathfrak{B}$. Note that its feasible set $\mathcal{U}(\mathbf x)$ is parameterized by $\mathbf x$ in its RHS, while bases are independent of  $\mathbf x$ and $\boldsymbol\beta$. According to the theory of linear programming, a basis corresponds to a basic solution (BS), and if feasible, a basic feasible solution (BFS) that is an extreme point of underlying polyhedron $\mathcal{U}(\mathbf x)$.  Following the literature on LP, we assume that all rows of $\mathbf F$ are linearly independent. Next, we present a result regarding a basis and $\mathcal{U}(\mathbf x)$.

\begin{lem}
	\label{lem_feasible_then_optimal}
	Consider $\mathbf{LP}(\mathbf x^0, \boldsymbol\beta)$ for a fixed $\boldsymbol\beta$, and suppose that basis $\mathfrak{B}^0$ is an optimal basis, i.e., its BS with respect to $\mathcal{U}(\mathbf x^0)$ is a BFS and an optimal solution. If $\mathfrak{B}^0$'s BS with respect to $\mathcal{U}(\mathbf x^1)$ is feasible, i.e., a BFS, it is also optimal to $\mathbf{LP}(\mathbf x^1, \boldsymbol\beta)$.  Moreover, if $\mathfrak{B}^0$ yields the unique optimal solution to $\mathbf{LP}(\mathbf x^0, \boldsymbol\beta)$, it also yields the unique one to $\mathbf{LP}(\mathbf x^1, \boldsymbol\beta)$.
\end{lem}
\begin{proof} See its proof in Appendix \ref{apd:proofS3}.
\end{proof}	

Then, as shown in the following, we can make use of bases to terminate Variant 2. Without loss of generality, we assume that the oracle computing $\mathbf{MP2}$  returns BFS  for sets $\mathcal{OU}(\mathbf x, \boldsymbol\pi)$ and $\mathcal{OV}(\mathbf x, \boldsymbol\gamma)$ (i.e., extreme point solutions of $\mathcal{U}(\mathbf x)$)  as well as the associated optimal bases  (or arbitrary one(s) when degeneracy exists). Additionally,  let  $\mathbb{B}^t$ denote the \textit{set} of bases obtained from all $\mathcal{OU}$ and $\mathcal{OV}$ sets in computing $\mathbf{MP2}$ in $t$-th iteration.


\begin{lem}
	\label{lem_repeated_bases_ccg}
	If $\mathbb{B}^{t_1}=\mathbb{B}^{t_2}$ with $t_1<t_2$, Variant 2 terminates, and $\mathbf x^1$, an optimal solution to $\mathbf{MP2}$ in $t_1$-th iteration, is optimal to $\mathbf{2-Stg \ RO}$.
\end{lem}
\begin{proof} See its proof in Appendix \ref{apd:proofS3}.
\end{proof}	

Note that the total number of bases of $\mathbf{LP}(\mathbf x, \boldsymbol\beta)$ is finite and bounded by $\binom{n_u+\mu_u}{\mu_u}$  with $\mu_{u}$ denoting the number of rows of matrix $\mathbf F$. Hence,  
by enumerating all combinations of bases, which is up to $2^{\binom{n_u+\mu_u}{\mu_u}}$,  the following result can be derived from Lemma \ref{lem_repeated_bases_ccg} easily.
\begin{prop}
\label{prop_repeatbases_CCG}
    Variant 2 either reports that $\mathbf{2-Stg \ RO}$ is infeasible, or converges to its optimal value and solution in a finite number of iterations, which is bounded by $2^{\binom{n_u+\mu_u}{\mu_u}}$. Hence, the algorithm is of ${O}(2^{\binom{n_u+\mu_u}{\mu_u}})$ iteration complexity. ~\hfill$\square$
\end{prop}

\begin{rem} $(i)$ Because of Corollary \ref{cor_V2_complexity} and Propositions \ref{prop_cardinality_repeat_CCG} and \ref{prop_repeatbases_CCG}, it is straightforward to conclude that Variant 2's iteration complexity is of ${O}(\min\{\binom{n_y+\mu_y}{\mu_y},2^{\binom{n_u+\mu_u}{\mu_u}}\})$ or \\ ${O}(\min\{\binom{n_y+\mu_y}{\mu_y},2^{\binom{n_u+\mu_u}{\mu_u}}, |\mathcal{X}|\})$ when $\mathcal{X}$ is a finite discrete set.
	\\
	$(ii)$  The significance of Proposition \ref{prop_repeatbases_CCG} lies in that it guarantees the finite convergence by $\mathcal{U}(\mathbf x)$, instead of the recourse problem (or its dual problem). It supports us to handle more complex recourse problems. For example, an SOC program, whose dual problem is again an SOC program, could have infinite number of extreme points. By Proposition \ref{prop_repeatbases_CCG}, Variant 2  solves $\mathbf{2-Stg \ RO}$ with a polyhedral DDU set and an SOC recourse problem in finite iterations.  
\end{rem}

Obviously, those iteration complexity results altogether establish the finite convergence of Variant 2 for a broad class of two-stage RO. Nevertheless,  neither Corollaries \ref{cor_V2_complexity} nor Proposition
\ref{prop_repeatbases_CCG} generalizes the iteration complexity result of  basic C\&CG when handling DIU sets. Because of the connection between basis and extreme point, we believe that basis is the right tool and deserves a deeper study in the context of two-stage RO. Indeed, it helps us modify Variant 2 a little bit so that it comprehensively generalizes  basic C\&CG in both operations and computational complexity.

\subsubsection{Modifying Variant 2 to Achieve A Stronger Performance}
 In the context of linear program, basis and extreme point are rather identical unless degeneracy exists. With the indeterminate structure of $\mathcal{U}(\mathbf x)$, they become very different and basis turns out to be a more powerful and flexible tool. In the following, we consider a special situation which perfectly demonstrates the strength of this concept in helping us understand Variant 2's theoretical performance. 
 
 
 To facilitate our derivation, we include the following additional operation in Variant 2 between $(i)$ and $(ii)$ for both cases in Step 4, and number it by $(i.a)$:\\
$(\mathbf{i.a})$ compute $\mathbf{LP}(\mathbf x^*, \boldsymbol\pi^*)$ (or $\mathbf{LP}(\mathbf x^*, \boldsymbol\gamma^*)$, respectively) with an optimal BFS and the associated basis $\mathfrak{B}^{*}$.

Note that $\mathbf u^*_s$ obtained from computing $\mathbf{SP2}$ is clearly optimal to this linear program but is not necessarily identical to this BFS. We mention that this additional operation does not change the behavior of the algorithm except for providing information on $\mathfrak{B}^{*}$ for our proof. Also, recall that Variant 2 adopts the unified cutting sets as its default implementation.

\begin{prop}
\label{prop_complexity_uniqueness}
Assume that $\mathbf{LP}(\mathbf x, \cdot)$ always has a unique optimal solution in the execution of Variant 2 (with operation $(i.a)$), which is referred to as the \textit{unique optimal solution} (or ``uniqueness'' for short) property. Then, the number of iterations before  termination is bounded by the number of bases of $\mathbf{LP}(\mathbf x, \cdot)$, which is of $O(\binom{n_u+\mu_u}{\mu_u})$.
\end{prop}
\begin{proof} See its proof in Appendix \ref{apd:proofS3}.
\end{proof}	

With the connection between bases and extreme points, this result clearly generalizes and is comparable to the previous one developed for basic CC\&G with respect to DIU. Certainly it would be of a great interest if the uniqueness assumption can be dropped. 
Note that this property is only necessary to ensure that a previously derived basis contributes its BFS  in $\mathbf{MP2}$ in the iteration whenever it is contained in $\mathcal{OU}$ or $\mathcal{OV}$ again, which therefore renders $LB=UB$. Indeed, by a very classical result stated in the following,  this property always holds if we make simple modifications on $\mathbf{LP}(\mathbf x, \cdot)$ --- for any extreme point of a polyhedron (e.g., $\mathcal{U}(\mathbf x)$), there is a linear objective function to which this point is uniquely  optimal~\citep{NemWol.1988}. 

Our basic idea is to slightly modify the objective function coefficients of $\mathbf{LP}(\mathbf x, \cdot)$ according to reduced costs to achieve the uniqueness property. Note that if some non-basic variable(s) has zero reduced cost (such information is always available after an optimal BFS and the associated optimal basis $\mathfrak{B}^{*}$ are provided), it indicates the existence of multiple optimal BFSs.  We can make simple changes in $\mathbf{LP}(\mathbf x, \cdot)$'s objective function  to eliminate occurrence of that situation, which ensures that $\mathfrak{B}^{*}$'s BS is uniquely optimal. The updated objective function will then be used to define set $\mathcal{OU}$ or $\mathcal{OV}$, respectively, which is denoted by $\widehat{\mathcal{OU}}$ or $\widehat{\mathcal{OV}}$ to highlight the difference.  
By Lemma \ref{lem_feasible_then_optimal},  $\mathfrak{B}^{*}$'s BS is the only solution in $\widehat{\mathcal{OU}}$ (or $\widehat{\mathcal{OV}}$, respectively) whenever it is a BFS. Detailed modifications are available in Appendix~\ref{Asect_modiunique}, and we refer to this modified version as the \textit{modified Variant~2}. 

Using an argument similar to that of  Proposition \ref{prop_complexity_uniqueness}, we can easily show that the modified Variant 2 has a desirable iteration complexity whenever $\mathcal{U}(\mathbf x)$ is with RHS dependence.
\begin{prop}
	\label{prop_complexity_mCCG_0}
	The number of iterations of the modified Variant 2 before  termination is bounded by the number of bases of $\mathbf{LP}(\mathbf x, \cdot)$, which is of $O(\binom{n_u+\mu_u}{\mu_u})$.
\end{prop}

%

 \begin{rem}
 	\label{rem_complexity_mCCG}
 	Proposition \ref{prop_complexity_mCCG_0} and associated derivations present several new insights. \\  
 	($i$) When the modifications on the objective function coefficients of $\mathbf{LP}(\mathbf x, \cdot)$ is sufficiently small, it can be seen that $\widehat{\mathcal{OU}}(\mathbf x, \cdot)\subseteq \mathcal{OU}(\mathbf x, \cdot)$ (or $\widehat{\mathcal{OV}}\subseteq \mathcal{OU}(\mathbf x, \cdot)$, respectively) for all $\mathbf x\in\mathcal{X}$, i.e., it does not  expand to include a new BFS. This result ensures that the modified Variant 2 also inherits Variant 2's complexity result described in Corollary \ref{cor_V2_complexity}.  So, the modified Variant~2's iteration complexity is of ${O}(\min\{\binom{n_y+\mu_y}{\mu_y},\binom{n_u+\mu_u}{\mu_u}\})$ or  ${O}(\min\{\binom{n_y+\mu_y}{\mu_y},\binom{n_u+\mu_u}{\mu_u}, |\mathcal{X}|\})$ when $\mathcal{X}$ is a finite discrete set. We can further conclude that, because of the connection between bases and extreme points,  the modified Variant 2 extends and generalizes  basic C\&CG in both aspects, i.e., the operations and computational complexity.\\
 	($ii$) It is worth highlighting the basis based analysis approach, which we believe is a rather novel application of the ``basis'' concept from LP theory. In addition to using it for analysis,  we are inspired to develop another variant of C\&CG in the following subsection, which explicitly makes use of bases within its  operations.  \\
 	($iii$) The reduced cost based modifications involve  some operation, i.e., obtaining the associated optimal basis after  solving an LP, that is not supported by current professional solvers. This situation clearly imposes a technical difficulty in implementation.  Also, those modifications require non-trivial operations. Nevertheless, in practice we often have structural insights that can help to implement the modified Variant 2 without performing those operations. A demonstration on the reliable p-median problem is presented in Section \ref{subsect_Uniqueness_Pareto}. We observe that on uncapacitated instances the modified version dominates the standard Variant~2 with a clear advantage.
 \end{rem}

\subsection{The Third Variant: Generating Bases by C\&CG}
\label{subsect_3rdvariant}
The concept of basis and related methodology lay the foundation of Simplex algorithm in LP, which is theoretically elegant and practically fast. In the previous subsection, we also adopt them  to develop a set of new convergence and complexity analyses for Variant~2.
It would be interesting to design some algorithm that directly makes use of bases and other powerful LP tools to achieve a strong performance. This idea is investigated in this subsection. Before presenting our results, we assume in this subsection that the oracle solving linear programs always returns both an optimal BFS and the associated basis.

\subsubsection{Another Single-Level Reformulation through Basis Enumeration}
Consider an extreme point and the associated basis $\mathfrak{B}$ of $\mathbf{LP}(\mathbf x,\boldsymbol\beta)$ defined in \eqref{eq_LP_parametric} for a given $\mathbf x\in \mathcal{X}$.
Based on $\mathfrak{B}$ (and its complement $\mathfrak N$) we can always rewrite and reorganize constraints of $\mathcal U(\mathbf x)$ in the following form, where equality constraints correspond to those whose slack variables are in $\mathfrak{N}$ and inequality constraints represent the remaining ones, and $\mathbf u_{\mathfrak B}$ and $\mathbf u_{\mathfrak N}$ denote the subvectors  consisting of $u_j$ with $i\in\mathfrak B$ and $i\in \mathfrak N$ respectively. For the purpose of simplicity, we do not explicitly list slack variables, noting that they  either are provided or can be inferred easily according to the context.  Let $\mathcal{BS}(\mathfrak B, \mathbf x)$ denote this singleton as in the following, noting that its only element is that extreme point. 
\begin{equation}
	\label{eq_basis_constraints}
 \mathcal{BS}(\mathfrak B, \mathbf x) =\left\{
  \begin{array}{l}
   \mathbf{F}_\mathfrak{N}(\mathbf x)\mathbf u =\mathbf{h}_{\mathfrak{N}}+\mathbf{G}_{\mathfrak{N}}\mathbf x\\
   \mathbf{F}_{\mathfrak{B}}(\mathbf x)\mathbf u \leq\mathbf{h}_{\mathfrak{B}}+\mathbf{G}_{\mathfrak{B}}\mathbf x\\
    \ \ \mathbf{u}_{\mathfrak{B}} \geq  \mathbf 0, \ \mathbf{u}_{\mathfrak{N}} = \mathbf 0
      \end{array}\right\}
 \end{equation}

 Because of the connection between extreme points and bases and Proposition \ref{prop_extremepoint},  it is anticipated that we can enumerate all such  bases of $\mathcal{U}(\mathbf x)$ to convert $\mathbf{2-Stg \ RO}$  into a single-level optimization model as those in \eqref{eqn_2stgRO_Pi_KKT} and \eqref{eq_2stgRO_PIOmega}. Previous approaches employ optimization problems (equivalently their optimality conditions) to elicit   critical scenarios, whose validity holds regardless of the particular realization of parameterized set $\mathcal U(\mathbf x)$. Nevertheless, directly utilizing an explicit and fixed basis in algorithm design is technically more challenging. Note that the change of $\mathbf x$ may alter $\mathcal U(\mathbf x)$ drastically, especially when LHS dependence presents, i.e., constraint matrix $\mathbf F$ is a function of $\mathbf x$. Next, we give a simple $\mathcal U(\mathbf x)$ as in the following to show this challenge, where $\mathbf u^s$ represents slack variables and $\mathbf F(\mathbf x)=\begin{bmatrix}
 	f_{11} & f_{12} \\
 	f_{21} & f_{22}
 \end{bmatrix}$ is a function of the first stage decision $\mathbf x$.
 \begin{align*}
 f_{11}u_1+f_{12}u_2+u^{s}_1= \ &3\\
  f_{21}u_1+f_{22}u_2+u^{s}_2=\ &5\\
  u_1\geq 0, u_2 \geq 0, u^s_1 \geq 0, \ & u^s_2 \geq 0  	
 \end{align*}
 Assume that for some $\mathbf x^0$ we have $\mathbf F(\mathbf x^0) =\begin{bmatrix}
 	2 & 1\\
 	1 & 2
 \end{bmatrix}$. Considering basis $\mathfrak{B}^0=\{u_1,u_2\}$, we have $\mathcal{BS}(\mathfrak B^0, \mathbf x^0)=\left\{\left(\frac 1 3, \frac 7 3\right)\right\}$, an extreme point solution. If $\mathbf x^1$ renders $\mathbf F(\mathbf x^1)=\begin{bmatrix}
 1 & 1\\
 1 & 1
\end{bmatrix}$, note that $\mathcal U(\mathbf x^1)$ is not empty. Nevertheless, if rewriting the whole system with respect to $\mathfrak{B}^0$, we have
 \begin{equation*}
 \label{eq_eg_basis}
     \mathcal{BS}(\mathfrak B^0, \mathbf x^1) = \left\{\begin{array}{l}
u_1+ u_2= \ 3\\	 u_1+u_2=\ 5\\
u_1\geq 0, \ u_2 \geq 0
     \end{array}\right\}
\end{equation*}
which is simply empty (i.e., infeasible). Hence, if \eqref{eq_basis_constraints} with $\mathfrak{B}^0$, i.e., $\{f_{11}u_1+f_{12}u_2=3,
f_{21}u_1+f_{22}2u_2=5, u_1\geq 0, u_2\geq 0\}$,  is imposed as constraints on $\mathcal{X}$, $\mathbf x^1$
will be eliminated from $\mathcal{X}$. However, it should not be the case given that $\mathcal U(\mathbf x^1)$ is meaningful.

The aforementioned situation actually can be detected  by applying one core result for the system of linear inequalities, i.e., Theorem of Alternatives \citep{boyd2004convex}. In the following, we present a variant of this theorem tailored for \eqref{eq_basis_constraints}, where $\mathbf F_{\mathfrak B}$ and $\mathbf F_{\mathfrak N}$ are simplified to 
$\hat{\mathbf F}_{\mathfrak B}$ and $\hat{\mathbf F}_{\mathfrak N}$, respectively, after removing columns associated with $\mathbf u_{\mathfrak N}$.

\begin{lem}[ Theorem of Alternatives]
	\label{lem_alternative_thm}
	One and only one  of the following two statements holds:
	either
	$0=\max\{0: \mathbf u\in \mathcal{BS}(\mathfrak B, \mathbf x)\}= \min\{(\mathbf h_{\mathfrak N}+\mathbf G_{\mathfrak N}\mathbf x)^\intercal\boldsymbol\lambda_{\mathfrak N}+(\mathbf h_{\mathfrak B}+\mathbf G_{\mathfrak B}\mathbf x)^\intercal\boldsymbol\lambda_{\mathfrak B}: \hat{\mathbf F}_{\mathfrak N}(\mathbf x)^\intercal\boldsymbol\lambda_{\mathfrak N}+\hat{\mathbf F}_{\mathfrak B}(\mathbf x)^\intercal\boldsymbol\lambda_{\mathfrak B}\geq \mathbf 0, \boldsymbol\lambda_{\mathfrak B}\geq \mathbf 0\}$, or 
	the preceding maximization problem is infeasible and the minimization problem is unbounded. 
	\hfill $\square$
\end{lem}

 Next, we discuss some changes to achieve the computational feasibility of directly using bases. Let $\mathbb{B}$ denote the collection of possible bases of $\mathbf{LP}(\mathbf x, \cdot)$, and for $\mathfrak{B}\in \mathbb{B}$, we extend $\mathcal{BS}$ to $\overline{\mathcal{BS}}$ as in the following. We mention that with additional variables $(\bar{\mathbf {u}}^1, \bar{\mathbf {u}}^2, \bar{\mathbf {u}}^3)$, which are 3 vectors of appropriate dimensions,  $\overline{\mathcal{BS}}(\mathfrak B, \mathbf x)$ is not empty for $\mathbf x\in \mathcal{X}$.

\begin{equation}
	\label{eq_basis_constraints_ex}
	\overline{\mathcal{BS}}(\mathfrak B, \mathbf x) = \left\{
	\begin{array}{l}
		\hat{\mathbf{F}}_{\mathfrak{N}}(\mathbf x)\mathbf u_{\mathfrak B} =\mathbf{h}_{\mathfrak{N}}+ \mathbf{G}_{\mathfrak{N}}\mathbf x+ \bar{\mathbf {u}}^{1}- \bar{\mathbf {u}}^{2}\\
		\hat{\mathbf{F}}_{\mathfrak{B}}(\mathbf x)\mathbf u_{\mathfrak B} \leq\mathbf{h}_{\mathfrak{B}}+\mathbf{G}_{\mathfrak{B}}\mathbf x+\bar{\mathbf {u}}^3\\
		\mathbf u_{\mathfrak B}\geq \mathbf 0, \  \bar{\mathbf {u}}^j\geq \mathbf 0,  j=1,2,3
	\end{array}\right\}
\end{equation}

By mainly taking advantage of Lemma \ref{lem_alternative_thm} and enumeration of bases in $\mathbb{B}$, we next show that $\mathbf{2-Stg \ RO}$ can be converted into a new single-level formulation.

\begin{thm}
	\label{thm_basis_enumeration}
   Formulation $\mathbf{2-Stg \ RO}$ in \eqref{eq_2RO} (and its equivalences) is equivalent to a sing-level optimization program as in the following.
   \begin{subequations}
   	\label{eq_2stgRO_Bases}
   	\begin{align}
   		\mathbf{2\!-\!Stg \ RO(\mathbb{B})}: \   w^*=\min \quad & \mathbf{c}_1\mathbf x+ \eta   \\
   		\mathrm{s.t.} \quad & \mathbf x  \in \mathcal{X}\\
   		 \bigg\{&\eta\geq   \mathbf{c_2}\mathbf y+ M(\mathbf1^\intercal\bar{\mathbf u}^1+\mathbf1^\intercal\bar{\mathbf u}^2+\mathbf1^\intercal\bar{\mathbf u}^3) \label{eq_enu3_basis}\\
   		& + (\mathbf h_{\mathfrak N}+\mathbf G_{\mathfrak N}\mathbf x)^\intercal\boldsymbol\lambda_{\mathfrak N}+(\mathbf h_{\mathfrak B}+\mathbf G_{\mathfrak B}\mathbf x)^\intercal\boldsymbol\lambda_{\mathfrak B} \notag \\
   		&   (\mathbf{u},\bar{\mathbf u}^1,\bar{\mathbf u}^2,\bar{\mathbf u}^3)\in \overline{\mathcal{BS}}(\mathfrak B, \mathbf x) \label{eq_enu3_basis2}\\
   		&\hat{\mathbf F}_{\mathfrak N}(\mathbf x)^\intercal\boldsymbol\lambda_{\mathfrak N}+\hat{\mathbf F}_{\mathfrak B}(\mathbf x)^\intercal\boldsymbol\lambda_{\mathfrak B}\geq \mathbf 0, \ \boldsymbol\lambda_{\mathfrak N}\geq \mathbf 0 \\
   		& \mathbf y\in \mathcal Y(\mathbf x, \mathbf u) 	\bigg\}  \ \forall \mathfrak B\in \mathbb{B} \label{eq_enu3_optimality_basis}
   	\end{align}
   \end{subequations}
\end{thm}
\begin{proof}
	See its proof in Appendix \ref{apd:proofS3}.
\end{proof}

\begin{rem}
	Compared to the previous two single-level reformulations, $\mathbf{2-Stg \ RO(\mathbb B)}$ seems to have a relatively simpler structure, given that there are no complex optimality conditions involved. As noted earlier, bases and extreme points of a polyhedron are directly connected. So,
	$\mathbf{2-Stg \ RO(\mathbb B)}$ should reduce to \eqref{eq_DIU_enumeration} of Corollary \ref{cor_DIU_Omega}  when $\mathcal U(\mathbf x)=\mathcal U^0$ for all $\mathbf x$. Actually, except the infeasible bases, it is the case according to the proof of Theorem \ref{thm_basis_enumeration}. Hence,  $\mathbf{2-Stg \ RO(\mathbb B)}$ can be treated as an algebraic generalization of \eqref{eq_DIU_enumeration}. 
\end{rem}
Similar to Corollaries \ref{cor_Pi_partial_KKT} and  \ref{cor_Pi_partial_KKTOmega}, we derive a relaxation to $\mathbf{2\!-\!Stg \ RO(\mathbb{B})}$ (and therefore to $\mathbf{2\!-\!Stg \ RO}$) based on partial enumeration.
\begin{cor}
	\label{cor_Bases_partial}
	Consider set $\hat {\mathbb B}\subseteq \mathbb B$, and denote formulation~\eqref{eq_2stgRO_Bases} defined with respect to $\hat{\mathbb B}$ by $\mathbf{2-Stg \ RO(\hat{\mathbb B})}$. Then,  $\mathbf{2-Stg \ RO(\hat{\mathbb B})}$ is a relaxation to $\mathbf{2-Stg \ RO(\mathbb B)}$ (and
	$\mathbf{2-Stg \ RO}$ in \eqref{eq_2RO}). And the optimal value of $\mathbf{2-Stg \ RO(\hat{\mathbb B})}$ is smaller than or equal to~$w^*$.  \hfill $\square$
\end{cor}

\subsubsection{Variant 3: Basis Based C\&CG and Analysis}
Theorem \ref{thm_basis_enumeration} and Corollary \ref{cor_Bases_partial} naturally yield a foundation to develop our third variant of C\&CG.  Similar to the case of Variant 2, we describe it by specifying modifications on top of Variant 1 to minimize repetition. Also, we simply adopt the unified cutting sets for all identified bases, in spite of they are derived due to optimality or feasibility reasons.


\begin{description}
	\item[$\bullet$ Master Problem] Throughout the algorithm, $\mathbf{MP1}$ is replaced by the next one. Recall that $\mathfrak N$ is uniquely defined for any fixed $\mathfrak B$.
	\begin{eqnarray*}
		\label{eq_master_1}
		\mathbf{MP3}: \ \underline w=\min && \mathbf{c}_1\mathbf{x}+ \eta \notag \\
		\mbox{s.t.} && \mathbf{x} \in \mathcal{X} \notag\\
		&& \eta\geq \mathbf{c_2}\mathbf y^{\mathfrak B}+ M(\mathbf1^\intercal\bar{\mathbf u}^1_{\mathfrak B}+\mathbf1^\intercal\bar{\mathbf u}^2_{\mathfrak B}+\mathbf1^\intercal\bar{\mathbf u}^3_{\mathfrak B}) + (\mathbf h_{\mathfrak N}+\mathbf G_{\mathfrak N}\mathbf x)^\intercal\boldsymbol\lambda_{\mathfrak N}\\
		&& +(\mathbf h_{\mathfrak B}+\mathbf G_{\mathfrak B}\mathbf x)^\intercal\boldsymbol\lambda_{\mathfrak B} \ \ \forall \mathfrak B\in \mathbb{\hat B} \\
        && \hat{\mathbf{F}}_{\mathfrak{N}}(\mathbf x)\mathbf u_{\mathfrak B} =\mathbf{h}_{\mathfrak{N}}+ \mathbf{G}_{\mathfrak{N}}\mathbf x+ \bar{\mathbf {u}}^{1}_{\mathfrak B}- \bar{\mathbf {u}}^{2}_{\mathfrak B} \ \ \forall \mathfrak B\in \mathbb{\hat B}\\
        &&\hat{\mathbf{F}}_{\mathfrak{B}}(\mathbf x)\mathbf u_{\mathfrak{B}} \leq\mathbf{h}_{\mathfrak{B}}+\mathbf{G}_{\mathfrak{B}}\mathbf x+\bar{\mathbf {u}}^3_{\mathfrak B} \ \ \forall \mathfrak B\in \mathbb{\hat B}\\
		&&\mathbf F_{\mathfrak N}(\mathbf x)^\intercal\boldsymbol\lambda_{\mathfrak N}+\mathbf F_{\mathfrak B}(\mathbf x)^\intercal\boldsymbol\lambda_{\mathfrak B}\geq \mathbf 0 \ \ \forall \mathfrak B\in \mathbb{\hat B} \\
		&& \mathbf{B}_2\mathbf y^{\mathfrak B}\geq \mathbf{d-\mathbf B}_1\mathbf x-\mathbf{Eu}_{\mathfrak B} \ \ \forall \mathfrak B\in \mathbb{\hat B}\\
		&& \mathbf u_{\mathfrak B}\geq  \mathbf 0, \  \bar{\mathbf {u}}^j_{\mathfrak B}\geq \mathbf 0,  j=1,2,3,\ \ \boldsymbol\lambda_{\mathfrak N}\geq \mathbf 0, \ \mathbf y^{\mathfrak B}\geq \mathbf 0 \ \ \forall \mathfrak B \in \mathbb{\hat B} \label{eq_enu3_basis2}
		\end{eqnarray*}
	\item[$\bullet$ Cutting Sets]: For every $\mathfrak B$ identified by computing subproblems, we generate and incorporate a replicate of recourse problem (including recourse variables and constraints) to $\mathbf{MP3}$. Specifically, we modify \textbf{Step 4 [Case A]} as in the following.
	
	$(i)$ compute  $\mathbf{SP2}$ in \eqref{eq_SP2} to derive $\eta_s(\mathbf x^*)$, $\mathbf u^*_s$ and corresponding extreme  point $\boldsymbol\pi^*$ of $\Pi$; $(i.a)$ compute $\mathbf{LP}(\mathbf x^*, \boldsymbol\pi^*)$ with an optimal BFS and its associated basis $\mathfrak{B}^*$; $(ii)$ update $\mathbb{\hat B}= \mathbb{\hat B}\cup\{\mathfrak B^*\}$ and augment master problem $\mathbf{MP3}$ accordingly, i.e., create variables $\mathbf u_{\mathfrak B^*}, \bar{\mathbf u}^j_{\mathfrak B^*}, j=1,2,3, \boldsymbol\lambda_{\mathfrak N^*}, \boldsymbol\lambda_{\mathfrak B^*}$ and $\mathbf y^{\mathfrak B^*}$ and add following constraints  to~$\mathbf{MP3}$.
\begin{align}
\label{eq_enu3_cuttingset}
\begin{split}
	 \eta&\geq \mathbf{c_2}\mathbf y^{\mathfrak B^*}+ M(\mathbf1^\intercal\bar{\mathbf u}^1_{\mathfrak B^*}+\mathbf1^\intercal\bar{\mathbf u}^2_{\mathfrak B^*}+\mathbf1^\intercal\bar{\mathbf u}^3_{\mathfrak B^*})+ (\mathbf h_{\mathfrak N^*}+\mathbf G_{\mathfrak N^*}\mathbf x)^\intercal\boldsymbol\lambda_{\mathfrak N^*}\\
	 & \ +(\mathbf h_{\mathfrak B^*}+\mathbf G_{\mathfrak B^*}\mathbf x)^\intercal\boldsymbol\lambda_{\mathfrak B^*} \\
	& \hat{\mathbf{F}}_{\mathfrak{N}^*}(\mathbf x)\mathbf u_{\mathfrak B^*} =\mathbf{h}_{\mathfrak{N}^*}+ \mathbf{G}_{\mathfrak{N}^*}\mathbf x+ \bar{\mathbf {u}}^{1}_{\mathfrak B^*}- \bar{\mathbf {u}}^{2}_{\mathfrak B^*}\\
	&\hat{\mathbf{F}}_{\mathfrak{B}^*}(\mathbf x)\mathbf u_{\mathfrak{B}^*} \leq\mathbf{h}_{\mathfrak{B}^*}+\mathbf{G}_{\mathfrak{B}^*}\mathbf x+\bar{\mathbf {u}}^3_{\mathfrak B^*}\\	
	& \hat{\mathbf F}_{\mathfrak N^*}(\mathbf x)^\intercal\boldsymbol\lambda_{\mathfrak N^*}+\mathbf F_{\mathfrak B^*}(\mathbf x)^\intercal\boldsymbol\lambda_{\mathfrak B^*}\geq \mathbf 0\\
	& \mathbf{B}_2\mathbf y^{\mathfrak B^*}\geq \mathbf{d-\mathbf B}_1\mathbf x-\mathbf{Eu}_{\mathfrak B^*}\\
	& \mathbf u_{\mathfrak B^*}\geq  \mathbf 0, \  \bar{\mathbf {u}}^j_{\mathfrak B^*}\geq \mathbf 0,  j=1,2,3\ \ \boldsymbol\lambda_{\mathfrak N^*}\geq \mathbf 0, \ \mathbf y^{\mathfrak B^*}\geq \mathbf 0
\end{split}
\end{align}

	\indent  For operations in \textbf{Step 4 [Case B]},  almost same modifications will be made, except that $(\mathbf{i.a})$ is changed to ``compute $\mathbf{LP}(\mathbf x^*, \boldsymbol\gamma^*)$ with an optimal solution and the associated basis $\mathfrak{B}^*$''.   Note that  the feasibility cutting set in \eqref{eq_BD_feasibility} is replaced by the unified cutting set \eqref{eq_enu3_cuttingset} defined with respect to this identified $\mathfrak B^*$.
\end{description}

In the following, we provide the convergence and iteration complexity results of Variant~3. Similar to Proposition \ref{prop_complexity_uniqueness}, they can be easily proven by arguing that the algorithm converges if a particular basis is identified more than once.

\begin{thm}
	\label{thm_V2_converge}
	Suppose $T\!O\!L=0$.  When terminates, Variant 3 either reports that $\mathbf{2-Stg \ RO}$ in \eqref{eq_2RO} is infeasible, or converges to its optimal value and exact solution. ~\hfill$\square$
\end{thm}

\begin{cor}
	\label{cor_V3_complexity}
	Variant 3 is  of $O(\binom{n_u+\mu_u}{\mu_u})$ iteration complexity.~\hfill$\square$
\end{cor}

\begin{rem}
	$(i)$ Since $\mathbf{2-Stg \ RO(\mathbb B)}$ is an algebraic generalization of the extreme point based reformulation in
	\eqref{eq_DIU_enumeration}, Variant 3 is clearly an algebraic generalization of  basic C\&CG. As it is built upon bases, we also refer to Variant 3 as \textit{basis based C\&CG} method.  Note that, compared to Variant 2, Variant 3 is rather a static generalization. \\
	$(ii)$ Variant 3 demonstrates a theoretically strong iteration complexity regardless of LHS or RHS. It is clearly desired as this iteration complexity actually is comparable to and generalizes that of basic C\&CG.   Nevertheless, in our numerical study, its computational performance is not as good as that of Variant 2.  One primary reason is the complex  constraints with big-M and bilinear terms  introduced in \eqref{eq_enu3_basis}.\\
	 $(iii)$  As noted earlier, an additional technical difficulty is that current professional MIP solvers do not provide optimal basis when they report optimal solution, although we might be able to infer this information according to variables' values in an optimal solution. Indeed, it is often the case that we cannot single out a specific basis, as multiple bases could be optimal. If it occurs in our algorithm implementation, we simply include all those bases (i.e., their associated cutting sets in \eqref{eq_enu3_cuttingset}) into $\mathbf{MP3}$.
\end{rem}

Given Variant 3's nice theoretical properties, we still believe that it deserves further studies and deep analyses, especially stronger and better integration with MIP solvers. Moreover, generalizing extreme points to either bases or sets of optimal solutions (as in Variant 2) provides a higher level abstraction and could inspire new strategies to handle problems with indeterminate structures. Hence, we anticipate that more sophisticated algorithm designs will be developed in this line of research.

\subsection{Computing  Complex Two-stage RO Formulations}
\label{subsect_complex_RO}
We would like to note that  basic C\&CG and variants developed in this section are rather general computing schemes.  When the utilized oracle for master and subproblems is powerful, they can be called to compute  complex nonlinear two-stage RO formulations beyond the linear one  in the form of \eqref{eq_2RO}. We next
discuss briefly the applicability of those algorithms to more complex two-stage RO formulations.

Regarding the feasible set for the first stage decision, i.e., $\mathcal{X}$, it does not need to be a linear MIP set for any of three variants. It actually can be more complex sets, as long as master problems can be solved exactly. For example, $\mathcal{X}$ can be a mixed integer SOC set or a mixed integer bilinear set, given that they are supported by commercial solvers. As for the recourse problem, we can employ convex programs whose optimality conditions are supported by some oracles. For example, the recourse problem can be an SOC formulation with the strong duality. By replacing the recourse problem with its dual problem, we then can convert a $\max-\min$ subproblem into a $\max$ formulation, which is directly computable by some commercial solvers.   As for the uncertainty set, as Variant 3 depends on the concept of basis and related mathematical representation,  it is restricted to a polyhedral DDU set. For the other two, the uncertainty set can be a set of convex program whose optimality conditions are friendly to existing oracles. Moreover, we can bound the iteration complexity of those algorithms according to the context and by using some of the following structural information, i.e.,  the cardinality of $\mathcal{X}$ if it is finite, the extreme points and rays of $\Pi$ if the recourse problem is an LP, and the bases of $\mathbf{LP}(\mathbf x, \cdot)$ if $\mathcal{U}(\mathbf x)$  is a polyhedron. 

For two-stage RO models with SOC recourse problem (and SOC uncertainty set, respectively), we present  detailed algorithm operations using the scheme of Variant 2  in Appendix~\ref{Asect:V2SOC}. 
In Section \ref{subsect_complexRO}, we consider 
two such formulations and perform a small computational study by taking advantage of a professional solver's solution capacity on mixed integer SOC programs. It is observed that the numbers of iterations for those complex two-stage RO models are roughly same as those for their linear counterparts. Nevertheless, much longer computational time is needed for those models, which indicates the great challenge imposed by their nonlinear structures.  

Before closing this section, we mention that, similar to the study in \citet{an2015exploring}, those new C\&CG variants can solve, with little change, risk constrained two-stage RO with DDU, where risk constraint(s) in the form of 
$\displaystyle\max_{\mathbf u\in \mathcal{U}(\mathbf x)}\min\{\mathbf{c}_2\mathbf y: \mathbf y\in\mathcal{Y}(\mathbf x,\mathbf u)\}\leq b_0$ is appended to $\mathcal{X}$. Also, the underlying sets $\mathcal X$, $\mathcal{U}(x)$ and $\mathcal{Y}(\mathbf x,\mathbf u)$ can take the form of the aforementioned complex structures, and can be exactly handled.

\section{Further Investigations on Solving Two-stage RO}
\label{sect_more}
In this section,  we present a few studies on extending and enhancing current  solution approaches for two-stage RO.
In particular, the first one deserves more attention. It seems  counterintuitive, but demonstrates a clear improvement through incorporating ``deep knowledge'' of the model. Based on our theoretical derivations and numerical results, we strongly believe that this strategy, i.e., making use of non-trivial structural properties and insights within solution algorithms, may pave the way to solve difficult two-stage RO instances.



\subsection{Converting DIU into DDU: Deep Knowledge and Implications}
\label{subsect_DIU-DDU}
As mentioned in Section \ref{subsect_specialDDU}, existing strategies often transform a DDU-based $\mathbf{2-Stg \ RO}$ into a 	DIU-based one that is solvable by existing methods. In this subsection, we investigate the reverse transformation, i.e., transforming a DIU-based RO into a DDU-based one, and computing the resulting formulation by our new algorithms. Apparently this strategy is counterintuitive. Nevertheless, it not only demonstrates a  significant improvement, but also yields a general and flexible scheme to obtain quantifiable approximate solutions, showing an unusual potential to tackle challenging instances.

Consider $\mathbf{2-Stg \ RO}$ with DIU set $\mathcal{U}^0$. Because of its structure or our sound insights on the real system, we often have some ``deep  knowledge'' (exact or approximate) on the connection between  $\mathbf x$ and worst case scenarios in $\mathcal{U}^0$. If such knowledge is represented analytically, i.e., converting $\mathcal{U}^0$ into $\mathcal{U}(\mathbf x)$ that is parameterized by $\mathbf x$, a DDU-based model could be formulated to bound or even to derive an optimal solution to the original DIU-based model. We formalize this result in  the following that can be proven easily.

\begin{prop}
	\label{prop_DIU_to_DDU}
Consider sets $\mathcal{X}^0\subseteq \mathcal{X}^1$, $\mathcal{U}^0$ and $\mathcal{U}(\mathbf x)$, and $\mathcal{Y}^0(\mathbf x,\mathbf u)$ and $\mathcal{Y}^1(\mathbf x,\mathbf u)$.
If the following inequality holds for $\mathbf x\in \mathcal{X}^0$,
\begin{eqnarray}
	\label{eq_max-min_DIUDDU}
\max_{\mathbf u\in \mathcal{U}^0}\min\{\mathbf{c}^0_2\mathbf y: \mathbf y\in \mathcal{Y}^0(\mathbf x,\mathbf u)\}\geq \max_{\mathbf u\in \mathcal{U}(\mathbf x)}\min\{\mathbf{c}^1_2\mathbf y: \mathbf y\in \mathcal{Y}^1(\mathbf x,\mathbf u)\}
\end{eqnarray}
we have
\begin{equation}
	   	\label{eq_DIUDDU}
   \begin{split}
   	&\mathit{w}(\mathcal X^0, \mathcal{U}^0, \mathcal Y^0(\mathbf x,\mathbf u)) = \min_{\mathbf x\in \mathcal{X}^0}\mathbf c_1\mathbf x+\max_{\mathbf u\in \mathcal{U}^0}\min\{\mathbf{c}^0_2\mathbf y: \mathbf y\in \mathcal{Y}^0(\mathbf x,\mathbf u)\}\geq \\  & \mathit{w}(\mathcal X^1, \mathcal{U}(\mathbf x), \mathcal Y^1(\mathbf x,\mathbf u)) = \min_{\mathbf x\in \mathcal{X}^1}\mathbf c_1\mathbf x+\max_{\mathbf u\in \mathcal{U}(\mathbf x)}\min\{\mathbf{c}^1_2\mathbf y: \mathbf y\in \mathcal{Y}^1(\mathbf x,\mathbf u)\}.
  \end{split}
\end{equation}
Moreover, if $\mathcal{X}^0=\mathcal{X}^1$ and  \eqref{eq_max-min_DIUDDU} achieves equality for $\mathbf x\in \mathcal{X}^0$, \eqref{eq_DIUDDU} reduces to an equation, and $\mathbf x^*$, which is  optimal to the RO with DDU, is also optimal to the RO with DIU.\qed
\end{prop}
\begin{rem}
	 $(i)$ We can have a similar result with sets $\mathcal{X}^0 \supseteq \mathcal{X}^1$, and ``$\geq$'' replaced by ``$\leq$'' for $\mathbf x\in \mathcal{X}^1$ in \eqref{eq_max-min_DIUDDU} and therefore in \eqref{eq_DIUDDU}. That means computing the associated DDU-based RO model provides an upper bound to $ \mathit{w}(\mathcal X^0, \mathcal{U}^0, \mathcal Y^0(\mathbf x,\mathbf u))$. Yet, constructing such a DDU-based model is less useful for our algorithm development. \\
	$(ii)$ Naturally if $\mathcal{U}(\mathbf x)\subseteq \mathcal U^0$ for $\mathbf x\in \mathcal{X}^0=\mathcal{X}^1$, it follows from Proposition \ref{prop_relax_1} directly that inequality \eqref{eq_DIUDDU} holds.
	Differently, Proposition \ref{prop_DIU_to_DDU} is more flexible and general. It actually is not necessarily true that $\mathcal{U}(\mathbf x)$ is a subset of $\mathcal U^0$, or $\mathcal{Y}^1$ and $\mathcal{Y}^2$ should be of the same structure. One essential research is to derive  $\mathcal{U}(\mathbf x)$ such that \eqref{eq_max-min_DIUDDU} achieves the equality for $\mathbf x \in \mathcal{X}^0$.  Then, some mathematical analysis  may fit $\mathcal{U}(\mathbf x)$ into the form presented in \eqref{eq_uncer_set}. If this is the case,  previously developed algorithms can be adopted to derive an optimal solution.
\end{rem}
In the remainder of this subsection, we make good use of the reliable p-median facility location model developed in \citet{an2014reliable} to illustrate the benefits from transforming DIU to DDU.
To make this paper self-contained, we first give the problem description and its conventional $\mathbf{2-Stg \ RO}$ formulation with a DIU set.

\subsubsection{The Classical Reliable Facility Location Model}
Let $I$ be the set of client sites and $J\subseteq I$ the set of potential facility sites. Each client site $i\in I$ has demand $d_i$ and the unit cost of serving demand of $i$ by the
facility at $j\in J$ is $c_{ij}\geq 0$ with $c_{ii} =0$. In the first stage, the decision maker determines the constructions of $p$ uncapacitated/capacitated (depending on capacity parameter $A_j$) facilities and allocations of demands under the normal situation. After disruptions, the recourse problem  re-allocates demands to survived facilities. The whole decision making problem seeks to minimize the weighted sum of operational costs in the normal disruption-free scenario and in the worst disruptive scenarios, which is formulated as the following.
\begin{equation}
\label{eq_reliable_pmedian}
\begin{split}
\mathbf{Reliable \ FL}: \ \mathit{w}_{R}(\mathcal{X}, \mathcal{U}^0, \mathcal{Y}(\mathbf z, \mathbf u))=\min_{(\mathbf x_c,\mathbf x_d)\in\mathcal X} & (1-\rho)\sum_{i\in I}\sum_{j\in J}c_{ij}x_{c,ij}+ \\
   \rho\max_{\mathbf u\in\mathcal U^0}\min_{(\mathbf y_1,\mathbf y_2)\in\mathcal{Y}(\mathbf x,\mathbf u)} & \ \sum_{i\in I}\sum_{j\in J}c_{ij}y_{1,ij}+\sum_{i\in I}\textsl{C}y_{2,i}
\end{split}
\end{equation}
with sets
\begin{align*}
	\mathcal X=\{(\mathbf x_c, \mathbf x_d)\in \mathbb{R}^{|I|\times|J|}_+\times\{0,1\}^{|J|}:\ \sum_{j\in J}x_{d,j}=p, \ \sum_{j\in J}x_{c,ij}\geq d_i \ \forall i, \ \sum_{i\in I}x_{c,ij}\leq A_jx_{d,j} \ \forall j\},
\end{align*}
\begin{align*}
	\mathcal U^0=\{u_i\in\{0,1\}^{|I|}: \sum_{i\in I}u_i\leq k\},
\end{align*}
and
\begin{align*}
	\mathcal Y(\mathbf x,\mathbf u)=\{(\mathbf y_1,\mathbf y_2)\mathbb{R}^{|I||J|}_+\times \mathbb{R}^{|I|}_+: & \sum_{j\in J}y_{1,ij}+y_{2,i}\geq (1+\theta_iu_i)d_i \ \forall i, \ \sum_{i\in I}y_{1,ij}\leq A_jx_{d,j} \ \forall j, \\
	&  \sum_{i\in I}y_{1,ij}\leq A_j(1-u_j) \ \forall j\}.
\end{align*}
Variables $\mathbf x=(\mathbf x_c, \mathbf x_d)$ are continuous and binary representing allocations of demands and yes-no construction decisions of facilities, respectively. Constraints in $\mathcal X$ capture the basic requirements in the normal disruption-free scenario, where $p$ facilities are to be installed with all client sites' demands being satisfied and facilities' capacities being observed.  DIU set $\mathcal U^0$ contains all disruptive scenarios with up to $k$ disruptions, where $u_i$ is binary with~$1$ denoting a disruption at site $i$ and $0$ otherwise. After disruptions,  variables $\mathbf y =(\mathbf y_{1},  \mathbf y_{2})$ represent re-allocations of demands to survived facilities and unmet demand penalized by coefficient \textsl{C}, respectively. We note in $\mathcal{Y}(\mathbf x, \mathbf u)$ that site $i$, after hit by a disruption, may have a different demand pattern, depending on parameter $\theta_i$. This phenomenon often happens in real life, where demand becomes significantly larger (e.g., demand for personal protective equipment after a pandemic)  or smaller (e.g.,  demand for luxury cars after an earthquake).

\subsubsection{Constructing DDU Set for Exact Solutions}
  An analysis in \citet{an2014reliable} shows that if disruptions do not incur more demands (i.e., $\theta_i\leq 0$ for all $i$), disruptions do not occur at non-facility sites. For the DIU-based RO model in \eqref{eq_reliable_pmedian}, this structural property is regarded as a pure theoretical insight, rather than an understanding benefiting computation. Nevertheless, by making use of this deep knowledge, we are able to construct a DDU set and then convert \eqref{eq_reliable_pmedian} into an equivalent DDU-based RO model as in  the following.
 \begin{cor}
 	\label{cor_reliable_pmedian}
     Let a DDU set be
     \begin{eqnarray}
     \label{eq_DDU_reliable_pmedian}
      \mathcal{U}^k(\mathbf x)=\{\mathbf {u}\in \mathbb{R}^{|I|}_+: \sum_{j\in J}u_j\leq k, u_j\leq x_{d,j} \forall j, \ u_i=0 \ \forall i\notin J\}.
      \end{eqnarray}
      If $\textsl{C}\geq \max_{ij}\{c_{ij}\}$ and $\theta_i\leq 0$ for all $i$,
   the two-stage RO in \eqref{eq_reliable_pmedian} is equivalent to
   \begin{equation}
   	\label{eq_reliable_pmedian_DDU}
   	\begin{split}
   	 \mathit{w}_{R}(\mathcal{X}, \mathcal{U}^k(\mathbf x), \mathcal{Y}(\mathbf x, \mathbf u)) = \min_{(\mathbf x_c,\mathbf x_d)\in\mathcal X} & (1-\rho)\sum_{i\in I}\sum_{j\in J}c_{ij}x_{c,ij}+\\
   	 & \rho\max_{\mathbf u\in\mathcal U^k(\mathbf x)}\min_{(\mathbf y_1,\mathbf y_2)\in\mathcal{Y}(\mathbf x,\mathbf u)} \ \sum_{i\in I}\sum_{j\in J}c_{ij}y_{1,ij}+\sum_{i\in I}\textsl{C}y_{2,i}.
   	\end{split}
   \end{equation}
 \end{cor}
 \begin{proof}
 		See its proof in Appendix \ref{apd:proofS4}.
 \end{proof}

 As a direct consequence of Corollary \ref{cor_reliable_pmedian} and the polytope structure of the underlying DDU set, we can solve the reliable p-median facility location problem exactly by computing the new DDU-based two-stage RO in \eqref{eq_reliable_pmedian_DDU}. We mention the  totally unimodularity property involved in the proof, noting that $\mathcal{U}^k(\mathbf x)$ does not belong to $\mathcal{U}^0$. 
  
  According to our computational study,  using parametric C\&CG to solve \eqref{eq_reliable_pmedian_DDU} demonstrates a significant improvement for the difficult capacitated instances over the traditional approach, i.e., computing \eqref{eq_reliable_pmedian} by   basic C\&CG. It either converges in fewer  iterations with reduced computational time or produces solutions with much smaller optimality gaps. Indeed, we do think that such improvement is  not surprising. Incorporating  deep knowledge  into $\mathcal{U}(\mathbf x)$ that reveals the hidden connection between $\mathbf x$ and worst case scenarios, it enables parametric C\&CG to generate cutting sets that maintain their strength even though $\mathbf x$ changes. On the contrary,  basic C\&CG generates static ones based on fixed scenarios, which are worst ones for particular $\mathbf x'$s generated in the execution of the algorithm. Nevertheless, those static cutting sets typically become weak when $\mathbf x$ takes a value different from them.

  \subsubsection{Constructing DDU Set for Approximations}
We believe that the improvement depends  on the strength of structural properties embedded in $\mathcal{U}(\mathbf x)$ on capturing the dependence between worst case scenarios and $\mathbf x$. For complex problems, deriving such important structural properties could be very demanding.   Nevertheless, practitioners often develop their own understanding based on rich experiences and historical data. Although it might not be proven theoretically, such heuristic understanding is likely to be useful or correct in practice. Indeed, through Proposition  \ref{prop_relax_1}, Corollary \ref{cor_3_level_relaxation} and Proposition \ref{prop_DIU_to_DDU},  it is feasible to incorporate heuristic understanding into the construction of $\mathcal{U}(\mathbf x)$ to obtain an approximation scheme in the context of parametric C\&CG. 


Consider a situation where constructed $\mathcal U(\mathbf x)$ and $\mathcal{Y}^1$ are in the same spaces as $\mathcal U^0$ and $\mathcal{Y}^0$, respectively. Also, \eqref{eq_max-min_DIUDDU} holds for $\mathbf x\in\mathcal{X}^0\subseteq \mathcal{X}^1$. Let $\Pi^0$ denote the feasible set of the dual problem for the recourse problem on $\mathcal{Y}^0$.  Our approximation scheme involves two parts. One part is that, for given $\mathbf x^*\in \mathcal{X}^0$, we solve subproblems with respect to $\mathbf x^*$ and $\mathcal U^0$, which results in $\eta_s(\mathbf x^*)$ and critical  $\pi^*\in \Pi^0$. Let $UB=\min\{UB,\mathbf c_1\mathbf x^*+\eta_s(\mathbf x^*)\}$.  Another one is to adopt $\pi^*$ and $\mathcal{U}(\mathbf x)$ to define sets  $\mathcal{OU}$ and $\mathcal{OV}$ and to build master problem $\mathbf{MP2}$. It is clear that the optimal value of $\mathbf{MP2}$ is less than or equal to $w(\mathcal X^1, \mathcal{U}(\mathbf x), \mathcal{Y}^1(\mathbf x, \mathbf u))$. We let $LB$ be this optimal value.  If $\mathbf{MP2}$'s optimal solution $\mathbf x^*\in \mathcal X^0$,  $\mathbf x^*$ can be directly supplied to compute subproblems. Otherwise, we can employ simple operations to convert $\mathbf x^*$ to a point in $\mathcal X^0$  to carry out Variant 2.  

Regarding the termination condition of this approximation scheme, there is no standard one, noting that the lower and upper bounds may never converge.  One termination condition is that computing $\mathbf{SP2}$ generates a $\boldsymbol\pi^*$ that has  been derived before. If this happens and all formulations have unique optimal solutions, it can be shown that the difference between $UB$ and $LB$ will not decrease anymore. Another one is that the algorithm terminates if the number of iterations reaches a pre-defined limit. In our computational experiments, in addition to the optimality tolerance,   the first one (along with the time limit) is adopted. The validity of these modifications and therefore the approximation scheme holds as follows.
  \begin{prop}
  	\label{prop_CCG_DDU_approximation}
  	For Variant 2 with the  aforementioned modifications, the values of $LB$ and $UB$ are valid lower and upper bounds to $\mathit{w}(\mathcal X, \mathcal U^0, \mathcal Y^0(\mathbf x,\mathbf u))$, respectively.  \hfill $\square$
  \end{prop}
  \begin{rem}
  	Note that the quality of the approximation, i.e., the (relative) optimality gap between $UB$ and $LB$, provides an empirical evaluation on the extent of correctness of our heuristic understanding. If not satisfactory, it can always be revised for a better quality. Actually, solutions to $\mathbf{SP1}$ and $\mathbf{SP2}$ provide a basis and are informative to  guide our revision.  
  \end{rem}

We next illustrate this approximation scheme and a revision to achieve a better quality. As shown in \citet{an2014reliable}, worst case disruptions could occur at non-facility sites if those sites' demands increase once disrupted, i.e.,  $\theta_i>0$ for some $i$. Nevertheless, we may subjectively believe that disruptions do not occur at non-facility sites regardless of $\theta_i$. This heuristic understanding, which  is reflected by $\mathcal{U}^k(\mathbf x)$ in \eqref{eq_DDU_reliable_pmedian},  certainly  can be used to derive approximate solutions according to Propositions \ref{prop_DIU_to_DDU} and \ref{prop_CCG_DDU_approximation}. Moreover, it is observed that an instance with a large optimality gap often has solutions to $\mathbf{SP2}$ that have worst case disruptions occurring at sites with large demands. Hence, we revise $\mathcal{U}^k(\mathbf x)$ with $D_q\subseteq I$, which is the set of  $q$ sites with largest demands, to $\mathcal{U}^{kq}(\mathbf x)$ as in the following.
 \begin{eqnarray}
 	\label{eq_DDU2_reliable_pmedian}
 	\begin{split}
 	\mathcal{U}^{kq}(\mathbf x)=\{\mathbf {u}\in \mathbb{R}^{|I|}_+: & \sum_{j\in J\cup{D}_q}u_j\leq k, \ u_j\leq x_{d,j} \ \forall j\in J\backslash D_q, \ u_i\leq 1 \ \forall i\in{D}_q,   \\
 	& u_i=0 \ \forall i\notin J\cup{D}_q\}
    \end{split}
 \end{eqnarray}
  Based on Propositions  \ref{prop_relax_1} and \ref{prop_DIU_to_DDU}, it is clear that
  	$$\mathit{w}_{R}(\mathcal{X}, \mathcal{U}^k(\mathbf x), \mathcal{Y}(\mathbf x, \mathbf u))\leq \mathit{w}_{R}(\mathcal{X}, \mathcal{U}^{kq}(\mathbf x), \mathcal{Y}(\mathbf x, \mathbf u))\leq  \mathit{w}_{R}(\mathcal{X}, \mathcal{U}^0, \mathcal{Y}(\mathbf x, \mathbf u)).$$
   Hence, the quality of approximation, especially $LB$, should be improved if  $\mathcal{U}^{kq}(\mathbf x)$ replaces $\mathcal{U}^k(\mathbf x)$. Actually in our numerical study with $\mathcal{U}^{kq}(\mathbf x)$, we derive not only optimal solutions for all instances solved by the DIU-based formulation, but also high quality approximation solutions with significantly reduced optimality gaps for all other unsolved ones.  
   
  \subsubsection{Constructing DDU Set with Deep Information}

   To construct $\mathcal{U}(\mathbf x)$ that reflects a sophisticated understanding, we often depend on hidden or complex information beyond what is represented by original variables in $\mathcal X$. This issue can be addressed by introducing auxiliary variables, including binary ones, and constraints to collect and represent that ``deep information'', provided that they do not affect set $\mathcal X$. Then,  those new variables can help us convert that sophisticated deep knowledge, either exact or heuristic, into decision dependence in $\mathcal{U} (\mathbf x)$.   
   
   For example, we may believe  the worst case disruptions should occur on facility sites of the highest service costs, which clearly requires facilities sorted according to their service costs. Nevertheless,  neither  $\mathbf x_c$ nor $\mathbf x_d$ carry the needed information. To this end, we introduce  a continuous variable $x^0_{r}$ and a set of binary variables $\mathbf x_r$ to sort out $q$ facility sites with the most service costs. Note that $x_{r,j}=1$ if site $j$ is one of those sites, $0$ otherwise. With the necessary information carried by $\mathbf x_{r}$, we can construct a DDU set such that disruptions  only occur among those $q$ facilities.  The updated sets $\mathcal X$ and $\mathcal U(\mathbf x)$ are:
   \begin{align}
   	\label{eq_app_X_sorting}
   	\begin{split}
   	\mathcal X^r=\{&(\mathbf x_c, \mathbf x_d, \mathbf x_r)\in \mathbb{R}^{|I|\times|J|+1}_+\times\{0,1\}^{2|J|}:  \sum_{j\in J}x_{d,j}=p, \ \sum_{j\in J}x_{c,ij}\geq d_i \ \forall i, \\ & \sum_{i\in I}x_{c,ij}\leq A_jx_{d,j} \ \forall j, \ x_{r,j}\leq x_{d,j} \ \forall j, \ \sum_{i\in I}c_{ij}x_{c,ij}\geq x^0_r-M(1-x_{r,j}) \ \forall j, \\ & \sum_{i\in I}c_{ij}x_{c,ij}\leq x^0_{r}+Mx_{r,j} \ \forall j, \ \sum_{j\in J}x_{r,j}=q\}
   	\end{split}
   	\end{align}
    \vspace{-10pt}
    \begin{align}
    \label{eq_app_DDU_sorting}   
   	\mathcal{U}^r(\mathbf x)&=\{\mathbf {u}\in \mathbb{R}^{|I|}_+: \sum_{j\in J}u_j\leq k, \ u_j\leq x_{r,j} \ \forall j, \  u_i=0 \ \forall i\notin J\}.
   \end{align}
    Note that new variables and constraints for sorting do not change the feasible set of the first stage decisions,  given that the projection of $\mathcal X^r$ onto the space hosting $\mathcal{X}$ is identical to $\mathcal{X}$.  It is easy to show that 
    $\mathit{w}_{R}(\mathcal{X}^r, \mathcal{U}^r(\mathbf x), \mathcal{Y}(\mathbf x, \mathbf u))\leq \mathit{w}_{R}(\mathcal{X}, \mathcal{U}^0, \mathcal{Y}(\mathbf x, \mathbf u)),$ ensuring the applicability of $\mathcal{U}^r(\mathbf x)$ for approximation. The logically stronger connection between $\mathcal X^r$ and  $\mathcal{U}^r(\mathbf x)$ might render the corresponding $\mathbf{2-Stg \ RO}$ easier to analyze. Indeed, for the extreme case where $q=k$ and there are no multiple facility sites with the same service cost, it can be further simplified to $\mathcal{U}^r(\mathbf x) = \{u_j=x_{r,j} \ \forall j,  u_i=0 \ \forall i\notin J\}$, directly converting this RO  into a computationally friendly deterministic model. Hence, it is reasonable to believe that $\mathcal{U}^r(\mathbf x)$ may bring a non-trivial computational advantage, e.g., fast computation with less number of iterations before termination. Nevertheless, because $x_{r,j}\leq x_{d,j}$, it follows that $\mathcal{U}^r(\mathbf x)\subseteq \mathcal{U}^k(\mathbf x)$. Hence, as showed in \eqref{eq_revision_fails}, the strength of  the approximation using $\mathcal{U}^r(\mathbf x)$ might not be strong, which is observed in our numerical study in Section \ref{subsect: DIUtoDDU}.     
   \begin{eqnarray}
    \label{eq_revision_fails}
    	\mathit{w}_{R}(\mathcal{X}^r, \mathcal{U}^r(\mathbf x), \mathcal{Y}(\mathbf x, \mathbf u))\leq \mathit{w}_{R}(\mathcal{X}, \mathcal{U}^k(\mathbf x), \mathcal{Y}(\mathbf x, \mathbf u))\leq  \mathit{w}_{R}(\mathcal{X}, \mathcal{U}^0, \mathcal{Y}(\mathbf x, \mathbf u))
    \end{eqnarray}
    Again, if the trade-off between the approximation quality and the computational efficiency is not satisfactory, revision on the decision dependence is desired.

    \begin{rem}
      $(i)$ Overall, the approximation scheme discussed  in this subsection is general in handling various applications, and is flexible in absorbing  different heuristic understanding.
      In particular, we mention that its implementation, including revision of $\mathcal U(\mathbf x)$ (and possibly of $\mathcal{X}$ and $\mathcal{Y}(\mathbf x, \mathbf u)$) according to  computational results and the reconstruction of master problem (MP) (and possibly of subproblems (SPs)), can be dynamic and iterative in a loop, which is depicted in the following diagram.  \vspace{-10pt}\\

    	\hspace{-15pt}\begin{tikzcd}[scale=0.6, column sep=scriptsize]
    		\vspace{0pt}\textrm{MP and SP reconstruction}\ar[r]  \ar[leftarrow,to path={ -- ([yshift=-5ex]\tikztostart.south) -| (\tikztotarget)},
    		rounded corners=12pt]{rr}& {\textrm{Computation \& Evaluation}}\ar[r]& {\textrm{Observation \& }  \textrm{Revision}}  
    	\end{tikzcd}\vspace{10pt} 
$(ii)$  By introducing  auxiliary (especially binary) variables and constraints, as demonstrated in $\mathcal{X}^r$, deep information can always be captured and represented. Thus, our primary task is to uncover deep knowledge regarding the connection between the first stage decision and worst case scenarios, and to identify necessary information for converting DIU to DDU. Note that strong modeling tools and techniques should be utilized to minimize the complexity and computational burden incurred by those auxiliary variables and constraints.  
      \\ $(iii)$ As suggested earlier, solutions to $\mathbf{SP1}$ and $\mathbf{SP2}$ are critical feedback for revising $\mathcal U(\mathbf x)$  (and possibly $\mathcal{X}$ and $\mathcal{Y}(\mathbf x, \mathbf u)$).
      Nevertheless, one central question is how to modify existing decision dependence or introduce new ones algebraically on the fly. So,  an interesting research direction is to design an intelligent subroutine within Variant 2 to  achieve automated revisions,  instead of depending on human decision maker's involvement. We believe that this direction may lead to powerful methodologies for  practical and large-scale applications.   \\
      $(iv)$ We are not restricted to employ a single DDU set to approximate $\mathcal {U}^0$.  For a general situation, it is unrealistic to precisely capture the connection between $\mathbf x$ and corresponding worst case scenarios in $\mathcal {U}^0$ by a simple DDU set. Hence, multiple DDU sets can be constructed to  build a stronger approximation.  Results in Proposition~\ref{prop_DIU_to_DDU} hold (with trivial changes) when multiple DDU sets are employed. One demonstration is presented in Section \ref{subsect: DIUtoDDU}. \\
      $(v)$ We also highlight that the basic ideas presented here are also applicable to converting the initial DDU set, denoted by $\mathcal U^0(\mathbf x)$, to a ``knowledge-richer'' set, $\mathcal U(\mathbf x)$, to achieve a better solution capacity. Note that Proposition \ref{prop_DIU_to_DDU} readily holds if $\mathcal U^0$ is replaced by $\mathcal U^0(\mathbf x)$. Hence, structural properties or strong heuristic understanding connecting $\mathbf x$ and the worst case scenarios in $\mathcal U^0(\mathbf x)$ are worth investigation.   As previously argued, this idea can be used  to develop fast exact or approximate algorithms for complex DDU-based RO.
    \end{rem}

\subsection{Dealing with Mixed Integer Recourse by Approximation}
\label{subsect_MIP}

Very often discrete decisions have to be modeled and taken into account in the second stage~\citep{bertsimas2010power,zhao2011exact,hanasusanto2015k}. Nevertheless, it is observed in the context of DIU-based $\mathbf{2-Stg \ RO}$  that MIP recourse problem brings a much greater technical challenge. In stead of relying on a nested implementation of basic C\&CG to derive exact solutions \citep{zhao2011exact}, approximation strategies, e.g., setting recourse decisions by a fixed decision rule or one of $K$ pre-defined policies \citep{bertsimas2010power,hanasusanto2015k}, have been investigated to compute practical instances. We mention an effective approximation strategy implemented with   basic C\&CG~\citep{zhang2018ambulance} that often generates solutions of very small optimality gaps within a reasonable time. In this subsection, we extend it in a way such that it handles DDU and derives solutions of a quantifiable quality for $\mathbf{2-Stg \ RO}$ with MIP recourse.

Similar to the previously discussed approximation scheme, this approximation approach  involves two ideas. The first one is to replace the MIP recourse problem by its LP relaxation. It yields a two-stage RO with an LP recourse problem, which helps develop a lower bound to $\mathbf{2-Stg \ RO}$. The second one is to fix discrete recourse decisions to particular values, which reduces the recourse problem to an LP (in a lower dimensional space). As the resulting 2-stage RO model is a restriction to the original $\mathbf{2-Stg \ RO}$, it yields an upper bound. Next, we describe operations to incorporate these ideas into Variant 2 assuming  that the relatively completely recourse property holds, i.e., the continuous portion of the recourse problem is feasible for any possible $(\mathbf x, \mathbf u, \mathbf y_d)$.
This property can be realized by penalizing the deficiency of recourse constraints, allowing us to skip \textbf{Step 3} and \textbf{Step 4 (Case B)}.

Recall that $\mathcal{Y}_r(\mathbf x,\mathbf u)$ denotes the LP relaxation of the feasible set of the original MIP recourse problem. As MIP  problem does not have the strong duality in general, we, with a slight abuse of notation, adopt $\Pi$ to denote the feasible set of the dual problem for this relaxation, and $\mathcal{P}_{\Pi}$ the set of associated extreme points. Also, we redefine $\mathbf{SP2}$ as
\begin{eqnarray}
	\label{eq_SP2_MIP}
	\mathbf{SP2}: \ \max_{\mathbf u\in \mathcal{U}(\mathbf x^*)} \min\{
	\mathbf c_2\mathbf y: \mathbf y\in \mathcal{Y}_r(\mathbf x^*, \mathbf u)\},
\end{eqnarray}
and introduce $\mathbf{SP4}$ with given $\mathbf x^*$ and $\mathbf y^*_d$ as
\begin{eqnarray}
	\label{eq_SP4}
	\mathbf{SP4}: \ \tilde \eta_s(\mathbf x^*, \mathbf y^*_d)=\max_{\mathbf u\in \mathcal{U}(\mathbf x^*)} \min\{
	\mathbf c_2\mathbf y: \mathbf y\in \mathcal{Y}(\mathbf x^*, \mathbf u), \mathbf y_d=\mathbf y^*_d\}.
\end{eqnarray}
Note that with $\mathbf y_d=\mathbf y^*_d$, the lower-level problem in $\mathbf{SP4}$ is an LP. Then, the original Variant~2 is  modified with the following changes.
\begin{itemize}
	\item In \textbf{Step 4 (Case A)}, instead of computing \eqref{eq_SP2}, we solve $\mathbf{SP2}$ in  \eqref{eq_SP2_MIP} to obtain optimal  $\mathbf u^*_s$ and corresponding $\boldsymbol\pi^*$, and perform remaining operations in this step. Note that the cutting set added to $\mathbf{MP2}$, i.e., \eqref{eq_CCG_optimality}, is defined with respect to new variables $\mathbf y^{\boldsymbol\pi^*}=(\mathbf y^{\boldsymbol\pi^*}_c, \mathbf y^{\boldsymbol\pi^*}_d)$. So, it is necessary to include $\mathbf y^{\boldsymbol\pi^*}_d\in \mathbb{Z}^{m_y}_+$ into \eqref{eq_CCG_optimality}.
	\item  In \textbf{Step 5}, $(i)$ compute the original MIP recourse problem for given $\mathbf x^*$ and $\mathbf u^*_s$ to obtain optimal $(\mathbf y^*_c,\mathbf y^*_d)$; $(ii)$ compute $\mathbf{SP4}$ in $\eqref{eq_SP4}$ with given $\mathbf x^*$ and $\mathbf y^*_d$ to obtain $\tilde \eta_s(\mathbf x^*, \mathbf y^*_d)$; $(iii)$ update $UB=\min \{UB, \mathbf c_1\mathbf x^*+ \tilde \eta_s(\mathbf x^*, \mathbf y^*_d)\}$.
\end{itemize}
Given that $LB$ is set to $\mathbf{MP2}$'s optimal value, the validity of this approximation scheme is rather clear.  
\begin{prop}
	\label{prop_CCG_MIP_approximation}
	For Variant 2 with the  aforementioned modifications, the values of $LB$ and $UB$, i.e., the optimal value of $\mathbf{MP2}$ and $\mathbf c_1\mathbf x^*+ \tilde \eta_s(\mathbf x^*, \mathbf y^*_d)$ respectively, are valid lower and upper bounds to the optimal value of $\mathbf{2-Stg \ RO}$ in \eqref{eq_2RO}.
\end{prop}
\begin{proof}
See its proof in Appendix \ref{apd:proofS4}.
\end{proof}

As for the termination of this approximation scheme, those conditions discussed before Proposition \ref{prop_CCG_DDU_approximation} are still applicable.

\begin{rem}
	$(i)$  Unlike other existing approximation strategies, we note this approximation scheme largely reserves the full recourse capacity. Also, it progressively updates lower and upper bounds, rendering a quantitative measure on the strength of the approximation. According to the numerical studies presented in Section \ref{subsect_computation_MIPrecourse} and those in~\citet{zhang2018ambulance}, this approximation scheme often achieves a highly desirable balance between solution quality and computational complexity, noting that  optimality gaps are generally small at termination. Indeed, if the MIP recourse problem has no integrality gap for    $(\mathbf x, \mathbf u) \in (\mathcal{X}, \mathcal U(\mathbf x))$, it becomes exact without optimality gap.   \\
	$(ii)$ This approximation idea is rather general. It can be easily incorporated into the other two variants so that they can also handle MIP recourse problem too. Nevertheless, as Variant 1 generates Benders cutting sets based on LP's strong duality, it entirely relies on the MIP recourse problem's LP relaxation to build $\mathbf{MP1}$ and to derive lower bounds.  Hence, the best lower bound cannot be better than the one obtained from directly computing a two-stage RO model with that LP relaxation as its recourse problem. 
\end{rem}

\subsection{Addressing Uniqueness and Pareto Optimality}
\label{subsect_Uniqueness_Pareto}
In Section \ref{subsect_Variant2}, it is shown that the modified Variant 2, based on its \textit{uniqueness} property,  enjoys a theoretically stronger iteration complexity than the standard one.  Nevertheless, we note the implementation challenges with current professional MIP solvers. In this subsection, we show that, by making use of specific structures of the underlying model, that uniqueness property can be achieved easily with trivial extra work. Another consideration for computational improvement is Pareto optimality in generating or selecting cutting sets. As shown in \citet{magnanti1981accelerating} in the context of Benders decomposition, cutting sets with Pareto optimality often have a significantly better performance. In this subsection, we also develop extensions of  Variant 2 to have this property.

\subsubsection{Utilizing Structural Information to Attain Uniqueness}
Instead of depending on MIP solvers to provide optimal bases and related information, it is often possible to analyze  specific structures of one application and then implement simple modifications to achieve the uniqueness property. We next present an illustration using the reliable p-median problem formulated in \eqref{eq_reliable_pmedian_DDU} with $\theta_i=0$ for $i\in I$ and $\mathcal{U}^k(\mathbf x)$ defined in~$\eqref{eq_DDU_reliable_pmedian}$.

Suppose that $\mathbf x^*$ is output by solving $\mathbf{MP2}$ in some iteration. As the recourse problem is always feasible, we can simply compute $\mathbf{SP2}$ in its bilinear form
\begin{align*}
	\max \quad  & \sum_{j\in J} A_j\left[(u_j-1)\pi_{3,j}-x_{d,j}^*\pi_{2,j}\right]+\sum_{i\in I}d_i\pi_{1,i}\\
	\mathrm{s.t.}\quad & \pi_{1,i}-\pi_{3,j}-\pi_{2,j}\leq c_{ij} \ \ \forall i,j,  \ \ \pi_{1,i}\leq \textsl{C} \ \ \forall i \\
	& \mathbf u\in\mathcal U^k(\mathbf x^*), \ \ \pi_{1,i}\geq 0 \ \ \forall i, \ \ \pi_{2,j}\geq 0 \ \ \forall j, \ \pi_{3,j}\geq 0 \ \ \forall j
\end{align*}
and obtain optimal $\mathbf u^*_s$ and $\boldsymbol\pi^*=(\boldsymbol{\pi}_{1}^*, \boldsymbol{\pi}_{2}^*, \pi_{3}^*)$. Hence, for $\boldsymbol\pi^*$, its associated set $\mathcal{OU}$ consists of optimal solutions to the following linear program with $\mathfrak c^{\pi^*}_j= A_j\pi^*_{3,j}$ for $j\in J$.
\begin{align}
\label{eq_OU_LP_pmedian}	
\max\large\{\sum_{j\in J} \mathfrak c^{\boldsymbol\pi^*}_j u_j: \sum_{j\in J} u_j\leq k, \ u_j\leq x_{d,j} \ \forall j, \ u_j\geq 0 \ \forall j, \ u_i=0 \ \forall i\notin J\large\}
\end{align}

Obviously, in an optimal solution, we have $u_j=1$ if $x_{d,j}=1$ and its coefficient $\mathfrak c^{\boldsymbol\pi^*}_j$ is among the largest $k$ values, and remaining $u_j$'s set to $0$. Nevertheless, multiple $\mathfrak c^{\boldsymbol\pi^*}_j$'s could be same, especially for the uncapacitated model where $A_j$ is set to $\sum_{i}d_i$ for $j\in J$. If this is the case, the aforementioned linear program is likely to have multiple optimal solutions, i.e., $\mathcal{OU}(\mathbf x, \boldsymbol\pi^*)$ is not a singleton.  Instead of performing operations described in Appendix
\ref{Asect_modiunique} to implement the modified Variant 2 , we can render $\mathbf u^*_s$ the unique optimal solution of \eqref{eq_OU_LP_pmedian} for $\mathbf x^*$  by the following simple manipulations. 

Given that $\mathfrak c^{\boldsymbol\pi^*}_j\geq 0$ for $j\in J$, it is without loss of generality to assume that $u^*_{s,j_1}=\dots=u^*_{s,j_k}=1$. As stated after Proposition \ref{prop_complexity_uniqueness}, the uniqueness property can be achieved by modifying the objective function coefficients of \eqref{eq_OU_LP_pmedian}. Specifically, let $\overline{\mathfrak c}^{\boldsymbol\pi^*}$ denote the updated objective function coefficients. After computing facilities' service costs, we set $\hat{\mathfrak c}^{\boldsymbol\pi^*}_j= \sum_{i}c_{ij}x^*_{c,ij}$ for $j\in J\backslash\{j_1,\dots, j_k\}$ and
$\hat{\mathfrak c}^{\boldsymbol\pi^*}_j= \max_{j\in J}\{\sum_{i}c_{ij}x^*_{c,ij}\}+\tau$ for $j\in \{j_1,\dots, j_k\}$, with $\tau>0$. Because of the structure of $\hat{\mathfrak c}^{\boldsymbol\pi^*}$, the  next result simply follows.
\begin{cor}	\label{cor_direct_manu_unique}
	The linear program $\max\{\sum_{j\in J} \hat{\mathfrak c}^{\boldsymbol\pi^*}_j u_j: \ \mathbf u\in \mathcal{U}(\mathbf x^*)\}$ has a unique optimal solution, i.e., its optimal solution set, denoted by $\widehat{\mathcal{OU}}(\mathbf x^*, \boldsymbol\pi^*)$, is $\{\mathbf u^*_s\}$. Moreover, as long as $\mathbf u^*_s\in \mathcal U(\mathbf x)$ for some $\mathbf x\in \mathcal{X}$, we have
	$\widehat{\mathcal{OU}}(\mathbf x, \boldsymbol\pi^*)=\{\mathbf u^*_s\}$.
\end{cor}
We are now ready to implement the modified Variant 2 to compute this reliable $p$-median problem. It is interesting to mention that $\mathbf u^*_s$ is a highly degenerate extreme point of $\mathcal U(\mathbf x^*)$. That means the iteration complexity based on the number of bases, which is $O(\binom{2|J|+1}{|J|+1})$ according to Proposition \ref{prop_complexity_mCCG_0}, is an overestimation. Indeed, because of the structure of $\mathcal{U}(\mathbf x)$, that iteration complexity can be strengthened as in the following. 

 \begin{cor}
 	The number of iterations of the modified Variant 2 to converge is bounded by $O(\binom{|J|}{k})$. \end{cor}

\begin{rem}
	$(i)$ As mentioned, those simple manipulations guarantee the uniqueness property and therefore support implementation of the modified Variant 2. As  neither special outputs from an MIP solver are required nor computationally-heavy operations are involved, we believe that this strategy, i.e., structure-based manipulation, is practically useful and can be applied to many real problems on which we have deep insights.\\
	$(ii)$ To ensure the uniqueness property, it is sufficient to   set $\hat{\mathfrak c}^{\boldsymbol\pi^*}_j$ to a positive number  for $j\in \{j_1,\dots, j_k\}$ only and set others to zeros. Nevertheless, once $x_j=0$ for some $j\in \{j_1,\dots,j_k\}$, the cutting set based on $\widehat{\mathcal{OU}}(\mathbf x, \boldsymbol\pi^*)$ becomes weaker. In the extreme case where $x_j=0$ for $j\in \{j_1,\dots,j_k\}$, that cutting set is actually useless. Note that if this is the case,  $\widehat{\mathcal{OU}}(\mathbf x, \boldsymbol\pi^*)$ has many optimal solutions, including the one without any disruptions such that $u_j=0$ for all $j$. Then, $\mathbf{MP2}$, to minimize its objective function value,  simply selects that solution from $\widehat{\mathcal{OU}}(\mathbf x, \boldsymbol\pi^*)$, which renders the associated cutting set no different from the constraint set of the normal disruption-free situation.  Hence, we would suggest to incorporate more structural information or insights, either exact or heuristic, when setting $\hat{\mathfrak c}^{\boldsymbol\pi^*}$ to have the uniqueness property. As shown in our illustration, the reflected understanding is that facility sites of higher service costs are more likely to be disrupted in worst case scenarios, in addition to achieving the uniqueness.    
\end{rem}

 \subsubsection{Obtaining Pareto Optimal Cutting Sets}
 Using a strategy similar to that of \citet{magnanti1981accelerating} and \citet{papadakos2008practical} for Benders decomposition, we can design and compute subproblems with Pareto optimality consideration to introduce strong cutting sets in Variant 2. To present the main ideas, we focus on deriving optimality cutting sets with that consideration. Feasibility counterparts can be derived by implementing similar operations.

   Assume that we have $\mathbf x^0\neq \mathbf x^*$ for which $\mathcal U(\mathbf x^0)$ is non-empty, given that
   $\mathbf x^*$ is output from $\mathbf{MP2}$ in some iteration. Also, assume $\mathbf{SP2}$ has been solved with $\eta_s(\mathbf x^*)$, $\mathbf{u}^*_s$ and $\boldsymbol\pi^*$.  To search for optimality cutting set that is Pareto optimal (with respect to $\mathbf x^0$ and its $\mathcal{U}(\mathbf x^0)$), we additionally solve the following subproblem.
   \begin{equation*}
   	 \begin{split}
   	\mathbf{SP2}_{PO}: \  \eta_{s}(\mathbf x^0)=\max\quad & (\mathbf d-\mathbf B_1\mathbf x^0-\mathbf u)^\intercal\boldsymbol\pi \label{eq_PO_obj}\\
   	\mathrm{s.t.} \quad& \mathbf u\in \mathcal{U}(\mathbf x^0), \ \mathbf B^\intercal_2\boldsymbol\pi\leq \mathbf{c}^\intercal_2, \ \boldsymbol\pi\geq \mathbf 0 \\
   	& \tilde{\mathbf u}\in \mathcal{U}(\mathbf x^*), \ (\mathbf d-\mathbf B_1\mathbf x^*-\tilde{\mathbf u})^\intercal\boldsymbol\pi\geq \eta_s(\mathbf x^*)
   	\end{split}
   \end{equation*}

    Clearly, $\mathbf{SP2}_{PO}$ is a challenging  non-convex quadratically constrained model. Once an optimal solution, denoted by $(\dot{\mathbf u}_s, \dot{\boldsymbol\pi})_s$, is derived, we can employ $\dot{\boldsymbol\pi}$ to create a Pareto optimal cutting set. Note that we should still use $\eta_s(\mathbf x^*)$ to update the algorithm's upper bound. 

    $\mathbf{SP2}_{PO}$  is much more complicated than the one used in Benders decomposition. Alternatively, instead of computing this complex formulation exactly,
     $\mathbf{SP2}_{PO}$ can be replaced by some simpler formulations for fast computation. One is to fix $\tilde{\mathbf{u}} = \mathbf u^*_s$ and penalize the deficiency to $\eta_s(\mathbf x^*)$,  which  reduces  $\mathbf{SP2}_{PO}$ to the following bilevel linear program.
   \begin{align*}
   \mathbf{SP2}_{POB}:  \ &  \max  -\textsl{C}\tilde\tau + \mathbf c_2\mathbf y \\
    \mbox{s.t.}  \  & \mathbf u\in \mathcal{U}(\mathbf x^0), \  (\mathbf d-\mathbf B_1\mathbf x^*-\mathbf u^*_s)^\intercal\boldsymbol\pi+\tilde\tau\geq \eta_s(\mathbf x^*) \\ 
    & (\mathbf y,\boldsymbol\pi)\in \arg\min\{\mathbf c_2\mathbf y: \mathbf y\in \mathcal{Y}(\mathbf x^0, \mathbf u)\} 
   \end{align*}
    The  second constraint is a little bit of abuse of notation, noting that optimal $\boldsymbol\pi$ is always obtainable if the minimization problem is replaced by KKT or strong duality based optimality conditions. Clearly, if $\tilde \tau=0$ in its optimal solution, the derived $\boldsymbol{\pi}$ is also optimal to $\mathbf{SP2}_{PO}$. 
    Indeed, we can further simplify $\mathbf{SP2}_{PO}$ to the next LP by fixing  $\mathbf u=\mathbf u'$ that is different from $\mathbf u^*_s$. For the sake of simplicity, it is not necessary to have $\mathbf u'\in \mathcal{U}(\mathbf x^0)$.  Actually,  solving  it yields a solution that is better than $\boldsymbol\pi^*$, if not same,  but suboptimal to $\mathbf{SP2}_{PO}$. 
    \begin{eqnarray}
    \label{eq_Pareto_LP}
    	\mathbf{SP2}_{POL}:  \ \max \{ (\mathbf d-\mathbf B_1\mathbf x^0-\mathbf u')^\intercal\boldsymbol\pi: \mathbf B^\intercal_2\boldsymbol\pi\leq \mathbf{c}^\intercal_2, \ \boldsymbol\pi\geq \mathbf 0, \ (\mathbf d-\mathbf B_1\mathbf x^*-\mathbf u^*_s)^\intercal\boldsymbol\pi\geq \eta_s(\mathbf x^*)\}
    \end{eqnarray}
    Note that the modified Variant 2 (which achieves the uniqueness property), especially the structural information based implementation, can be applied in conjunction with Pareto optimal cutting sets. We actually believe that the uniqueness property and Pareto optimality are different but connected. More systematic studies are worth doing in future research.

\section{Numerical Studies}
\label{sect_numerical}
In this section, we report and discuss numerical results obtained from computational experiments. Our focus is on testing, evaluating and analyzing three C\&CG variants developed in Section \ref{sect_algorithms}, as well as their modified versions and extensions described in Section \ref{sect_more}. We do not consider reformulation techniques presented in Section \ref{subsect_specialDDU} as they just convert specially structured DDU-based models to those that are computable by existing  methods.

We employ  two types of robust facility location problems as the test bed, given that they have been popular applications or demonstration platforms for many  optimization models and algorithms. One is the reliable facility location problem described in Section~\ref{subsect_DIU-DDU} that considers uncertain site disruptions. Another type is the robust facility location problem with uncertain demands. For both types, the first stage decision is to determine the locations  of service facilities. We also need to decide the capacities of established facilities if they are not uncapacitated. Then, once the random scenario is revealed, we make recourse decisions that serve demands by (available) facilities according to their capacities.  The overall objective is to minimize the total cost across both stages. Data regarding sites' locations, distances and  basic demands are adopted from \citet{snyder2005reliability}, and parameter $M$ is set to 10,000. All solution methods are implemented by \texttt{Julia} with
\texttt{JuMP} and professional MIP solver \texttt{Gurobi~9.1} on a Windows PC with E5-1620 CPU and 32G RAM. Unless noted otherwise, the relative optimality tolerances of any algorithm  and the solver are set to $.1\%$ and $.01\%$ (the default value), respectively. Also, the time limit is set to 3,600 seconds for any instance.

\subsection{Computational Performances of Three Variants}
\label{subsect_compu_3ccg}
As noted earlier, all three C\&CG variants developed in Section \ref{sect_algorithms} are exact algorithms with different theoretical complexities. It is critical to understand their practical performances in computing actual instances. Also, we would like to recognize their strength and limitations, and, if applicable, identify the most efficient one. In this subsection, we adopt the robust facility location problem with two types of DDU demands to perform our experiments.

\subsubsection{Robust Facility Location Model and DDU Sets}
\label{subsubsect:RobustFL}

Recall that $I$ and $J\subseteq I$ are sets of client sites and potential facility sites. Parameters $f_j$ and $a_j$ are the fixed cost of building a facility at $j$ and the unit capacity cost at $j$, respectively. Moreover, $c_{ij}$ is the unit service cost incurred by serving client $i$'s demand by the facility at $j$, and $p_i$ is the unit profit received after satisfying client $i$'s demand. The complete formulation is as the following.
\begin{equation}
	\label{eq_robust_facility}
	\begin{split}
	\mathbf{Robust \ FL}: \	\mathit{w}_{D}(\mathcal{X}, \mathcal{U}(\mathbf x), \mathcal{Y}(\mathbf x, \mathbf u))=\min_{(\mathbf x_c,\mathbf x_d)\in\mathcal X} & \sum_{j\in J} (f_jx_{d,j}+a_jx_{c,j})+ \\
		& \max_{\mathbf u\in\mathcal U(\mathbf x)}\min_{\mathbf y\in\mathcal{Y}(\mathbf x,\mathbf u)} \ \sum_{i\in I}\sum_{j\in J}(c_{ij}-p_i)y_{ij}
	\end{split}
\end{equation}
with sets
\begin{align}
\label{eq_RFL_X}
	\mathcal X=\{(\mathbf x_c, \mathbf x_d)\in \mathbb{R}^{|J|}_+\times\{0,1\}^{|J|}: \underline A_jx_{d,j}\leq x_{c,j}\leq \overline A_jx_{d,j} \ \forall j\in J\},
\end{align}
and
\begin{align}
	\mathcal Y(\mathbf x,\mathbf u)=\{\mathbf y\in \mathbb{R}^{|I||J|}_+: & \sum_{j\in J}y_{ij}\geq u_i \ \forall i, \  \sum_{i\in I}y_{ij}\leq x_{c,j} \ \forall j\}.
\end{align}
In $\mathbf{Robust \ FL}$, $x_{d,j}$ and $x_{c,j}$ are binary and continuous variables representing 0-1 construction and capacity decisions on site $j$, respectively, and $y_{ij}$ is service allocation decision between the client at $i$ and the  facility on $j$.  
As argued in Section \ref{sect_algorithms}, the developed algorithms may have different behaviors or complexities depending on the underlying DDU's RHS or LHS dependence. So, regarding $\mathcal U(\mathbf x)$ that captures uncertain demands, we consider the following two types DDU sets. Note that those DDU sets are designed rather for the demonstration~purpose.\\

\noindent\textbf{$(i)$ DDU with RHS Dependence}\\
Assume that prior to the constructions of facilities, there is a constant basic demand, denoted by $\underline u_i$, for the client at site $i$. Let $ u^0=\sum_{i\in I}\underline u_i$. After facilities are established, it is expected that they will induce additional demand $\tilde u_i$ (expressed in terms of $\underline u_i$). The complete DDU demand set with RHS dependence is
\begin{subequations}
\label{ex_eq_DDU_RHS}
\begin{align}
	\mathcal U^R(\mathbf x)=\{(\mathbf u, \mathbf d): & \ u_i=\underline u_i(1+\tilde u_i) \ \ \forall i \label{eq_DDU_RHS0}\\
	& \sum_{j\in J(i)}\underline{\xi}^i_jx_{d,j}\leq \tilde u_i\leq \sum_{j\in  J(i)}\overline{\xi}^i_jx_{d,j} \ \ \forall i \label{eq_DDU_RHS1}\\
	& |\tilde u_i-\frac{1}{2}\sum_{j\in J(i)}(\underline{\xi}^i_j+\overline{\xi}^i_j)x_{d,j}|\leq d_i \ \ \forall i \label{eq_DDU_RHS2}\\
	&\sum_{i\in I}d_i\leq \alpha \frac{\sum_{j\in J} x_{c,j}}{ u^0} \label{eq_DDU_RHS3}
	\}.
\end{align}
\end{subequations}
As in \eqref{eq_DDU_RHS0}, the actual demand from the client at site $i$, denoted by $u_i$,  includes both $\underline u_i$ and uncertain  (relative) increase
$\tilde u_i$. The latter one is influenced by the construction of facilities in its neighborhood $J(i)$  as shown in \eqref{eq_DDU_RHS1}. Inequalities in \eqref{eq_DDU_RHS2} measure the deviation, denoted by $d_i$, of that uncertain increase from the nominal one for each site $i$, which can be linearized easily as $-d_i\leq \tilde u_i-\frac{1}{2}\sum_{j\in J(i)}(\underline{\xi}^i_j+\overline{\xi}^i_j)x_{d,j}\leq d_i \ \forall i$. 
The overall deviation is bounded by a multiple of the quotient between the total capacity and $ u^0$ as in~\eqref{eq_DDU_RHS3}.  \newline

\noindent\textbf{$(ii)$ DDU with LHS (and RHS) Dependence} \\
Similar to \eqref{ex_eq_DDU_RHS}, we assume that site $i$'s actual demand again consists two portions. Nevertheless, they are both influenced by the first stage decisions, showing LHS (and RHS) dependence.  Indeed, they may demonstrate very different relationships with regard to the installed capacities. The complete DDU demand set is
\begin{subequations}
	\begin{align}
		\mathcal U^{L\!R}(\mathbf x)=\{\mathbf u:\quad&u_i=\hat u_i+\tilde u_i \quad\forall i\label{FLDDDL_TolDem}\\
		&\hat u_i\leq \underline u_i+k_1\sum_{j\in J(i)}x_{c,j} \quad\forall i\label{FLDDDL_RegDem}\\
		&\tilde u_i\leq\gamma\underline u_i  \quad\forall i\label{FLDDDL_SigUp}\\
		&\sum_{i\in I}\tilde u_i\left(\underline u_i+k_2\sum_{j\in J(i)}x_{c,j}\right)\leq\gamma\sum_{j\in J}\underline u_i^2\label{FLDDDL_TolUp}\\
		&\hat u_i\geq 0, \ \tilde u_i\geq 0 \ \forall i\}.
	\end{align}
\end{subequations}
As in \eqref{FLDDDL_TolDem}, $u_i$ has two portions, i.e., $\hat u_i$ and $\tilde u_i$. Note that  $\hat u_i$ is rather predictable, which is positively related to and bounded by a multiple of the installed capacities in its neighborhood, in addition to basic demand $\underline u_i$. The other one, $\tilde u_i$, is bounded by a multiple of the basic demand as  in \eqref{FLDDDL_SigUp}. Yet, constraint \eqref{FLDDDL_TolUp}  reflects a system-level negative correlation between $\tilde u_i$s and the installed capacities.

\subsubsection{Computational Results and Analyses}
\begin{figure}[htp]
\hspace{15pt}\begin{subfigure}[t]{0.45\textwidth}
\begin{tikzpicture}
	\pgfplotsset{tick label style={font=\small},
		label style={font=\small},
		legend style={font=\tiny}
	}
\begin{axis} [height=6cm,width=7.5cm,
    legend style={at={(1.1,1.2)},anchor=north,legend columns=-1,/tikz/every even column/.append style={column sep=5pt}},
    ybar=0.5pt,
    bar width = 6pt,
    ymin = 0,
    ymax = 20,
    yticklabel={\pgfmathparse{\tick}\pgfmathprintnumber{\pgfmathresult}\%},
    symbolic x coords={25L,25H,40L,40H},
    xtick = data,
    enlarge x limits = 0.2,
    nodes near coords={\pgfmathprintnumber\pgfplotspointmeta\%},
    every node near coord/.append style={rotate=90, anchor=west, font=\tiny}
]
\pgfkeys{/pgf/number format/.cd,fixed,precision=2,};
\addplot[draw = blue!50, fill=blue!50] coordinates {(25H,0) (25L,0) (40H,3.39) (40L,15.01)};
\addplot[draw = red!50, fill=red!50] coordinates {(25H,0.03) (25L,0.03) (40H,0) (40L,0)};
\addplot[draw = green!50, fill=green!50] coordinates {(25H,0) (25L,0.02) (40H,0.02) (40L,0.01)};
\legend {Benders C\&CG, Parametric C\&CG, Basis Based C\&CG};
\end{axis}
\end{tikzpicture}
\caption{Relative Gap \label{fig:alg_compa_RHS_gap}}
\end{subfigure}
\begin{subfigure}[t]{0.45\textwidth}
\center
\begin{tikzpicture}
	\pgfplotsset{tick label style={font=\small},
		label style={font=\small},
		legend style={font=\tiny}
	}
\begin{axis} [height=6cm,width=7.5cm,
    ybar=0.5pt,
    bar width = 6pt,
    ymin = 0,
    ymax = 2100,
    symbolic x coords={25L,25H,40L,40H},
    xtick = data,
    enlarge x limits = 0.2,
    nodes near coords,
    every node near coord/.append style={rotate=90, anchor=west, font=\tiny}
]
\pgfkeys{/pgf/number format/.cd,fixed,precision=2,set thousands separator={}};
\addplot[draw = blue!50, fill=blue!50, point meta=explicit symbolic]
table[meta=label]{
	x		y		label
	25H		25		{25 (8)}
	25L		41.31	{41.31 (8.6)}
	40H		665.49	{665.49 (14.33)}
	40L		1299.35	{1299.35 (18)}
};
\addplot[draw = red!50, fill=red!50, point meta=explicit symbolic]
table[meta=label]{
	x		y		label
	25H		5.42	{5.42 (3.4)}
	25L		7.33	{7.33 (3.4)}
	40H		10.4	{10.4 (3)}
	40L		25.43	{25.43 (3.6)}
};
\addplot[draw = green!50, fill=green!50, point meta=explicit symbolic]
table[meta=label]{
	x		y		label
	25H		10.24	{10.24 (3.6)}
	25L		9.92	{9.92 (3.6)}
	40H		385.78	{385.78 (4)}
	40L		198.64	{198.64 (4.2)}
};
\end{axis}
\end{tikzpicture}
\caption{Computation Time (s) \label{fig:alg_compa_RHS_time}}
\end{subfigure}
\caption{Computational Results of Algorithms for $\mathbf{Robust \ FL}$ with $\mathcal U^{R}(\mathbf x)$ \label{fig:alg_compa_RHS}}
\end{figure}
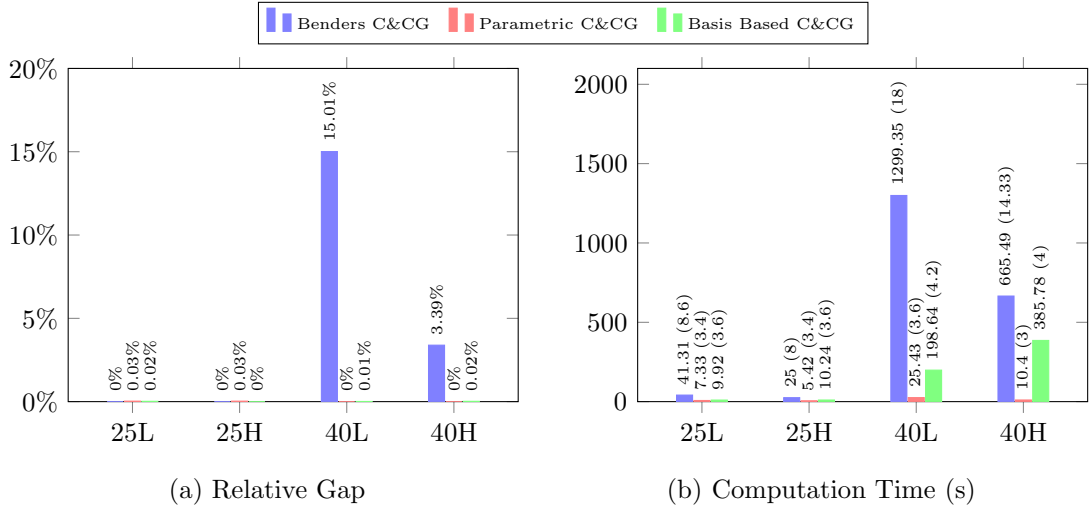

By setting $\mathcal{U}(x)=\mathcal U^R(\mathbf x)$, we solve $\mathbf{Robust \ FL}$ by three algorithms developed in Section~\ref{sect_algorithms}: Benders C\&CG (i.e., Variant 1), parametric C\&CG (i.e., Variant 2), and basis based C\&CG (i.e., Variant 3). We consider two groups of instances for 25 and 40 sites respectively, resulting in 4 groups altogether.  Instances are derived by modifying parameters of $\mathcal U^R(\mathbf x)$ according to values specified in Table \ref{tbl_FLR25-40-RHS}  presented in Appendix \ref{Asect_tables}, which also contains all detailed computational results.  We summarize and compare three C\&CG variants' performances for every group in Figure \ref{fig:alg_compa_RHS}. Note that $25$ or $40$ indicates the number of sites for instances in that group. $L$ and $H$ are used to denote instances with low and high fixed costs, respectively, where $H$ instances have fixed costs one third larger than those of $L$ instances. The bar chart in Figure \ref{fig:alg_compa_RHS}.(a) shows the average relative gaps upon termination, where the numerical values are also displayed on top of bars. Figure \ref{fig:alg_compa_RHS}.(b) presents the average computational time in seconds across instances solved to optimality for each group. The numerical values are displayed on top of bars, along with the average numbers of iterations displayed within parentheses. We mention that such arrangements basically hold for all figures presented in this section, unless otherwise noted. 

Based on those numerical results, we have the following two observations.
\begin{enumerate}
	\item Among all three algorithms, parametric C\&CG is undoubtedly the strongest one. It solves all instances in several iterations with no more than 40 seconds computational time. For basis based C\&CG, although the number of iterations is generally small and comparable to that of parametric C\&CG, its computational time is significantly longer. For some instance, the computational time could be up to two orders of magnitude more than that of parametric C\&CG. Benders C\&CG obviously has the worst performance, in both  the number of iterations and the computational time. For large instances, it often fails to close the gap between lower and upper bounds,  resulting in solutions with poor qualities. Even if it generates an optimal solution, its  solution time could be hundreds of times more than that of parametric C\&CG.
  \item Both parametric and basis based C\&CGs demonstrate a great scalability in terms of the number of iterations. Especially for parametric one, the number of iterations basically does not change when the instance's size increases from 25 to 40 sites. Different from them, Benders C\&CG computes many more iterations when the size of instance increases. Because of the nature of C\&CG, all those three algorithms introduce many new variables and constraints in each iteration, resulting in ever-increasing master problems. Hence, if many iterations are involved before termination, the most recent master problems are  very
  substantial and may need extremely long computational time. This is clearly demonstrated in Benders C\&CG's results, if we compare its performance between instances of 25 and 40 sites. Note that the increase in the computational time is unproportionally larger than that in the number of iterations. Therefore,  the aforementioned scalability ensures that parametric C\&CG maintains a strong solution capacity towards complex instances.
\end{enumerate}

\begin{figure}[htp]
	\begin{subfigure}[t]{0.47\textwidth}
		\pgfplotstableread{Sp2.dat}{\Sp}
		\pgfplotsset{tick label style={font=\small\bfseries},
			label style={font=\small},
			legend style={font=\tiny}
		}
		\begin{tikzpicture}[scale=1]
			\begin{axis}[ymajorgrids=true,legend style={at={(1.1,1.3)},anchor=north,legend columns=-1,/tikz/every even column/.append style={column sep=5pt}},height=5cm,width=7.5cm,
				xmin=0,xmax=16,ymin=0,ymax=6000,
				xtick={0,2,4,6,8,10,12,14,16},
				ytick={0,1000,2000,3000,4000,5000,6000},
				xlabel= Iterations]
				\addplot [blue!50,very thick] table [x={BeIt}, y={BeB}] {\Sp};
				\addlegendentry{Benders C\&CG}
				\addplot [red!50,very thick] table [x={PaIt}, y={PaB}] {\Sp};
				\addlegendentry{Parametric C\&CG}
				\addplot [green!50,very thick] table [x={BaIt}, y={BaB}] {\Sp};
				\addlegendentry{Basis Based C\&CG}
			\end{axis}
		\end{tikzpicture}
		\caption{Convergence over Iterations}
		\label{ConIter}
	\end{subfigure}
	\begin{subfigure}[t]{0.47\textwidth}
		\pgfplotstableread{Sp1.dat}{\Sp}
		\pgfplotsset{tick label style={font=\small\bfseries},
			label style={font=\small},
			legend style={font=\tiny}
		}
		\begin{tikzpicture}[scale=1]
			\begin{axis}[ymajorgrids=true,height=5cm,width=7.5cm,
				xmin=0,xmax=1350,ymin=0,ymax=6000,
				xtick={0,200,400,600,800,1000,1200},
				ytick={0,1000,2000,3000,4000,5000,6000},
				xlabel= Time (s)]
				\addplot [blue!50,very thick] table [x={BeTi}, y={BeB}] {\Sp};
				\addplot [red!50,very thick] table [x={PaTi}, y={PaB}] {\Sp};
				\addplot [green!50,very thick] table [x={BaTi}, y={BaB}] {\Sp};
			\end{axis}
		\end{tikzpicture}
		\caption{Convergence over Time}
		\label{ConTime}
	\end{subfigure}
	\caption{Convergence of Algorithms for $\mathbf{Robust \ FL}$ with $\mathcal U^{R}(\mathbf x)$}
	\label{ConFL}
\end{figure}
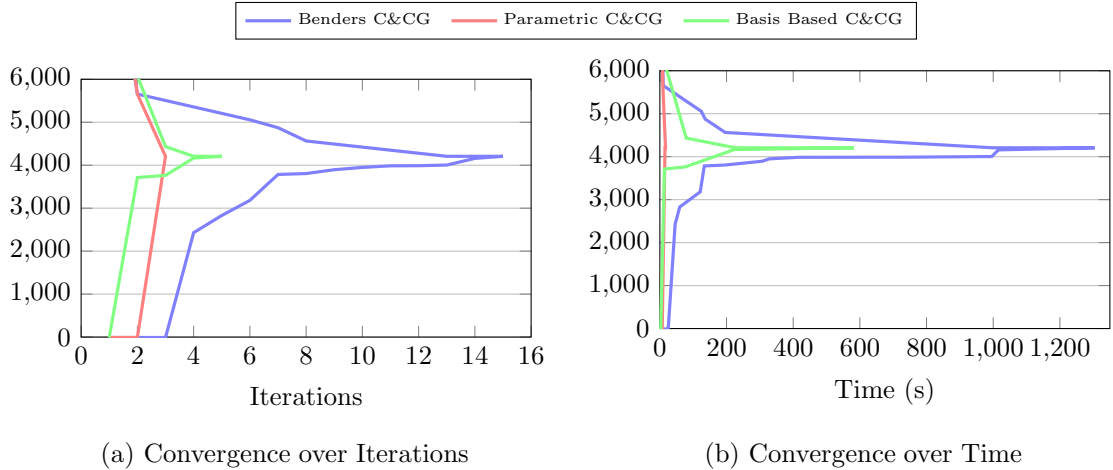

To further understand those algorithms' computational behaviors,  we plot their convergent progresses over iterations/time for a particular instance in Figure \ref{ConFL}. The instance is with 40 sites, high fixed cost, $\boldsymbol{\underline \xi}=0.08$, $\boldsymbol{\overline \xi} = 0.12$, and $\alpha=0.08$. It is clear that parametric C\&CG, i.e., Variant 2, quickly improves lower and upper bounds to reduce their gap. It actually converges in just a couple of iterations. On the contrary,  lower and upper bound curves of Benders C\&CG increase and decrease in a much slower pace, resulting in a much larger number of iterations to converge. Correspondingly, this algorithm takes a very long time to derive an optimal solution. As for basis based C\&CG, it has a performance between those of the other two algorithms. Nevertheless, due to the heavy burden associated with computing $\mathbf{MP3}$,  its progress over time lags behind that over iterations. We mention an interesting observation that lower bound from Variant 2 (also from Variant 1) does not improve until the second iteration (or the third iteration for Variant 1), while lower bound of Variant 3 simply increases after a cutting set is introduced and then is dominated by that of Variant 2 in following iterations. On the one hand, such an observation indicates that basis based cutting sets are effective, but less flexible in handling the changing first stage decision. On the other hand, parametric cutting sets become particularly strong whenever a non-trivial number of them are generated.

Also, we investigate those algorithms in handling DDU with LHS (and RHS) dependence by computing instances of $\mathbf{Robust \ FL}$ with $\mathcal{U}(x)=\mathcal U^{L\!R}(\mathbf x)$. Because of the huge computational burden, we employ a rather small set of instances of 25, 30, 35 and 40 sites to perform experiments. The complete results are reported in Table \ref{tbl_FLLR} in Appendix \ref{Asect_tables}, and the average performances are summarized in Figure \ref{Fig:alg_compa_LRH}. As none of instances is solved to optimality by Benders or basis based C\&CG before the time limit, we use bars of dotted lines labeling with $T$ for indication, together with the average numbers of iterations completed within the time limit displayed inside parentheses. A few observations are made as in the following.
\begin{enumerate}
	\item Comparing results for $\mathbf{Robust \ FL}$ with $\mathcal U^{L\!R}(\mathbf x)$ and $\mathcal U^{R}(\mathbf x)$, we note that the dominance of parametric C\&CG over the other two C\&CG variants is more remarkable. Given the existence of large gaps upon termination for the other two  methods and the very short computational time for parametric C\&CG, it is reasonably to believe parametric C\&CG  performs 3 or more orders of magnitude faster than other methods. Hence, it is the only practically feasible method to compute this type of problems.
	\item DDU with LHS dependence clearly brings a great computational challenge. Note that $\mathcal{U}^{L\!R}$ only has a single constraint, i.e., the one in \eqref{FLDDDL_TolUp},  with LHS dependence. Nevertheless,  this single constraint drastically changes performances of both Benders and basis based C\&CGs. For the former one, it generally takes a few hundred iterations with little improvement on reducing optimality gap, suggesting that computing $\mathbf{MP1}$ is not very demanding while cutting sets are very weak. On the contrary, $\mathbf{MP3}$ of the latter one is very difficult to compute. As mentioned in the remark after  Corollary~\ref{cor_V3_complexity}, the complex structure of $\textbf{MP3}$ is probably the primary reason behind the slow computation. We also observe that
	a nonzero gap often exists after parametric C\&CG terminates, which is less common for DDU with RHS dependence. One explanation is that  the numerical issue is more significant when LHS dependence appears.
	\item Although  basis based C\&CG  has a theoretically  desirable  iteration complexity even if LHS dependence is involved, we realize that its actual computational performance is not satisfactory. Hence, one future research direction is to investigate computationally more friendly reformulations and specialized algorithms for the master problem so that Variant 3's practical solution capacity matches with its theoretical strength.
\end{enumerate}

\begin{figure}
\begin{subfigure}[t]{0.45\textwidth}
	\pgfplotsset{tick label style={font=\small},
		label style={font=\small},
		legend style={font=\tiny}
	}
\begin{tikzpicture}
\begin{axis} [height=6cm,width=7.5cm,
    legend style={at={(1.1,1.2)},anchor=north,legend columns=-1,/tikz/every even column/.append style={column sep=5pt}},
    ybar=0.5pt,
    bar width = 6pt,
    ymin = 0,
    ymax = 100,
    yticklabel={\pgfmathparse{\tick}\pgfmathprintnumber{\pgfmathresult}\%},
    symbolic x coords={25,30,35,40},
    xtick = data,
    enlarge x limits = 0.2,
    every node near coord/.append style={rotate=90, anchor=west, font=\tiny}
]
\pgfkeys{/pgf/number format/.cd,fixed,precision=2};
\addplot[draw = blue!50, fill=blue!50, nodes near coords={\pgfmathprintnumber\pgfplotspointmeta\%}] coordinates {(25,34.86) (30,44.45) (35, 44.68) (40, 39.99)};
\addplot[draw = red!50, fill=red!50, nodes near coords={\pgfmathprintnumber\pgfplotspointmeta\%}] coordinates {(25,0.08) (30,0.01) (35, 0.07) (40, 0.07)};
\addplot[draw = green!50, fill=green!50, nodes near coords={\pgfmathprintnumber\pgfplotspointmeta\%}] coordinates {(25,66.98) (30,69.69) (35, 39.5) (40, 72.31)};
\legend {Benders C\&CG, Parametric C\&CG, Basis Based C\&CG};
\end{axis}
\end{tikzpicture}
\caption{Relative Gap}
\end{subfigure}
\begin{subfigure}[t]{0.45\textwidth}
\center
\pgfplotsset{tick label style={font=\small},
	label style={font=\small},
	legend style={font=\tiny}
}
\begin{tikzpicture}
\begin{axis} [height=6cm,width=7.5cm,
    ybar=0.5pt,
    bar width = 6pt,
    ymin = 0,
    ymax = 5000,
    symbolic x coords={25,30,35,40},
    xtick = data,
    enlarge x limits = 0.2,
]
\pgfkeys{/pgf/number format/.cd,fixed,precision=2}
\addplot[draw = blue!50, pattern = dots, pattern color = blue!50,point meta=explicit symbolic,nodes near coords,every node near coord/.append style={rotate=90, anchor=west, font=\tiny}]
table[meta=label]{
	x		y		label
	25		3600	{T (477)}
	30		3600	{T (422.5)}
	35		3600	{T (377)}
	40		3600	{T (337.5)}
};
\addplot[draw = red!50, fill=red!50, point meta=explicit symbolic, nodes near coords, every node near coord/.append style={rotate=90, anchor=west, font=\tiny}]
table[meta=label]{
	x		y		label
	25		14.82	{14.82 (5)}
	30		67.98	{67.98 (5)}
	35		60.02	{60.02 (5)}
	40		57.65	{57.65 (7)}
};
\addplot[draw = green!50, pattern = dots, pattern color = green!50,point meta=explicit symbolic,nodes near coords,every node near coord/.append style={rotate=90, anchor=west, font=\tiny}]
table[meta=label]{
	x		y		label
	25		3600	{T (8.5)}
	30		3600	{T (8)}
	35		3600	{T (9.5)}
	40		3600	{T (8)}
};
\end{axis}
\end{tikzpicture}
\caption{Computation Time (s)}
\end{subfigure}
\caption{Computational Results of Algorithms for $\mathbf{Robust \ FL}$ with $\mathcal U^{L\!R}(\mathbf x)$
\label{Fig:alg_compa_LRH}}
\end{figure}

 \subsection{Utilizing Deep Knowledge to Solve RO: from DIU to DDU}
 \label{subsect: DIUtoDDU}
 In Section \ref{subsect_DIU-DDU}, we have shown  that the classical DIU-based RO formulation for the reliable p-median facility location problem can be converted into a DDU-based formulation by using deep knowledge. Those two formulations are equivalent if the deep knowledge can be proven, and otherwise the DDU-based one can be used to derive quantifiable approximate solutions.   In this subsection, we perform a computational benchmark on them  to gain a deeper understanding on employing DDU-based reformulations for computation.
 As for test instances, they are generated by varying $\rho$, $p$ and $k$ in $\mathbf{Reliable \ FL}$ with 40 sites adopted from the previous subsection. As for computational algorithms, basic and parametric  C\&CGs are employed for the DIU-based and the DDU-based formulations, respectively. Note that parametric C\&CG reduces to the basic one if the DDU set is actually DIU. 

 \subsubsection{Computing Exact Solutions by Equivalent DDU Reformulation}
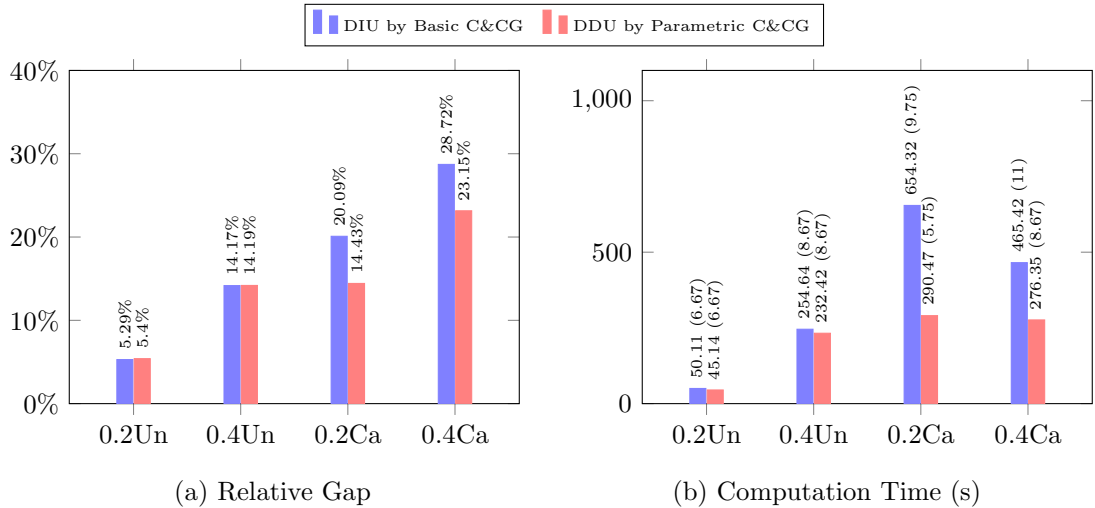
\begin{figure}[!htp]
	\begin{subfigure}[t]{0.45\textwidth}
		\pgfplotsset{tick label style={font=\small},
			label style={font=\small},
			legend style={font=\tiny}
		}
		\begin{tikzpicture}
			\begin{axis} [height=6cm,width=7.5cm,
				legend style={at={(1.1,1.2)},anchor=north,legend columns=-1,/tikz/every even column/.append style={column sep=5pt}},
				ybar=0.5pt,
				bar width = 6pt,
				ymin = 0,
				ymax = 40,
				yticklabel={\pgfmathparse{\tick}\pgfmathprintnumber{\pgfmathresult}\%},
				symbolic x coords={0.2Un,0.4Un,0.2Ca,0.4Ca},
				xtick = data,
				enlarge x limits = 0.2,
				nodes near coords={\pgfmathprintnumber\pgfplotspointmeta\%},
				every node near coord/.append style={rotate=90, anchor=west, font=\tiny}
				]
				\pgfkeys{/pgf/number format/.cd,fixed,precision=2};
				\addplot[draw = blue!50, fill=blue!50] coordinates {(0.2Un,5.29) (0.4Un,14.17) (0.2Ca, 20.09) (0.4Ca, 28.72)};
				\addplot[draw = red!50, fill=red!50] coordinates {(0.2Un,5.4) (0.4Un,14.19) (0.2Ca, 14.43) (0.4Ca, 23.15)};
				\legend {DIU by Basic C\&CG, DDU by Parametric C\&CG};
			\end{axis}
		\end{tikzpicture}
		\caption{Relative Gap}
	\end{subfigure}
	\begin{subfigure}[t]{0.45\textwidth}
		\center
		\pgfplotsset{tick label style={font=\small},
			label style={font=\small},
			legend style={font=\tiny}
		}
		\begin{tikzpicture}
			\begin{axis} [height=6cm,width=7.5cm,
				ybar=0.5pt,
				bar width = 6pt,
				ymin = 0,
				ymax = 1100,
				symbolic x coords={0.2Un,0.4Un,0.2Ca,0.4Ca},
				xtick = data,
				enlarge x limits = 0.2,
				nodes near coords,
				every node near coord/.append style={rotate=90, anchor=west, font=\tiny}
				]
				\pgfkeys{/pgf/number format/.cd,fixed,precision=2};
				\addplot[draw = blue!50, fill=blue!50, point meta=explicit symbolic]
				table[meta=label]{
					x		y		label
					0.2Un	50.11	{50.11 (6.67)}
					0.4Un	245.64	{254.64 (8.67)}
					0.2Ca	654.32	{654.32 (9.75)}
					0.4Ca	465.42	{465.42 (11)}
				};
				\addplot[draw = red!50, fill=red!50, point meta=explicit symbolic]
				table[meta=label]{
					x		y		label
					0.2Un	45.14	{45.14 (6.67)}
					0.4Un	232.42	{232.42 (8.67)}
					0.2Ca	290.47	{290.47 (5.75)}
					0.4Ca	276.35	{276.35 (8.67)}
				};
			\end{axis}
		\end{tikzpicture}
		\caption{Computation Time (s)}
	\end{subfigure}
	\caption{Computational Results of the Reliable P-Median Problem\label{fig:p_median_equiv}}
\end{figure}

 As noted in Corollary \ref{cor_reliable_pmedian}, the DIU-based reliable p-median formulation in \eqref{eq_reliable_pmedian} using $\mathcal{U}^0$ and the DDU-based reformulation  in  \eqref{eq_reliable_pmedian_DDU} using $\mathcal{U}^k(\mathbf x)$ are equivalent  for both uncapacitated and capacitated cases, if disruptions do not increase demand. Given that, we let demands equal to nominal ones, and present their computational results in Table~\ref{tbl:PMD} in Appendix \ref{Asect_tables}. Their average performances are displayed in Figure \ref{fig:p_median_equiv}, where $0.2$ and $0.4$ are values of parameter $\rho$, $Un$ and $Ca$ denote uncapacitated and capacitated cases respectively.
 	
 For the uncapacitated case,  we observe that those two formulations have comparable performances. Computing the DIU-based formulation has a slightly better optimality gap across all instances, while the DDU-based formulation can be computed a little bit faster for instances that are solved exactly. For the capacitated case, nevertheless, we highlight that the DDU-based formulation has a significantly better performance.  For the average gap across all instances, a reduction around 28\% to 56\%  is observed. For instances that can be solved exactly,  a  larger reduction, averagely around 41\% to 55\%, can be achieved in computational time. For some instances, this improvement could be much more substantial. We note that for the instance with $\rho=0.2,  p=6, k=1$,  the computational efficiency of the DDU-based formulation is almost 3 times faster than that of the DIU-based one.

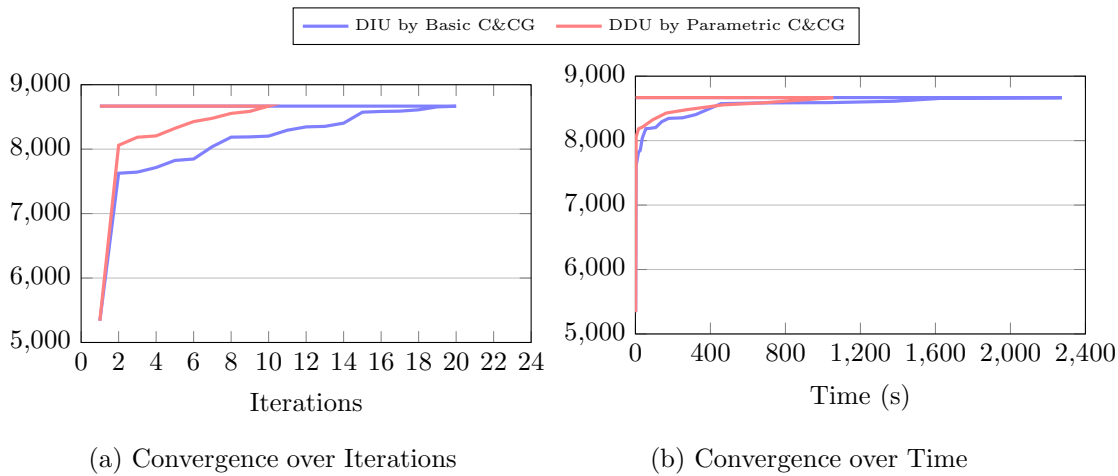
\begin{figure}[h]
\begin{subfigure}[t]{0.45\textwidth}
	\pgfplotstableread{Sp4.dat}{\Sp}
	\pgfplotsset{tick label style={font=\small\bfseries},
		label style={font=\small},
		legend style={font=\tiny}
	}
	\begin{tikzpicture}[scale=1]
		\begin{axis}[ymajorgrids=true,legend style={at={(1.1,1.3)},anchor=north,legend columns=-1,/tikz/every even column/.append style={column sep=5pt}},height=5cm,width=7.5cm,
			xmin=0,xmax=24,ymin=5000,ymax=9000,
			xtick={0,2,4,6,8,10,12,14,16,18,20,22,24},
			ytick={5000,6000,7000,8000,9000},
			xlabel= Iterations]
			\addplot [blue!50,very thick] table [x={DIUIt}, y={DIUB}] {\Sp};
			\addlegendentry{DIU by Basic C\&CG}
			\addplot [red!50,very thick] table [x={DDUIt}, y={DDUB}] {\Sp};
			\addlegendentry{DDU  by Parametric C\&CG}
		\end{axis}
	\end{tikzpicture}
		\caption{Convergence over Iterations}
\label{ConPMIt}
\end{subfigure}
\begin{subfigure}[t]{0.42\textwidth}
	\pgfplotstableread{Sp3.dat}{\Sp}
	\pgfplotsset{tick label style={font=\small\bfseries},
		label style={font=\small},
		legend style={font=\tiny}
	}
	\begin{tikzpicture}[scale=1]
		\begin{axis}[ymajorgrids=true,height=5cm,width=7.5cm,
			xmin=0,xmax=2400,ymin=5000,ymax=9000,
			xtick={0,400,800,1200,1600,2000,2400},
			ytick={5000,6000,7000,8000,9000},
			xlabel= Time (s)]
			\addplot [blue!50,very thick] table [x={DIUTi}, y={DIUB}] {\Sp};
			\addplot [red!50,very thick] table [x={DDUTi}, y={DDUB}] {\Sp};
		\end{axis}
	\end{tikzpicture}
	\caption{Convergence over Time}
\label{ConPMTi}
\end{subfigure}
\caption{Convergence of DIU- and DDU-based Reliable P-median Formulations \label{ConTFs}}
\end{figure}

 To further understand algorithms' dynamic behaviors on these two formulations,  we again plot their convergent progresses for the instance with  $\rho=0.2,  p=2, k=2$ in Figure~\ref{ConTFs}. Besides their obvious differences in the number of iterations and computational time,  it is interesting to observe that the progress curves of the DDU-based formulation is much smoother than that of the DIU-based one. The reason we believe is that the connection reflected in DDU ensures the parametric cutting sets directly capture the associated recourse cost of the first stage decision, without depending on any particular scenarios. They hence force the master problem to search for effective first stage decisions that are substantially different from  one iteration to another one.  On the contrary, each cutting set generated by basic C\&CG depends on a specific scenario. The master problem can often slightly modify an existing first stage decision from one iteration to its following iteration to render a recently identified scenario less disruptive. In other words, unless cutting sets from multiple scenarios are imposed, a substantially different first stage decision will not be produced. As a result, the lower bound curve increases rather slowly with a few sizable jumps.  It is also interesting to note that lower bound curve of the DIU-based formulation crosses over that of the DDU-based one a little bit in Figure \ref{ConTFs}.(b), which is not the case in Figure \ref{ConTFs}.(a). The reason is that the master problem of basic C\&CG can be solved faster. So,  on a particular time point,  basic  C\&CG may complete many iterations and produce a solution that is better than a solution derived by parametric C\&CG, which only completes less iterations up to that point.

Clearly, results for the capacitated case  demonstrate the great benefit of incorporating deep knowledge  on improving our solution capacity. Regarding the ineffectiveness of this idea on the  uncapacitated case, one explanation is that  the associated $\mathcal{OU}$ and $\mathcal{OV}$  are not singletons in general, i.e., they do not have the uniqueness property. As argued in Section \ref{subsect_pCCG_complexity}, the computational performance should be improved if this property can be ensured. In the latter of this section, we adopt a strategy described in  Section \ref{subsect_Uniqueness_Pareto} to achieve this property and to implement the modified parametric C\&CG with a better performance.

\vspace{0pt}
\subsubsection{Computing Approximation Solutions by Various DDU Sets}
\label{subsect_DDU_app}

As noted in Section \ref{subsect_DIU-DDU}, for the reliable p-median problem,  if one site's demand could increase after disruption, the DDU-based formulation in  \eqref{eq_reliable_pmedian_DDU} does not represent the problem exactly. Nevertheless, it follows from the descriptions before and in Proposition~\ref{prop_CCG_DDU_approximation} that \eqref{eq_reliable_pmedian_DDU} with $\mathcal{U}^k(\mathbf x)$ can be used to derive quantifiable approximation solutions. Assuming that one site's demand will be doubled after disruption, we perform a set of experiments on capacitated instances to evaluate this approximation. We also consider $\mathcal{U}^{kq}(\mathbf x)$ with $q=3$ defined in \eqref{eq_DDU2_reliable_pmedian}. It extends $\mathcal{U}^k(\mathbf x)$ by considering 3 more sites with the largest demands, which reflects an empirical understanding that worst case disruptions often occur on sites with large demands.

\begin{figure}[!h]
\begin{subfigure}[t]{0.45\textwidth}
	\pgfplotsset{tick label style={font=\small},
		label style={font=\small},
		legend style={font=\tiny}
	}
\begin{tikzpicture}
\begin{axis} [height=6cm,width=7.5cm,
    legend style={at={(1.1,1.2)},anchor=north,legend columns=-1,/tikz/every even column/.append style={column sep=5pt}},
    ybar=0.5pt,
    bar width = 6pt,
    ymin = 0,
    ymax = 40,
    yticklabel={\pgfmathparse{\tick}\pgfmathprintnumber{\pgfmathresult}\%},
    symbolic x coords={0.2,0.4},
    xtick = data,
    enlarge x limits = 0.5,
    nodes near coords={\pgfmathprintnumber\pgfplotspointmeta\%},
    every node near coord/.append style={rotate=90, anchor=west, font=\tiny}
]
\pgfkeys{/pgf/number format/.cd,fixed,precision=2};
\addplot[draw = blue!50, fill=blue!50] coordinates {(0.2,21.43) (0.4,29.92)};
\addplot[draw = red!50, fill=red!50] coordinates {(0.2,19.27) (0.4,29.73)};
\addplot[draw = green!50, fill=green!50] coordinates {(0.2,17.91) (0.4,23.64)};
\legend {DIU by basic C\&CG, $\mathcal{U}^k(\mathbf x)$ by Parametric C\&CG, $\mathcal{U}^{kq}(\mathbf x)$ by Parametric C\&CG};
\end{axis}
\end{tikzpicture}
\caption{Average Gap of All Instances}
\end{subfigure}
\begin{subfigure}[t]{0.45\textwidth}
\center
\pgfplotsset{tick label style={font=\small},
	label style={font=\small},
	legend style={font=\tiny}
}
\begin{tikzpicture}
\begin{axis} [height=6cm,width=7.5cm,
    ybar=0.5pt,
    bar width = 6pt,
    ymin = 0,
    ymax = 60,
    yticklabel={\pgfmathparse{\tick}\pgfmathprintnumber{\pgfmathresult}\%},
    symbolic x coords={0.2,0.4},
    xtick = data,
    enlarge x limits = 0.5,
    nodes near coords={\pgfmathprintnumber\pgfplotspointmeta\%},
    every node near coord/.append style={rotate=90, anchor=west, font=\tiny}
]
\pgfkeys{/pgf/number format/.cd,fixed,precision=2,set thousands separator={}};
\addplot[draw = blue!50, fill=blue!50] coordinates {(0.2,38.57) (0.4,44.88)};
\addplot[draw = red!50, fill=red!50] coordinates {(0.2,32.72) (0.4,42.41)};
\addplot[draw = green!50, fill=green!50] coordinates {(0.2,32.24) (0.4,35.46)};
\end{axis}
\end{tikzpicture}
\caption{Average Gap of Unsolved Instances}
\end{subfigure}
\caption{DDU Approximation when Demands Doubled after Disruptions\label{fig:Pmedian_double_gap}}
\end{figure}
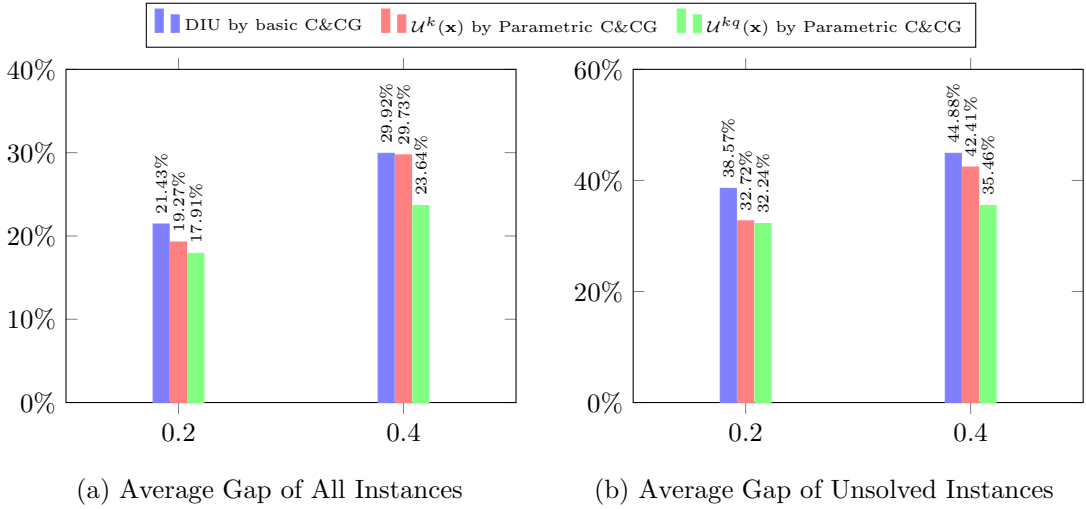

Detailed results are reported in Table \ref{tbl:CapPMDD} in Appendix \ref{Asect_tables}. The average optimality gaps on  all instances, and on unsolved instances within the time limit are summarized in Figure \ref{fig:Pmedian_double_gap} according to $\rho=0.2$ and $0.4$, respectively. Note that the first termination condition  discussed before Proposition \ref{prop_CCG_DDU_approximation} is adopted to stop the algorithm execution. On the one hand,  the approximation quality from using $\mathcal{U}^k(\mathbf x)$ is not bad. The average gap upon termination for unsolved instances is clearly better than that of the DIU-based formulation. Nevertheless, it often fails to generate optimal solutions for easy instances.  
On the other hand,  the approximation quality from using $\mathcal{U}^{kq}(\mathbf x)$ is significantly improved.  The average gap upon termination is reduced by  almost 20\%. Especially for difficult instances where the DIU-based formulation terminates with large gaps, the improvement could be more substantial, up to 47\% reduction. Actually, for all instances that are solved to optimality by basic C\&CG, the implementation using $\mathcal{U}^{kq}(\mathbf x)$ produces exact solutions as well. We mention that such an improvement from $\mathcal{U}^{kq}(\mathbf x)$ demonstrates the value of incorporating valid heuristic understanding and information feedback from computing subproblems in constructing more effective DDU set.

%

\vspace{-2pt}

We also implement and evaluate DDU sets constructed with help of auxiliary variables to reflect sophisticated understanding. One is DDU set $\mathcal U^{r}(\mathbf x)$ described in \eqref{eq_app_DDU_sorting} (and updated set $\mathcal X$ defined in \eqref{eq_app_X_sorting}). It reflects a heuristic understanding that facilities with heavy first stage service costs should be disrupted in worst case scenarios. We also employ two DDU sets simultaneously in the purpose to better capture the connection between the first stage decision and worse case scenarios. So, we further augment $\mathcal{X}$ and consider another DDU set as in the following. It reflects a heuristic understanding that facility sites with large recourse costs should be disrupted in worst case scenarios. Similar to  $\mathbf x_{r}$ in \eqref{eq_app_X_sorting}, auxiliary variable $\mathbf x_{s}$ is introduced to $\mathcal{X}$ that sorts facilities according to their recourse costs and provides necessary information for the construction of the expected DDU set. Specifically, $x_{s,j}=1$ if  $j$ is one of those sites with largest recourse costs  (measured by the product between the demand served by $j$ and the total distance to other facilities), $0$ otherwise. We have

\begin{align}
	\label{eq_app_X_sorting2}
	\begin{split}
		\mathcal X^{rs}=\{&(\mathbf x_c, \mathbf x_d, \mathbf x_r, \mathbf x_s)\in \mathbb{R}^{|I|\times|J|+2}_+\times\{0,1\}^{3|J|}:  \sum_{j\in J}x_{d,j}=p, \ \sum_{j\in J}x_{c,ij}\geq d_i \ \forall i, \\ & \sum_{i\in I}x_{c,ij}\leq A_jx_{d,j} \ \forall j, \ x_{r,j}\leq x_{d,j} \ \forall j, \ \sum_{i\in I}c_{ij}x_{c,ij}\geq x^0_r-M(1-x_{r,j}) \ \forall j, \\ & \sum_{i\in I}c_{ij}x_{c,ij}\leq x^0_{r}+Mx_{r,j} \ \ \forall j, \ \sum_{j\in J}x_{r,j}=q_1, \\
		& x_{s,j}\leq x_{d,j} \ \forall j, \ (\sum_{i\in I}x_{c,ij})(\sum_{j'\in J}c_{jj'}x_{d,j}x_{d,j'})\geq x^0_s-M(1-x_{s,j}) \ \forall j, \\
		& (\sum_{i\in I}x_{c,ij})(\sum_{j'\in J}c_{jj'}x_{d,j}x_{d,j'})\leq x^0_{s}+Mx_{s,j} \ \forall j, \ \sum_{j\in J}x_{s,j}=q_2\}.
	\end{split}
\end{align}
In addition to identifying $q_1$ facilities with the most first stage service costs,
the last two rows of constraints sort out the top $q_2$ facilities according to
$$(\sum_{i\in I}x_{c,ij})(\sum_{j'\in J}c_{jj'}x_{d,j}x_{d,j'}).$$ The first term represents the demands served by the facility site $j$, and the second term computes its overall distances to other facilities, representing the recourse cost after $j$ is disrupted. Note that the nonlinear expressions in those constraints can be easily linearized, given that $x_{d,j}$ and $x_{d,j'}$ are binary. We next define the second DDU set accordingly. 
\begin{align*}
	\label{eq_app_DDU_sorting2}
	\mathcal{U}^s(\mathbf x)&=\{\mathbf {u}\in \mathbb{R}^{|I|}_+: \sum_{j\in J}u_j\leq k, \ u_j\leq x_{s,j} \ \forall j, \  u_i=0 \ \forall i\notin J\}
\end{align*}


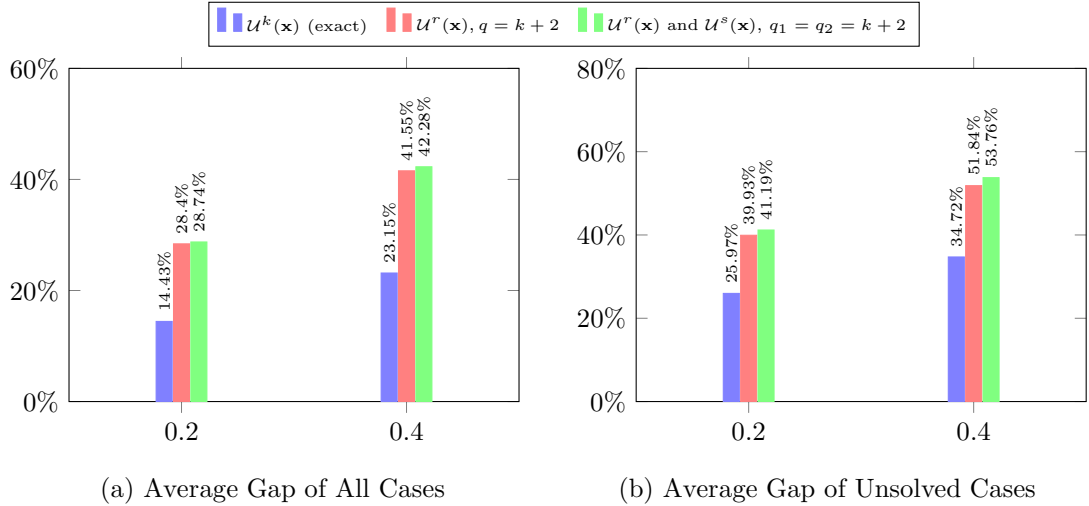
\begin{figure}
	\begin{subfigure}[t]{0.45\textwidth}
		\pgfplotsset{tick label style={font=\small},
			label style={font=\small},
			legend style={font=\tiny}
		}
		\begin{tikzpicture}
			\begin{axis} [height=6cm,width=7.5cm,
				legend style={at={(1.1,1.2)},anchor=north,legend columns=-1,/tikz/every even column/.append style={column sep=5pt}},
				ybar=0.5pt,
				bar width = 6pt,
				ymin = 0,
				ymax = 60,
				yticklabel={\pgfmathparse{\tick}\pgfmathprintnumber{\pgfmathresult}\%},
				symbolic x coords={0.2,0.4},
				xtick = data,
				enlarge x limits = 0.5,
				nodes near coords={\pgfmathprintnumber\pgfplotspointmeta\%},
				every node near coord/.append style={rotate=90, anchor=west, font=\tiny}
				]
				\pgfkeys{/pgf/number format/.cd,fixed,precision=2};
				\addplot[draw = blue!50, fill=blue!50] coordinates {(0.2,14.43) (0.4,23.15)};
				\addplot[draw = red!50, fill=red!50] coordinates {(0.2,28.4) (0.4,41.55)};
				\addplot[draw = green!50, fill=green!50] coordinates {(0.2,28.74) (0.4,42.28)};
				\legend {$\mathcal{U}^k(\mathbf x)$ (exact), {$\mathcal{U}^r(\mathbf x), q=k+2$}, {  $\mathcal{U}^r(\mathbf x)$ and $\mathcal{U}^s(\mathbf x)$, $q_1=q_2=k+2$}};
			\end{axis}
		\end{tikzpicture}
		\caption{Average Gap of All Cases}
	\end{subfigure}
	\begin{subfigure}[t]{0.45\textwidth}
		\center
		\pgfplotsset{tick label style={font=\small},
			label style={font=\small},
			legend style={font=\tiny}
		}
		\begin{tikzpicture}
			\begin{axis} [height=6cm,width=7.5cm,
				ybar=0.5pt,
				bar width = 6pt,
				ymin = 0,
				ymax = 80,
				yticklabel={\pgfmathparse{\tick}\pgfmathprintnumber{\pgfmathresult}\%},
				symbolic x coords={0.2,0.4},
				xtick = data,
				enlarge x limits = 0.5,
				nodes near coords={\pgfmathprintnumber\pgfplotspointmeta\%},
				every node near coord/.append style={rotate=90, anchor=west, font=\tiny}
				]
				\pgfkeys{/pgf/number format/.cd,fixed,precision=2,set thousands separator={}};
				\addplot[draw = blue!50, fill=blue!50] coordinates {(0.2,25.97) (0.4,34.72)};
				\addplot[draw = red!50, fill=red!50] coordinates {(0.2,39.93) (0.4,51.84)};
				\addplot[draw = green!50, fill=green!50] coordinates {(0.2,41.19) (0.4,53.76)};
			\end{axis}
		\end{tikzpicture}
		\caption{Average Gap of Unsolved Cases}
	\end{subfigure}
	\caption{DDU Approximation with Sorting \label{fig:P_median_appr_sort}}
\end{figure}

Detailed computational results of parametric C\&CG with those sorting based DDU sets, along with those from the exact reformulation based on $\mathcal{U}^k(\mathbf x)$, are reported in Table \ref{tbl:AppCapPM_sorting} in Appendix \ref{Asect_tables} for the  capacitated case with nominal demands. Again, the average optimality gaps on all instances, and on unsolved instances within the time limit are summarized respectively in Figure \ref{fig:P_median_appr_sort}.   It can been seen that either the implementation with single set $\mathcal{U}^r(\mathbf x)$,  or the one with two sets, i.e.,  $\mathcal{U}^{r}(\mathbf x)$ and $\mathcal{U}^s(\mathbf x)$, cannot provide a satisfactory approximation quality with respect to results obtained by using $\mathcal{U}^k(\mathbf x)$. For the majority of instances, their associated gaps and computational time are worse than those from $\mathcal{U}^k(\mathbf x)$. Note that lower bounds from approximations are often small, which suggests that they might not represent the most disruptive scenarios. Comparing implementations with a single DDU set and with two DDU sets jointly, the latter one fails to yield  any improvement, while requires significantly longer computational time. Hence, it indicates $\mathcal U^{s}(\mathbf x)$ may not carry much new information on the actual worse case scenarios, and $\mathcal X^{rs}$ slows down our computation.

Although DDU sets $\mathcal U^r(\mathbf x)$ and $\mathcal U^s(\mathbf x)$ do not lead to strong approximations in general, we observe that for some large instances, e.g., those with $p=8$ and $k=3$, both of them, especially $\mathcal{U}^{r}(\mathbf x)$,  generate feasible solutions better than those produced by using $\mathcal{U}^k(\mathbf x)$ in a very short time. It shows that the reflected heuristic understandings still capture some nontrivial connection between the first stage decision and critical disruptive scenarios. 


Before closing this subsection, we mention that this strategy, i.e., employing DDU sets to derive exact and quantifiable approximation solutions, is critical to compute practical instances with a strong scalability, which certainly requires deeper analyses and study. 


\subsection{Computing RO with Mixed Integer Recourse}
\label{subsect_computation_MIPrecourse}
In this subsection, we study the empirical performance of the approximate scheme designed to compute  $\mathbf{2-Stg \ RO}$ with mixed integer recourse problems. Regarding the testing platform, we extend $\mathbf{Robust \ FL}$ presented in Section \ref{subsubsect:RobustFL} to build an RO model with mixed integer recourse, where binary variables $\mathbf z$ are introduced to represent establishment decisions of temporary facilities with fixed capacities. Let $z_j=1$ denote the establishment of a temporary facility on site $j$ and $0$ otherwise. The new objective function and the feasible set for the recourse problem are updated as in the following.
\begin{equation*}
	\begin{split}
	\mathbf{Robust \ FL-MIP}: & \min_{(\mathbf x_c,\mathbf x_d)\in\mathcal X} \sum_{j\in J} (f_jx_{d,j}+a_jx_{c,j})  \\
	  & + \max_{\mathbf u\in\mathcal U(\mathbf x)} \ \min_{(\mathbf y_1,\mathbf y_2\mathbf, \mathbf z)\in\mathcal{Y}(\mathbf x,\mathbf u)} \ \sum_{i\in I}\sum_{j\in J}(c_{ij}-p_i)y_{1,ij}+\sum_{j\in J} \hat f_jz_j+\sum_{i\in I} \textsl{C}y_{2,i}
	\end{split}
\end{equation*}
\begin{eqnarray*}
		\mathcal Y(\mathbf x,\mathbf u)=\{(\mathbf y_1, \mathbf y_2, \mathbf z):\quad & \sum_{j\in J}y_{1,ij}+ y_{2,i} \geq u_i \ \forall i,
		\ \sum_{i\in I}y_{1,ij}\leq x_{c,j}+h_jz_j\quad\forall j,\\
		&y_{1,ij}\geq 0, \ y_{2,i}\geq 0, \ z_j\in\{0,1\}\}
\end{eqnarray*}

In $\mathbf{Robust \ FL-MIP}$, $\hat f_j$ is the fixed cost of establishing a temporary facility with fixed capacity $h_j$ at site $j$. Variable $y_{1,ij}$ represents site $i$'s demand served by the facility(s) at~$j$, where the total capacity may include that from the pre-established and from the temporary ones, and  $y_{2,i}$ represents its unmet demand, which is penalized with \textsl{C} per unit.


\begin{figure}[h]
\begin{subfigure}[t]{0.45\textwidth}
	\pgfplotsset{tick label style={font=\small},
		label style={font=\small},
		legend style={font=\tiny}
	}
\begin{tikzpicture}
\begin{axis} [height=6cm,width=7.5cm,
    legend style={at={(1.1,1.2)},anchor=north,legend columns=-1,/tikz/every even column/.append style={column sep=5pt}},
    ybar=0.5pt,
    bar width = 6pt,
    ymin = 0,
    ymax = 0.15,
    yticklabel={\pgfmathparse{\tick}{\pgfmathresult}\%},
    symbolic x coords={25L,25H,40L,40H},
    xtick = data,
    enlarge x limits = 0.2,
    nodes near coords={\pgfmathprintnumber\pgfplotspointmeta\%},
    every node near coord/.append style={rotate=90, anchor=west, font=\tiny}
]
\pgfkeys{/pgf/number format/.cd,fixed,precision=2};
\addplot[draw = blue!50, fill=blue!50] coordinates {(25H,0.04) (25L,0) (40H, 0) (40L, 0.01)};
\addplot[draw = red!50, fill=red!50] coordinates {(25H,0.04) (25L,0) (40H, 0.04) (40L, 0)};
\legend {w/ Temp. Facilities (Ave), w/o Temp.  Facilities (Ave)};
\end{axis}
\end{tikzpicture}
\caption{Relative Gap}
\end{subfigure}
\begin{subfigure}[t]{0.45\textwidth}
\center
\pgfplotsset{tick label style={font=\small},
	label style={font=\small},
	legend style={font=\tiny}
}
\begin{tikzpicture}
\begin{axis} [height=6cm,width=7.5cm,
    ybar=0.5pt,
    bar width = 6pt,
    ymin = 0,
    ymax = 25,
    symbolic x coords={25L,25H,40L,40H},
    xtick = data,
    enlarge x limits = 0.2,
    nodes near coords,
    every node near coord/.append style={rotate=90, anchor=west, font=\tiny}
]
\pgfkeys{/pgf/number format/.cd,fixed,precision=2,set thousands separator={}};
\addplot[draw = blue!50, fill=blue!50, point meta=explicit symbolic]
table[meta=label]{
	x		y		label
	25H		3.59	{3.59 (2.2)}
	25L		6.2		{6.2 (3)}
	40H		10.7	{10.7 (2)}
	40L		16.96	{16.96 (2)}
};
\addplot[draw = red!50, fill=red!50, point meta=explicit symbolic]
table[meta=label]{
	x		y		label
	25H		2.16	{2.16 (2)}
	25L		2.9		{2.9 (3)}
	40H		4.54	{4.54 (2)}
	40L		4.87	{4.87 (3.2)}
};
\end{axis}
\end{tikzpicture}
\caption{Computation Time (s)}
\end{subfigure}
\caption{Computational Results for Mixed Integer Recourse\label{fig:appr_MIPrecourse}}
\end{figure}
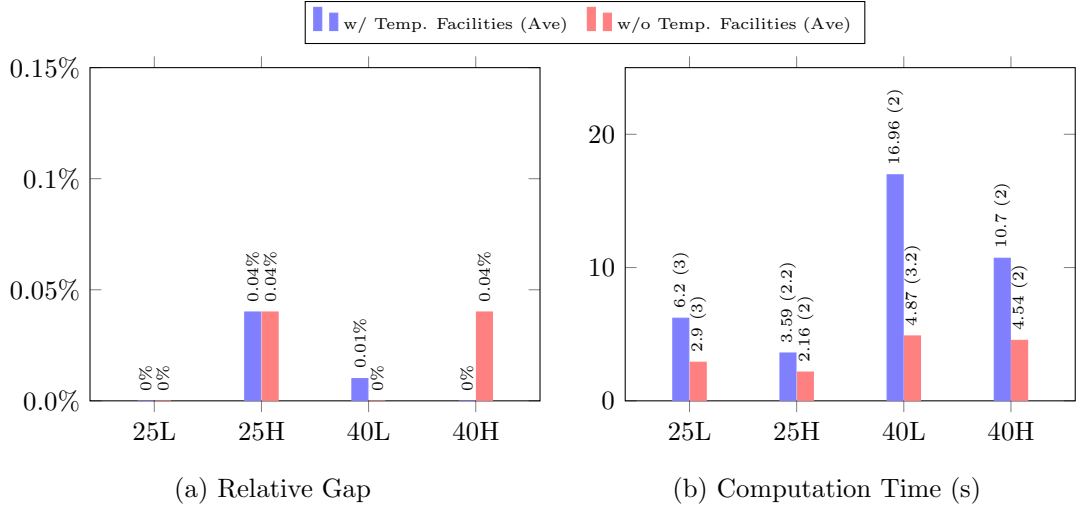

Testing instances are obtained by modifying those adopted in Section \ref{subsubsect:RobustFL} with $\mathcal{U}(\mathbf x)=\mathcal{U}^R(\mathbf x)$. Changes include setting $h_j$ to a random number in $[30,80]$, $\hat f_j=5*h_j*\displaystyle\max_{j'\in J}\{a_{j'}\}$ and $\textsl{C}=1.5*\displaystyle\max_{i\in I,j\in J}\{c_{ij}\} $.  Detailed computational results are reported in Table \ref{tbl:MIPRC} in Appendix~\ref{Asect_tables}, along with those obtained for $\mathbf{Robust \ FL-MIP}$ with $\mathbf z=\mathbf 0$ (which has an LP recourse problem with no temporary facilities). In Figure \ref{fig:appr_MIPrecourse}, we present the summary of average performances. Numerical results show that this approximation scheme is very capable in computing high quality solutions. For the majority of cases, the optimality gaps upon termination are zeros or well below the optimality tolerance, which agrees with one observation made in \citet{zhang2018ambulance} for a DIU-based RO model with MIP recourse. Also, the number of iterations before termination is comparable to that of the case with LP recourse. This is particularly useful, indicating that little extra complexity is involved except for the computational burden associated the MIP formulation of the recourse problem.

\subsection{Enhancements with Uniqueness and Pareto Optimality}
\begin{figure}[!htp]
\begin{subfigure}[t]{0.45\textwidth}
	\pgfplotsset{tick label style={font=\small},
		label style={font=\small},
		legend style={font=\tiny}
	}
\begin{tikzpicture}
\begin{axis} [height=6cm,width=7.5cm,
    legend style={at={(1.1,1.2)},anchor=north,legend columns=-1,/tikz/every even column/.append style={column sep=5pt}},
    ybar=0.5pt,
    bar width = 6pt,
    ymin = 0,
    ymax = 40,
    yticklabel={\pgfmathparse{\tick}\pgfmathprintnumber{\pgfmathresult}\%},
    symbolic x coords={0.2,0.4},
    xtick = data,
    enlarge x limits = 0.5,
    nodes near coords={\pgfmathprintnumber\pgfplotspointmeta\%},
    every node near coord/.append style={rotate=90, anchor=west, font=\tiny}
]
\pgfkeys{/pgf/number format/.cd,fixed,precision=2};
\addplot[draw = blue!50, fill=blue!50] coordinates {(0.2,14.43) (0.4,23.15)};
\addplot[draw = red!50, fill=red!50] coordinates {(0.2,15.18) (0.4,23.22)};
\addplot[draw = green!50, fill=green!50] coordinates {(0.2,16.36) (0.4,24.18)};
\legend {Parametric C\&CG,Parametric C\&CG w/ Pareto,Parametric C\&CG w/ Uniqueness};
\end{axis}
\end{tikzpicture}
\caption{Relative Gap}
\end{subfigure}
\begin{subfigure}[t]{0.45\textwidth}
\center
\pgfplotsset{tick label style={font=\small},
	label style={font=\small},
	legend style={font=\tiny}
}
\begin{tikzpicture}
\begin{axis} [height=6cm,width=7.5cm,
    ybar=0.5pt,
    bar width = 6pt,
    ymin = 0,
    ymax = 520,
    symbolic x coords={0.2,0.4},
    xtick = data,
    enlarge x limits = 0.5,
    nodes near coords,
    every node near coord/.append style={rotate=90, anchor=west, font=\tiny}
]
\pgfkeys{/pgf/number format/.cd,fixed,precision=2,set thousands separator={}};
\addplot[draw = blue!50, fill=blue!50, point meta=explicit symbolic]
table[meta=label]{
	x		y		label
	0.2		290.47	{290.47 (5.75)}
	0.4		276.35	{276.35 (8.67)}
};
\addplot[draw = red!50, fill=red!50, point meta=explicit symbolic]
table[meta=label]{
	x		y		label
	0.2		324.98	{324.98 (5.5)}
	0.4		271.58	{271.58 (8.33)}
};
\addplot[draw = green!50, fill=green!50, point meta=explicit symbolic]
table[meta=label]{
	x		y		label
	0.2		315.16	{315.16 (5.25)}
	0.4		247.77	{247.77 (7.67)}
};
\end{axis}
\end{tikzpicture}
\caption{Computation Time (s)}
\end{subfigure}
\caption{Capacitated Reliable P-Median Problem (with Enhancements)\label{fig:capacitated_enhancement}}
\end{figure}
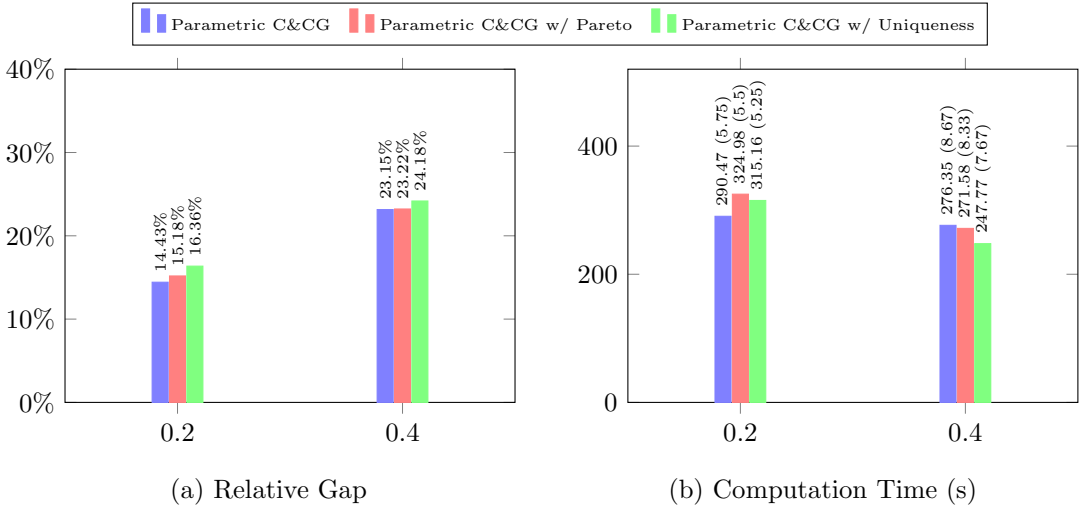

In this subsection, we implement and evaluate two enhancement strategies introduced in Section \ref{subsect_Uniqueness_Pareto}. One is to directly modify  the objective function that defines set $\mathcal{OU}$ based on structural information, which ensures $\mathcal{OU}$  with the uniqueness property. Detailed operations have been described before Corollary \ref{cor_direct_manu_unique} on the DDU-based model for the reliable $p$-median problem. Another one is to select cutting sets with Pareto optimality by computing an extra optimization problem. To reduce the additional computational burden, we adopt the simplified LP version, i.e., $\mathbf{SP2}_{POL}$ defined in \eqref{eq_Pareto_LP}, in our implementation. We also set $\mathbf x^0$ to an optimal solution of the classical p-median problem, and $u'_j=1$ if $j$ is chosen  as a facility site according to $\mathbf x^0$  and $u'_j=0$ otherwise. Numerical results of those two enhancements, along with those from standard parametric C\&CG,  are reported in Tables \ref{tbl:UnCapPM} and \ref{tbl:CapPM} for uncapacitated and capacitated cases, respectively, in Appendix \ref{Asect_tables}. Also, Figures~\ref{fig:uncap_enhancement} and \ref{fig:capacitated_enhancement} display the overall average performances.

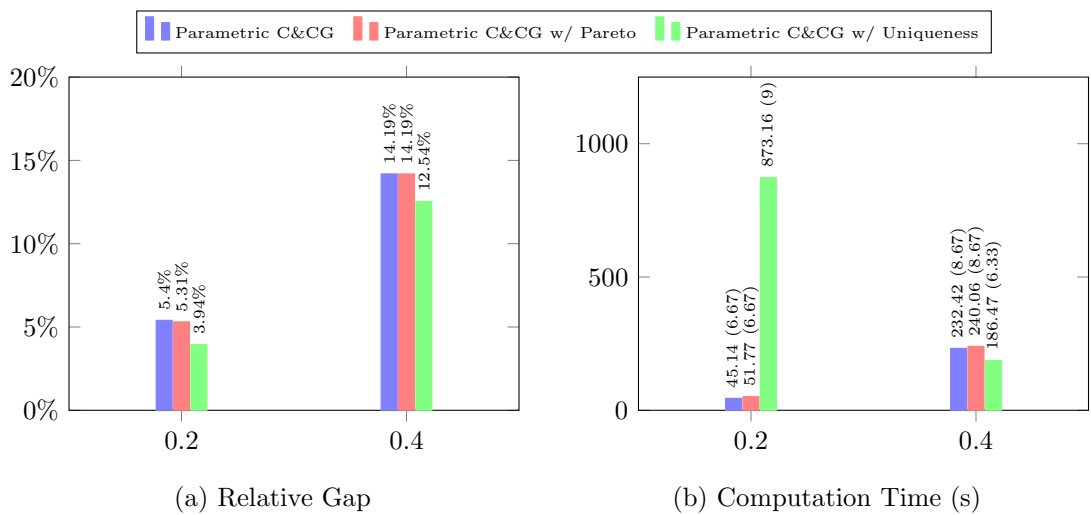
\begin{figure}[!htp]
	\begin{subfigure}[t]{0.45\textwidth}
		\pgfplotsset{tick label style={font=\small},
			label style={font=\small},
			legend style={font=\tiny}
		}
		\begin{tikzpicture}
			\begin{axis} [height=6cm,width=7.5cm,
				legend style={at={(1.1,1.2)},anchor=north,legend columns=-1,/tikz/every even column/.append style={column sep=5pt}},
				ybar=0.5pt,
				bar width = 6pt,
				ymin = 0,
				ymax = 20,
				yticklabel={\pgfmathparse{\tick}\pgfmathprintnumber{\pgfmathresult}\%},
				symbolic x coords={0.2,0.4},
				xtick = data,
				enlarge x limits = 0.5,
				nodes near coords={\pgfmathprintnumber\pgfplotspointmeta\%},
				every node near coord/.append style={rotate=90, anchor=west, font=\tiny}
				]
				\pgfkeys{/pgf/number format/.cd,fixed,precision=2};
				\addplot[draw = blue!50, fill=blue!50] coordinates {(0.2,5.4) (0.4,14.19)};
				\addplot[draw = red!50, fill=red!50] coordinates {(0.2,5.31) (0.4,14.19)};
				\addplot[draw = green!50, fill=green!50] coordinates {(0.2,3.94) (0.4,12.54)};
				\legend {Parametric C\&CG,Parametric C\&CG w/ Pareto,Parametric C\&CG w/ Uniqueness};
			\end{axis}
		\end{tikzpicture}
		\caption{Relative Gap}
	\end{subfigure}
	\begin{subfigure}[t]{0.45\textwidth}
		\center
		\pgfplotsset{tick label style={font=\small},
			label style={font=\small},
			legend style={font=\tiny}
		}
		\begin{tikzpicture}
			\begin{axis} [height=6cm,width=7.5cm,
				ybar=0.5pt,
				bar width = 6pt,
				ymin = 0,
				ymax = 1250,
				symbolic x coords={0.2,0.4},
				xtick = data,
				enlarge x limits = 0.5,
				nodes near coords,
				every node near coord/.append style={rotate=90, anchor=west, font=\tiny}
				]
				\pgfkeys{/pgf/number format/.cd,fixed,precision=2,set thousands separator={}};
				\addplot[draw = blue!50, fill=blue!50, point meta=explicit symbolic]
				table[meta=label]{
					x		y		label
					0.2		45.14	{45.14 (6.67)}
					0.4		232.42	{232.42 (8.67)}
				};
				\addplot[draw = red!50, fill=red!50, point meta=explicit symbolic]
				table[meta=label]{
					x		y		label
					0.2		51.77	{51.77 (6.67)}
					0.4		240.06	{240.06 (8.67)}
				};
				\addplot[draw = green!50, fill=green!50, point meta=explicit symbolic]
				table[meta=label]{
					x		y		label
					0.2		873.16	{873.16 (9)}
					0.4		186.47	{186.47 (6.33)}
				};
			\end{axis}
		\end{tikzpicture}
		\caption{Computation Time (s)}
	\end{subfigure}
	\caption{Uncapacitated Reliable P-Median Problem (with Enhancements)\label{fig:uncap_enhancement}}
\end{figure}

Based on results presented in those figures and tables, we note that the improvement from utilizing Pareto optimal cutting sets over  standard parametric C\&CG is very marginal, if exists, for both  uncapacitated and capacitated cases. For the majority of instances, these two implementations are roughly same in terms of gap and number of iterations. Also, because of the additional computational burden associated with the extra optimization problem,  more computational time is often necessary.  Hence, for this RO problem, the enhancement with Pareto optimal cutting sets is not recommended.

As for another enhancement that ensures the uniqueness property, we observe that it yields a non-trivial improvement for the uncapacitated case. On average, about 11\% to 27\% gap reduction is achieved, compared to standard parametric C\&CG. Especially for difficult instances, the reduction can be up to 50\%. Similar reduction in computational time is also observed when $\rho=0.4$. As for the large average computation time for $\rho=0.2$ in Figure \ref{fig:uncap_enhancement}.(b), the reason is that one more instance is exactly solved by this enhancement  before the time limit.  As for the capacitated case, which often naturally has the uniqueness property, this enhancement does not help much. Those observations confirm our previous understanding that the lack of uniqueness property causes parametric C\&CG less effective. Hence, research to attain this property, especially  structural property based strategies, is worth  further investigation.

\vspace{-10pt}
\subsection{Compute Complex Two-Stage RO}
\label{subsect_complexRO}
As we mentioned in Section \ref{subsect_complex_RO}, the developed C\&CG algorithms can be applied to solve two-stage RO formulations with complex DDU set or recourse problem, if powerful solvers are available to compute master problems and subproblems. In this subsection, we employ parametric C\&CG and make use of \texttt{Gurobi}'s exact solution capacity on mixed integer SOC and (small) bilinear programs to provide an illustration. We first consider an extension of $\mathbf{Robust \ FL}$ defined in \eqref{eq_robust_facility} with the following modified objective function.
\begin{equation*}
	\begin{split}
	\ \min_{(\mathbf x_c,\mathbf x_d)\in\mathcal X}  \sum_{j\in J} (f_jx_{d,j}+a_jx_{c,j}) + \max_{\mathbf u\in\mathcal U(\mathbf x)} \ \min_{(\mathbf y_1,\mathbf y_2,\mathbf z)\in\mathcal{Y}(\mathbf x,\mathbf u)} \sum_{i\in I}\sum_{j\in J}(c_{ij}-p_i)y_{ij}+\sum_{j\in J}s_j\left(\sum_{i\in I}y_{ij}\right)^2
	\end{split}
\end{equation*}
The quadratic term in the recourse problem represents a complex cost component, e.g., cost incurred by penalizing on-site waiting time that is captured by quadratic of the service volume.  Through replacing it by $\tilde \tau$ and augmenting $\mathcal{Y}(\mathbf x,\mathbf u)$ with constraint $\tilde \tau\geq \displaystyle\sum_{j\in J}s_j\left(\sum_{i\in I}y_{ij}\right)^2$, the recourse problem is converted into a typical SOC problem. As $\tilde \tau$ is unrestricted, it can be easily verified that this recourse problem satisfies Slater's condition and the strong duality property holds.

We also modify $\mathcal U^R(\mathbf x)$ defined in \eqref{ex_eq_DDU_RHS} to obtain the following DDU set, i.e., replacing the last two sets of constraints by an SOC constraint.
	\begin{align*}
		\mathcal U^{R^2}(\mathbf x)=\{(\mathbf u, \mathbf d): & \ u_i=\underline u_i(1+\tilde u_i) \ \forall i\\
		& \sum_{j\in J(i)}\underline{\xi}^i_jx_{d,j}\leq \tilde u_i\leq \sum_{j\in J(i)}\overline{\xi}^i_jx_{d,j} \ \forall i \\
		&\sum_{i\in I}\tilde u^2_i\leq \alpha \frac{\sum_{j\in J} x_{c,j}}{ u^0}
		\}
	\end{align*}
 For this SOC set, by setting parameters $\alpha> 0$ and $\underline \xi^i_j=0$ for $i\in I$ and $j\in J$, it is guaranteed that $\mathcal U^{R^2}(\mathbf x)$ has an interior point to ensure Slater's condition. So, the strong duality holds for any linear objective function.

 By adopting the aforementioned SOC recourse problem and uncertainty set respectively, we have two variants of $\mathbf{Robust \ FL}$, referred to as 	$\mathbf{Robust \ FL-RC}$ and $\mathbf{Robust \ FL-UC}$ respectively. Note that $\mathbf{Robust \ FL-RC}$ and $	\mathbf{Robust \ FL-UC}$ are naturally more challenging than  $\mathbf{Robust \ FL}$. As noted earlier, we present detailed modifications on top of  parametric C\&CG in Appendix \ref{Asect:V2SOC} to handle general two-stage RO formulation with such SOC recourse problem or DDU set.

\begin{figure}
\begin{subfigure}[t]{0.45\textwidth}
	\pgfplotsset{tick label style={font=\small},
		label style={font=\small},
		legend style={font=\tiny}
	}
\begin{tikzpicture}
\begin{axis} [height=6cm,width=7.5cm,
    legend style={at={(0.5,1.2)},anchor=north,legend columns=-1,/tikz/every even column/.append style={column sep=5pt}},
    ybar=0.5pt,
    bar width = 6pt,
    ymin = 0,
    ymax = 700,
    symbolic x coords={25L,25H},
    xtick = data,
    enlarge x limits = 0.5,
    nodes near coords,
    every node near coord/.append style={rotate=90, anchor=west, font=\tiny}
]
\pgfkeys{/pgf/number format/.cd,fixed,precision=2};
\addplot[draw = blue!50, fill=blue!50, point meta=explicit symbolic]
table[meta=label]{
	x		y		label
	25H		269.57	{269.57 (2)}
	25L		468.44	{468.44 (2)}
};
\addplot[draw = red!50, fill=red!50, point meta=explicit symbolic]
table[meta=label]{
	x		y		label
	25H		2.79	{2.79 (3)}
	25L		3.45	{3.45 (3)}
};
\legend { $\mathbf{Robust \ FL-RC}$, \ $\mathbf{Robust \ FL}$};
\end{axis}
\end{tikzpicture}
\caption{Different Recourse Problems}
\end{subfigure}
\begin{subfigure}[t]{0.45\textwidth}
\center
\pgfplotsset{tick label style={font=\small},
	label style={font=\small},
	legend style={font=\tiny}
}
\begin{tikzpicture}
\begin{axis} [height=6cm,width=7.5cm,
    legend style={at={(0.5,1.2)},anchor=north,legend columns=-1,/tikz/every even column/.append style={column sep=5pt}},
    ybar=0.5pt,
    bar width = 6pt,
    ymin = 0,
    ymax = 36,
    symbolic x coords={25L,25H},
    xtick = data,
    enlarge x limits = 0.5,
    nodes near coords,
    every node near coord/.append style={rotate=90, anchor=west, font=\tiny}
]
\pgfkeys{/pgf/number format/.cd,fixed,precision=2,set thousands separator={}};
\addplot[draw = blue!50, fill=blue!50, point meta=explicit symbolic]
table[meta=label]{
	x		y		label
	25H		10.47	{10.47 (2.6)}
	25L		24.03	{24.03 (2.8)}
};
\addplot[draw = red!50, fill=red!50, point meta=explicit symbolic]
table[meta=label]{
	x		y		label
	25H		4.66	{4.66 (3.6)}
	25L		6.62	{6.62 (3.6)}
};
\legend {$\mathbf{Robust \ FL-UC}$, \ $\mathbf{Robust \ FL}$};
\end{axis}
\end{tikzpicture}
\caption{Different Uncertainty Sets}
\end{subfigure}
\caption{Average Computational Time (s) of Complex 2-Stage RO
\label{fig_complexRO}}
\end{figure}
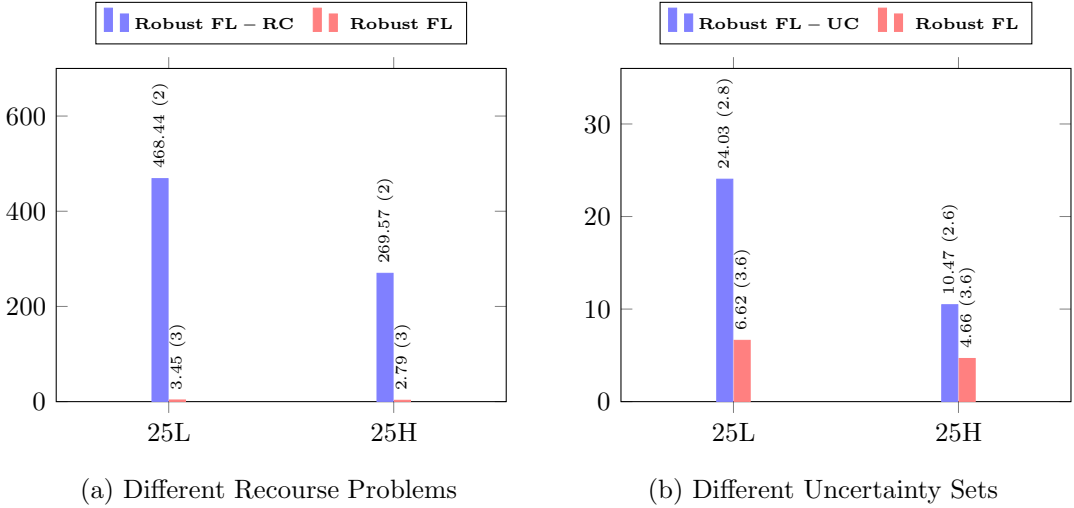

 We observe that quadratic term $\displaystyle\sum_{j\in J}s_j\left(\sum_{i\in I}y_{ij}\right)^2$ in $\mathbf{Robust \ FL-RC}$ may cause numerical challenges to the solver in our empirical  study. To alleviate this issue, we
  modify testing instances of 25 sites used in Section \ref{subsect_compu_3ccg} for experiments on $\mathbf{Robust \ FL-RC}$. Specifically, both upper and lower bounds in \eqref{eq_RFL_X} and basic demand in \eqref{eq_DDU_RHS0} are set to $1/10$ of their original values. Testing instances of $\mathbf{Robust \ FL-UC}$ are simply those of 25 sites with $\underline{\xi}^i_j$ set to $0$ for all $i$ and $j$. Detailed computational results on these two variants are reported in Tables \ref{tbl:SOCPR} and \ref{tbl:SOCPU} in Appendix \ref{Asect_tables},  and the average computational times are displayed in Figure \ref{fig_complexRO}, along with those from computing original $\mathbf{Robust \ FL}$ for comparison.


 Based on numerical results,  we note that the SOC recourse problem does not require more C\&CG iterations. Indeed,  less number of iterations are typically involved. Similar observations are also made on $\mathbf{Robust \ FL-UC}$.
It clearly confirms the strength of the algorithm on dealing with SOC recourse problem and SOC uncertainty set.
Nevertheless, as noted earlier, the computational burden imposed by the SOC structure is drastically heavier. On average, the computational time of $\mathbf{Robust \ FL-RC}$ is about two orders of magnitude more than that of  $\mathbf{Robust \ FL}$. Comparatively, the computational burden associated with SOC uncertainty set is not so demanding, which still is much more than that of polyhedron uncertainty set.  Moreover,  among $\mathbf{Robust \ FL-RC}$,  $\mathbf{Robust \ FL-UC}$ and $\mathbf{Robust \ FL}$, the first two are clearly very sensitive to parameters. It indicates that existing methods for linear structures, e.g., polyhedral DDU set and LP recourse problem, demonstrate a stronger stability and scalability, which, however, is not the case for nonlinear structures.

\section{Conclusion}
In this paper, we present a systematic study on two-stage robust optimization subject to DDU,
whose general form basically remains unsolved in the current literature. Our main results include three sophisticated  variants of C\&CG method to exactly compute DDU-based two-stage RO, along with rigorous and novel analyses on their convergence and iteration complexity. We highlight those analyses as they present a novel showcase of  a few core linear programming concepts in understanding and solving more complex optimization paradigms. Also, it is interesting to note that DDU-based two-stage RO theoretically might not be more demanding than its DIU-based counterpart in terms of the iteration complexity. Actually, a counterintuitive  discovery is that converting a DIU set into a DDU set by making use of ``deep knowledge'' and then computing the resulting DDU-based formulation may lead to a significant improvement.
Other noteworthy results include basic structural properties of two-stage RO with DDU, an approximation scheme to deal with mixed integer recourse, and a couple of enhancement techniques for the developed algorithms. We also conduct an organized  numerical study on all developed algorithms and techniques to evaluate their performances and appreciate their strength.
	
Currently, we have observed that two-stage RO has become a major decision making model in many real systems, and basic C\&CG is a mainstream solution method for its DIU-based instances. The presented study naturally extends and complements previous research work in a unified manner, and remarkably broadens our capability to handle difficult uncertainty issues in decision making.

We have noted a few possible research tasks throughout the paper. Among them, probably the most critical one is to  advance the methodology on utilizing ``deep knowledge'' to improve our solution capacity. New strategies to derive exact or highly effective knowledge and to represent such knowledge analytically are fundamental topics to work on. We mention that intelligent subroutines taking advantage of feedback to  generate, revise and represent deep knowledge within the iterative procedure is of a particular interest, which may help us achieve automated fast computation for large scale instances in the future.


\setlength{\bibsep}{0pt plus 0.3ex}
\bibliographystyle{plain}
\bibliography{robustUC}

\newpage
\begin{center}
	\textbf{\Large{Appendix}}
\end{center}
\setcounter{section}{0}

\setcounter{equation}{0}
\def\theequation{A.\arabic{equation}}

\renewcommand{\thesection}{A.\arabic{section}}
\section{Proofs of Section \ref{sect_general_formulation}}
\label{apd:proofS2}

%

\noindent\textbf{Proposition \ref{prop_reform_DIU}.} 
Formulation $\mathbf{2-Stg \ RO}$ with  the DDU set defined in \eqref{eq_US_binary} is equivalent to the following DIU-based 2-stage RO
$$\min_{\mathbf x\in \mathcal{X}}\mathbf c_1\mathbf x+\max_{\mathbf u\in \mathcal{U}^0}\ \min_{\mathbf y\in\bar{\mathcal Y}(\mathbf x,\mathbf u)}\mathbf c_2\mathbf y$$
where $\mathcal{U}^0=\{\mathbf {u}\in \{0,1\}^{m_u}: \mathbf{Fu}\leq \mathbf{h}\}$ and
\begin{eqnarray}
\label{eq_reformulation_ybar1}
	\bar{\mathcal Y}(\mathbf x,\mathbf u)=\{\mathbf y\in\mathbb Z^{m_y}_+\times \mathbb{R}^{n_y}_+:\mathbf B_2\mathbf y\geq\mathbf d-\mathbf B_1\mathbf x-\mathbf E(\mathbf u\circ \mathbf x')\},
\end{eqnarray}
or equivalently
\begin{eqnarray}
\begin{split}
	\label{eq_reformulation_ybar2}
	(\mathbf y, \mathbf v)\in\bar{\mathcal Y}(\mathbf x,\mathbf u)&=&\{\mathbf y\in\mathbb Z^{m_y}_+\times \mathbb{R}^{n_y}_+, \mathbf v\in \mathbb{R}^{m_u}_+:\mathbf B_2\mathbf y\geq\mathbf d-\mathbf B_1\mathbf x-\mathbf E\mathbf v\\
	&& \ \mathbf v\leq \mathbf x', \mathbf v\leq \mathbf u, \mathbf {v}\geq
	\mathbf {x}'+\mathbf u-\mathbf{1}\}.
\end{split}
\end{eqnarray}

\begin{proof}
	We prove first that the original two-stage RO with that DDU is equivalent to
	\begin{eqnarray*}
		\label{eq_neutralization_nolinearization}
		\min_{\mathbf x\in \mathcal{X}}\mathbf c_1\mathbf x+\max_{\mathbf u\in \mathcal{U}^0}\ \min_{\mathbf y\in \bar{\mathcal Y}(\mathbf x,\mathbf u)}\mathbf c_2\mathbf y
	\end{eqnarray*}
with $\bar{\mathcal Y}(\mathbf x,\mathbf u)$ defined in \eqref{eq_reformulation_ybar1}. 
	
	Claim 1: If $\mathbf u\in \mathcal U(\mathbf x)$, we have $\mathbf u\circ\mathbf x'=\mathbf u$.
	\begin{proof}[Proof of Claim 1:] If $ x'_i=1$, $ u_ix'_i= u_i$. Otherwise, i.e., $x'_i=0$, it follows from \eqref{eq_US_binary} that $u_i=0=u_ix'_i$. Hence, the claim is valid.
	\end{proof}
	Claim 2: If $\mathbf u\in \mathcal{U}^0$ and 
	$\mathbf u\notin \mathcal U(\mathbf x)$,  we have $\mathbf u\circ\mathbf x'\in \mathcal U(\mathbf x)$.
	\begin{proof}[Proof of Claim 2:] Given $u_i\in \{0,1\}$, $u_ix'_i\leq x'_i$ for all $i$. Similarly, we have $u_ix'_i\leq u_i$ for all $i$. Then, it follows from the downward-closedness property that $\mathbf F(\mathbf u\circ\mathbf x')\leq \mathbf h$. Hence, $\mathbf u\circ\mathbf x'\in \mathcal{U}(\mathbf x)$.
	\end{proof}
	Claim 3: For a given $\mathbf x\in \mathcal{X}$, we have
	$$\max_{\mathbf u\in \mathcal{U}(\mathbf x)}\ \min_{\mathbf y\in \mathcal Y(\mathbf x,\mathbf u)}\mathbf c_2\mathbf y = \max_{\mathbf u\in \mathcal{U}^0}\ \min_{\mathbf y\in \bar{\mathcal Y}(\mathbf x,\mathbf u)}\mathbf c_2\mathbf y.$$
	\begin{proof}[Proof of Claim 3:] Given that Claim 1 and Claim 2 are valid, this claim follows directly.
	\end{proof}
	As a result of Claim 3, the aforementioned equivalence to $\mathbf{2-Stg \ RO}$ is valid. Then, using standard linearization technique to linearize $\mathbf u\circ \mathbf x'$, i.e., introducing $ v_i$ to replace $u_i x'_i$ for all $i$ and including related constraints, \eqref{eq_reformulation_ybar1} is converted into \eqref{eq_reformulation_ybar2} and the expected result is derived.
\end{proof}

\noindent\textbf{Lemma \ref{lem_switching}.} 
For $\mathbf x\in \mathcal{X}$ such that $\mathcal Y(\mathbf x)$ is non-empty, we have
$$\max_{\mathbf u \in \mathcal U(\mathbf x)} \min_{\mathbf y\in \mathcal Y(\mathbf x)} \mathbf{ \hat c}_2(\mathbf u) \mathbf y = \min_{\mathbf y \in \mathcal Y(\mathbf x)} \max_{\mathbf u\in \mathcal U(\mathbf x)} \mathbf {\hat c}_2(\mathbf u) \mathbf y. $$
\begin{proof}
	By the definition of $\mathbf {\hat c}_2(\mathbf u)$,  objective function $\mathbf{\hat c}_2(\mathbf u) \mathbf y$ is linear in both $\mathbf u$ and $\mathbf y$. Also, both $\mathcal{Y}(\mathbf x)$ of the recourse problem, if not empty, and uncertainty set $\mathcal U(\mathbf x)$ are polyhedra. Moreover, according to assumption (\textit{A2}), $\mathcal U(\mathbf x)$  is bounded. Hence, the expected equality follows from the minimmax theorem \citep{sion1958general,rockafellar2015convex}.
\end{proof}

\noindent\textbf{Proposition \ref{prop_RO_switching}.}
	Two-stage RO with the DDU set defined in \eqref{eq_uncer_set} and the recourse problem defined in \eqref{eq_switching_recourse} is equivalent to
	\begin{eqnarray}
		\label{A:eq_equiv_switching}
		\min_{\mathbf x\in \mathcal{X}, \mathbf y\in \mathcal Y(\mathbf x)} \mathbf c_1\mathbf x + \mathbf c_2\mathbf y+\max_{\mathbf u\in \mathcal{U}(\mathbf x)} (\hat{\mathbf{E}}\mathbf u)^\intercal\mathbf y,
	\end{eqnarray}
	which is further equivalent to the following single-level optimization problem.
	\begin{eqnarray*}
		\min && \mathbf c_1\mathbf x+\mathbf c_2\mathbf y+ \mathbf{(h+Gx)}^\intercal\boldsymbol\lambda\\
		\mathrm{s.t.} && \mathbf x\in \mathcal{X}, \mathbf y\in \mathcal Y(\mathbf x)\\
		&& \mathbf F(\mathbf x)^\intercal\boldsymbol\lambda  - \mathbf{\hat E}^\intercal\mathbf y \geq \mathbf 0\\
		&& \boldsymbol\lambda \geq \mathbf{0}
	\end{eqnarray*}
\begin{proof}
	It is clear that the min-max bilevel equivalence follows directly from Lemma~\ref{lem_switching}. Note that for a fixed $(\mathbf x, \mathbf y)$, the maximization problem in that equivalence is actually a linear program with respect to $\mathbf u$. Specifically, it is
	$$\max \{ ((\mathbf{\hat E u})^\intercal+\mathbf c_2)\mathbf y: \mathbf F(\mathbf x)\mathbf u\leq \mathbf h+\mathbf{Gx}, \mathbf u\geq \mathbf 0\}=\mathbf c_2\mathbf y + \max \{ (\mathbf{\hat E u})^\intercal\mathbf y: \mathbf F(\mathbf x)\mathbf u\leq \mathbf h+\mathbf{Gx}, \mathbf u\geq \mathbf 0\}.$$
	Let $\boldsymbol\lambda$ denote dual variables. That linear program has the following dual problem.
	$$\min \{ (\mathbf h+\mathbf{Gx})^\intercal\boldsymbol{\lambda}: \mathbf F(\mathbf x)^\intercal\boldsymbol{\lambda} \geq \mathbf{\hat E}^\intercal\mathbf y, \boldsymbol {\lambda}\geq \mathbf 0\}$$
	Then, by replacing that maximization problem in the min-max formulation by its dual problem (and simply dropping the minimization sign) and reorganizing slightly, we have the equivalent single-level optimization.
\end{proof}

\section{Proofs and Reformulations of Section \ref{sect_algorithms}}
\label{apd:proofS3}

\noindent\textbf{Theorem \ref{thm_PI_enumeration}.}
	Let $\Sigma_{\Pi}=(\mathcal{P}_{\Pi}, \mathcal{R}_{\Pi})$ with $\mathcal{P}_{\Pi}$ and $\mathcal{R}_{\Pi}$ being the sets of extreme points and extreme rays of $\Pi$, respectively.  The two-stage RO in \eqref{eq_2RO} (and  the equivalence in \eqref{eq_2RO_NC_bilevel}) is equivalent to a bilevel linear optimization program as in the following.
	\begin{subequations}
		\label{A:eq_2stgRO_PI}
		\begin{align}
			\mathbf{2-Stg \ RO(\Sigma_{\Pi})}: \   w^*=\min \ & \ \mathbf{c}_1\mathbf x+ \eta\\
			\mbox{s.t.} \ & \ \mathbf x\in \mathcal{X}\\
			&  \bigg\{\eta\geq  (\mathbf d-\mathbf B_1\mathbf x)^\intercal\boldsymbol\pi + \max_{\mathbf u\in \mathcal{U}(\mathbf x)} \{(-\mathbf E\mathbf u)^\intercal\boldsymbol\pi\}\bigg\} \forall \boldsymbol\pi\in \mathcal{P}_{\Pi} \label{A:eq_enu_optimality}\\
			&  \bigg\{(\mathbf d-\mathbf B_1\mathbf x)^\intercal\boldsymbol\gamma+\max_{\mathbf v\in \mathcal{U}(\mathbf x)}(\mathbf E\mathbf v)^\intercal\boldsymbol\gamma \leq 0\bigg\} \forall \boldsymbol\gamma\in \mathcal{R}_{\Pi}. \label{A:eq_enu_feasibility}
		\end{align}
	\end{subequations}
\begin{proof}
	Note that it is sufficient to show the equivalence between the epigraph reformulation of \eqref{eq_2RO_NC_bilevel}, i.e.,
	\begin{eqnarray}
		\label{eq_bilinear_epi}
		\begin{split}
			\min & \ \mathbf{c}_1\mathbf x+ \eta\\
			\mbox{s.t.} & \ \mathbf x\in \mathcal{X}\\
			& \eta\geq \max\{(\mathbf d-\mathbf B_1\mathbf x-\mathbf E\mathbf u)^\intercal\boldsymbol\pi: \mathbf{u}\in \mathcal{U}(\mathbf x), \mathbf B_2^\intercal\boldsymbol\pi \leq \mathbf c_2, \boldsymbol\pi\geq \mathbf 0\}
		\end{split}
	\end{eqnarray}
	and $\mathbf{2-Stg \ RO(\Sigma_{\Pi})} $.
	
	When $\mathbf x$ renders the recourse problem feasible, according to Corollary \ref{cor_maxmax_dual}, we have \begin{equation}
		\begin{split}
			&\max\{(\mathbf d-\mathbf B_1\mathbf x-\mathbf E\mathbf u)^\intercal\boldsymbol\pi: \mathbf{u}\in
			\mathcal{U}(\mathbf x), \mathbf B_2^\intercal\boldsymbol\pi \leq \mathbf c_2, \boldsymbol\pi\geq \mathbf{0}\}\\
			&=\max\{(\mathbf d-\mathbf B_1\mathbf x-\mathbf E\mathbf u)^\intercal\boldsymbol\pi: \mathbf{u}\in
			\mathcal{U}(\mathbf x),  \boldsymbol\pi\in \mathcal{P}_{\Pi}\}\\
			& = \max_{\boldsymbol\pi\in \mathcal{P}_{\Pi}}\{\max\{(\mathbf d-\mathbf B_1\mathbf x-\mathbf E\mathbf u)^\intercal\boldsymbol\pi: \mathbf{u}\in
			\mathcal{U}(\mathbf x)\}\}
		\end{split}
	\end{equation}
	
	Similarly, by Corollary \ref{cor_maxmax_dual}, we have that $\mathbf x$ renders the recourse problem feasible if and only if all constraints in \eqref{A:eq_enu_feasibility}  are satisfied.
	
	With those two arguments being held, it is clear that  \eqref{eq_bilinear_epi} is equivalent to \eqref{eq_2stgRO_PI}.
\end{proof}

\noindent\textbf{Primal-dual Based Representations of $\mathcal{OU}$ and $\mathcal{OV}$}\\
In addition to KKT conditions based representation for $\mathcal{OU}(\mathbf x,\boldsymbol\pi^k)$ and $\mathcal{OV}(\mathbf x,\boldsymbol\gamma^l)$ in \eqref{eq_KKT_point}, we can employ both primal and dual formulations of $\mathbf{LP}(\mathbf x,\boldsymbol\beta)$ and the strong duality property to define those optimal solution sets. Specifically, we have
\begin{eqnarray}
	\label{eq_dual_point}
	\mathcal{OU}(\mathbf x, \boldsymbol\pi^k) =  \left\{\begin{array}{l}
		\mathbf F(\mathbf x)\mathbf u^k\leq \mathbf h+\mathbf{Gx}\\
		\mathbf F(\mathbf x)^\intercal\boldsymbol\lambda^k\geq - \mathbf E^\intercal\boldsymbol\pi^k \\
		(-\mathbf{Eu})^\intercal\boldsymbol\pi^k\geq (\mathbf h+\mathbf{Gx})^\intercal\boldsymbol\lambda^k\\
		\mathbf u^k\geq \mathbf{0}, \boldsymbol\lambda^k\geq \mathbf{0}
	\end{array}\right\}
\end{eqnarray}
As for $\mathcal{OV}(\mathbf x,\boldsymbol\gamma^l)$, it can be obtained through replacing $\mathbf u^k$ by $\mathbf v^l$ and $\boldsymbol\lambda^k$ by $\boldsymbol\zeta^l$ in \eqref{eq_dual_point}. 

We note that if those primal-dual based representations of $\mathcal{OU}$ and $\mathcal{OV}$ are adopted in the master problem of Variant 1 or 2, the nonlinear terms in the forms of 
$\mathbf x^\intercal\boldsymbol{\lambda}^k$ or $\mathbf x^\intercal\boldsymbol{\zeta}^l$ can be easily linearized for computational benefit when $\mathbf x$ is binary or can be represented by its binary expansion.  \hfill$\square$
\\

\noindent\textbf{Theorem 	\ref{thm_V1_convergence}.}
	Suppose $T\!O\!L=0$. When terminates, Variant 1 either reports that $\mathbf{2-Stg \ RO}$ in \eqref{eq_2RO} is infeasible,  or converges to its optimal value and exact solution.
\begin{proof} From Corollaries \ref{cor_KKT_BendersRe} and \ref{cor_Pi_partial_KKT}, it is clear that if $\mathbf{MP1}$ is infeasible, so is $\mathbf{2-Stg \ RO}$ in \eqref{eq_2RO}. Moreover, according to Corollary  \ref{cor_maxmax_dual} and $\mathbf {SP3}$ in \eqref{eq_SP3}, if computing $\mathbf {MP1}$ derives an $\mathbf x$ that is infeasible to $\mathbf{2-Stg \ RO}$, an extreme ray of $\Pi$ will be identified.
	
	\noindent Claim 1: If an extreme ray is output by solving $\mathbf{SP3}$ in some iteration, it will not be derived  in any following iterations.
	\begin{proof}[Proof of Claim 1:]
		Assume that $\mathbf x^*$ is one solution derived from computing $\mathbf{MP1}$ in some iteration that is infeasible to $\mathbf{2-Stg \ RO}$, and $\boldsymbol\gamma^*\in \mathcal{R}_{\Pi}$ is the corresponding extreme ray from solving $\mathbf{SP3}$.
		
		Note that $\mathcal{OV}(\mathbf x, \boldsymbol\gamma^*)$ is the KKT reformulation of  $\mathbf{LP}(\mathbf x,\boldsymbol\gamma^*)$. After the feasibility cutting set \eqref{eq_BD_feasibility} defined for $\boldsymbol\gamma^*$ is included into master problem $\mathbf {MP1}$, for any feasible $\mathbf x$ we have
		$$(\mathbf d -\mathbf B_1\mathbf x-\mathbf {Ev})^\intercal\boldsymbol\gamma^*\leq 0 \ \ \forall  \mathbf v\in \mathcal U(\mathbf x).$$
		It means $\boldsymbol\gamma^*$ is no longer an extreme ray along which $\mathbf {SP3}$ becomes unbounded, which proves Claim 1.
	\end{proof}
	If $\mathbf{2-Stg \ RO}$ is infeasible, it is clear that $\mathbf{SP3}$ is always called before  master problem $\mathbf{MP1}$ becomes infeasible. Given the fact that extreme ray set $\mathcal{R}_{\Pi}$ is fixed and finite, it follows from Claim 1 that $\mathbf {MP1}$ will certainly become infeasible after being augmented by a finite number of feasibility cutting sets.
	
	Next, we prove the algorithm converges to an exact solution of $\mathbf{2-Stg \ RO}$ if it is feasible.  To simplify our proof, let $\mathbf{MP}^{t}$ denote master problem $\mathbf{MP}$ in $t$-th iteration throughout this appendix.  
\\
	Claim 2: If an extreme point of $\Pi$ appears twice in an optimal solution to $\mathbf{SP2}$ in \eqref{eq_SP2}, we have $UB\leq LB$.
	\noindent \begin{proof}[Proof of Claim 2:] Assume that $\boldsymbol\pi^*\in \mathcal{P}_{\Pi}$ is derived as an optimal solution to $\mathbf{SP2}$ in iterations $t_1$ and $t_2$ with $t_1 < t_2$, and $\mathbf x^*$ is optimal to $\mathbf{MP1}^{t_2}$. Accordingly, we have the following result for $\mathbf{SP2}$ in iteration $t_2$.
		$$\eta_s(\mathbf x^*) = \max_{\mathbf u\in \mathcal{U}(\mathbf x^*)} \{(\mathbf d-\mathbf B_1\mathbf x^*-\mathbf E\mathbf u)^\intercal\boldsymbol\pi^*\}=(\mathbf d-\mathbf B_1\mathbf x^*)^\intercal\boldsymbol\pi^*+ \max_{\mathbf u\in \mathcal{U}(\mathbf x^*)}\{-(\mathbf E\mathbf u)^\intercal\boldsymbol\pi^*\}.$$

		On the one hand, according to Step 5 of the algorithm, we have $UB\leq  \mathbf{c}_1\mathbf x^*+\eta_s(\mathbf x^*)$. On the other hand, the optimality cutting set \eqref{eq_BD_optimality} defined with respect to $\boldsymbol\pi^*$ has been a part of constraints of $\mathbf{MP1}^{t}$ for $t\geq t_1+1$, after  $\boldsymbol\pi^*$ is derived 
		by computing $\mathbf{SP2}$ in iteration $t_1$. Hence, we have
		in iteration $t_2$ that
		\begin{align*}
			LB=\underline w = \mathbf{c}_1\mathbf x^*+\eta&\geq  \mathbf{c}_1\mathbf x^*+ (\boldsymbol\pi^*)^\intercal\mathbf d-(\boldsymbol\pi^*)^\intercal\mathbf B_1\mathbf x^* - (\boldsymbol\pi^*)^\intercal\mathbf E\mathbf u^{\boldsymbol\pi^*}\\
			&\quad\quad (\mathbf u^{\boldsymbol\pi^*}, \boldsymbol\lambda^{\boldsymbol\pi^*}) \in \mathcal{OU}(\mathbf x^*, \boldsymbol\pi^*)\\
			&= \mathbf{c}_1\mathbf x^*+(\mathbf d-\mathbf B_1\mathbf x^*)^\intercal\boldsymbol\pi^* + \max_{\mathbf u\in \mathcal{U}(\mathbf x^*)} \{-(\mathbf E\mathbf u)^\intercal\boldsymbol\pi^*\}\\
			&= \mathbf{c}_1\mathbf x^*+\eta_s(\mathbf x^*)\\
			&\geq UB \qedhere
		\end{align*}
	\end{proof}

Given that $\mathbf{MP1}$ is a relaxation to $\mathbf{2-Stg \ RO}$ and its optimal value is always  less than or equal to $UB$, we have $LB=UB$ and $\mathbf x^*$ is optimal to $\mathbf{2-Stg \ RO}$. Together with the fact that extreme point set $\mathcal{P}_{\Pi}$ is fixed and finite, it is straightforward to conclude that the algorithm converges to an exact solution after
a finite number of iterations. 
\end{proof}

\noindent\textbf{Lemma \ref{lem:parametric:maxmin}.}
For a given $\mathbf x$, we have
\begin{align*}
	\max_{\mathbf u\in \mathcal{U}(\mathbf x)}\min\{\mathbf c_2\mathbf y: \mathbf y\in \mathcal{Y}(\mathbf x, \mathbf u)\}=\max_{\mathbf u\in \ \mathcal{U^*}(\mathbf x)\bigcup \mathcal{V^*}(\mathbf x)}  \min\{\mathbf c_2\mathbf y: \mathbf y\in \mathcal{Y}(\mathbf x, \mathbf u)\} 	
\end{align*}
where $\mathcal{U^*}(\mathbf x)=\bigcup \limits_{k=1}^{K_p} \mathcal{OU}_{\mathbf u}(\mathbf x, \boldsymbol\pi^k)$ and $\mathcal{V^*}(\mathbf x)=\bigcup\limits_{l=1}^{K_r} \mathcal{OV}_{\mathbf u}(\mathbf x, \boldsymbol\gamma^l)$.

\begin{proof}
	To simplify our proof, let $\eta_l(\mathbf x)$ and $\eta_r(\mathbf x)$ represent the optimal values of the problems on LHS and RHS of the equal sign, respectively.
	
	Given the fact that $\mathcal{OU}_{\mathbf u}(\mathbf x, \boldsymbol\pi^k)\in \mathcal{U}(\mathbf x)$, $k=1,\dots, K_p$ and $\mathcal{OV}_{\mathbf u}(\mathbf x, \boldsymbol\gamma^l)\in \mathcal{U}(\mathbf x)$, $l=1,\dots, K_r$, it is straightforward that
	$\mathcal{U^*}(\mathbf x)\bigcup \mathcal{V^*}(\mathbf x)\subseteq \mathcal{U} (\mathbf x)$, which results in
	$\eta_l(\mathbf x)\geq \eta_r(\mathbf x)$.
	
	Next, we consider the other direction.
	For the case where $\eta_l(\mathbf x)<+\infty$, i.e., the recourse minimization problem is feasible for any $\mathbf u\in \mathcal{U}(\mathbf x)$,  by Proposition \ref{prop_extremepoint}, Corollary \ref{cor_maxmax_dual} and their remarks,  there exists $\mathbf u^*\in \mathcal{U}(\mathbf x)$ optimal to the left-hand-side $\max-\min$ problem that is also optimal to $\mathbf{LP} (\mathbf x, \boldsymbol\pi^*)$ in \eqref{eq_LP_parametric}
	with $\boldsymbol\pi^*$ being one extreme point of $\Pi$. 	 Given the definition of KKT conditions based set $\mathcal{OU}(\mathbf x, \boldsymbol\pi^*)$ in \eqref{eq_KKT_point}, there naturally exists  $\boldsymbol\lambda^*$
	such that $(\mathbf u^*, \boldsymbol\lambda^*) \in \mathcal{OU}(\mathbf x, \boldsymbol\pi^*)$. Hence, by the definition of projection, $\mathbf u^* \in \mathcal{U^*}(\mathbf x)$.
	
	Similarly, for the case where $\eta_l(\mathbf x)=+\infty$, i.e., there exists some $\mathbf u^*\in \mathcal{U}(\mathbf x)$ for which the  recourse minimization problem is infeasible,  there exists $\mathbf u^*$ that is optimal to  $\mathbf{LP} (\mathbf x, \boldsymbol\gamma^*)$ with $\boldsymbol\gamma^*$ being one extreme ray of $\Pi$, and we have $(\mathbf d-\mathbf B_1\mathbf x-\mathbf{Eu}^*)^\intercal\boldsymbol\gamma^* > 0$. Again, there exists optimal dual solution $\boldsymbol\zeta^*$ such that $(\mathbf u^*, \boldsymbol\zeta^*) \in \mathcal{OV}(\mathbf x, \boldsymbol\gamma^*)$. Hence, $\mathbf u^* \in \mathcal{V^*}(\mathbf x)$.
	
	Consequently, for both discussed cases,  we have $\eta_l(\mathbf x)\leq \eta_r(\mathbf x)$. Together with $\eta_l (\mathbf x)\geq \eta_r(\mathbf x)$, the expected conclusion follows.
\end{proof}

\noindent\textbf{Proposition \ref{prop_V2betterLB}.} 
Assume that the same sets $\mathcal{\hat P}_{\Pi}\subseteq \mathcal{P}_{\Pi}$ and $\mathcal{\hat R}_{\Pi}\subseteq \mathcal{R}_{\Pi}$ have been included for defining sets $\mathcal{OU}$ and $\mathcal{OV}$  in both $\mathbf{MP1}$ and $\mathbf{MP2}$. Then, the optimal value of $\mathbf{MP1}$ is an underestimation of that of $\mathbf{MP2}$.
\begin{proof}
	We first show that any first stage decision  feasible to $\mathbf{MP2}$ is also feasible to $\mathbf{MP1}$. Then, we prove that for any first stage decision that is feasible to $\mathbf{MP1}$, its corresponding value in $\mathbf{MP1}$ is less than or equal to its corresponding value in $\mathbf{MP2}$. Clearly, if both statements are valid, the expected conclusion follows.
	
	Claim 1: If $\mathbf x'\in \mathcal{X}$ is a first stage decision that is feasible to $\mathbf{MP2}$, it is also feasible to $\mathbf{MP1}$.	
	\begin{proof}[Proof of Claim 1:]
		Comparing $\mathbf{MP1}$ and $\mathbf{MP2}$, it is sufficient to show that
		$$\boldsymbol\gamma^\intercal\mathbf d-\boldsymbol\gamma^\intercal\mathbf B_1\mathbf x' - \boldsymbol\gamma^\intercal\mathbf E\mathbf v^{\boldsymbol\gamma} \leq 0$$
		with $(\mathbf v^{\boldsymbol\gamma}, \boldsymbol\zeta^{\boldsymbol\gamma})\in \mathcal{OV}(\mathbf x', \boldsymbol\gamma)$ for all $\boldsymbol\gamma\in \mathcal{\hat R}_{\Pi}$.
		
		Consider an arbitrary $\boldsymbol\gamma\in \mathcal{\hat R}_{\Pi}$. Because there exits $\mathbf y^{\boldsymbol\gamma}\geq 0$ such that $\mathbf{B}_2\mathbf y^{\boldsymbol\gamma}\geq \mathbf{d-\mathbf B}_1\mathbf x'-\mathbf{Ev}^{\boldsymbol\gamma,0}$ with $\mathbf{v}^{\boldsymbol\gamma,0}\in \mathcal{OV}(\mathbf x', \boldsymbol\gamma)$ in $\mathbf{MP2}$, the following optimization problem is feasible.
		$$\min\{\mathbf c_2\mathbf y: \mathbf{B}_2\mathbf y\geq \mathbf{d-\mathbf B}_1\mathbf x'-\mathbf{Ev}^{\boldsymbol\gamma,0},  \ \mathbf y\geq 0\}$$
		Hence its dual problem is bounded, i.e.,
		$$\max\{\boldsymbol\pi^\intercal(\mathbf{d-\mathbf B}_1\mathbf x'-\mathbf{Ev}^{\boldsymbol\gamma,0}): \boldsymbol\pi\in \mathcal{R}_{\Pi} \}\leq 0,$$
		recalling that $\mathcal{R}_{\Pi}$ is the extreme ray set of $\Pi=\{\mathbf{B}^\intercal_2\boldsymbol\pi\leq \mathbf c^\intercal_2, \boldsymbol\pi\geq \mathbf 0\}$.
		By setting $\boldsymbol\pi=\boldsymbol\gamma$, we have
		$$\boldsymbol\gamma^\intercal(\mathbf{d-\mathbf B}_1\mathbf x'-\mathbf{Ev}^{\boldsymbol\gamma,0})\leq 0,$$
		which justifies the claim.
	\end{proof}
	
	Claim 2: Let $\mathbf x'\in \mathcal{X}$ denote an arbitrary first stage decision that is feasible to $\mathbf{MP1}$.  Also, let $\underline w_{1}(\mathbf x')$ and $\underline w_{2} (\mathbf x')$  be the optimal values of $\mathbf{MP1}$ and $\mathbf{MP2}$, respectively, when $\mathbf x=\mathbf x'$. We have $\underline w_{1}(\mathbf x') \leq \underline w_{2}(\mathbf x').$
	
	\begin{proof}[Proof of Claim 2:]
		Note that $\underline w_{2} (\mathbf x') = \infty$ when  $\mathbf x'$ is infeasible to $\mathbf{MP2}$. Hence, it is only needed to study the non-trivial case where  $\mathbf x'$ is also feasible to $\mathbf{MP2}$.

		Consider constraints defined with respect to $\boldsymbol\pi \in \mathcal{\hat P}_{\Pi}$, i.e., optimality cutting sets, in $\mathbf{MP1}$ and $\mathbf{MP2}$. They can be equivalently reformulated as
		\begin{eqnarray}
			\label{eq_v1v2_LB1}
			\eta \geq \min \{\boldsymbol\pi^\intercal(\mathbf{d-\mathbf B}_1\mathbf x'-\mathbf{Eu}^{\boldsymbol\pi}), \ (\mathbf{u}^{\boldsymbol\pi},\boldsymbol\lambda^{\boldsymbol\pi})\in \mathcal{OU}(\mathbf x',\boldsymbol\pi)\}
		\end{eqnarray}
		for  $\mathbf{MP1}$ and
		\begin{eqnarray}
			\label{eq_v1v2_LB2}
			\eta \geq \min \{\mathbf{c_2}\mathbf y^{\boldsymbol\pi}: \ \mathbf{B}_2\mathbf y^{\boldsymbol\pi}\geq \mathbf{d-\mathbf B}_1\mathbf x'-\mathbf{Eu}^{\boldsymbol\pi},  \ \mathbf y^{\boldsymbol\pi} \geq \mathbf 0, \  (\mathbf{u}^{\boldsymbol\pi},\boldsymbol\lambda^{\boldsymbol\pi})\in \mathcal{OU}(\mathbf x',\boldsymbol\pi)\}
		\end{eqnarray}
		for $\mathbf{MP2}$, respectively, for every $\boldsymbol\pi\in \mathcal{\hat P}_{\Pi}$. Indeed, \eqref{eq_v1v2_LB2} can be further reformulated in the following.
		\begin{eqnarray*}
			\eta&\geq &\min \{ \ \min\{\mathbf{c_2}\mathbf y^{\boldsymbol\pi}: \ \mathbf{B}_2\mathbf y^{\boldsymbol\pi}\geq \mathbf{d-\mathbf B}_1\mathbf x'-\mathbf{Eu}^{\boldsymbol\pi},  \ \mathbf y^{\boldsymbol\pi} \geq \mathbf 0\}, \  (\mathbf{u}^{\boldsymbol\pi},\boldsymbol\lambda^{\boldsymbol\pi})\in \mathcal{OU}(\mathbf x',\boldsymbol\pi)\}.
		\end{eqnarray*}
		Since $\mathbf x'$ is feasible to $\mathbf{MP2}$, the RHS of the above formulation has a finite optimal value. We assume without loss of generality that it achieves the minimum at $(\mathbf u^{\boldsymbol\pi,0}, \boldsymbol\lambda^{\boldsymbol\pi,0})\in \mathcal{OU}(\mathbf x',\boldsymbol\pi)$. We have
		\begin{eqnarray*}
			\eta &\geq & \min\{\mathbf{c_2}\mathbf y^{\boldsymbol\pi}: \ \mathbf{B}_2\mathbf y^{\boldsymbol\pi}\geq \mathbf{d-\mathbf B}_1\mathbf x'-\mathbf{Eu}^{\boldsymbol\pi,0},  \ \mathbf y^{\boldsymbol\pi} \geq \mathbf 0\}\\
			&=&  \max\{\tilde {\boldsymbol\pi}^\intercal(\mathbf{d-\mathbf B}_1\mathbf x'-\mathbf{Eu}^{\boldsymbol{\boldsymbol\pi},0}):  \mathbf{B}^\intercal_2\tilde {\boldsymbol\pi}\leq \mathbf c^\intercal_2, \tilde {\boldsymbol\pi}\geq \mathbf 0\}\\
			&=& \max\{\tilde {\boldsymbol\pi}^\intercal(\mathbf{d-\mathbf B}_1\mathbf x'-\mathbf{Eu}^{\boldsymbol{\boldsymbol\pi},0}): \tilde {\boldsymbol\pi}\in \mathcal{P}_{\Pi}\}\\
			&\geq & \min \{\ \max\{\tilde {\boldsymbol\pi}^\intercal(\mathbf{d-\mathbf B}_1\mathbf x'-\mathbf{Eu}^{\boldsymbol\pi}): \tilde {\boldsymbol\pi}\in \mathcal{P}_{\Pi}\}, \  (\mathbf{u}^{\boldsymbol\pi},\boldsymbol\lambda^{\boldsymbol\pi})\in \mathcal{OU}(\mathbf x',\boldsymbol\pi)\}\\
			&\geq & \min \{\boldsymbol\pi^\intercal(\mathbf{d-\mathbf B}_1\mathbf x'-\mathbf{Eu}^{\boldsymbol\pi}), \ (\mathbf{u}^{\boldsymbol\pi},\boldsymbol\lambda^{\boldsymbol\pi})\in \mathcal{OU}(\mathbf x',\boldsymbol\pi)\}.
		\end{eqnarray*}
		The first equality is obtained by the strong duality of LP, the second equality is valid because of the definition of $\mathcal{P}_{\Pi}$, the second inequality is valid due to $(\mathbf u^{\boldsymbol\pi,0}, \boldsymbol\lambda^{\boldsymbol\pi,0})\in \mathcal{OU}(\mathbf x',\pi)$, and the last inequality is valid just by fixing $\tilde{\boldsymbol{\pi}} =\boldsymbol{\pi}$.
		
		Noting that the last inequality is simply the RHS of \eqref{eq_v1v2_LB1}, it follows that \eqref{eq_v1v2_LB2} dominates \eqref{eq_v1v2_LB1}.
		Since it holds for every $\boldsymbol\pi\in \mathcal{\hat P}_{\Pi}$, we conclude that $\underline w_2(\mathbf x')\geq \underline w_1(\mathbf x')$.
	\end{proof}	
	
	Because $\mathbf x'$ is an arbitrary solution that is feasible to $\mathbf {MP1}$, we have the desired conclusion.
\end{proof}

\noindent\textbf{Lemma	\ref{lem_feasible_then_optimal}.}
Consider $\mathbf{LP}(\mathbf x^0, \boldsymbol\beta)$ for a fixed $\boldsymbol\beta$, and suppose that basis $\mathfrak{B}^0$ is an optimal basis, i.e., its BS with respect to $\mathcal{U}(\mathbf x^0)$ is a BFS and an optimal solution. If $\mathfrak{B}^0$'s BS with respect to $\mathcal{U}(\mathbf x^1)$ is feasible, i.e., a BFS, it is also optimal to $\mathbf{LP}(\mathbf x^1, \boldsymbol\beta)$.  Moreover, if $\mathfrak{B}^0$ yields the unique optimal solution to $\mathbf{LP}(\mathbf x^0, \boldsymbol\beta)$, it also yields the unique one to $\mathbf{LP}(\mathbf x^1, \boldsymbol\beta)$.
\begin{proof}
	To support our reasoning, we simplify the representation of $\mathbf{LP}(\mathbf x^0, \boldsymbol\beta)$ in \eqref{eq_LP_parametric} as
	\begin{equation}
		\label{eq_LP_parametric_simplified}
		\max\{\mathfrak c^{\boldsymbol\beta}\mathbf u: \mathbf F\mathbf u\leq \mathbf {h}^0, \mathbf u\geq \mathbf 0\}
	\end{equation}
	where $\mathfrak c^{\boldsymbol\beta}=-\boldsymbol\beta^\intercal\mathbf{E}$ and $\mathbf h^0=\mathbf h+\mathbf{Gx}^0$. 	
	According to basis $\mathfrak B^0$, its complement $\mathfrak N^0$, and the classical Simplex algorithm, we can reorder variables and represent $\mathbf u$ as $(\mathbf u_{\mathfrak B^0}, \mathbf u_{\mathfrak N^0})$, and columns of $\mathbf F$ and represent it as $[\mathbf F_{\mathfrak B^0}|\mathbf F_{\mathfrak N^0}]$.
	
	Given that $\mathfrak B^0$ is an optimal basis, we have
	$\mathbf u^0_{\mathfrak B^0}=\mathbf F_{\mathfrak B^0}^{-1}\mathbf {h}^0\geq \mathbf 0$
	and the optimal value is 
	$$\mathfrak c^{\boldsymbol\beta}_{\mathfrak B^0}\mathbf F_{\mathfrak B^0}^{-1}\mathbf {h}^0 + (\mathfrak c^{\boldsymbol{\beta}}_{\mathfrak N^0}-\mathfrak c^{\boldsymbol{\beta}}_{\mathfrak B^0}\mathbf F_{\mathfrak B^0}^{-1}\mathbf F_{\mathfrak N^0})\mathbf u_{\mathfrak N^0}$$
	with $(\mathfrak c^{\boldsymbol{\beta}}_{\mathfrak N^0}-\mathfrak c^{\boldsymbol\beta}_{\mathfrak B^0}\mathbf F_{\mathfrak B^0}^{-1}\mathbf F_{\mathfrak N^0})_j\leq 0$, i.e., the reduced cost of $u_j$ is non-positive for $j\in \mathfrak N^0$.  Note that  those reduced costs are independent of $\mathbf h^0$. Hence, even if $\mathbf x^0$ is replaced by $\mathbf x^1$, they remain non-positive. Let $\mathbf h^1=\mathbf h+\mathbf{Gx}^1$.  As long as $\mathbf u^1_{\mathfrak B^0}=\mathbf F_{\mathfrak B^0}^{-1}\mathbf {h}^1\geq \mathbf 0$, we have that $\mathfrak B^0$ is still an optimal basis for $\mathbf x^1$, i.e.,  the associated BFS solves $\mathbf{LP}(\mathbf x^1, \boldsymbol\beta)$.
	
	When the BFS of $\mathfrak B^0$ is the unique optimal solution for $\mathbf{LP}(\mathbf x^0, \boldsymbol\beta)$, the reduced cost of $u_j$ is strictly negative for $j\in \mathfrak N^0$. Again, because of the aforementioned independence, those reduced costs remain strictly negative  when $\mathbf x^0$ is replaced by $\mathbf x^1$.  Hence, the BFS of $\mathfrak B^0$ is also the unique optimal solution for $\mathbf{LP}(\mathbf x^1, \boldsymbol\beta)$, provided that $\mathfrak B^0$ is a feasible basis.
\end{proof}

\noindent\textbf{Lemma \ref{lem_repeated_bases_ccg}.}
If $\mathbb{B}^{t_1}=\mathbb{B}^{t_2}$ with $t_1<t_2$, Variant 2 terminates, and $\mathbf x^1$, an optimal solution to $\mathbf{MP2}$ in $t_1$-th iteration, is optimal to $\mathbf{2-Stg \ RO}$.
\begin{proof}
	We consider the case where $\mathbf{MP2}$ is with $\mathcal{OU}$ sets only as other cases can be proven similarly.  Assume that $(\mathbf x^1, \eta^1)$ is optimal to $\mathbf{MP2}^{t_1}$, and $(\mathbf x^2, \eta^2)$ is optimal to $\mathbf{MP2}^{t_2}$. If $\mathbf x^1=\mathbf x^2$, Proposition \ref{prop_cardinality_repeat_CCG} ensures that $\mathbf x^1$ is optimal to $\mathbf{2-Stg \ RO}$ and the algorithm terminates. Hence, we assume that they are different.  According to the construction of $\mathbf{MP2}$, it is clear that $\mathbf{MP2}^{t_1}$ is a relaxation to $\mathbf{MP2}^{t_2}$.  Obviously, once $(\mathbf x^1,\eta^1)$ is feasible to $\mathbf{MP2}^{t_2}$, it is also optimal in this iteration. Consequently, by Proposition \ref{prop_cardinality_repeat_CCG} again, the desired statement holds.

	Let $\hat{\mathcal{P}}_{\Pi}=\{\boldsymbol\pi^1,\dots, \boldsymbol\pi^{|t_1|}\}$ in iteration $t_1$, $\hat{\mathcal{P}}_{\Pi} =\{\boldsymbol\pi^1,\dots, \boldsymbol\pi^{|t_1|},\dots,\pi^{|t_2|}\}$ in iteration $t_2$, and
	$\mathbb{B}^{t_1}=\mathbb{B}^{t_2}=\{\mathfrak{B}^1,\dots, \mathfrak{B}^k\}$. Consider an arbitrary $\boldsymbol\pi'$ that belongs to $\hat{\mathcal{P}}_{\Pi}$ in iteration $t_2$, and assume that computing $\mathbf{MP2}^{t_2}$ yields the BFS of basis $\mathfrak{B}'$ for $\mathcal{OU}(\mathbf x^2, \boldsymbol\pi')$. Note that $\boldsymbol\pi'$ might not appear in $\hat{\mathcal{P}}_{\Pi}$ in iteration $t_1$.
	
	Because of $\mathfrak{B}'\in \mathbb{B}^{t_2}=\mathbb{B}^{t_1}$, the BFS of $\mathfrak{B}'$ is in $\mathcal{U}(\mathbf x^1)$, i.e., it is a feasible extreme point solution to $\mathbf{LP}(\mathbf x^1, \boldsymbol\pi')$. Following from Lemma \ref{lem_feasible_then_optimal}, we have that the BFS of $\mathfrak{B}'$ belongs to $\mathcal{OU}(\mathbf x^1, \boldsymbol\pi')$.  Given that $\boldsymbol\pi'$ is an arbitrary one in $\hat{\mathcal{P}}_{\Pi}$ in iteration $t_2$, every $\mathcal{OU}$ of $\mathbf{MP2}^{t_2}$ with $\mathbf x=\mathbf x^1$ has a BFS that has already appeared in the solution of $\mathbf{MP2}^{t_1}$.
	Therefore,  $(\mathbf x^1,\eta^1)$ is feasible to $\mathbf{MP2}^{t_2}$. As mentioned, the desired result readily follows.
\end{proof}

\noindent\textbf{Proposition \ref{prop_complexity_uniqueness}.}
Assume that $\mathbf{LP}(\mathbf x, \cdot)$ always has a unique optimal solution in the execution of Variant 2 (with operation $(i.a)$), which is referred to as the \textit{unique optimal solution} (or ``uniqueness'' for short) property. Then, the number of iterations before  termination is bounded by the number of bases of $\mathbf{LP}(\mathbf x, \cdot)$, which is of $O(\binom{n_u+\mu_u}{\mu_u})$.
\begin{proof}
	As noted earlier, the number of bases of $\mathbf{LP}(\mathbf x, \boldsymbol\beta)$ in \eqref{eq_LP_parametric} is independent of $\mathbf x$ and $\boldsymbol\beta$ and is completely determined by  the dimensions of $\mathbf F$. Next, we show that the number of iterations before termination is bounded by that of bases.
	
	Claim: If the basis output from Step 4 has been derived in a previous iteration, Variant~2  terminates with $LB= UB$.
	\begin{proof}[Proof of Claim:]
		Suppose that basis $\mathfrak{B}'$ has been derived in Step 4 in both iterations $t_1$  and $t_2$ with $t_1 < t_2$, and denote the involved extreme point or ray in Step 4 of those two iterations by $\boldsymbol\pi^1$ and $\boldsymbol\pi^2$. Also, $\mathbf x^1$ and $\mathbf x^2$ denote the optimal first stage solutions obtained by computing $\mathbf{MP2}^{t_1}$ and $\mathbf{MP2}^{t_2}$, respectively.

		By definition, the BFS associated with basis $\mathfrak{B}'$, denoted by $\mathbf u'$, belongs to
		$\mathcal{U} (\mathbf x^2)$  and is optimal to $\mathbf{LP}(\mathbf x^2, \boldsymbol\pi^2)$. Also, when $\boldsymbol\pi^2$ is an extreme point of $\Pi$, it follows from \eqref{eq_SP2} that
		\begin{eqnarray*}
			\eta_s(\mathbf x^2) &=&\max_{\mathbf u\in \mathcal{U}(\mathbf x^2)} \min\{\mathbf c_2\mathbf y:
			\mathbf{y}\in \mathcal{Y}(\mathbf x^2, \mathbf u)\}\\ 
			&=& \min\{\mathbf c_2\mathbf y: \mathbf{y}\in \mathcal{Y}(\mathbf x^2, \mathbf u')\}
		\end{eqnarray*}
		When $\boldsymbol\pi^2$ is an extreme ray of $\Pi$, which is derived by
		computing \eqref{eq_SP3},  we have $\eta_s(\mathbf x^2)=+\infty$. As the optimal value of $\mathbf{LP}(\mathbf x^2, \boldsymbol\pi^2)$ is larger than 0 (otherwise \eqref{eq_SP3} will not be unbounded), it follows that the recourse problem for $\mathbf u'$ is infeasible. Hence, when $\boldsymbol\pi^2$ is an extreme ray, we have
		\begin{eqnarray*}
			\eta_s(\mathbf x^2) = +\infty = \min\{\mathbf c_2\mathbf y: \mathbf{y}\in \mathcal{Y}(\mathbf x^2, \mathbf u')\}
		\end{eqnarray*}

		Since $\mathfrak{B}'$ is also derived in iteration $t_1$, the associated  BS is certainly optimal to  $\mathbf{LP}(\mathbf x^1, \boldsymbol\pi^1)$. Given that $\mathfrak{B}'$ yields  $\mathbf u'$ in
		$\mathcal{U}(\mathbf x^2)$, it follows from Lemma \ref{lem_feasible_then_optimal} that $\mathbf u'$ is also optimal to $\mathbf{LP}(\mathbf x^2, \boldsymbol\pi^1)$.  Moreover, with the unique optimal solution assumption being held,  $\mathbf u'$  remains the unique one for $\mathbf{LP}(\mathbf x^2, \boldsymbol\pi^1)$.
		Hence, we have
		$\mathcal{OU}_{\mathbf u}(\mathbf x^2, \boldsymbol\pi^1)=\{\mathbf u'\}$
		or $\mathcal{OV}_{\mathbf u}(\mathbf x^2, \boldsymbol\pi^1)=\{\mathbf u'\}$.
		Indeed, the unified cutting set \eqref{eq_CCG_optimality}, based on either $\mathcal{OU}(\mathbf x, \boldsymbol\pi^1)$ or
		$\mathcal{OV}(\mathbf x, \boldsymbol\pi^1)$, is part of
		constraints of $\mathbf{MP2}^{t_2}$.
		So, for given $\mathbf x^2$ we have the following from $\mathbf{MP2}^{t_2}$.
		\begin{eqnarray*}
			\label{eq_convergence_V2}
			\eta &\geq & \min\{\mathbf{c_2}\mathbf y^{\boldsymbol\pi^1}: \      \mathbf{B}_2\mathbf y^{\boldsymbol\pi^1}\geq \mathbf{d-\mathbf B}_1\mathbf x^2-\mathbf{Eu}', \ \mathbf y^{\boldsymbol\pi^1} \geq \mathbf 0\} = \eta_s(\mathbf x^2)
		\end{eqnarray*}
		
		For Step 5 in iteration $t_2$, we have
		$$UB\leq \mathbf{c}_1\mathbf x^2+\eta_s(\mathbf x^2)\leq LB,$$
		which  is the desired result.
	\end{proof}
	The aforementioned claim directly indicates that a repeated basis leads to the algorithm termination. Together with the fact that the number of bases is bounded by
	$\binom{n_u+\mu_u}{\mu_u}$, it can be concluded that the iteration complexity of Variant 2 is $O(\binom{n_u+\mu_u}{\mu_u})$.
\end{proof}

\noindent\textbf{Theorem \ref{thm_basis_enumeration}.}
Formulation $\mathbf{2-Stg \ RO}$ in \eqref{eq_2RO} (and its equivalences) is equivalent to a sing-level optimization program as in the following.
\begin{subequations}
	\label{A:eq_2stgRO_Bases}
	\begin{align}
		\mathbf{2\!-\!Stg \ RO(\mathbb{B})}: \   w^*=\min \quad & \mathbf{c}_1\mathbf x+ \eta   \\
		\mathrm{s.t.} \quad & \mathbf x  \in \mathcal{X}\\
		\bigg\{&\eta\geq   \mathbf{c_2}\mathbf y+ M(\mathbf1^\intercal\bar{\mathbf u}^1+\mathbf1^\intercal\bar{\mathbf u}^2+\mathbf1^\intercal\bar{\mathbf u}^3) \label{A:eq_enu3_basis}\\
		& + (\mathbf h_{\mathfrak N}+\mathbf G_{\mathfrak N}\mathbf x)^\intercal\boldsymbol\lambda_{\mathfrak N}+(\mathbf h_{\mathfrak B}+\mathbf G_{\mathfrak B}\mathbf x)^\intercal\boldsymbol\lambda_{\mathfrak B} \notag \\
		&   (\mathbf{u},\bar{\mathbf u}^1,\bar{\mathbf u}^2,\bar{\mathbf u}^3)\in \overline{\mathcal{BS}}(\mathfrak B, \mathbf x) \label{A:eq_enu3_basis2}\\
		&\hat{\mathbf F}_{\mathfrak N}(\mathbf x)^\intercal\boldsymbol\lambda_{\mathfrak N}+\hat{\mathbf F}_{\mathfrak B}(\mathbf x)^\intercal\boldsymbol\lambda_{\mathfrak B}\geq \mathbf 0, \ \boldsymbol\lambda_{\mathfrak N}\geq \mathbf 0 \\
		& \mathbf y\in \mathcal Y(\mathbf x, \mathbf u) 	\bigg\}  \ \forall \mathfrak B\in \mathbb{B} \label{A:eq_enu3_optimality_basis}
	\end{align}
\end{subequations}

\begin{proof}
	Note that it is sufficient to show that if $\mathfrak{B}$ is a feasible basis for $\mathbf x$, $\overline{\mathcal{BS}}(\mathfrak B, \mathbf x)$ reduces to $\mathcal{BS}(\mathfrak B, \mathbf x)$. Otherwise, \eqref{A:eq_enu3_basis} generates no effect on $\eta$.
	
	Consider the case where  $\mathfrak{B}$ is a feasible basis for $\mathbf x$, i.e., $\mathcal{BS}(\mathfrak B, \mathbf x)$ is not empty. Because of the role of big M in \eqref{A:eq_enu3_basis}, it can be easily argued that in an optimal solution to $\mathbf{2\!-\!Stg \ RO(\mathbb{B})}$ we have $\bar{\mathbf {u}}^j=0, j=1,2,3$ and therefore $\overline{\mathcal{BS}}(\mathfrak B, \mathbf x)=\mathcal{BS}(\mathfrak B, \mathbf x)$. Furthermore, by Lemma \ref{lem_alternative_thm}, we have $(\mathbf h_{\mathfrak N}+\mathbf G_{\mathfrak N}\mathbf x)^\intercal\boldsymbol\lambda_{\mathfrak N} +(\mathbf h_{\mathfrak B}+\mathbf G_{\mathfrak B}\mathbf x)^\intercal\boldsymbol\lambda_{\mathfrak B}=0$ in that optimal solution. So, \eqref{A:eq_enu3_basis} becomes $\eta\geq \mathbf{c}_2\mathbf y$.
	
	Consider the other case where  $\mathfrak{B}$ is an infeasible basis for $\mathbf x$. In an optimal solution to $\mathbf{2\!-\!Stg \ RO(\mathbb{B})}$, one $\mathbf u$ will be selected, along with non-zero $\bar{\mathbf {u}}^j$ for some $j$ penalized by $M$. Nevertheless, by Lemma \ref{lem_alternative_thm} again, it follows that $(\mathbf h_{\mathfrak N}+\mathbf G_{\mathfrak N}\mathbf x)^\intercal\boldsymbol\lambda_{\mathfrak N} +(\mathbf h_{\mathfrak B}+\mathbf G_{\mathfrak B}\mathbf x)^\intercal\boldsymbol\lambda_{\mathfrak B}$ can be  arbitrarily negative. Hence, the RHS of \eqref{A:eq_enu3_basis} imposes no restriction on $\eta$.
	
	With arguments of those two cases, the desired conclusion follows.
\end{proof}

\section{Modifications on Variant 2 to Attain The Uniqueness Property}
\label{Asect_modiunique}
On the one hand, according to the classical Simplex method, a maximization linear program has a unique optimal solution  if reduced costs of all non-basic variables are negative. On the other hand, as noted in Remark \ref{rem_advantageCCG},  instead of restricting $\boldsymbol\beta$ of $\mathbf{LP}(\mathbf x, \boldsymbol\beta)$ to be $\boldsymbol\pi$ or $\boldsymbol\gamma$, i.e., extreme points or rays of $\Pi$, we can creatively  define the objective function of $\mathbf{LP}(\mathbf x, \cdot)$ (consequently its optimal solution sets  $\mathcal{OU}$ and $\mathcal{OV}$). Hence, if some non-basic variables have zero reduced costs in an optimal solution to $\mathbf{LP}(\mathbf x, \cdot)$, one strategy is to reduce their coefficients in the objective function. Subsequently, for the updated objective function, that optimal solution becomes the unique one. Because our treatments on $\boldsymbol\pi^*$-related and $\boldsymbol\gamma^*$-related terms are basically same, we next simply take the first one for illustration.

Specifically, suppose that the reduced cost, denoted by $\mathbf{r}^*$, is output from  solving $\mathbf{LP}(\mathbf x^*, \boldsymbol\pi^*)$, along with the optimal basis $\mathfrak{B}^{*}$, in \textbf{Step 4 $(i.a)$}.  As shown in \eqref{eq_LP_parametric}, the objective function coefficients of $\mathbf{LP}(\mathbf x^*, \boldsymbol\pi^*)$ are $\mathfrak c^{\boldsymbol\pi^*}=-(\mathbf{E}^\intercal\boldsymbol\pi^*)^\intercal$. To address the uniqueness issue, we introduce $\hat{\mathfrak c}^{\boldsymbol\pi^*}$ defined as
\begin{eqnarray}
	\label{eq_unique_property_manu}
	\hat{\mathfrak c}^{\boldsymbol\pi^*}_j = \left\{\begin{array}{l}
		\mathfrak{c}^{\boldsymbol\pi^*}_j, \ \ \  \ \ \ \mbox{if} \ j\in \mathfrak{B}^{*}\\
		\mathfrak{c}^{\boldsymbol\pi^*}_j, \ \ \ \ \ \ \mbox{if} \ j\notin \mathfrak{B}^{*} \ \mbox{and} \ \mathbf r^*_j<0	\\	
		\mathfrak{c}^{\boldsymbol\pi^*}_j-\epsilon, \ \mbox{if} \ j\notin \mathfrak{B}^{*} \ \mbox{and} \ \mathbf r^*_j=0	
	\end{array}\right\}
\end{eqnarray}
where $\epsilon$ is a positive value. Then, the next result can be proven easily as $\mathfrak{B}^{*}$ remains an optimal basis and reduced costs of non-basic variables are all negative. 

\noindent\textbf{Lemma A.1}	\label{lem_modifiedV2_unique}
	The linear program $\max\{\hat{\mathfrak c}^{\boldsymbol\pi^*}\mathbf u: \ \mathbf u\in \mathcal{U}(\mathbf x^*)\}$ has a unique optimal solution, i.e., its optimal solution set
	\begin{eqnarray}
		\label{eq_KKT_point_unique}
		\widehat{\mathcal{OU}}(\mathbf x^*, \boldsymbol\pi^*)= \left\{ \begin{array}{l}
			\mathbf F\mathbf u\leq \mathbf h+\mathbf{Gx}^*\\
			\mathbf F^\intercal\boldsymbol\lambda\geq \hat{\mathfrak c}^{\boldsymbol\pi^*}\\
			\boldsymbol\lambda \circ (\mathbf h+\mathbf{Gx}^*-\mathbf F\mathbf u)=\mathbf{0}\\
			\mathbf u \circ (\mathbf F^\intercal\boldsymbol\lambda- \hat{\mathfrak c}^{\boldsymbol\pi^*}) = \mathbf 0\\
			\mathbf u\geq \mathbf{0}, \boldsymbol\lambda\geq \mathbf{0}.
		\end{array}\right\}
	\end{eqnarray}
	is a singleton. 	$\hfill\square$
	
By performing the same set of operations, we can obtain  counterparts for $\boldsymbol\gamma^*$, i.e.,  $\hat{\mathbf c}^{\boldsymbol\gamma^*}$ and $\widehat{\mathcal{OV}}(\mathbf x^*,\boldsymbol\gamma^*)$. The latter one is also guaranteed to be a singleton. We then describe particular modifications to the original Variant 2 in the following:
\begin{description}
	\item[(M1)] In \textbf{Step 4}, for \textbf{Case A}, modify $(i.a)$ to ``compute $\mathbf{LP}(\mathbf x^*, \boldsymbol\pi^*)$ with an optimal BFS and the associated basis $\mathfrak{B}^{*}$, and the reduced cost $\mathbf r^*$.'' Similarly, for \textbf{Case B}, modify $(i.a)$ to ``compute $\mathbf{LP}(\mathbf x^*, \boldsymbol\gamma^*)$ with ...''
	\item[(M2)] In \textbf{Step 4}, for \textbf{Case A} after $(i.a)$, perform ``$(i.b)$ compute $\hat{\mathfrak c}^{\boldsymbol\pi^*}$''. Similarly, for \textbf{Case B}, perform ``$(i.b)$ compute $\hat{\mathfrak c}^{\boldsymbol\gamma^*}$''.
	\item [(M3)] In cutting set generation, replace $\mathcal{OU}(\mathbf x, \boldsymbol\pi^*)$ in  \eqref{eq_CCG_optimality}  and $\mathcal{OV}(\mathbf x, \boldsymbol\gamma^*)$ in \eqref{eq_CCG_feasibility} by $\widehat{\mathcal{OU}}(\mathbf x, \boldsymbol\pi^*)$ and $\widehat{\mathcal{OV}}(\mathbf x, \boldsymbol\gamma^*)$, respectively. Consequently, sets $\mathcal{OU}$ and $\mathcal{OV}$ in $\mathbf{MP2}$ are replaced by $\widehat{\mathcal{OU}}$ and $\widehat{\mathcal{OV}}$.
\end{description}

\section{Algorithm Details for Two-Stage RO with SOC Structures}
\label{Asect:V2SOC}

We consider first two-stage RO with an SOC recourse problem, and then a DDU set represented as an SOC set, with the overall formulation and other sets defined in  \eqref{eq_2RO}-\eqref{eq_recourse_set} when applicable. We also assume that both SOC structures satisfy Slater's condition.

\subsection{Two-Stage RO with SOC Recourse Problem}
Different from $\mathbf{2-Stg \ RO}$ in \eqref{eq_2RO}, the recourse problem, $\min\{\mathbf{c}_2\mathbf y: \mathbf y\in \mathcal{Y}(\mathbf x,\mathbf u)\}$, is defined with respect to the following SOC set 
\begin{align}
\hspace{-5pt}	\mathcal Y(\mathbf x,\mathbf u)=\left\{\mathbf y\in\mathbb{R}^{n_y}:||
	\mathbf Q_i\mathbf y+\mathbf f_i||_2\leq  \mathbf B_{2,\bar i}\mathbf y-d_i+\mathbf{B}_{1,\bar i}\mathbf x+\mathbf E_{\bar i}\mathbf u, \quad i=1,2,\dots,\mu_y\right\},
\end{align}
where a matrix  with subscript $\bar i$, e.g., $\mathbf{B}_{2,\bar i}$, indicates its $i$-th row. Clearly, if $\mathbf Q_i$ and $\mathbf f_i$ are zero matrix and vector respectively  for all $i$, $\mathcal Y(\mathbf x,\mathbf u)$ reduces to the original polyhedral set in \eqref{eq_recourse_set}.   
Note that vector $\mathbf f_i$ collectively represents components that are given as parameters, which might be  functions of $\mathbf x$ and/or $\mathbf u$.  

In the following, we described detailed subproblems and the master problem using parametric C\&CG scheme. To minimize repetition, we do not explicitly describe cutting sets, noting that they are (partially) enumerated components in the master problem. 
\begin{align*}
\hspace{-5pt}	\mathbf{SP1}: \ \eta_f (\mathbf x^*)  = \max_{\mathbf u\in \mathcal U(\mathbf x^*)} & \min
	\Big\{\mathbf{1}^\intercal\tilde{\mathbf y}: || \mathbf Q_i\mathbf y+\mathbf f_i||_2\leq  \mathbf B_{2,\bar i}\mathbf y-d_i+\mathbf{B}_{1,\bar i}\mathbf x +\mathbf E_{\bar i}\mathbf u+\tilde y_i \ \ \forall i, \ \tilde{\mathbf y}\geq
	\mathbf 0 \Big\}
\end{align*}
By taking the dual of the lower-level problem, $\mathbf{SP1}$ can be reformulated into the following one that can be directly solved by some professional solver. 
\begin{align*}
	\mathbf{SP1}:\quad\max\quad&\sum_i \mathbf f_i^\intercal\boldsymbol\pi_{i}+\mathbf1^\intercal(\mathbf d- \mathbf B_1 \mathbf x^*- \mathbf {Eu})\\
	\mathrm{s.t.}\quad&\sum_i \mathbf Q_i^\intercal\boldsymbol\pi_{i}=\mathbf B_2^\intercal\\
	&||\boldsymbol\pi_{i}||_2\leq 1 \ \ \forall i\\
	&\mathbf F(\mathbf x^*)\mathbf u \leq\mathbf h+\mathbf{Gx^*}\\
	&\mathbf u\in\mathbb Z^{m_u}_+\times \mathbb{R}^{n_u}_+
\end{align*}

Similarly, we can deal with $\mathbf{SP2}$ defined as the following  
\begin{eqnarray*}
	\mathbf{SP2}: \ \eta_s(\mathbf x^*)=\max_{\mathbf u\in \mathcal{U}(\mathbf x^*)} \min\{
	\mathbf c_2\mathbf y: \mathbf y\in \mathcal{Y}(\mathbf x^*, \mathbf u)\},
\end{eqnarray*}
and solve its single-level reformulation listed below. 
\begin{align*}
	\mathbf{SP2}:\quad\max\quad&\sum_{i}\mathbf f_i^\intercal\boldsymbol\pi_{1,i}+(\mathbf d- \mathbf B_1 \mathbf x^*- \mathbf {Eu})^\intercal\boldsymbol\pi_{2}\\
	\mathrm{s.t.}\quad&\sum_i \mathbf Q_i^\intercal\boldsymbol\pi_{1,i}-\mathbf B_2^\intercal\boldsymbol\pi_2+\mathbf c_2^\intercal=0\\
	&||\boldsymbol\pi_{1,i}||_2\leq\pi_{2,i} \ \  \forall i\\
	&\mathbf F(\mathbf x^*)\mathbf u \leq\mathbf h+\mathbf{Gx^*}\\
	&\mathbf u\in\mathbb Z^{m_u}_+\times \mathbb{R}^{n_u}_+
\end{align*}
Regarding $\mathbf{SP3}$, it is defined in the following
\begin{eqnarray*}
	\mathbf{SP3}: \ \max\left\{\sum_i \mathbf f_i^\intercal\boldsymbol\pi_{1,i}+(\mathbf d- \mathbf B_1 \mathbf x^*- \mathbf {Eu^*_f})^\intercal\boldsymbol\pi_{2}: (\boldsymbol\pi_{1,1},\boldsymbol\pi_{1,2},\cdots,\boldsymbol\pi_{1,\mu_y},\boldsymbol\pi_2)\in \Pi\right\}
\end{eqnarray*}
with $
	\Pi=\left\{\boldsymbol\pi:\sum_i \mathbf Q_i^\intercal\boldsymbol\pi_{1,i}-\mathbf B_2^\intercal\boldsymbol\pi_2+\mathbf c_2^\intercal=0, \ \ ||\boldsymbol\pi_{1,i}||_2\leq\pi_{2,i}\ \ \forall i \right\}$.
	
We note that professional solvers might not generate an extreme ray of $\Pi$. If this is the case, it can be augmented with $\pi_{2,i}\leq M$ for all $i$ to achieve the computational feasibility. Also, if the recourse problem is actually a convex quadratic program (with linear constraints), e.g., the recourse problem in $\mathbf{Robust \ FL-RC}$, $\mathbf{SP3}$ can be built with respect to the linear constraints so that an extreme ray will be provided by an LP solver. 

Next, we provide the general form of $\mathbf{MP2}$, along with set $\mathcal{OU}$. Note that $\mathcal{OV}$ is in the same form as $\mathcal{OU}$. We mention that if $\mathbf f_i$ include components that depend on $\mathbf u$, additional operations are needed to include $\boldsymbol{\pi}_{1,i}$ for some $i$ into the definition of set $\mathcal{OU}$ or $\mathcal{OV}$. 
\begin{eqnarray*}
	\mathbf{MP2}: \ \underline w=\min && \mathbf{c}_1\mathbf{x}+ \eta \notag \\
	\mbox{s.t.} && \mathbf{x} \in \mathcal{X} \notag\\
	&& \eta\geq \mathbf{c_2}\mathbf y^{\boldsymbol\pi} \ \ \forall \boldsymbol\pi\in \mathcal{\hat P}_{\Pi}\\
	&&||\mathbf Q_i\mathbf y^{\boldsymbol\pi}+\mathbf f_i||_2\leq  \mathbf B_{2,\bar i}\mathbf y^{\boldsymbol\pi}-d_i+\mathbf{B}_{1,\bar i}\mathbf x+\mathbf E_{\bar i}\mathbf u^{\boldsymbol\pi} \quad   \forall i \ \forall \boldsymbol\pi\in \mathcal{\hat P}_{\Pi}\\
	&& (\mathbf{u}^{\boldsymbol\pi},\boldsymbol\lambda^{\boldsymbol\pi})\in\mathcal{OU}(\mathbf x,\boldsymbol\pi) \  \  \forall \boldsymbol\pi\in \mathcal{\hat P}_{\Pi}\\
	&& \eta\geq \mathbf{c_2}\mathbf y^{\boldsymbol\gamma} \ \ \forall \boldsymbol\gamma\in \mathcal{\hat R}_{\Pi}\\
	&&||\mathbf Q_i\mathbf y^{\boldsymbol\gamma}+\mathbf f_i||_2\leq  \mathbf B_{2,\bar i}\mathbf y^{\boldsymbol\gamma}-d_i+\mathbf{B}_{1,\bar i}\mathbf x+\mathbf E_{\bar i}\mathbf u^{\boldsymbol\gamma} \quad  \forall i \ \forall \boldsymbol\gamma\in \mathcal{\hat R}_{\Pi}\\
	&&  (\mathbf{v}^{\boldsymbol\gamma},\boldsymbol\zeta^{\boldsymbol\gamma})\in \mathcal{OV}(\mathbf x, \boldsymbol\gamma) \ \ \forall \boldsymbol\gamma\in \mathcal{\hat R}_{\Pi}
\end{eqnarray*}
\begin{eqnarray*}
	\label{A:eq_OU_RC}
	\mathcal{OU}(\mathbf x, \boldsymbol\pi^k_2) =  \left\{\begin{array}{l}
		\mathbf F(\mathbf x)\mathbf u^k\leq \mathbf h+\mathbf{Gx}\\
		\mathbf F(\mathbf x)^\intercal\boldsymbol\lambda^k\geq - \mathbf E^\intercal\boldsymbol\pi^k_2 \\
		\boldsymbol\lambda^k \circ (\mathbf h+\mathbf{Gx}-\mathbf F(\mathbf x)\mathbf u^k)=\mathbf{0}\\
		\mathbf u^k \circ (\mathbf F(\mathbf x)^\intercal\boldsymbol\lambda^k+ \mathbf E^\intercal\boldsymbol\pi^k_2) = \mathbf 0\\
		\mathbf u^k\geq \mathbf{0}, \boldsymbol\lambda^k\geq \mathbf{0},
	\end{array}\right\}
\end{eqnarray*}

\subsection{Two-Stage RO with SOC DDU}
Different from $\mathbf{2-Stg \ RO}$ in \eqref{eq_2RO}, the DDU set is the following SOC set
\begin{align}
	\label{A:eq_uncer_set}
	\mathcal U(\mathbf x)=\left\{\mathbf u\in\mathbb{R}^{n_u}:||\mathbf W_i\mathbf u+ \mathbf g_i||_2 \leq h_i+\mathbf G_{\bar i}\mathbf x-\mathbf F_{\bar i}(\mathbf x)\mathbf u ,\quad i=1,2,\dots,\mu_u\right\}.
\end{align}
Again, if $\mathbf W_i$ and $\mathbf g_i$ are zero matrix and vector respectively  for all $i$, $\mathcal U(\mathbf x)$ reduces to the  polyhedral DDU set defined in \eqref{eq_uncer_set}. Also,  
 vector $\mathbf g_i$ collectively represents components that are given as parameters, which might be functions of $\mathbf x$.  

In the following, we introduce three subproblems. Since the recourse problem is an LP, they can be handled by the same approach used for original subproblems presented in Section \ref{subsubsec_Variant1}, except for a more sophisticated solver for the SOC structure.  
\begin{eqnarray*}
	\mathbf{SP1}: \ \eta_f (\mathbf x^*) = \max_{\mathbf u\in \mathcal U(\mathbf x^*)}  & \min
	\{\mathbf{1}^\intercal\tilde{\mathbf y}: \mathbf{B}_2\mathbf y+\tilde{\mathbf y}\geq \mathbf d-\mathbf B_1\mathbf x^*-\mathbf E\mathbf u, \mathbf y\geq \mathbf 0, \tilde{\mathbf y}\geq
	\mathbf 0 \}.
\end{eqnarray*}
\begin{eqnarray*}
	\mathbf{SP2}: \ \eta_s(\mathbf x^*)=\max_{\mathbf u\in \mathcal{U}(\mathbf x^*)} \min\{
	\mathbf c_2\mathbf y: \mathbf y\in \mathcal{Y}(\mathbf x^*, \mathbf u)\}
\end{eqnarray*}
\begin{eqnarray*}
	\mathbf{SP3}: \ \max\{(\mathbf d- \mathbf B_1 \mathbf x^*- \mathbf {Eu^*_f})^\intercal\boldsymbol\pi: \boldsymbol\pi\in \Pi\}
\end{eqnarray*}
with $\Pi=\{\mathbf B_2^\intercal\boldsymbol\pi \leq \mathbf c_2, \  \boldsymbol\pi\geq \mathbf{0}\}$.

Next, we provide the general form of $\mathbf{MP2}$, along with set $\mathcal{OU}$ (noting that $\mathcal{OV}$ can be defined similarly). We mention that $\mathcal{OU}$ (and $\mathcal{OV}$, respectively) is in the form of the primal-dual based representation, which allows us to take advantage of the SOC structure of the dual problem.  
\begin{eqnarray*}
	\mathbf{MP2}: \ \underline w=\min && \mathbf{c}_1\mathbf{x}+ \eta \notag \\
	\mbox{s.t.} && \mathbf{x} \in \mathcal{X} \notag\\
	&& \eta\geq \mathbf{c_2}\mathbf y^{\boldsymbol\pi} \ \ \forall \boldsymbol\pi\in \mathcal{\hat P}_{\Pi}\\
	&& \mathbf{B}_2\mathbf y^{\boldsymbol\pi}\geq \mathbf{d-\mathbf B}_1\mathbf x-\mathbf{Eu}^{\boldsymbol\pi} \ \  \forall \boldsymbol\pi\in \mathcal{\hat P}_{\Pi}\\
	&& (\mathbf{u}^{\boldsymbol\pi},\boldsymbol\lambda^{\boldsymbol\pi})\in
	\mathcal{OU}(\mathbf x,\boldsymbol\pi), \  \mathbf y^{\boldsymbol\pi} \geq \mathbf 0 \ \forall \boldsymbol\pi\in \mathcal{\hat P}_{\Pi}\\
	&& \eta\geq \mathbf{c_2}\mathbf y^{\boldsymbol\gamma} \ \ \forall \boldsymbol\gamma\in \mathcal{\hat R}_{\Pi}\\
	&&\mathbf{B}_2\mathbf y^{\boldsymbol\gamma}\geq \mathbf{d-\mathbf B}_1\mathbf x-\mathbf{Ev}^{\boldsymbol\gamma} \ \ \forall \boldsymbol\gamma\in \mathcal{\hat R}_{\Pi}\\
	&&  (\mathbf{v}^{\boldsymbol\gamma},\boldsymbol\zeta^{\boldsymbol\gamma})\in \mathcal{OV}(\mathbf x, \boldsymbol\gamma), \ \mathbf y^{\boldsymbol\gamma}\geq \mathbf 0 \ \ \forall \boldsymbol\gamma\in \mathcal{\hat R}_{\Pi}
\end{eqnarray*}
\begin{eqnarray*}
	\mathcal{OU}(\mathbf x, \boldsymbol\pi^k) =  \left\{\begin{array}{l}
		\boldsymbol\pi^{k\intercal}\mathbf E\mathbf u^k \leq \sum_{i=1}^{n_u}\mathbf g_i^\intercal\boldsymbol\lambda_{1,i}-(\mathbf h+\mathbf G\mathbf x)^\intercal\boldsymbol\lambda_2\\
		||\mathbf W_i\mathbf u^k+ \mathbf g_i||_2 \leq h_i+\mathbf G_{\bar i}\mathbf x-\mathbf F_{\bar i}(\mathbf x)\mathbf u^k \quad \forall i\\
		\sum_{i=1}^{n_u}W_i^\intercal\boldsymbol\lambda_{1,i}+\mathbf F^\intercal(\mathbf x)\boldsymbol\lambda_2+\mathbf E^\intercal\boldsymbol\pi=0\\
		||\boldsymbol\lambda_{1,i}||_2\leq\lambda_{2,i}\quad \forall i
	\end{array}\right\}
\end{eqnarray*}

\section{Proofs of Section \ref{sect_more}}
\label{apd:proofS4}
\noindent\textbf{Corollary \ref{cor_reliable_pmedian}.}
	Let a DDU set be
	\begin{eqnarray}
		\label{Aeq_DDU_reliable_pmedian}
		\mathcal{U}^k(\mathbf x)=\{\mathbf {u}\in \mathbb{R}^{|I|}_+: \sum_{j\in J}u_j\leq k, u_j\leq x_{d,j} \forall j, \ u_i=0 \ \forall i\notin J\}.
	\end{eqnarray}
	If $\textsl{C}\geq \max_{ij}\{c_{ij}\}$ and $\theta_i\leq 0$ for all $i$,
	the two-stage RO in \eqref{eq_reliable_pmedian} is equivalent to
	\begin{equation*}
		\begin{split}
			\mathit{w}_{R}(\mathcal{X}, \mathcal{U}^k(\mathbf x), \mathcal{Y}(\mathbf x, \mathbf u)) = \min_{(\mathbf x_c,\mathbf x_d)\in\mathcal X} & (1-\rho)\sum_{i\in I}\sum_{j\in J}c_{ij}x_{c,ij}+\\
			& \rho\max_{\mathbf u\in\mathcal U^k(\mathbf x)}\min_{(\mathbf y_1,\mathbf y_2)\in\mathcal{Y}(\mathbf x,\mathbf u)} \ \sum_{i\in I}\sum_{j\in J}c_{ij}y_{1,ij}+\sum_{i\in I}\textsl{C}y_{2,i}.
		\end{split}
	\end{equation*}
\begin{proof}
    It is sufficient to consider the non-trivial case where $p\geq k$.  If the listed conditions hold,  it is proven in \citet{an2014reliable} that for any fixed $\mathbf x$ we have
	$$\max_{\mathbf u\in\mathcal U^0}\min_{(\mathbf y_1,\mathbf y_2)\in\mathcal{Y}(\mathbf x,\mathbf u)} \ \sum_{i\in I}\sum_{j\in J}c_{ij}y_{1,ij}+\sum_{i\in I}\textsl{C}y_{2,i} = \max_{\mathbf u\in\mathcal {U'}(\mathbf x)}\min_{(\mathbf y_1,\mathbf y_2)\in\mathcal{Y}(\mathbf x,\mathbf u)} \ \sum_{i\in I}\sum_{j\in J}c_{ij}y_{1,ij}+\sum_{i\in I}\textsl{C}y_{2,i}$$
	with $\mathcal{U'}(\mathbf x)=\{\mathbf u\in\{0,1\}^{|I|}: \sum_{j}u_j\leq k, u_j\leq x_{d,j} \ \forall j, \ u_i=0 \ \forall i\notin J\}$.
	As $u_i$ is fixed to 0 for all $i\notin J$, it is without loss of generality to assume that $I=J$ in the following argument.
	
	Consider the continuous relaxation of $\mathcal{U'}(\mathbf x)$, which is
	$$\mathcal{U}'_r(\mathbf x)=\{\mathbf u \in \mathbb{R}^{|I|}_+: \sum_{i\in I}u_i\leq k, u_i\leq x_{d,i} \ \forall i\},$$
	parameterized by binary variables $\mathbf x_d$. \\
	Claim: $\mathcal{U}'_r(\mathbf x)$'s  constraint matrix is a totally unimodular matrix.
	\begin{proof}[Proof of Claim:]
		For this constraint matrix, we have \\
		$(i)$ every entry is either 0 or 1; \\
		$(ii)$ every column contains only two ``1'' entries;\\
		$(iii)$ there is a partition such that the first constraint is in one subset and the remaining constrains are in the other subset, which ensures that each column has exactly a ``1'' entry in one of those two subsets.\\
		Hence, according to Proposition 3.2 of \citet{wolsey2020integer}, this constraint matrix is totally unimodular.
	\end{proof}
With this claim proved, it follows that
$$co(\mathcal{U'}(\mathbf x))=\mathcal{U}'_r(\mathbf x)=\{\mathbf u \in \mathbb{R}^{|I|}_+: \sum_{i\in I}u_i\leq k, u_i\leq x_{d,i}, \ \forall i\}.$$ Then, according to Propositions \ref{prop_extremepoint},  we have
$$\max_{\mathbf u\in\mathcal U^0}\min_{(\mathbf y_1,\mathbf y_2)\in\mathcal{Y}(\mathbf x,\mathbf u)} \ \sum_{i\in I}\sum_{j\in J}c_{ij}y_{1,ij}+\sum_{i\in I}\textsl{C}y_{2,i} = \max_{\mathbf u\in\mathcal {U}'_r(\mathbf x)}\min_{(\mathbf y_1,\mathbf y_2)\in\mathcal{Y}(\mathbf x,\mathbf u)} \ \sum_{i\in I}\sum_{j\in J}c_{ij}y_{1,ij}+\sum_{i\in I}\textsl{C}y_{2,i}.$$
Finally, because of the arbitrarity of $\mathbf x$, $\mathcal{U}^k(\mathbf x)= \mathcal {U}'_r(\mathbf x)$, and Proposition \ref{prop_DIU_to_DDU}, the desired result follows.
\end{proof}	

\noindent\textbf{Proposition \ref{prop_CCG_MIP_approximation}.}
	For Variant 2 with the  aforementioned modifications, the values of $LB$ and $UB$, i.e., the optimal value of $\mathbf{MP2}$ and $\mathbf c_1\mathbf x^*+ \tilde \eta_s(\mathbf x^*, \mathbf y^*_d)$ respectively, are valid lower and upper bounds to the optimal value of $\mathbf{2-Stg \ RO}$ in \eqref{eq_2RO}.
\begin{proof}
	Note that we have $\mathcal{OU}(\mathbf x, \boldsymbol\beta)\in \mathcal{U}(\mathbf x)$ for an arbitrary $\boldsymbol\beta$. According to $\mathbf{2-Stg \ RO}$'s epigraph reformulation, it is straightforward to conclude that
	$\mathbf{MP2}$ is a relaxation to $\mathbf{2-Stg \ RO}$, and its optimal value is a lower bound to that of $\mathbf{2-Stg \ RO}$.
	
	For given $\mathbf x^*$, $\mathbf y^*_d$ and an arbitrary $\mathbf u\in \mathcal{U}(\mathbf x^*)$, we have
	$$\min\{\mathbf c_2\mathbf y: \mathbf y\in \mathcal{Y}(\mathbf x^*, \mathbf u)\}\leq \min\{
	\mathbf c_2\mathbf y: \mathbf y\in \mathcal{Y}(\mathbf x^*, \mathbf u), \mathbf y^d=\mathbf y^*_d\}.$$
	Because of $\mathbf u$'s arbitrarity, it follows that
	$$ \max_{\mathbf u\in \mathcal{U}(\mathbf x^*)} \min\{
	\mathbf c_2\mathbf y: \mathbf y\in \mathcal{Y}(\mathbf x^*, \mathbf u)\}\leq \max_{\mathbf u\in \mathcal{U}(\mathbf x^*)} \min\{
	\mathbf c_2\mathbf y: \mathbf y\in \mathcal{Y}(\mathbf x^*, \mathbf u), \mathbf y_d=\mathbf y^*_d\} =\tilde{\eta}_s(\mathbf x^*, \mathbf y^*_d).$$
	Note that the sum of $\mathbf c_1\mathbf x^*$ and the LHS of the above inequality provides an upper bound to the optimal value of $\mathbf{2-Stg \ RO}$. So does $\mathbf c_1\mathbf x^*+ \tilde{\eta}_s(\mathbf x^*, \mathbf y^*_d)$.
\end{proof}

\section{Detailed Computational Results}
\label{Asect_tables}

In all presented tables, columns ``LB'' and ``UB'' give the lower and upper bounds obtained when algorithms terminate, respectively. Note that if no (first stage solution) feasible solution is derived for some instance, its upper bound is marked by ``NA''. Columns ``Gap'', ``Iter'', and ``Time(s)'' display the relative gap between bounds, the number of iterations performed before termination, and the running time in seconds, respectively. For most tables, we show in ``Average'' row the averages of ``Iter'' and ``Time(s)'' across instances solved to optimality, and the average of ``gap'' across all instances. For tables reporting approximation results, including Tables \ref{tbl:CapPMDD}-\ref{tbl:MIPRC}, that ``Average'' row simply reports the average performances across all instances.  

\newpage
\begin{landscape}
	\begin{table}[htbp]
		\centering
		\caption{Computational Results of $\mathbf{Robust \ FL}$ with $\mathcal U^{R}(\mathbf x)$}
		\resizebox{!}{0.66\height}{
				\begin{tabular}{|c|cccc|ccccc|ccccc|ccccc|}
					\hline
					\multirow{2}[4]{*}{\# of Sites} & \multirow{2}[4]{*}{Fixed Cost} & \multirow{2}[4]{*}{$\underline{\boldsymbol\xi}$} & \multirow{2}[4]{*}{$\overline{\boldsymbol\xi}$} & \multirow{2}[4]{*}{$\alpha$} & \multicolumn{5}{c|}{Benders C\&CG}      & \multicolumn{5}{c|}{Parametric C\&CG}   & \multicolumn{5}{c|}{Basis Based C\&CG} \\
					\cline{6-20}          &       &       &       &       & LB    & UB    & Gap   & Iter.  & Time(s) & LB    & UB    & Gap   & Iter.  & Time(s) & LB    & UB    & Gap   & Iter.  & Time(s) \\
					\hline\hline
					\multirow{11}[6]{*}{25} & \multirow{5}[2]{*}{High} & 0.01  & 0.02  & 0.01  & 3887.87 & 3887.87 & 0.00\% & 8     & 32.09 & 3887.00 & 3887.87 & 0.02\% & 3     & 14.04 & 3887.87 & 3887.87 & 0.00\% & 3     & 9.65 \\
					&       & 0.02  & 0.05  & 0.02  & 3915.67 & 3915.67 & 0.00\% & 8     & 12.33 & 3913.83 & 3915.67 & 0.05\% & 3     & 2.07  & 3915.67 & 3915.67 & 0.00\% & 5     & 16.45 \\
					&       & 0.05  & 0.08  & 0.05  & 3955.06 & 3955.06 & 0.00\% & 8     & 22.22 & 3952.21 & 3955.06 & 0.07\% & 3     & 2.94  & 3955.06 & 3955.06 & 0.00\% & 4     & 13.75 \\
					&       & 0.08  & 0.12  & 0.08  & 4001.03 & 4001.03 & 0.00\% & 8     & 24.97 & 4001.03 & 4001.03 & 0.00\% & 4     & 5.34  & 4000.74 & 4001.03 & 0.01\% & 3     & 5.22 \\
					&       & 0.1   & 0.15  & 0.1   & 4034.51 & 4034.51 & 0.00\% & 8     & 33.41 & 4034.51 & 4034.51 & 0.00\% & 4     & 2.73  & 4034.14 & 4034.51 & 0.01\% & 3     & 6.13 \\
					\cline{2-20}          & \multirow{5}[2]{*}{Low} & 0.01  & 0.02  & 0.01  & 512.70 & 512.70 & 0.00\% & 9     & 54.39 & 512.33 & 512.70 & 0.07\% & 3     & 13.49 & 512.70 & 512.70 & 0.00\% & 3     & 10.65 \\
					&       & 0.02  & 0.05  & 0.02  & 543.32 & 543.32 & 0.00\% & 9     & 40.78 & 543.06 & 543.32 & 0.05\% & 3     & 6.94  & 543.32 & 543.32 & 0.00\% & 4     & 6.87 \\
					&       & 0.05  & 0.08  & 0.05  & 585.81 & 585.81 & 0.00\% & 9     & 32.33 & 585.51 & 585.81 & 0.05\% & 3     & 5.28  & 585.81 & 585.81 & 0.00\% & 4     & 8.25 \\
					&       & 0.08  & 0.12  & 0.08  & 635.16 & 635.16 & 0.00\% & 8     & 37.35 & 635.16 & 635.16 & 0.00\% & 4     & 4.39  & 634.57 & 635.16 & 0.09\% & 3     & 6.36 \\
					&       & 0.1   & 0.15  & 0.1   & 671.43 & 671.43 & 0.00\% & 8     & 41.69 & 671.43 & 671.43 & 0.00\% & 4     & 6.53  & 671.43 & 671.43 & 0.00\% & 4     & 17.48 \\
					\cline{2-20}          & \multicolumn{4}{c|}{\textbf{Average}}  &       &       & \textbf{0.00\%} & \textbf{8.3} & \textbf{33.16} &       &       & \textbf{0.03\%} & \textbf{3.4} & \textbf{6.37} &       &       & \textbf{0.01\%} & \textbf{3.6} & \textbf{10.08} \\
					\hline
					\hline
					\multirow{11}[6]{*}{40} & \multirow{5}[2]{*}{High} & 0.01  & 0.02  & 0.01  & 4086.89 & 4086.89 & 0.00\% & 14    & 380.10 & 4086.89 & 4086.89 & 0.00\% & 3     & 14.23 & 4086.87 & 4086.89 & 0.00\% & 3     & 15.82 \\
					&       & 0.02  & 0.05  & 0.02  & 4116.12 & 4116.12 & 0.00\% & 14    & 312.31 & 4116.12 & 4116.12 & 0.00\% & 3     & 3.76  & 4115.35 & 4116.12 & 0.02\% & 3     & 17.41 \\
					&       & 0.05  & 0.08  & 0.05  & 3962.99 & 4158.20 & 4.69\% & 13    & T     & 4158.20 & 4158.20 & 0.00\% & 3     & 5.05  & 4158.20 & 4158.20 & 0.00\% & 4     & 122.50 \\
					&       & 0.08  & 0.12  & 0.08  & 4207.12 & 4207.12 & 0.00\% & 15    & 1304.05 & 4207.12 & 4207.12 & 0.00\% & 3     & 16.97 & 4207.12 & 4207.12 & 0.00\% & 5     & 582.24 \\
					&       & 0.1   & 0.15  & 0.1   & 4020.62 & 4583.18 & 12.27\% & 12    & T     & 4242.67 & 4242.67 & 0.00\% & 3     & 12.01 & 4238.92 & 4242.67 & 0.09\% & 5     & 1190.94 \\
					\cline{2-20}          & \multirow{5}[2]{*}{Low} & 0.01  & 0.02  & 0.01  & -923.69 & -923.69 & 0.00\% & 18    & 1533.37 & -923.69 & -923.69 & 0.00\% & 3     & 20.50 & -924.09 & -923.69 & 0.04\% & 3     & 27.94 \\
					&       & 0.02  & 0.05  & 0.02  & -873.71 & -873.71 & 0.00\% & 18    & 767.24 & -873.71 & -873.71 & 0.00\% & 3     & 19.55 & -873.71 & -873.71 & 0.00\% & 5     & 131.18 \\
					&       & 0.05  & 0.08  & 0.05  & -802.62 & -802.62 & 0.00\% & 18    & 1597.44 & -802.74 & -802.62 & 0.02\% & 4     & 22.65 & -802.62 & -802.62 & 0.00\% & 5     & 261.41 \\
					&       & 0.08  & 0.12  & 0.08  & -1246.14 & -724.21 & 72.07\% & 10    & T     & -724.21 & -724.21 & 0.00\% & 4     & 34.35 & -724.21 & -724.21 & 0.00\% & 4     & 226.10 \\
					&       & 0.1   & 0.15  & 0.1   & -686.71 & -666.94 & 2.96\% & 17    & T     & -666.94 & -666.94 & 0.00\% & 4     & 30.10 & -666.94 & -666.94 & 0.00\% & 4     & 346.58 \\
					\cline{2-20}          & \multicolumn{4}{c|}{\textbf{Average}}  &       &       & \textbf{9.20\%} & \textbf{16.17} & \textbf{982.42} &       &       & \textbf{0.00\%} & \textbf{3.3} & \textbf{17.92} &       &       & \textbf{0.02\%} & \textbf{4.1} & \textbf{292.21} \\
					\hline
				\end{tabular}%
				\label{tbl_FLR25-40-RHS}%
			}
			\vspace{5pt}
			\centering
			\caption{Computational Results of $\mathbf{Robust \ FL}$ with $\mathcal U^{L\!R}(\mathbf x)$}
			\resizebox{!}{0.66\height}{
					\begin{tabular}{|c|c|ccccc|ccccc|ccccc|}
						\hline
						\multirow{2}[4]{*}{\# of Sites} & \multirow{2}[4]{*}{Fixed Cost} & \multicolumn{5}{c|}{Benders C\&CG}      & \multicolumn{5}{c|}{Parametric C\&CG}   & \multicolumn{5}{c|}{Basis Based C\&CG} \\
						\cline{3-17}          &       & LB    & UB    & Gap   & Iter.  & Time(s) & LB    & UB    & Gap   & Iter.  & Time(s) & LB    & UB    & Gap   & Iter.  & Time(s) \\
						\hline
						25    & High  & 2127.93 & 3367.59 & 36.81\% & 456   & T & 3142.91 & 3145.18 & 0.07\% & 5     & 15.06 & 1271.68 & 3596.40 & 64.64\% & 9     & T \\
						& Low   & 1904.66 & 2839.00 & 32.91\% & 498   & T & 2699.72 & 2701.94 & 0.08\% & 5     & 14.59 & 1046.68 & 3411.67 & 69.32\% & 8     & T \\
						30    & High  & 2036.12 & 3765.39 & 45.93\% & 406   & T & 3336.37 & 3336.37 & 0.00\% & 5     & 21.32 & 1257.66 & 4399.70 & 71.41\% & 8     & T \\
						& Low   & 1798.91 & 3154.77 & 42.98\% & 439   & T & 2959.10 & 2959.59 & 0.02\% & 5     & 114.64 & 1114.51 & 3478.55 & 67.96\% & 8     & T \\
						35    & High  & 2144.29 & 3678.60 & 41.71\% & 370   & T & 3464.75 & 3468.13 & 0.10\% & 5     & 44.13 & 1294.57 & 4326.70 & 70.08\% & 9     & T \\
						& Low   & 1789.55 & 3417.84 & 47.64\% & 384   & T & 3068.14 & 3069.28 & 0.04\% & 5     & 75.91 & 2977.67 & 3269.51 & 8.93\% & 10    & T \\
						40    & High  & 2206.46 & 3593.39 & 38.60\% & 339   & T & 3334.31 & 3335.91 & 0.05\% & 8     & 65.95 & 1285.74 & 4116.04 & 68.76\% & 8     & T \\
						& Low   & 1869.44 & 3189.43 & 41.39\% & 336   & T & 2973.24 & 2975.57 & 0.08\% & 6     & 49.35 & 1116.47 & 4623.80 & 75.85\% & 8     & T \\
						\hline
					\end{tabular}%
					\label{tbl_FLLR}%
				}
			\end{table}%
		\end{landscape}

\newpage
\begin{landscape}
	\begin{table}[htbp]
		\centering
		\caption{Computational Results of the Reliable P-Median Problem}
		\resizebox{!}{0.65\height}{
			\begin{tabular}{|ccc|ccccc|ccccc|ccccc|ccccc|}
				\hline
				\multirow{3}[6]{*}{$\rho$} & \multirow{3}[6]{*}{p} & \multirow{3}[6]{*}{k} & \multicolumn{10}{c|}{Uncapacitated}                                           & \multicolumn{10}{c|}{Capacitated} \\
				\cline{4-23}          &       &       & \multicolumn{5}{c|}{DIU by Basic C\&CG} & \multicolumn{5}{c|}{DDU by Parametric C\&CG} & \multicolumn{5}{c|}{DIU by Basic C\&CG} & \multicolumn{5}{c|}{DDU by Parametric C\&CG} \\
				\cline{4-23}          &       &       & LB    & UB    & Gap   & Iter.  & Time(s) & LB    & UB    & Gap   & Iter.  & Time(s) & LB    & UB    & Gap   & Iter.  & Time(s) & LB    & UB    & Gap   & Iter.  & Time(s) \\
				\hline
				\multirow{9}[2]{*}{0.2} & 6     & 1     & 8378.25 & 8378.25 & 0.00\% & 8     & 98.15 & 8378.25 & 8378.25 & 0.00\% & 8     & 80.90 & 8979.84 & 8979.84 & 0.00\% & 11    & 325.76 & 8979.84 & 8979.84 & 0.00\% & 6     & 86.08 \\
				& 6     & 2     & 9137.37 & 9339.90 & 2.17\% & 24    & T     & 9104.73 & 9339.90 & 2.52\% & 23    & T     & 8972.98 & 16553.29 & 45.79\% & 25    & T     & 10329.96 & 13867.41 & 25.51\% & 9     & T \\
				& 6     & 3     & 9047.58 & 10749.89 & 15.84\% & 18    & T     & 9047.58 & 10749.89 & 15.84\% & 18    & T     & 9546.12 & 28113.18 & 66.04\% & 20    & T     & 10972.16 & 25527.78 & 57.02\% & 11    & T \\
				& 7     & 1     & 6853.43 & 6853.43 & 0.00\% & 9     & 50.14 & 6853.43 & 6853.43 & 0.00\% & 9     & 52.55 & 7475.97 & 7475.97 & 0.00\% & 3     & 3.79  & 7475.97 & 7475.97 & 0.00\% & 2     & 4.06 \\
				& 7     & 2     & 8018.31 & 8146.95 & 1.58\% & 24    & T     & 8019.73 & 8146.95 & 1.56\% & 25    & T     & 8665.60 & 8667.66 & 0.02\% & 20    & 2272.87 & 8667.66 & 8667.66 & 0.00\% & 10    & 1054.88 \\
				& 7     & 3     & 8233.66 & 9097.02 & 9.49\% & 23    & T     & 8233.66 & 9097.02 & 9.49\% & 23    & T     & 7876.11 & 15628.32 & 49.60\% & 29    & T     & 9141.63 & 14005.55 & 34.73\% & 8     & T \\
				& 8     & 1     & 5718.85 & 5718.85 & 0.00\% & 3     & 2.06  & 5718.85 & 5718.85 & 0.00\% & 3     & 1.97  & 6530.64 & 6533.51 & 0.04\% & 5     & 14.85 & 6530.64 & 6533.51 & 0.04\% & 5     & 16.85 \\
				& 8     & 2     & 6792.38 & 7204.49 & 5.72\% & 30    & T     & 6748.79 & 7204.49 & 6.33\% & 29    & T     & 7406.77 & 7696.92 & 3.77\% & 27    & T     & 7444.13 & 7696.92 & 3.28\% & 20    & T \\
				& 8     & 3     & 7055.41 & 8089.93 & 12.79\% & 27    & T     & 7045.60 & 8089.93 & 12.91\% & 26    & T     & 7392.81 & 8751.71 & 15.53\% & 25    & T     & 7874.48 & 8681.52 & 9.30\% & 11    & T \\
				\hline
				\multicolumn{3}{|c|}{\textbf{Average}} &       &       & \textbf{5.29\%} & \textbf{6.67}  & \textbf{50.11} &       &       & \textbf{5.40\%} & \textbf{6.67}  & \textbf{45.14} &       &       & \textbf{20.09\%} & \textbf{9.75}  & \textbf{654.32} &       &       & \textbf{14.43\%} & \textbf{5.75}  & \textbf{290.47} \\
				\hline
				\hline
				\multirow{9}[2]{*}{0.4} & 6     & 1     & 9001.85 & 9001.85 & 0.00\% & 12    & 612.62 & 9001.85 & 9001.85 & 0.00\% & 12    & 573.70 & 9822.07 & 9822.07 & 0.00\% & 14    & 978.75 & 9822.07 & 9822.07 & 0.00\% & 8     & 405.79 \\
				& 6     & 2     & 10176.89 & 11801.40 & 13.77\% & 21    & T     & 10178.93 & 11801.40 & 13.75\% & 22    & T     & 9514.23 & 23361.76 & 59.27\% & 25    & T     & 11434.54 & 18314.39 & 37.57\% & 8     & T \\
				& 6     & 3     & 9865.02 & 14827.25 & 33.47\% & 17    & T     & 9865.02 & 14827.25 & 33.47\% & 17    & T     & 10272.99 & 45001.75 & 77.17\% & 17    & T     & 11889.96 & 42561.90 & 72.06\% & 11    & T \\
				& 7     & 1     & 7463.87 & 7463.87 & 0.00\% & 10    & 120.04 & 7463.87 & 7463.87 & 0.00\% & 10    & 119.46 & 8279.68 & 8279.68 & 0.00\% & 8     & 104.65 & 8279.68 & 8279.68 & 0.00\% & 7     & 110.77 \\
				& 7     & 2     & 8969.83 & 10499.90 & 14.57\% & 23    & T     & 8969.83 & 10499.90 & 14.57\% & 23    & T     & 9498.56 & 10663.06 & 10.92\% & 22    & T     & 9853.03 & 10663.06 & 7.60\% & 11    & T \\
				& 7     & 3     & 9091.07 & 11919.62 & 23.73\% & 21    & T     & 9091.07 & 11919.62 & 23.73\% & 21    & T     & 8520.38 & 22732.53 & 62.52\% & 24    & T     & 10453.59 & 22707.14 & 53.96\% & 8     & T \\
				& 8     & 1     & 6200.78 & 6200.78 & 0.00\% & 4     & 4.28  & 6200.78 & 6200.78 & 0.00\% & 4     & 4.12  & 7308.85 & 7308.85 & 0.00\% & 11    & 312.85 & 7308.85 & 7308.85 & 0.00\% & 11    & 312.50 \\
				& 8     & 2     & 7496.31 & 9122.22 & 17.82\% & 28    & T     & 7484.67 & 9122.22 & 17.95\% & 27    & T     & 7710.40 & 9377.39 & 17.78\% & 26    & T     & 8299.00 & 9498.59 & 12.63\% & 16    & T \\
				& 8     & 3     & 7880.40 & 10397.08 & 24.21\% & 22    & T     & 7880.40 & 10397.08 & 24.21\% & 22    & T     & 8143.25 & 11773.64 & 30.83\% & 21    & T     & 8888.18 & 11773.64 & 24.51\% & 8     & T \\
				\hline
				\multicolumn{3}{|c|}{\textbf{Average}} &       &       & \textbf{14.17\%} & \textbf{8.67}  & \textbf{245.64} &       &       & \textbf{14.19\%} & \textbf{8.67}  & \textbf{232.42} &       &       & \textbf{28.72\%} & \textbf{11.00} & \textbf{465.42} &       &       & \textbf{23.15\%} & \textbf{8.67}  & \textbf{276.35} \\
				\hline
			\end{tabular}%
			\label{tbl:PMD}%
		}
	\end{table}%
\end{landscape}

\begin{landscape}
	\begin{table}[htbp]
		\centering
		\caption{DDU Approximation for  Capacitated Reliable P-Median Problem with Double Post-Disruption Demand}
		\resizebox{!}{0.85\height}{
			\begin{tabular}{|ccc|ccccc|ccccc|ccccc|}
				\hline
				\multirow{2}[4]{*}{$\rho$} & \multirow{2}[4]{*}{p} & \multirow{2}[4]{*}{k} & \multicolumn{5}{c|}{DIU by basic C\&CG} & \multicolumn{5}{c|}{$\mathcal{U}(\mathbf x)$ by Parametric C\&CG} & \multicolumn{5}{c|}{ $\mathcal{U}^{kq}(\mathbf x)$ by Parametric C\&CG} \\
				\cline{4-18}          &       &       & LB    & UB    & Gap   & Iter  & Time(s) & LB    & UB    & Gap   & Iter  & Time(s) & LB    & UB    & Gap   & Iter  & Time(s) \\
				\cline{4-18}    \multirow{9}[2]{*}{0.2} & 6     & 1     & 9725.62 & 9725.62 & 0.00\% & 5     & 24.17 & 9698.27 & 9770.88 & 0.74\% & 8     & 117.41 & 9725.62 & 9725.62 & 0.00\% & 5     & 27.12 \\
				& 6     & 2     & 9482.75 & 17408.62 & 45.53\% & 23    & T     & 11239.04 & 15597.74 & 27.94\% & 9     & T     & 11551.33 & 16383.55 & 29.49\% & 9     & T \\
				& 6     & 3     & 10057.39 & 30528.11 & 67.06\% & 18    & T     & 11684.33 & 29332.53 & 60.17\% & 11    & T     & 12232.39 & 29808.22 & 58.96\% & 11    & T \\
				& 7     & 1     & 8394.50 & 8394.50 & 0.00\% & 5     & 20.52 & 8096.03 & 8394.50 & 3.56\% & 5     & 23.64 & 8394.50 & 8394.50 & 0.00\% & 6     & 39.41 \\
				& 7     & 2     & 9742.34 & 9868.46 & 1.28\% & 26    & T     & 9350.13 & 9757.97 & 4.18\% & 10    & 1400.76 & 9488.19 & 10189.11 & 6.88\% & 7     & 371.66 \\
				& 7     & 3     & 8245.26 & 18491.38 & 55.41\% & 24    & T     & 9977.56 & 20492.54 & 51.31\% & 10    & T     & 10315.71 & 19090.51 & 45.96\% & 8     & T \\
				& 8     & 1     & 7284.81 & 7284.81 & 0.00\% & 5     & 17.75 & 7153.57 & 7284.81 & 1.80\% & 5     & 15.87 & 7284.81 & 7284.81 & 0.00\% & 5     & 18.76 \\
				& 8     & 2     & 8268.59 & 8268.59 & 0.00\% & 18    & 702.26 & 8182.99 & 8497.55 & 3.70\% & 17    & T     & 8268.59 & 8268.59 & 0.00\% & 18    & 1478.81 \\
				& 8     & 3     & 8727.38 & 11420.92 & 23.58\% & 35    & T     & 8699.43 & 10876.21 & 20.01\% & 11    & T     & 9149.72 & 11420.92 & 19.89\% & 12    & T \\
				\hline
				\multicolumn{3}{|c|}{\textbf{Average}} &       &       & \textbf{21.43\%} & \textbf{17.67}&  &       &       & \textbf{19.27\%} & \textbf{9.56}  &  &       &       & \textbf{17.91\%} & \textbf{9.00}  &  \\
				\hline
				\hline
				\multirow{9}[2]{*}{0.4} & 6     & 1     & 11001.45 & 11001.45 & 0.00\% & 10    & 239.17 & 10699.08 & 11296.02 & 5.28\% & 7     & 280.96 & 11001.45 & 11001.45 & 0.00\% & 8     & 168.81 \\
				& 6     & 2     & 10178.13 & 26603.99 & 61.74\% & 21    & T     & 12397.74 & 22903.99 & 45.87\% & 8     & T     & 13236.18 & 19671.50 & 32.71\% & 10    & T \\
				& 6     & 3     & 10963.58 & 50161.33 & 78.14\% & 16    & T     & 13050.76 & 50446.90 & 74.13\% & 10    & T     & 14334.97 & 45766.24 & 68.68\% & 10    & T \\
				& 7     & 1     & 9492.39 & 9492.39 & 0.00\% & 7     & 91.77 & 9026.23 & 9623.16 & 6.20\% & 7     & 166.03 & 9492.39 & 9492.39 & 0.00\% & 7     & 119.38 \\
				& 7     & 2     & 11181.10 & 12647.93 & 11.60\% & 29    & T     & 10737.95 & 13111.66 & 18.10\% & 9     & T     & 10940.08 & 12571.25 & 12.98\% & 13    & T \\
				& 7     & 3     & 9125.29 & 29058.40 & 68.60\% & 23    & T     & 11258.55 & 29305.95 & 61.58\% & 8     & T     & 12214.81 & 25282.52 & 51.69\% & 9     & T \\
				& 8     & 1     & 8155.80 & 8155.80 & 0.00\% & 6     & 24.80 & 8026.10 & 8155.80 & 1.59\% & 8     & 238.27 & 8155.80 & 8155.80 & 0.00\% & 6     & 34.18 \\
				& 8     & 2     & 9381.25 & 10385.98 & 9.67\% & 29    & T     & 9221.47 & 10483.92 & 12.04\% & 18    & T     & 9380.34 & 10500.46 & 10.67\% & 16    & 2221.66 \\
				& 8     & 3     & 10352.95 & 17112.06 & 39.50\% & 25    & T     & 9795.56 & 17112.06 & 42.76\% & 10    & T     & 10943.25 & 17112.06 & 36.05\% & 11    & T \\
				\hline
				\multicolumn{3}{|c|}{\textbf{Average}} &       &       & \textbf{29.92\%} & \textbf{18.44} &  &       &       & \textbf{29.73\%} & \textbf{9.44}  &   &       &       & \textbf{23.64\%} & \textbf{10.00}  &  \\
				\hline
			\end{tabular}%
			\label{tbl:CapPMDD}%
		}
	\end{table}%
\end{landscape}

\begin{landscape}
	\begin{table}[htbp]
		\centering
		\caption{DDU Approximation for Capacitated Reliable P-Median Problem with Sorting}
		\resizebox{!}{0.85\height}{
			\begin{tabular}{|ccc|ccccc|ccccc|ccccc|}
				\hline
				\multirow{2}[4]{*}{$\rho$} & \multirow{2}[4]{*}{p} & \multirow{2}[4]{*}{k} & \multicolumn{5}{c|}{$\mathcal{U}(\mathbf x)$ (exact)} & \multicolumn{5}{c|}{$\mathcal{X}^r$ and $\mathcal{U}^r(\mathbf x), q_1=k+2$} & \multicolumn{5}{c|}{$\mathcal{X}^{rs}$,  $\mathcal{U}^r(\mathbf x)$ and $\mathcal{U}^s(\mathbf x)$, $q_1=q_2=k+2$} \\
				\cline{4-18}          &       &       & LB    & UB    & Gap   & Iter.  & Time(s) & LB    & UB    & Gap   & Iter.  & Time(s) & LB    & UB    & Gap   & Iter.  & Time(s) \\
				\hline
				\multirow{9}[2]{*}{0.2} & 6     & 1     & 8979.84 & 8979.84 & 0.00\% & 6     & 86.08 & 8068.97 & 9066.46 & 11.00\% & 4     & 21.67 & 8357.43 & 9066.46 & 7.82\% & 6     & 836.51 \\
				& 6     & 2     & 10329.96 & 13867.41 & 25.51\% & 9     & T     & 9385.68 & 15840.84 & 40.75\% & 10    & T     & 8920.30 & 15857.32 & 43.75\% & 6     & T \\
				& 6     & 3     & 10972.16 & 25527.78 & 57.02\% & 11    & T     & 9980.93 & 26843.26 & 62.82\% & 10    & T     & 9740.68 & 26859.07 & 63.73\% & 8     & T \\
				& 7     & 1     & 7475.97 & 7475.97 & 0.00\% & 2     & 4.06  & 6699.67 & 7475.97 & 10.38\% & 2     & 1.65  & 6699.67 & 7475.97 & 10.38\% & 2     & 8.05 \\
				& 7     & 2     & 8667.66 & 8667.66 & 0.00\% & 10    & 1054.88 & 6739.85 & 8667.66 & 22.24\% & 2     & 1.98  & 6739.85 & 8667.66 & 22.24\% & 2     & 10.08 \\
				& 7     & 3     & 9141.63 & 14005.55 & 34.73\% & 8     & T     & 8146.49 & 15364.39 & 46.98\% & 9     & T     & 7909.98 & 15613.91 & 49.34\% & 6     & T \\
				& 8     & 1     & 6530.64 & 6533.51 & 0.04\% & 5     & 16.85 & 5729.79 & 6533.51 & 12.30\% & 2     & 1.19  & 5729.79 & 6533.51 & 12.30\% & 2     & 7.91 \\
				& 8     & 2     & 7444.13 & 7696.92 & 3.28\% & 20    & T     & 5871.52 & 7725.20 & 24.00\% & 3     & 5.03  & 5871.52 & 7725.20 & 24.00\% & 3     & 20.78 \\
				& 8     & 3     & 7874.48 & 8681.52 & 9.30\% & 11    & T     & 6467.46 & 8636.84 & 25.12\% & 6     & 171.18 & 6467.46 & 8636.84 & 25.12\% & 6     & 801.15 \\
				\hline
				\multicolumn{3}{|c|}{\textbf{Average}} &       &       & \textbf{14.43\%} & \textbf{9.11}  &       &       &       & \textbf{28.40\%} & \textbf{5.33}  &       &       &       & \textbf{28.74\%} & \textbf{4.56}  &  \\
				\hline
				\hline
				\multirow{9}[2]{*}{0.4} & 6     & 1     & 9822.07 & 9822.07 & 0.00\% & 8     & 405.79 & 8026.28 & 10268.74 & 21.84\% & 4     & 20.73 & 8302.48 & 9986.93 & 16.87\% & 6     & 906.09 \\
				& 6     & 2     & 11434.54 & 18314.39 & 37.57\% & 8     & T     & 9856.85 & 19410.38 & 49.22\% & 8     & T     & 9221.99 & 22502.02 & 59.02\% & 6     & T \\
				& 6     & 3     & 11889.96 & 42561.90 & 72.06\% & 11    & T     & 10741.62 & 42743.94 & 74.87\% & 10    & T     & 10380.24 & 42743.94 & 75.72\% & 7     & T \\
				& 7     & 1     & 8279.68 & 8279.68 & 0.00\% & 7     & 110.77 & 6694.58 & 8279.68 & 19.14\% & 2     & 1.88  & 6694.58 & 8279.68 & 19.14\% & 2     & 10.12 \\
				& 7     & 2     & 9853.03 & 10663.06 & 7.60\% & 11    & T     & 6722.95 & 10663.06 & 36.95\% & 2     & 1.95  & 6722.95 & 10663.06 & 36.95\% & 2     & 10.91 \\
				& 7     & 3     & 10453.59 & 22707.14 & 53.96\% & 8     & T     & 8652.43 & 24207.33 & 64.26\% & 8     & T     & 8437.74 & 24207.33 & 65.14\% & 6     & T \\
				& 8     & 1     & 7308.85 & 7308.85 & 0.00\% & 11    & 312.50 & 5729.79 & 7337.22 & 21.91\% & 2     & 1.01  & 5729.79 & 7337.22 & 21.91\% & 2     & 7.18 \\
				& 8     & 2     & 8299.00 & 9498.59 & 12.63\% & 16    & T     & 5836.09 & 9720.60 & 39.96\% & 2     & 1.45  & 5836.09 & 9720.60 & 39.96\% & 2     & 10.46 \\
				& 8     & 3     & 8888.18 & 11773.64 & 24.51\% & 8     & T     & 6174.89 & 11386.81 & 45.77\% & 3     & 8.24  & 6174.89 & 11386.81 & 45.77\% & 3     & 47.13 \\
				\hline
				\multicolumn{3}{|c|}{\textbf{Average}} &       &       & \textbf{23.15\%} & \textbf{9.78}  &       &       &       & \textbf{41.55\%} & \textbf{4.56}  &       &       &       & \textbf{42.28\%} & \textbf{4.00}  &  \\
				\hline
			\end{tabular}%
			\label{tbl:AppCapPM_sorting}%
		}
	\end{table}%
\end{landscape}

\begin{landscape}
	\begin{table}[htbp]
		\centering
		\caption{Computational Results of the Approximation Scheme for Mixed Integer Recourse}
		\resizebox{!}{0.85\height}{
			\begin{tabular}{|c|cccc|ccccc|ccccc|cc|}
				\hline
				\multirow{2}[4]{*}{\# of Sites} & \multirow{2}[4]{*}{Fixed Cost} & \multirow{2}[4]{*}{$\underline{\boldsymbol\xi}$} & \multirow{2}[4]{*}{$\overline{\boldsymbol\xi}$} & \multirow{2}[4]{*}{$\alpha$} & \multicolumn{5}{c|}{with Temp. Facilities} & \multicolumn{5}{c|}{without Temp.  Facilities} & \multicolumn{2}{c|}{Difference} \\
				\cline{6-17}          &       &       &       &       & LB    & UB    & Gap   & Iter  & Time(s) & LB    & UB    & Gap   & Iter  & Time(s) & Iter  & Time(s) \\
				\hline
				\multirow{11}[6]{*}{25} & \multirow{5}[2]{*}{High} & 0.01  & 0.02  & 0.01  & 3264.38 & 3264.84 & 0.01\% & 2     & 9.78  & 3887.13 & 3887.87 & 0.02\% & 2     & 6.96  & 0     & 2.82 \\
				&       & 0.02  & 0.05  & 0.02  & 3278.34 & 3279.80 & 0.04\% & 2     & 1.71  & 3914.89 & 3915.67 & 0.02\% & 2     & 0.69  & 0     & 1.02 \\
				&       & 0.05  & 0.08  & 0.05  & 3300.30 & 3302.52 & 0.07\% & 2     & 1.68  & 3952.98 & 3955.06 & 0.05\% & 2     & 0.75  & 0     & 0.93 \\
				&       & 0.08  & 0.12  & 0.08  & 3325.19 & 3328.45 & 0.11\% & 2     & 2.14  & 3999.18 & 4001.03 & 0.05\% & 2     & 1.51  & 0     & 0.63 \\
				&       & 0.1   & 0.15  & 0.1   & 3347.20 & 3347.20 & 0.00\% & 3     & 2.64  & 4032.20 & 4034.51 & 0.06\% & 2     & 0.90  & 1     & 1.74 \\
				\cline{2-17}          & \multirow{5}[2]{*}{Low} & 0.01  & 0.02  & 0.01  & 101.70 & 101.70 & 0.00\% & 3     & 14.79 & 512.70 & 512.70 & 0.00\% & 3     & 7.04  & 0     & 7.75 \\
				&       & 0.02  & 0.05  & 0.02  & 132.32 & 132.32 & 0.00\% & 3     & 5.00  & 543.32 & 543.32 & 0.00\% & 3     & 1.02  & 0     & 3.98 \\
				&       & 0.05  & 0.08  & 0.05  & 174.80 & 174.80 & 0.00\% & 3     & 4.62  & 585.81 & 585.81 & 0.00\% & 3     & 2.28  & 0     & 2.34 \\
				&       & 0.08  & 0.12  & 0.08  & 224.16 & 224.16 & 0.00\% & 3     & 3.29  & 635.16 & 635.16 & 0.00\% & 3     & 2.84  & 0     & 0.46 \\
				&       & 0.1   & 0.15  & 0.1   & 260.43 & 260.43 & 0.00\% & 3     & 3.30  & 671.43 & 671.43 & 0.00\% & 3     & 1.30  & 0     & 1.99 \\
				\cline{2-17}          & \multicolumn{4}{c|}{\textbf{Average}}  &       &       & \textbf{0.02\%} & \textbf{2.6}   & \textbf{4.89}  &       &       & \textbf{0.02\%} & \textbf{2.5}   & \textbf{2.53}  & \textbf{0.1}   & \textbf{2.37} \\
				\hline
				\hline
				\multirow{11}[6]{*}{40} & \multirow{5}[2]{*}{High} & 0.01  & 0.02  & 0.01  & 3095.24 & 3095.24 & 0.00\% & 2     & 12.79 & 4085.99 & 4086.89 & 0.02\% & 2     & 6.31  & 0     & 6.48 \\
				&       & 0.02  & 0.05  & 0.02  & 3124.46 & 3124.46 & 0.00\% & 2     & 5.46  & 4115.87 & 4116.12 & 0.01\% & 2     & 1.47  & 0     & 3.99 \\
				&       & 0.05  & 0.08  & 0.05  & 3166.55 & 3166.55 & 0.00\% & 2     & 10.39 & 4156.27 & 4158.20 & 0.05\% & 2     & 4.33  & 0     & 6.06 \\
				&       & 0.08  & 0.12  & 0.08  & 3204.31 & 3204.63 & 0.01\% & 2     & 14.11 & 4205.14 & 4207.12 & 0.05\% & 2     & 5.60  & 0     & 8.50 \\
				&       & 0.1   & 0.15  & 0.1   & 3223.36 & 3223.76 & 0.01\% & 2     & 10.73 & 4240.19 & 4242.67 & 0.06\% & 2     & 4.97  & 0     & 5.76 \\
				\cline{2-17}          & \multirow{5}[2]{*}{Low} & 0.01  & 0.02  & 0.01  & -1720.08 & -1720.08 & 0.00\% & 2     & 27.73 & -923.69 & -923.69 & 0.00\% & 3     & 8.38  & -1    & 19.35 \\
				&       & 0.02  & 0.05  & 0.02  & -1670.09 & -1670.09 & 0.00\% & 2     & 8.11  & -873.71 & -873.71 & 0.00\% & 3     & 2.54  & -1    & 5.58 \\
				&       & 0.05  & 0.08  & 0.05  & -1599.13 & -1598.97 & 0.01\% & 2     & 12.17 & -802.58 & -802.62 & 0.00\% & 4     & 7.12  & -2    & 5.05 \\
				&       & 0.08  & 0.12  & 0.08  & -1516.43 & -1516.34 & 0.01\% & 2     & 7.31  & -724.21 & -724.21 & 0.00\% & 3     & 3.60  & -1    & 3.71 \\
				&       & 0.1   & 0.15  & 0.1   & -1456.08 & -1455.97 & 0.01\% & 2     & 29.49 & -666.94 & -666.94 & 0.00\% & 3     & 2.71  & -1    & 26.78 \\
				\cline{2-17}          & \multicolumn{4}{c|}{\textbf{Average}}  &       &       & \textbf{0.00\%} & \textbf{2.00}  & \textbf{13.83} &       &       & \textbf{0.02\%} & \textbf{2.60}  & \textbf{4.70}  & \textbf{-0.60} & \textbf{9.13} \\
				\hline
			\end{tabular}%
			\label{tbl:MIPRC}%
		}
	\end{table}%
\end{landscape}

\newpage
\begin{landscape}
	\begin{table}[htbp]
		\centering
		\caption{Uncapacitated Reliable P-Median Problem (with Enhancements)}
		\resizebox{!}{0.85\height}{
			\begin{tabular}{|ccc|ccccc|ccccc|ccccc|}
				\hline
				\multirow{2}[4]{*}{$\rho$} & \multirow{2}[4]{*}{p} & \multirow{2}[4]{*}{k} & \multicolumn{5}{c|}{Parametric C\&CG} & \multicolumn{5}{c|}{Parametric C\&CG with Pareto}   & \multicolumn{5}{c|}{Parametric C\&CG with Uniqueness} \\
				\cline{4-18}          &       &       & LB    & UB    & Gap   & Iter  & Time(s) & LB    & UB    & Gap   & Iter  & Time(s) & LB    & UB    & Gap   & Iter  & Time(s) \\
				\hline
    \multirow{9}[2]{*}{0.2} & 6     & 1     & 8378.25 & 8378.25 & 0.00\% & 8     & 80.90 & 8378.25 & 8378.25 & 0.00\% & 8     & 100.27 & 8378.25 & 8378.25 & 0.00\% & 8     & 174.35 \\
& 6     & 2     & 9104.73 & 9339.90 & 2.52\% & 23    & T     & 9137.37 & 9339.90 & 2.17\% & 24    & T     & 9230.12 & 9405.91 & 1.87\% & 15    & T \\
& 6     & 3     & 9047.58 & 10749.89 & 15.84\% & 18    & T     & 9047.58 & 10749.89 & 15.84\% & 18    & T     & 9470.48 & 10749.89 & 11.90\% & 11    & T \\
& 7     & 1     & 6853.43 & 6853.43 & 0.00\% & 9     & 52.55 & 6853.43 & 6853.43 & 0.00\% & 9     & 52.95 & 6853.43 & 6853.43 & 0.00\% & 7     & 36.43 \\
& 7     & 2     & 8019.73 & 8146.95 & 1.56\% & 25    & T     & 8019.73 & 8146.95 & 1.56\% & 25    & T     & 8146.95 & 8146.95 & 0.00\% & 18    & 3279.72 \\
& 7     & 3     & 8233.66 & 9097.02 & 9.49\% & 23    & T     & 8218.92 & 9097.02 & 9.65\% & 22    & T     & 8406.53 & 9225.67 & 8.88\% & 16    & T \\
& 8     & 1     & 5718.85 & 5718.85 & 0.00\% & 3     & 1.97  & 5718.85 & 5718.85 & 0.00\% & 3     & 2.09  & 5718.85 & 5718.85 & 0.00\% & 3     & 2.16 \\
& 8     & 2     & 6748.79 & 7204.49 & 6.33\% & 29    & T     & 6778.88 & 7204.49 & 5.91\% & 29    & T     & 6913.56 & 7148.04 & 3.28\% & 19    & T \\
& 8     & 3     & 7045.60 & 8089.93 & 12.91\% & 26    & T     & 7064.58 & 8089.93 & 12.67\% & 28    & T     & 7207.98 & 7965.62 & 9.51\% & 16    & T \\
\hline
\multicolumn{3}{|c|}{\textbf{Average}} &       &       & \textbf{5.40\%} & \textbf{6.67} & \textbf{45.14} &       &       & \textbf{5.31\%} & \textbf{6.67} & \textbf{51.77} &       &       & \textbf{3.94\%} & \textbf{9.00} & \textbf{873.16} \\
\hline
\hline
\multirow{9}[2]{*}{0.4} & 6     & 1     & 9001.85 & 9001.85 & 0.00\% & 12    & 573.70 & 9001.85 & 9001.85 & 0.00\% & 12    & 595.01 & 9001.85 & 9001.85 & 0.00\% & 8     & 414.51 \\
& 6     & 2     & 10178.93 & 11801.40 & 13.75\% & 22    & T     & 10178.93 & 11801.40 & 13.75\% & 22    & T     & 10235.21 & 11772.54 & 13.06\% & 11    & T \\
& 6     & 3     & 9865.02 & 14827.25 & 33.47\% & 17    & T     & 9865.02 & 14827.25 & 33.47\% & 17    & T     & 10798.42 & 14827.25 & 27.17\% & 11    & T \\
& 7     & 1     & 7463.87 & 7463.87 & 0.00\% & 10    & 119.46 & 7463.87 & 7463.87 & 0.00\% & 10    & 120.79 & 7463.87 & 7463.87 & 0.00\% & 8     & 141.61 \\
& 7     & 2     & 8969.83 & 10499.90 & 14.57\% & 23    & T     & 8969.83 & 10499.90 & 14.57\% & 23    & T     & 9033.49 & 10427.17 & 13.37\% & 15    & T \\
& 7     & 3     & 9091.07 & 11919.62 & 23.73\% & 21    & T     & 9091.07 & 11919.62 & 23.73\% & 21    & T     & 9327.13 & 12878.77 & 27.58\% & 9     & T \\
& 8     & 1     & 6200.78 & 6200.78 & 0.00\% & 4     & 4.12  & 6200.78 & 6200.78 & 0.00\% & 4     & 4.36  & 6200.78 & 6200.78 & 0.00\% & 3     & 3.28 \\
& 8     & 2     & 7484.67 & 9122.22 & 17.95\% & 27    & T     & 7484.67 & 9122.22 & 17.95\% & 27    & T     & 7880.31 & 8598.00 & 8.35\% & 15    & T \\
& 8     & 3     & 7880.40 & 10397.08 & 24.21\% & 22    & T     & 7880.40 & 10397.08 & 24.21\% & 22    & T     & 8305.78 & 10830.24 & 23.31\% & 12    & T \\
\hline
\multicolumn{3}{|c|}{\textbf{Average}} &       &       & \textbf{14.19\%} & \textbf{8.67} & \textbf{232.42} &       &       & \textbf{14.19\%} & \textbf{8.67} & \textbf{240.06} &       &       & \textbf{12.54\%} & \textbf{6.33} & \textbf{186.47} \\
				\hline
			\end{tabular}%
			\label{tbl:UnCapPM}%
		}
	\end{table}%
\end{landscape}

\newpage
\begin{landscape}
	\begin{table}[htbp]
		\centering
		\caption{Capacitated Reliable P-Median Problem (with Enhancements)}
		\resizebox{!}{0.85\height}{
			\begin{tabular}{|ccc|ccccc|ccccc|ccccc|}
				\hline
				\multirow{2}[4]{*}{$\rho$} & \multirow{2}[4]{*}{p} & \multirow{2}[4]{*}{k} & \multicolumn{5}{c|}{Parametric C\&CG} & \multicolumn{5}{c|}{Parametric C\&CG with Pareto}   & \multicolumn{5}{c|}{Parametric C\&CG with Uniqueness} \\
				\cline{4-18}          &       &       & LB    & UB    & Gap   & Iter  & Time(s) & LB    & UB    & Gap   & Iter  & Time(s) & LB    & UB    & Gap   & Iter  & Time(s) \\
				\hline
    \multirow{9}[2]{*}{0.2} & 6     & 1     & 8979.84 & 8979.84 & 0.00\% & 6     & 86.08 & 8979.84 & 8979.84 & 0.00\% & 6     & 112.57 & 8979.84 & 8979.84 & 0.00\% & 6     & 126.62 \\
& 6     & 2     & 10329.96 & 13867.41 & 25.51\% & 9     & T     & 10421.02 & 14774.25 & 29.46\% & 9     & T     & 9980.93 & 14497.82 & 31.16\% & 11    & T \\
& 6     & 3     & 10972.16 & 25527.78 & 57.02\% & 11    & T     & 10896.57 & 27409.29 & 60.24\% & 10    & T     & 10085.86 & 25274.49 & 60.09\% & 13    & T \\
& 7     & 1     & 7475.97 & 7475.97 & 0.00\% & 2     & 4.06  & 7475.97 & 7475.97 & 0.00\% & 2     & 1.74  & 7475.97 & 7475.97 & 0.00\% & 2     & 2.18 \\
& 7     & 2     & 8667.66 & 8667.66 & 0.00\% & 10    & 1054.88 & 8667.66 & 8667.66 & 0.00\% & 10    & 1174.19 & 8664.13 & 8667.66 & 0.04\% & 9     & 1118.13 \\
& 7     & 3     & 9141.63 & 14005.55 & 34.73\% & 8     & T     & 9141.63 & 14005.55 & 34.73\% & 8     & T     & 8666.96 & 15540.92 & 44.23\% & 11    & T \\
& 8     & 1     & 6530.64 & 6533.51 & 0.04\% & 5     & 16.85 & 6530.64 & 6533.51 & 0.04\% & 4     & 11.44 & 6533.51 & 6533.51 & 0.00\% & 4     & 13.69 \\
& 8     & 2     & 7444.13 & 7696.92 & 3.28\% & 20    & T     & 7442.05 & 7696.92 & 3.31\% & 19    & T     & 7474.26 & 7628.99 & 2.03\% & 15    & T \\
& 8     & 3     & 7874.48 & 8681.52 & 9.30\% & 11    & T     & 7912.82 & 8681.52 & 8.85\% & 11    & T     & 7686.06 & 8511.53 & 9.70\% & 13    & T \\
\hline
\multicolumn{3}{|c|}{\textbf{Average}} &       &       & \textbf{14.43\%} & \textbf{5.75} & \textbf{290.47} &       &       & \textbf{15.18\%} & \textbf{5.50} & \textbf{324.98} &       &       & \textbf{16.36\%} & \textbf{5.25} & \textbf{315.16} \\
\hline
\hline
\multirow{9}[2]{*}{0.4} & 6     & 1     & 9822.07 & 9822.07 & 0.00\% & 8     & 405.79 & 9822.07 & 9822.07 & 0.00\% & 8     & 436.02 & 9822.07 & 9822.07 & 0.00\% & 9     & 350.70 \\
& 6     & 2     & 11434.54 & 18314.39 & 37.57\% & 8     & T     & 11434.54 & 18314.39 & 37.57\% & 8     & T     & 10741.62 & 18149.99 & 40.82\% & 10    & T \\
& 6     & 3     & 11889.96 & 42561.90 & 72.06\% & 11    & T     & 12453.75 & 42599.75 & 70.77\% & 11    & T     & 10761.23 & 42743.94 & 74.82\% & 11    & T \\
& 7     & 1     & 8279.68 & 8279.68 & 0.00\% & 7     & 110.77 & 8279.68 & 8279.68 & 0.00\% & 7     & 99.41 & 8279.68 & 8279.68 & 0.00\% & 5     & 50.54 \\
& 7     & 2     & 9853.03 & 10663.06 & 7.60\% & 11    & T     & 9729.87 & 10663.06 & 8.75\% & 10    & T     & 9781.99 & 10663.06 & 8.26\% & 10    & T \\
& 7     & 3     & 10453.59 & 22707.14 & 53.96\% & 8     & T     & 10475.66 & 23675.55 & 55.75\% & 8     & T     & 9519.82 & 23763.45 & 59.94\% & 10    & T \\
& 8     & 1     & 7308.85 & 7308.85 & 0.00\% & 11    & 312.50 & 7308.85 & 7308.85 & 0.00\% & 10    & 279.30 & 7308.85 & 7308.85 & 0.00\% & 9     & 342.08 \\
& 8     & 2     & 8299.00 & 9498.59 & 12.63\% & 16    & T     & 8402.62 & 9489.03 & 11.45\% & 14    & T     & 8264.31 & 9067.48 & 8.86\% & 13    & T \\
& 8     & 3     & 8888.18 & 11773.64 & 24.51\% & 8     & T     & 8870.46 & 11773.64 & 24.66\% & 9     & T     & 8840.23 & 11773.64 & 24.92\% & 10    & T \\
\hline
\multicolumn{3}{|c|}{\textbf{Average}} &       &       & \textbf{23.15\%} & \textbf{8.67} & \textbf{276.35} &       &       & \textbf{23.22\%} & \textbf{8.33} & \textbf{271.58} &       &       & \textbf{24.18\%} & \textbf{7.67} & \textbf{247.77} \\
				\hline
			\end{tabular}%
			\label{tbl:CapPM}%
		}
	\end{table}%
\end{landscape}
\begin{landscape}
	\begin{table}[htbp]
		\centering
		\caption{Computational Results for SOC Recourse Problem}
		\resizebox{!}{0.85\height}{
			\begin{tabular}{|c|cccc|ccccc|ccccc|cc|}
				\hline
				\multirow{2}[4]{*}{\# of Sites} & \multirow{2}[4]{*}{Fixed Cost} & \multirow{2}[4]{*}{$\underline{\boldsymbol\xi}$} & \multirow{2}[4]{*}{$\overline{\boldsymbol\xi}$} & \multirow{2}[4]{*}{$\alpha$} & \multicolumn{5}{c|}{SOC Recourse}    & \multicolumn{5}{c|}{LP Recourse}      & \multicolumn{2}{c|}{Difference} \\
				\cline{6-17}          &       &       &       &       & LB    & UB    & Gap   & Iter  & Time(s) & LB    & UB    & Gap   & Iter  & Time(s) & Iter  & Time(s) \\
				\hline
				\multirow{11}[6]{*}{25} & \multirow{5}[2]{*}{High} & 0.01  & 0.02  & 0.01  & 8942.26 & 8942.26 & 0.00\% & 2     & 68.50 & 3766.50 & 3766.50 & 0.00\% & 3     & 8.82  & -1    & 59.68 \\
				&       & 0.02  & 0.05  & 0.02  & 8968.31 & 8968.31 & 0.00\% & 2     & 116.71 & 3767.34 & 3767.34 & 0.00\% & 3     & 1.23  & -1    & 115.48 \\
				&       & 0.05  & 0.08  & 0.05  & 9006.82 & 9006.83 & 0.00\% & 2     & 115.61 & 3768.49 & 3768.52 & 0.00\% & 3     & 2.15  & -1    & 113.46 \\
				&       & 0.08  & 0.12  & 0.08  & 9052.32 & 9052.32 & 0.00\% & 2     & 344.76 & 3769.89 & 3769.89 & 0.00\% & 3     & 0.67  & -1    & 344.09 \\
				&       & 0.1   & 0.15  & 0.1   & 9085.50 & 9085.51 & 0.00\% & 2     & 702.25 & 3770.90 & 3770.90 & 0.00\% & 3     & 1.10  & -1    & 701.16 \\
				\cline{2-17}          & \multirow{5}[2]{*}{Low} & 0.01  & 0.02  & 0.01  & 6638.06 & 6638.06 & 0.00\% & 2     & 155.64 & 2485.25 & 2485.25 & 0.00\% & 3     & 9.51  & -1    & 146.12 \\
				&       & 0.02  & 0.05  & 0.02  & 6664.11 & 6664.11 & 0.00\% & 2     & 123.36 & 2486.09 & 2486.09 & 0.00\% & 3     & 0.97  & -1    & 122.39 \\
				&       & 0.05  & 0.08  & 0.05  & 6702.62 & 6702.63 & 0.00\% & 2     & 440.53 & 2487.24 & 2487.27 & 0.00\% & 3     & 1.55  & -1    & 438.98 \\
				&       & 0.08  & 0.12  & 0.08  & 6748.12 & 6748.13 & 0.00\% & 2     & 711.16 & 2488.64 & 2488.64 & 0.00\% & 3     & 3.66  & -1    & 707.50 \\
				&       & 0.1   & 0.15  & 0.1   & 6781.30 & 6781.31 & 0.00\% & 2     & 911.51 & 2489.65 & 2489.65 & 0.00\% & 3     & 1.56  & -1    & 909.95 \\
				\cline{2-17}          & \multicolumn{4}{c|}{\textbf{Average}}  &       &       & \textbf{0.00\%} & \textbf{2} & \textbf{369.00} &       &       & \textbf{0.00\%} & \textbf{3} & \textbf{3.12} & \textbf{-1}    & \textbf{365.88} \\
				\hline
			\end{tabular}%
			\label{tbl:SOCPR}%
		}
		
		\vspace{20pt} \centering
		\caption{Computational Results for SOC Uncertainty Set}
		\resizebox{!}{0.85\height}{
			\begin{tabular}{|c|cccc|ccccc|ccccc|cc|}
				\hline
				\multirow{2}[4]{*}{\# of Sites} & \multirow{2}[4]{*}{Fixed Cost} & \multirow{2}[4]{*}{$\underline{\boldsymbol\xi}$} & \multirow{2}[4]{*}{$\overline{\boldsymbol\xi}$} & \multirow{2}[4]{*}{$\alpha$} & \multicolumn{5}{c|}{SOC DDU} & \multicolumn{5}{c|}{Polyhedral DDU}   & \multicolumn{2}{c|}{Difference} \\
				\cline{6-17}          &       &       &       &       & LB    & UB    & Gap   & Iter  & Time(s) & LB    & UB    & Gap   & Iter  & Time(s) & Iter  & Time(s) \\
				\hline
    \multirow{11}[6]{*}{25} & \multirow{5}[2]{*}{High} & 0     & 0.02  & 0.01  & 3889.66 & 3889.67 & 0.00\% & 2     & 5.22  & 3882.65 & 3883.32 & 0.02\% & 3     & 9.70  & -1    & -4.48 \\
&       & 0     & 0.05  & 0.02  & 3921.19 & 3923.49 & 0.06\% & 2     & 1.37  & 3904.72 & 3905.52 & 0.02\% & 3     & 1.15  & -1    & 0.22 \\
&       & 0     & 0.08  & 0.05  & 3957.32 & 3957.32 & 0.00\% & 3     & 13.26 & 3935.41 & 3935.41 & 0.00\% & 4     & 2.25  & -1    & 11.01 \\
&       & 0     & 0.12  & 0.08  & 3996.67 & 3996.67 & 0.00\% & 3     & 19.57 & 3971.11 & 3971.11 & 0.00\% & 4     & 2.15  & -1    & 17.43 \\
&       & 0     & 0.15  & 0.1   & 4021.77 & 4021.77 & 0.00\% & 3     & 12.90 & 3997.11 & 3997.11 & 0.00\% & 4     & 8.06  & -1    & 4.84 \\
\cline{2-17}          & \multirow{5}[2]{*}{Low} & 0     & 0.02  & 0.01  & 514.39 & 514.41 & 0.00\% & 2     & 6.34  & 508.04 & 508.37 & 0.06\% & 3     & 11.36 & -1    & -5.02 \\
&       & 0     & 0.05  & 0.02  & 549.44 & 549.45 & 0.00\% & 3     & 12.24 & 533.68 & 533.83 & 0.03\% & 3     & 5.17  & 0     & 7.07 \\
&       & 0     & 0.08  & 0.05  & 584.72 & 584.73 & 0.00\% & 3     & 43.63 & 566.28 & 566.28 & 0.00\% & 4     & 3.28  & -1    & 40.35 \\
&       & 0     & 0.12  & 0.08  & 622.99 & 623.00 & 0.00\% & 3     & 41.79 & 605.75 & 605.75 & 0.00\% & 4     & 4.66  & -1    & 37.14 \\
&       & 0     & 0.15  & 0.1   & 647.29 & 647.30 & 0.00\% & 3     & 16.15 & 634.66 & 634.66 & 0.00\% & 4     & 8.65  & -1    & 7.51 \\
\cline{2-17}          & \multicolumn{4}{c|}{\textbf{Average}}  &       &       & \textbf{0.01\%} & \textbf{2.7}   & \textbf{17.25} &       &       & \textbf{0.01\%} & \textbf{3.6}   & \textbf{5.64}  & \textbf{-0.9}  & \textbf{11.61} \\
				\hline
			\end{tabular}%
			\label{tbl:SOCPU}%
		}
	\end{table}%
\end{landscape}

\end{document}